# INTERVAL LINEAR ALGEBRA

**W. B. Vasantha Kandasamy**
**Florentin Smarandache**

**2010**

# INTERVAL LINEAR ALGEBRA

# CONTENTS

~ DEDICATED TO ~

Dr C.N Deivanayagam

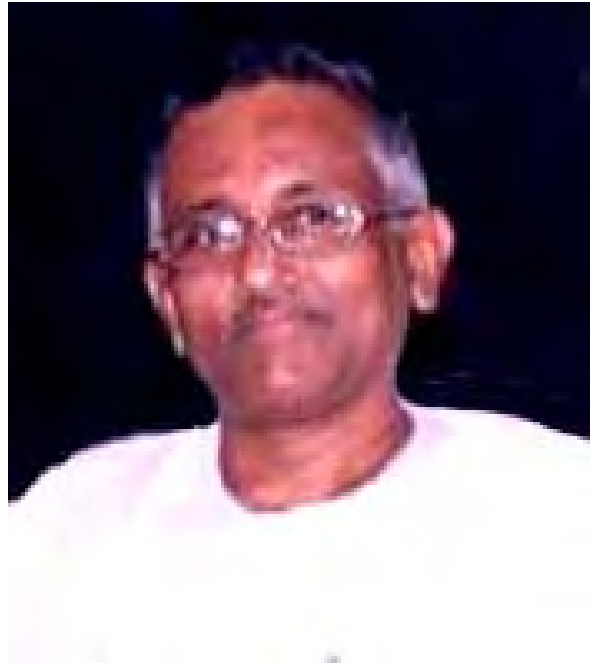

*This book is dedicated to Dr C.N Deivanayagam, Founder, Health India Foundation for his unostentatious service to all patients, especially those who are economically impoverished and socially marginalized. He was a pioneer in serving people living with HIV/AIDS at the Government Hospital of Thoracic Medicine (Tambaram Chennai). When the first author of this book had an opportunity of interacting with the patients, she learnt of his tireless service. His innovative practice of combining traditional Siddha Medicine alongside Allopathic remedies and his advocacy of ancient systems co-existing with modern health care distinguishes him. This dedication is a mere token of appreciation for his humanitarian service.*

# PREFACE

This Interval arithmetic or interval mathematics developed in 1950's and 1960's by mathematicians as an approach to putting bounds on rounding errors and measurement error in mathematical computations. However no proper interval algebraic structures have been defined or studies. In this book we for the first time introduce several types of interval linear algebras and study them.

This structure has become indispensable for these concepts will find applications in numerical optimization and validation of structural designs.

In this book we use only special types of intervals and introduce the notion of different types of interval linear algebras and interval vector spaces using the intervals of the form [0, a] where the intervals are from $Z_n$ or $Z^+ \cup \{0\}$ or $Q^+ \cup \{0\}$ or $R^+ \cup \{0\}$.

A systematic development is made starting from set interval vector spaces to group interval vector spaces. Vector spaces are taken as interval polynomials or interval matrices or just intervals over suitable sets or semigroups or groups. Main feature of this book is the authors have given over 350 examples.

This book has six chapters. Chapter one is introductory in nature. Chapter two introduces the notion of set interval linear algebras of type one and two. Set fuzzy interval linear algebras and their algebras and their properties are discussed in chapter three.

Chapter four introduces several types of interval linear bialgebras and bivector spaces and studies them. The possible applications are given in chapter five. Chapter six suggests nearly 110 problems of all levels.

The authors deeply acknowledge Dr. Kandasamy for the proof reading and Meena and Kama for the formatting and designing of the book.

<div style="text-align: right;">
W.B.VASANTHA KANDASAMY  
FLORENTIN SMARANDACHE
</div>

**Chapter One**

# INTRODUCTION

In this chapter we just define some basic properties of intervals used in this book. Throughout this book [a, b] denotes an interval $a \leq b$. If $a = b$ we say the interval degenerates to a point a. We assume the intervals [a, b] is such that $0 \leq a \leq b$. We just give the notations.

*Notations:* Let
$Z_I^+ = \{[a, b] \mid a, b \in Z^+ \cup \{0\}, a \leq b\}$
$Q_I^+ = \{[a, b] \mid a, b \in Q^+ \cup \{0\}, a \leq b\}$
$R_I^+ = \{[a, b] \mid a, b \in R^+ \cup \{0\}, a \leq b\}$.

Clearly $Z_I^+ \subseteq Q_I^+ \subseteq R_I^+$. Consider $Z_n^I = \{[0, r] \mid r \in Z_n\}$ is the set of intervals in $Z_n$.

However from the context one can easily follow from which set the intervals are taken.

While working we further refrain and use mainly intervals of the form [0, a] where a $\in Z_n$ or $Z^+ \cup \{0\}$ or $Q^+ \cup \{0\}$ or $R^+ \cup \{0\}$. We add intervals as [[a, b] + [c, d] = [ac, bd]

In case of [0, a] type of intervals [0, a] + [0, b] = [0, a + b] and [0, a]. [0, b] =[0, ab] for a, b in $Z_n$ or $Z^+\{0\}$ or so on. We use only interval of the form [a, b] where a < b for in our collection of intervals we do not accept the degenerate intervals except 0. When we say A = ($a_{ij}$) is an interval matrix the entries $a_{ij}$ are intervals.

For example

$$\begin{bmatrix} [0,5] & [0,3] \\ [0,1] & [0,4] \\ [0,2] & [0,7] \end{bmatrix}$$

is a 3 × 2 interval matrix.

For more about these concepts please refer [52].

**Chapter Two**

# SET INTERVAL LINEAR ALGEBRAS OF TYPE I AND THEIR GENERALIZATIONS

In this chapter we for the first time introduce the new notion of set interval linear algebras of type I and their fuzzy analogue. This chapter has two sections.

## 2.1 Set Interval Linear Algebras of Type I

In this section we define two classes of set interval linear algebras one built using just subsets from Z or Q or R or C or $Z_n$ (n < ∞) and the other built using intervals from $Z^+ \cup \{0\}$ or $Q^+ \cup \{0\}$ or $R^+ \cup \{0\}$ or $Z_n$ discussed in chapter one of this book.

**DEFINITION 2.1.1**: *Let S denote a collection of intervals of the form $\{[x_i, y_i]; y_i, x_i \in Z; 1 \leq i \leq n\}$ (This set S need not be closed under any operation just an arbitrary collection of intervals). Let F be a subset of $Z^+ \cup \{0\}$ If for every $c \in F$ and $s = [x_i, y_i] \in S$, we have $cs = [cx_i, cy_i] \in S$; then we define S to be a set interval integer vector space over the subset F. If the number of distinct elements in S is finite we call S to be a finite set interval integer vector space; if $|S| = \infty$ we say S is an set integer interval vector space of infinite order.*

We will illustrate this situation by some examples.

***Example 2.1.1:*** Let $S = \{[2n, 2m], n < m \mid m, n \in Z^+\} \subseteq Z_I^+$, take $F = \{2, 4, 8, \ldots, 2^{12}\} \subseteq Z$. S is a set integer interval vector space of infinite order over the set F.

***Example 2.1.2:*** Let $S = \{[1, 2], [0, 0], [4, 7], [-2, 3], [4, 21], [-45, 37] [3, 7], [147, 2011]\} \subseteq Z_I$ be a subset of integer intervals. Take $F = \{0, 1\} \subseteq Z$. We see S is a set integer interval vector space over the set F. Clearly S is of finite cardinality and $o(S) = |S| = $ eight.

Now having seen the structure of set integer interval vector space of finite and infinite dimension we now proceed on to define set rational interval vector space.

**DEFINITION 2.1.2:** *Let $S \subseteq Q_I$ or $Q_I^+$ be a subset of intervals of $Q_I$ or $Q_I^+$. Let $F \subseteq Z^+ \cup \{0\}$ or $Q^+ \cup \{0\}$ be a subset of $Z^+$ or $Q^+$; (according as S is from $Q_I^+$ or $Q_I$).*

*If for every $c \in F$ and $s = [x, y] \in S$ we have $c s = [cx, cy]$ and $sc = [xc, yc]$ is in S then we define S to be a set rational interval vector space over F.*

If the number of distinct elements in S is finite we say the cardinality is finite otherwise infinite. Throughout this chapter unless otherwise stated the set F over which the vector spaces are defined is assumed to be subsets of $Z^+ \cup \{0\}$ or $Q^+ \cup \{0\}$ or $R^+ \cup \{0\}$ that is $F \subseteq Z^+ \cup \{0\}$ (or $Q^+ \cup \{0\}$ or $R^+ \cup \{0\}$).

We shall illustrate this situation by some examples.

*Example 2.1.3*: Let
$$S = \left\{ \left[\frac{1}{n}, \frac{1}{n-2}\right] \,\middle|\, 3 \le n \le \infty \right\} \subseteq Q_I^+,$$
be a subset of intervals. Take $F = \{0, 1\} \subseteq Q$. Clearly S is a set rational interval vector space over the set S of infinite cardinality.

*Example 2.1.4:* Let S =

$$\left\{ \left[\frac{7}{2}, 9\right], \left[\frac{5}{3}, 4\right], \left[\frac{17}{5}, 19\right], \left[\frac{22}{7}, 40\right], [0,0], \left[\frac{121}{2}, 149\right], \left[\frac{231}{4}, 504\right] \right\}$$

$\subseteq Q_I^+$ be an interval subset of $Q_I^+$.

It is easily verified S is a set rational vector space over the set $F = \{0, 1\}$ and the cardinality of S is seven.

**DEFINITION 2.1.3:** *Let $S \subseteq R_I^+$ (or $R_I$) be the subset of intervals of reals. Let $F \subseteq Z^+$ or $Q^+$ or $R^+$ (Z or Q or R). If for all $s \in S$ and $c \in F$, sc and cs is in S then we define S to be a set real interval vector space over F. If the number of elements in S is finite we say S is of finite order otherwise S is of infinite order.*

We shall illustrate both the situations by some examples.

*Example 2.1.5*: Let
$$S = \left\{ \left[n\sqrt{5}, n\sqrt{13}\right] \,\middle|\, n \in Z^+ \cup \{0\} \right\} = R_I^+,$$
be the subset of intervals from the positive reals. Take $F = \{1, 2, 3, 4, \ldots, 256\} \subseteq Z^+$. S is an infinite set real interval vector space over the set S.

***Example 2.1.6***: Let
$$S = \left\{ \left[ \frac{n}{\sqrt{7}}, \frac{n}{\sqrt{2}} \right] \middle| 1 \leq n \leq \infty \right\} \subseteq R_I^+$$
be a subset of intervals. Take $F = Z^+ \subseteq R^+$. Clearly S is an infinite set real vector space over F.

***Example 2.1.7***: Let $S = \{[0, 0], [0, 1], [\sqrt{2}, \sqrt{7}], [-\sqrt{3}, 4], [-\sqrt{13}, \sqrt{43}], [5, 8], [\sqrt{17}, 41]\} \subseteq R_I$ subset of real intervals. Take $F = \{0, 1\} \subseteq R$. We see S is a set real interval vector space over the set F. S is of finite dimension or cardinality and the number of elements in S is 7.

Now we will define the concept of set modulo integer interval vector spaces.

**DEFINITION 2.1.4:** *Let $S = \{[x, y] / x, y \in Z_n, x < y\} \subseteq Z_n^I$ be a subset of intervals from the modulo integers. Take $F \subseteq Z_n$ to be proper subset of $Z_n$. If for every $c \in F$ and all $s = [x, y] \in S$, [cx (mod n), cy (mod n)] and [xc (mod n), yc (mod n)] $\in S$ then we say S is a set modulo integer interval vector space over a subset $\{0, 1\} \subseteq Z_n$ ($n < \infty$) any other subset $S_1 \subseteq Z_n$ is choosen provided if $x < y$ implies $sx < sy$ $\forall s \in S_1$ and $\forall [x, y]$ in S.*

We will illustrate this situation by some examples.

***Example 2.1.8***: Let $S = \{[0, 0], [0, 1], [0, 2], [1, 1], [2, 2]\} \subseteq Z_3^I$ be the subset of intervals of $Z_3$. Take $F = Z_3$ it is easily verified that S is a set modulo integer interval vector space over $Z_3 = F$.

***Example 2.1.9***: Let $S = \{[0, 0], [2, 4], [4, 6], [6, 8], [8, 10]\} \subseteq Z_{12}^I$. Take $F = \{0, 1\} \subseteq Z_{12}$. It is easily verified that S is a set modulo integer interval vector space over F.

It is pertinent to mention here that all set modulo integer interval vector spaces are only of finite dimension.

Thus it is convenient to use these structures when our need is just finite.

Now we proceed onto define the notion of set complex interval vector spaces.

**DEFINITION 2.1.5:** *Let $S \subseteq C_I$ subset of intervals of complex numbers. Take F to be a subset of $Z^+ \cup \{0\}$ or $R^+ \cup \{0\}$ or $Q^+ \cup \{0\}$. If for be the every c in F and for every s = [x, y] in S sc, cs $\in$ S then we call S to be a set complex interval vector space over the set F.*

We will illustrate this situation by some examples.

*Example 2.1.10*: Let S = {[2i, 4i + 2], [7, 3i + 13], [0, 0], [14i + 1, 27i + 4]} $\subseteq C_I$ be a subset of intervals from $C_I$. Choose F = {0, 1}; we see S is a set complex interval vector space of cardinality four over the set F = {0, 1}.

*Example 2.1.11*: Let S = {[ni, ni + n] | n $\in Z^+$} $\subseteq C_I$ be a subset of intervals from $C_I$. Choose F = {1, 2, …, 24} we see S is an infinite set complex interval vector space over F.

We now proceed onto describe substructures in these algebraic structures.

**DEFINITION 2.1.6**: *Let $S \subseteq Z_I^+ \subseteq (Z_I)$ be a set integer interval vector space over the set $F \subseteq Z^+$ we say a proper integer interval subset $P \subseteq S$ to be a set integer interval vector subspace of S over F if P itself is a set integer interval vector space over F.*

We will illustrate this by some examples.

*Example 2.1.12*: Let S = {[0, 0], [3, 9], [4, 14], [-5, 17], [13, 19], [41, 53]} $\subseteq Z_I$ be an integer interval subset of $Z_I$. Take F = {0, 1} $\subseteq Z$, S is a set integer interval vector space over F. Take

P = {[0, 0], [41, 53], [–5, 17], [3, 9]} ⊆ S, P is a set integer interval vector subspace of S over the set F.

***Example 2.1.13***: Let S = {[0, 0], [$(2m)^n$, $(2m)^{n+1}$]; $1 \leq n, m \leq \infty$} ⊆ $Z_I^+$; S is a set integer interval vector space over the set F = {0, 2, $2^2$, ..., $2^{40}$} ⊆ $Z^+$. Choose P = {[0, 0], [$(4m)^n$, $(4m)^{n+1}$] | $1 \leq n, m \leq \infty$} ⊆ S; P is a set integer interval vector subspace of S over F.

Now we can as in case of set integer interval vector spaces define for set real (complex, rational, modulo integers) set real interval vector spaces (complex, rational, modulo integer) interval vector subspaces with appropriate simple changes.

We shall however illustrate this situation by some examples.

***Example 2.1.14***: Let S = {[0, 0], [2, 2], [1, 1], [0, 1], [0, 2], [3, 3], [0, 3]} ⊆ $Z_4^I$ be a set modulo integer interval vector space built using $Z_4$. Take F = {0, 1, 2, 3} ⊆ $Z_4$. We see S is a set modulo integer interval vector space over F.

Take P = {[0, 0], [1, 1], [2, 2], [3, 3], [0, 2]} ⊆ S; P is a set modulo integer interval vector subspace of S over F.

***Example 2.1.15***: Let S = {[0, 0], [1, $\sqrt{2}$], [1, $\sqrt{3}$], [$\sqrt{2}$, $\sqrt{3}$], [$\sqrt{17}, \sqrt{23}$]} ⊆ $R_I$ be a set real interval vector space over the set F = {0, 1}. Choose P = {[0, 0], [1, $\sqrt{3}$], [$\sqrt{2}$, $\sqrt{3}$]} ⊆ S; P is a set real interval vector subspace of S over F.

***Example 2.1.16***: Let
$$S = \{[0, 0], \left[\frac{3}{2}, \frac{5}{2}\right], \left[\frac{5}{2}, \frac{7}{2}\right], ..., \left[\frac{43}{2}, \frac{45}{2}\right]\} \subseteq Q_I$$
be a set rational interval vector space over the set F = {0, 1}.

Take

$P = \{[0, 0], \left[\dfrac{5}{2}, \dfrac{7}{2}\right], \left[\dfrac{9}{2}, \dfrac{11}{2}\right], \left[\dfrac{23}{2}, \dfrac{25}{2}\right], \left[\dfrac{35}{2}, \dfrac{37}{2}\right], \left[\dfrac{41}{2}, \dfrac{43}{2}\right]\}$

$\subseteq$ S; it is easily verified P is a set rational interval vector subspace of S over the set F = {0, 1}.

***Example 2.1.17***: Let S = {[n $\sqrt{2}$, n $\sqrt{23}$], [0, 0] | n $\in$ $Z^+$} $\subseteq$ $R_I$ be a set real interval vector space over the set F = {0, 1}. Choose P = {[3n $\sqrt{2}$, 3n $\sqrt{23}$], [0, 0]} $\subseteq$ S; P is a set real interval vector subspace of S over the set F = {0, 1}.

***Example 2.1.18***: Let S = {[mi, (m + 3) + (m + 3)i], [0, 0] | m $\in$ $Z^+$} $\subseteq$ $C_I$ be a set complex interval vector space over the set F = {0, 1}. Choose P = {[5mi, [5(m + 3) + 5(m + 3)i] | m $\in$ $Z^+$} $\subseteq$ S $\subseteq$ $C_I$. P is a set complex interval vector subspace of S over F.

Now we call a set integer (real or complex or rational or modulo integer) interval vector space S to be a simple set integer (real or complex or rational or modulo integer) interval vector space if it has no proper set integer (real or complex or rational or modulo integer) interval vector subspace P; where P $\neq$ [0, 0] or S over F.

We will illustrate by some simple examples the notion of set integer (real or complex or rational or complex or modulo integer) simple vector space.

***Example 2.1.19***: Let S = {[0, 0], [5, 7]} be a set integer interval vector space over the set F = {0, 1}. We see S is a simple set integer interval vector space over F.

***Example 2.1.20***: Let S = {[0, 0], [0, 1], [0, 2], [0, 3], [0, 4], [0, 5], [0, 6]} $\subseteq$ $Z_7^I$, be a set modulo integer interval vector space over the set F = {0, 1, 2, 3, 4, 5, 6}. We see S is a simple set modulo integer interval vector space over F.

***Example 2.1.21***: Let S = {[0, 0], [$5/2, 7/2$]} $\subseteq Q_I$ be a rational interval vector space over the set F = {0, 1}. S is a simple set rational interval vector space over F = {0, 1}.

***Example 2.1.22***: Let S = {[0, 0], [1, 3 + i]} be a set complex interval vector space over the set F = {0, 1}. S is a simple set complex interval vector space over F = {0, 1}.

***Example 2.1.23***: Let S = {[0, 0], [$\sqrt{7}, \sqrt{3}+40$]} be a set real interval vector space over the set F = {0, 1}. Clearly S is a simple set real interval vector space over F.

We now proceed onto define the new notion of subset integer (real or complex or rational or modulo integer) interval vector subspace defined over a subset T $\subseteq$ F of a set integer (real or complex or rational or modulo integer) interval vector space defined over F.

**DEFINITION 2.1.7**: *Let S $\subseteq Z_I$ be a set integer interval vector space defined over the set F $\subseteq Z^+ \cup \{0\}$. Suppose P $\subseteq$ S (P a proper subset of S, P $\neq$ [0, 0] or P $\neq$ S) is a set integer interval vector space over the subset T $\subseteq$ F (T $\neq$ (0) or T $\neq$ P and |T| > 1) then we define P to be a subset integer interval vector subspace of S over the subset T of F. Similar definition can be made in case of set real or complex or rational or modulo integer interval vector spaces with suitable modifications.*

However we will illustrate this situation by some examples.

***Example 2.1.24***: Let S = {[0, 0] [0, 1] [0, 2], …, [0, n] | n < ∞} $\subseteq Z_I$ be a set integer interval vector space over the set F = {0, 1, 2, 3, 4}. Choose P = {[0, 0], [0, 2], [0, 4], …, [0, 2n]} $\subseteq$ S. P is a subset integer interval vector subspace of S over the subset T = {0, 2} $\subseteq$ F.

***Example 2.1.25***: Let S = {[0, 0], [0, 1], [0, 2], [0, 3], [0, 4], [1, 1], [2, 2], [3, 3], [4, 4]} be a set modulo integer interval vector space over the set F = {0, 1, 2, 3, 4} $\subseteq Z_5$. Choose P = {[0, 0], [0, 1], [0, 2], [0, 3], [0, 4]} ⊆ S, P is a subset modulo integer interval vector subspace of S over the subset T = {0, 1} ⊆ F.

We now proceed onto define the new notion of pseudo simple set integer (real or rational or complex modulo integer) simple interval vector spaces.

**DEFINITION 2.1.8:** *Let $S \subseteq Z_I$ (or $Q_I$ or $Z_n^I$ or $R_I$ or $C_I$) be a set integer (rational or modulo integer or real $\subseteq Z^+ \cup \{0\}$ or complex) interval vector space over the subset $F \subseteq Z^+ \cup \{0\}$. Suppose S has no proper subset integer (rational or modulo integer or real) interval vector subspace over a proper subset T of F then we define S to be a pseudo simple set integer (rational or modulo integer or real) interval vector space over F.*

If S is both a simple set interval vector space as well as pseudo simple set interval vector space over F then we define S to be a doubly simple set interval integer (real or rational or modulo integer) vector space.

We will give some illustrations before we proceed onto prove some properties.

*Example 2.1.26*: Let S = {[0, 0], [0, 1], [0, 2], …, [0, 415]} ⊆ $Z_I^+$ be a set integer interval vector space over the set F = {0, 1}. Clearly S is a pseudo simple set integer interval vector space over F. However S is not a simple integer interval vector space as S has several set integer interval vector subspaces over F = {0, 1}. Thus S is not a doubly simple set integer interval vector space over F = {0, 1}.

In view of this we have the following theorem.

**THEOREM 2.1.1**: *Let $S \subseteq Z_I$ (or $Q_I$ or $Z_n^I$ or $R_I$ or $C_I$) be a set integer (rational or modulo integer or real or complex) interval vector space over the set F = {0, 1}. Then S is a pseudo simple set integer (rational or modulo integer or real or complex) interval vector space over the set F = {0, 1}.*

*Proof:* The result follows from the fact that the set $F = \{0, 1\}$ has no proper subset T of order greater than or equal to two. Thus we cannot have any subset integer (rational or modulo integer or real or complex) vector space over $F = \{0, 1\}$. Hence the theorem.

**THEOREM 2.1.2:** *Let $S = \{[0, 0], [x, y]\} \subseteq Z_I$ {or QI or $Z_n^I$ or $R_I$ or $C_I$) be a set integer (rational or modulo integer or real or complex) interval vector space over the set $F = \{0, 1\}$. Then S is a doubly simple set integer (rational or modulo integer or real or complex) interval vector space over the set $F = \{0, 1\}$.*

*Proof:* Obvious from the very definition and the cardinalities of S and F. S is a doubly simple set integer (rational or modulo integer or real or complex) interval vector space over F.

Now we will give an example of a doubly simple set interval integer vector space.

***Example 2.1.27***: Let $S = \{[0, 0], [\sqrt{7}, 3\sqrt{19}]\} \subseteq R_I$ be a set real interval vector space over the set $F = \{0, 1\}$. Clearly S is a doubly simple set real interval vector space over F.

We now proceed onto define the notion of set interval vector space interval linear transformations.

**DEFINITION 2.1.9**: *Let S and T be any two set integer (rational or modulo integer or real or complex) interval vector spaces defined over the same set F. We call a map $T_I : S \rightarrow T$ which maps intervals of S into intervals of T and $T_I(cs) = cT_I(s)$ for all $s \in S$ and $c \in F$ to be a interval linear transformations of S to T.*
   *The collection of such interval linear transformations of S to T is denoted by $IHom_F (S, T)$.*

We will give some illustrations of this definition.

***Example 2.1.28***: Let S = {[0, 0], [0, 2], …, [0, 45]} and T = {[0, 0], [1, 2], [1, 3], …, [1, 45]} be two set integer interval vector spaces defined over the set F = {0, 1}.
Define $T_I : S \to T$ by
$T_I \{[0, 0]\}$ = $\{[0, 0]\}$
$T_I \{[0, n]\}$ = $\{[1, n]\}$     $2 \le n \le 45$.

$T_I$ is an interval linear transformation of S to T.

***Example 2.1.29***: Let
S = {[0, 0], [$\sqrt{2},\sqrt{7}$], [$\sqrt{7},\sqrt{11}$], [$\sqrt{11},\sqrt{43}$], [$\sqrt{43},20\sqrt{53}$]}
and T = {[0, 0], [7, 9], [3, 11], [24, 45], [10, 29]} be two set real interval vector spaces defined over the set F = {0, 1}.
Define $T_I$ ([0, 0]) = [0, 0].
$T_I$ ([$\sqrt{2},\sqrt{7}$])     =     [7, 9]
$T_I$ ([$\sqrt{7},\sqrt{11}$])     =     [3, 11]
$T_I$ ([$\sqrt{11},\sqrt{43}$])     =     [3, 11] and
$T_I$ ([$\sqrt{43},20\sqrt{45}$])     =     [24, 45].
$T_I$ is an interval linear transformation of S to T.

It is important to mention here that S and T can be any type of set interval vector space built using integers or reals or complex or so on but only criteria we need is that both should be defined over the same set F. This is evident from the following example.

As we do not demand any thing from the set map $T_I$ except $T_I$ (cs) = c$T_I$ (s) for every c $\in$ F and s $\in$ S.

As in case of usual vector spaces we say a interval linear transformation is an interval linear operator if S = T in the definition 2.1.9.

Now having seen interval linear transformation $T_I$ we can define kernel of $T_I$ only if [0, 0] $\in$ S otherwise the notion of kernel of $T_I$ remains undefined.

We will illustrate this by some examples.

***Example 2.1.30:*** Let S = {[$2^n$, $2^{n+4}$] | n = 1, 2, …, ∞} and T = {[$4^n$, $4^{n+4}$] | n = 1, 2, …, ∞} be two set integer interval vector

spaces defined over the set $F = \{4, 4^2, 4^3, 4^4, 4^5\}$. Define $T_I : S \to T$ by $T_I [2^n, 2^{n+4}] = [4^n, 4^{n+4}]$, $n = 1, 2, \ldots, \infty$.

It is easily verified that $T_I$ is a interval linear transformation of S to T.

We see the notion of kernel $T_I$ has no meaning as $[0, 0] \notin S$.

Now we proceed onto give one example of a linear interval operator (interval linear operator) on a set interval vector space.

***Example 2.1.31***: Let $S = \{[3^n, 3^{n+3}] \mid n = 1, 2, \ldots, \infty\}$ be a set integer interval vector space over the set $F = \{0, 1\}$. Define $T_I : V \to V$ by $T_I [3^n, 3^{n+3}] \to [3^{2n}, 3^{2n+3}]$, $n = 1, 2, \ldots, \infty$.

It is easily verified that $T_I$ is a interval linear operator on V. Further $T_I$ has kernel.

Next we proceed onto define set interval linear algebras built using integer intervals, real intervals and so on.

**DEFINITION 2.1.10**: *Let $S_1, S_2, \ldots, S_k$ be a collection of subset integer (real, complex, rational or modulo integer) interval vector subspaces of S defined over the subsets $T_1, \ldots, T_k$ of F respectively (that is each $S_i$ is a subset interval vector subspace of S over the subset $T_i$ of F; $i=1, 2, \ldots, k$). If $W = \cap S_i \neq \phi$ and $T = \cap T_i \neq \phi$ then we call W to be a sectional subset interval vector sectional subspace of S over T.*

We will illustrate this situation by an example.

***Example 2.1.32***: Let $S = \{[0, 2n], [0, 6n], [0, 5n], [0, 11n], [0, 14n] / n = 0, 1, 2, \ldots, \infty\}$ be a set integer interval vector space over the set $F = Z^+$.

Take $S_1 = \{[0, 2n] / n = 0, 1, 2, \ldots, \infty\}$ $S_2 = \{[0, 6n] \mid n = 0, 1, 2, \ldots, \infty\}$, $S_3 = \{[0, 5n] / n = 0, 1, 2, \ldots, \infty\}$, $S_4 = \{[0, 14n] / n = 0, 1, 2, \ldots, \infty\}$ and $S_5 = \{[0, 11n] / n = 0, 1, 2, \ldots, \infty\}$ be subset integer interval vector subspaces of S over the subsets $T_1 = 2Z$, $T_2 = 3Z$, $T_3 = 5Z$, $7Z = T_4$ and $T_5 = 11Z$ respectively. Clearly $W = \cap S_i \neq \phi$ and $T = \cap T_i \neq \phi$. Hence W is a sectional subset interval vector sectional subspace of S over $T \subseteq F$.

We have the following interesting theorem the proof of which is left as an exercise for the reader.

**THEOREM 2.1.3:** *Every sectional subset interval sectional vector subspace W of the set interval vector space S over the set F is a subset interval vector subspace of a subset F, but not conversely.*

We can as in case of set vector spaces define the generating interval set of a set interval vector space.

**DEFINITION 2.1.11**: *Let S be a set interval vector space built using interval integers or reals or rationals or complex or modulo integers over the set F. We say a subset of intervals B of S generates S if every interval s of S can be got as cs some $c \in F$ $s_j \neq cs_i$ and $s_i \neq cs_j$ for $s_i \neq s_j$; $s_i, s_j \in B$ and $c \in F$. B is called the generating interval set of S over F.*

We will illustrate this by some simple examples.

*Example 2.1.33*: Let S = {[0, 2n], [0, 3n], [0, 5n], [0, 7n] | n = 0, 1, 2, …, ∞} be a set integer interval vector space over the set F = {0, 1, 2, …, ∞}.
 Take B = {[0, 2], [0, 3], [0, 5], [0, 7]} ⊆ S; B is the generating interval subset of S over F.

*Example 2.1.34*: Let S = {[2n, 3n], [5n, 7n], [11n, 13n], [15n, 29n], [12n, 31n] | n = 0, 1, 2, …, ∞} be a set integer interval vector space over the set F = $Z^+ \cup \{0\}$. Take B = {[2, 3], [5, 7], [11, 13], [15, 29], [12, 31]} ⊆ S, B is the interval generating subset of S over the set F.

We as in case of set vector spaces say a proper interval subset B of a set interval vector space to be linearly independent interval set if x, y ∈ B then x ≠ cy or y ≠ dx ; c, d ∈ F. If the interval set B is not linearly independent then we say B is a linearly dependent interval set.
 We see in the example 2.1.34, B is a linearly independent subset of S. If we take D = {[2, 3], [8, 12], [5, 7], [10, 14]} ⊆ S.

D is not a linearly independent interval set as [8, 12] = 4 [2, 3] and [10, 14] = 2 [5, 7] for 4, 2 $\in Z^+ \cup \{0\}$.

It is left as an exercise for the reader to prove the following theorem.

**THEOREM 2.1.4:** *Let S be a set interval vector space over the set F. Let B $\subseteq$ S be a generating interval set of S over F then B is a linearly independent interval set of S over F. Further if S, P and V be any three set interval vector spaces over the set F (S, P and V may be integer interval or real interval or complex interval or rational interval or modulo integer interval) such that if $T_I$ and $M_I$ be interval linear transformations where*

$$T_I : S \to P$$
*and* $\quad M_I : P \to V.$
*Then*
$$T_I \, o \, M_I : S \to V.$$
*That is* $\quad (T_I \, o \, M_I)(s)$ *(for $s \in S$)*
$$= M_I(T_I(s))$$
$$= M_I(p) \, (p \in P)$$
$$= v; \, v \in V;$$
*is a interval linear transformation for S to V.*

*We can define invertible interval linear transformation of $T_I$ where $T_I: S \to P$ then $T_I^{-1} : P \to S$ and derive related properties.*

It is pertinent to mention that we cannot define for these set interval vector spaces set interval linear functional; as F the set over which S is defined is not an interval set.

Now we proceed onto define the notion of set interval linear algebras using integer intervals or real intervals or rational intervals or complex intervals or modulo integer intervals.

**DEFINITION 2.1.12:** *Let S be a set integer (real or complex or rational or modulo integer) vector space defined over the set F. If S is closed under the operation '+' of interval addition i.e., if $s = [x, y]$ and $s_1 = [a, b]$, $s + s_1 = [c, d]$ is in S for every s, $s_1 \in$ S and $c \, (s + s_1) = cs + cs_1$ for all s, $s_1 \in S$ and $c \in F$ then we call*

S to be a *set integer (real or complex or rational or modulo integer) interval linear algebra over F.*

We will first illustrate this situation by some examples.

***Example 2.1.35***: Let $S = \{[0, 2n] \mid n = 0, 1, 2, \ldots, \infty\} \subseteq Z_I^+$ be a set interval linear algebra over the set $F = \{0, 1\}$. Clearly S is closed under interval addition. For if $x = [0, 2n]$ and $y = [0, 2m]$ are in S then

$$\begin{aligned} x + y &= [0, 2n] + [0, 2m] \\ &= [0 + 0, 2n + 2m] \\ &= [0, 2(n + m)] \in S. \end{aligned}$$

***Example 2.1.36***: Let $S = \{[0, \sqrt{5}\,n] \mid n = \{0, 1, 2, \ldots, \infty\}\} \subseteq R_I^+$ be a set interval linear algebra over the set $F = \{0, 1, 2, \ldots, \infty\}$.

***Example 2.1.37***: Let $S = \{[5n, 9n] \mid n = 0, 1, 2, \ldots, \infty\}$ is a set interval linear algebra over the set $F = \{0, 1\}$. For if $x = [5, 9]$ and $y = [20, 36]$ then $x + y = [5, 9] + [20, 36] = [25, 45] = [5.5, 9.5] \in S$.

Now having seen examples of set interval linear algebras defined using real intervals or integer intervals or rational intervals or modulo integer intervals or complex intervals, now we proceed on to define set real (or complex or integer or rational or modulo integer) interval linear subalgebra.

**DEFINITION 2.1.13:** *Let S be a set interval linear algebra using (integer intervals or real intervals or complex intervals or rational intervals or modulo integer intervals) over the set F. Suppose $P \subseteq S$ is a proper subset of S and P itself is a set interval linear algebra over F then we define P to be a set interval linear subalgebra of S over F.*

We will illustrate this by some examples.

***Example 2.1.38***: Let $S = \{[n(1 + i), n(20 + 20i)] \mid n \in Z^+\}$ be a set complex interval linear algebra over the set $F = \{1, 2, \ldots, \infty\}$. Take $P = \{[4n(1 + i), 4n(20 + 20i)] \mid n \in Z^+\} \subseteq S$; P is a set complex interval linear subalgebra of S over F.

***Example 2.1.39***: Let $S = \{[21n, 43n] \mid n = 0, 1, 2, \ldots, \infty\}$ be a set interval linear algebra over the set $F = Z^+$. Let $P = \{[21 \times 5n, 43 \times 5n] \mid n = 0, 1, 2, \ldots, \infty\} \subseteq S$. P is a set interval linear subalgebra of S over the set $F = Z^+$.

We illustrate this situation by some examples.

***Example 2.1.40***: Let $S = \{[0, \sqrt{7}\, n] / n = 0, 1, 2, \ldots, \infty\}$ be a set real interval linear algebra over the set $F = \{0, 1\}$. Take $P = \{[0, \sqrt{7} \times 5n] / n = 0, 1, 2, \ldots, \infty\} \subseteq S$; P is a set real interval linear subalgebra of S over F.

***Example 2.1.41***: Let $S = \{[n(2 + 3i), n(12 + 17i)] \mid n = 0, 1, 2, \ldots, \infty\}$ be a set complex interval linear algebra over the set $F = \{0, 1\}$. Choose $P = \{[6n(2 + 3i), 6n(12 + 17i)] \mid n = 0, 1, 2, \ldots, \infty\} \subseteq S$, P is a set complex interval linear subalgebra of S over F.

Now we proceed onto define subset interval linear subalgebra built using integer intervals or complex intervals or real intervals or rational intervals or modulo integer intervals.

**DEFINITION 2.1.14**: *Let S be a set integer (real or complex or rational or modulo integer) interval linear algebra over the set F. Let $P \subseteq S$ be a proper subset of S ($P \neq \phi$ and $P \neq S$); if P is a set integer (real or complex or rational or modulo integer) interval linear algebra over a proper subset T of F ($T \neq F$) then we define P to be subset integer (real or complex or rational or modulo integer) interval linear subalgebra of S over the subset T of F.*

We will illustrate this situation by some examples.

***Example 2.1.42***: Let S = {[0, (3 + $\sqrt{17}$ )n] be such that n = 0, 1, 2, …, ∞} be a set real interval linear algebra over the set F = {0, 1, 2, …, n = ∞}. Choose P = {[0, (3 +$\sqrt{17}$ )n] | n = 0, 2, 4, 6, 8, …, ∞; that is n is even} ⊆ S; P is a subset real interval linear subalgebra of S over the subset T = {4n | n = 0, 1, 2, …, ∞} ⊆ F.

***Example 2.1.43***: Let S = {[0, 0], [0, 1], [0, 2], [0, 3], [0, 4], [0, 5]} be a set modulo integer linear algebra over the set F = {0, 1, 2, 3, 4, 5}. Choose P = {[0, 0], [0, 2], [0, 4]} ⊆ S, P is a subset modulo integer $Z_6$ interval linear subalgebra of S over the subset T = {0, 2, 4} ⊆ F.

Now if we have a set interval linear algebra S built over the (real intervals or rational intervals or complex intervals or modulo integer intervals) over the set F and if S has no proper set interval linear subalgebra over F then we define S to be simple set interval linear algebra over F. If S has no subset interval linear subalgebra over any proper subset T of F then we define S to be pseudo a simple set interval linear algebra. If S is both a simple set interval linear algebra and a pseudo simple set interval linear algebra then we define S to be a doubly simple set interval linear algebra.

We will illustrate this by some simple examples.

***Example 2.1.44***: Let S = {[0, 0], [0, 1], [0, 2] [0, 3], [0, 4]} ⊆ $Z_5$I be a set modulo integer 5 interval linear algebra over the set F = {0, 1, 2, 3, 4} then S is a simple set modulo integer 5 interval linear algebra. Infact S is also a pseudo simple set modulo integer 5 interval linear algebra. Thus S is a doubly simple set modulo integer 5 interval linear algebra.

Consequent of this we give a class of doubly simple set interval linear algebras.

**THEOREM 2.1.5:** *Let S = {[0, 0], [0, 1], …, [0, p-1] / p is a prime and intervals are from $Z_p^I$} be a set modulo integer p interval linear algebra over F = {0, 1} then S is a doubly simple* *set modulo integer p interval linear algebra. S is a doubly simple set modulo p integer interval linear algebra.*

The proof is left as an exercise for the reader.

Now as in case of set interval vector spaces we in case of set interval linear algebras define interval linear transformations.

**DEFINITION 2.1.15**: *Let S and M to two set integer (real or complex or rational or modulo integer) interval linear algebras over a set F. Suppose $T_I$ is a map from S to M, $T_I$ is called a interval linear transformation if the following condition holds;*
$$T_I(cs + s_1) = cT_I(s) + T_I(s_1)$$
*for all intervals s, $s_1$ in S and for all c in F.*

It is important to mention that interval linear transformation is defined if and only if both the set linear algebras are defined over the same set F. Further set linear interval transformations of set interval vector spaces are different from set interval linear algebras.

If in the definition 2.1.15, M is replaced by S then we call the set interval linear transformation to be a set interval linear operator on S. As in case of set interval vector spaces we define the notion of generating set linearly independent elements and set linearly dependent elements.

We see in case of set interval linear algebra S over F; a subset of intervals B ⊆ S is said to be a linearly independent interval subset if there is no s ∈ B such that s can be written as

$$s = \sum_i c_i s_i \; ; \; s_i \in B \text{ and } c_i \in F;$$

otherwise we say the set B is a linearly dependent interval subset. We say a linearly independent interval subset B of S to be a generating interval subset if every s ∈ S can be written as

$$s = \sum_i c_i b_i \; ; \; c_i \in F \text{ and } b_i \in B.$$

We will illustrate this situation by some examples.

***Example 2.1.45***: Let $S = \{[0, n\sqrt{2}\,] \mid n = 0, 1, 2, \ldots, \infty\} \subseteq R_I$ be a set real interval linear algebra over the set $F = \{0, 1\}$. Take $B = \{[0, \sqrt{2}\,]\} \subseteq S$, B is the generating interval subset of S. Consider $\{[0, \sqrt{2}\,], [0, 5\sqrt{2}\,]\} = C \subseteq S$, C is a linearly dependent interval of S as $[0, 5\sqrt{2}\,] = [0, \sqrt{2}\,] + [0, \sqrt{2}\,] + [0, \sqrt{2}\,] + [0, \sqrt{2}\,] + [0, \sqrt{2}\,]$.

We call the set interval linear algebra S over F to be finite dimensional if B is a generating interval subset of S over F and the number of elements in B is finite; otherwise we say S is an infinite dimensional set interval linear algebra of over F. The dimension of S given in example 2.1.45 is finite and is one.

Interested reader can construct and study more about the dimension of set interval linear algebras.

Now having seen only class of set interval linear algebras we now proceed onto define another new class of interval linear algebras.

## 2.2 Semigroup Interval Vector Spaces

In this section we proceed on to define a new class of semigroup interval vector spaces and discuss a few of their properties. However every semigroup interval vector space is a set interval vector space and not vice versa.

**DEFINITION 2.2.1**: *Let S be a subset of intervals from $Z_I$ or $R_I$ or $Z_n^I$ or $Q_I$ or $C_I$. F be any additive semigroup with zero. We call S a semigroup interval vector space over the semigroup F if the following conditions hold good.*
   1. *$cs \in S$ for all $c \in F$ and for all $s \in S$.*
   2. *$0s = 0 \in S$ for all $s \in S$ and $0 \in F$.*
   3. *$(c_1 + c_2) s = c_1 s + c_2 s$ for all $c_1, c_2 \in F$ and $s \in S$.*

We will illustrate this situation by some examples.

*Example 2.2.1*: Let $S = \{[0, 2n] \mid n = 0, 1, 2, \ldots, \infty \}$ be a semigroup interval vector space over the semigroup $F = Z^+ \cup \{0\}$ under addition.

*Example 2.2.2*: Let $S = \{[(1 + i)n, n (2 + 2i)n] \mid n = 0, 1, 2, \ldots, \infty\}$ be a semigroup interval vector space over the semigroup $F = 3Z^+ \cup \{0\}$ under addition.

*Example 2.2.3*: Let $S = \{[0, 0], [0, 2], [0, 4], [0, 8], [0, 6], [0, 10], [0, 12], [0, 14], [0, 16], [0, 18]\}$ be a semigroup interval vector space over the semigroup $F = Z_{20}$ (semigroup under addition modulo 20).

*Example 2.2.4*: Let $S = \{[0, 0], [0, 3]\} \subseteq Z_9^I$ be a semigroup interval vector space over the semigroup $F = \{0, 3, 6\}$ addition modulo 9.

Now we proceed on to define semigroup interval vector subspace of S.

**DEFINITION 2.2.2**: *Let S be a semigroup interval vector space over the semigroup F. Suppose $\phi \neq P \subseteq S$ ($P \neq S$ a proper subset S) is a semigroup interval vector space over the semigroup F then we define P to be a semigroup interval vector subspace of S over the semigroup F.*

We will illustrate this situation by some examples.

*Example 2.2.5*: Let $S = \{[0, 0], [0, 1], [0, 2] [0, 3], [0, 4], [0, 5], [0, 6], [0, 7], [0, 8], [0, 9], [0, 10], [0, 11]\}$ be a semigroup interval vector space over the semigroup $F = \{0, 2, 4, 6, 8, 10\} \subseteq Z_{12}$ addition modulo 12.

Take $P = \{[0, 0], [0, 4], [0, 8]\} \subseteq S$, P is a semigroup interval vector subspace of S over the semigroup F.

*Example 2.2.6*: Let $S = \{[0, n \sqrt{17}] \mid n = 0, 1, 2, \ldots, \infty\}$ be a semigroup interval vector space over the semigroup $F = 3Z^+ \cup$

{0} under addition. Take $P = \{[0, 4n \sqrt{17}\,]$ such that $n = 0, 1, 2,$ …, $\infty\} \subseteq S$; P is a semigroup interval vector subspace of S over the semigroup F.

If a semigroup interval vector space S over the semigroup S has no proper semigroup interval vector subspace over F other than $P = \{[0, 0]\}$, then we call S to be a simple semigroup interval vector space over F.
We will illustrate this situation by some examples.

*Example 2.2.7*: Let $S = \{[0, 0], [0, 1], [0, 2], [0, 3], [0, 4]\}$ be the semigroup interval vector space over the semigroup $F = Z_5$ under addition modulo 5. S is a simple semigroup interval vector space over F.

*Example 2.2.8*: Let $S = \{[0, 0], [0, n] / n = 1, 2, …, 22\} \subseteq Z_{23}^I$ be a semigroup interval vector space over the semigroup $F = Z_{23}$ under addition modulo 23. S is a simple semigroup interval vector space over $Z_{23}$.
In view of these examples we have the following theorem which guarantees the existence of a class of simple semigroup interval vector spaces.

**THEOREM 2.2.1:** *Let $S = \{[0, n] / n = 0, 1, 2, …, p-1\} \subseteq Z_p^I$, p a prime. $F = Z_p$ a semigroup under addition modulo p. S is a simple semigroup interval vector space over F.*

*Proof:* Follows from the fact that no proper interval subset P of S ($P \neq [0, 0]$ or $P \neq S$) can be a semigroup interval vector spaces over F. Hence the claim.

Thus we have a class of infinite number of simple semigroup interval vector space over the semigroup F.

*Example 2.2.9*: Let $S = \{[0, 0], [0, n] \mid n \in Z_m$; m a non prime integer $m < \infty\}$ be a semigroup interval vector space over the semigroup $Z_m = F$. We see S is not a simple semigroup interval vector space over $Z_m = F$.

In view of this we have the following theorem.

**THEOREM 2.2.2:** *Let $Z_m = \{0, 1, 2, ..., m – 1\}$; $m = p_1^{\alpha_1}...p_t^{\alpha_t}$ where $p_1, ..., p_t$ are t distinct primes and $\alpha_i \geq 1$, $1 \leq i \leq t$, be the set of integers modulo m. $S = \{[0, n] \mid n \in Z_m\} \subseteq Z_m^I$. S is a semigroup interval vector space over the semigroup $F = Z_m$. Infact S is not a simple semigroup interval vector space and $P_i = \{[0, np_i] / p_i$ a prime such that $p_i^{\alpha_i} / m$ and $p_i^{\alpha_i+1} + /m$; $1 \leq i \leq t$, $n, p_i \in Z_m\} \subseteq S$ are semigroup interval vector subspaces of S over $F = Z_m$.*

The proof is straight forward and left as an exercise for the reader.

We will illustrate the above theorem by some examples.

***Example 2.2.10***: Let $Z_{30} = \{0, 1, 2, ..., 29\}$ be the modulo integer 30 and $30 = 2.3.5$.

$S = \{[0, n] \mid n \in Z_{30}\}$ be a semigroup interval vector space over the semigroup $F = Z_{30}$. Take $P_1 = \{[0, 0] [0, 2], [0, 4], [0, 6], ..., [0, 28]\} = \{[0, 2n] \mid 2, n \in Z_{30}\} \subseteq S$.

It is easily verified $P_1$ is a semigroup interval vector subspace of S over F.

Take $P_2 = \{[0, 3n] \mid 3, n \in Z_{30}\} \subseteq S$, $P_2$ is a semigroup interval vector subspace of S over F.

$P_3 = \{[0, 5n] \mid 5, n \in Z_{30}\} \subseteq S$; $P_3$ is a semigroup interval vector subspace of S.

***Example 2.2.11***: Let $Z_{36} = \{0, 1, 2, ..., 35\}$ modulo 36, integers. $36 = 2^2 . 3^2$. Let $S = \{[0, n] \mid n \in Z_{36}\}$ be a semigroup interval vector space of S over the semigroup $Z_{36} = F$.

Choose $P_1 = \{[0, 2n] \mid n \in Z_{36}\} = \{[0, 0], [0, 2], [0, 4], [0, 6], ..., [0, 34]\} \subseteq Z_{36}^I$; $P_1$ is a semigroup interval vector subspace of S over the semigroup $F = Z_{36}$. $P_2 = \{[0, 4n] \mid n \in Z_{36}\} = \{[0, 0], [0, 4], [0, 8], [0, 12], [0, 16], [0, 20], [0, 24], [0, 28]\} \subseteq S$ is a semigroup interval vector subspace of S over the semigroup F

= $Z_{36}$. $P_3$ = {[0, 3n] / n ∈ $Z_{36}$} ⊆ S is a semigroup interval vector subspace of S over F = $Z_{36}$.

$P_4$ = {[0, 0], [0, 9], [0, 18], [0, 27]} ⊆ S is a semigroup interval vector subspace of S over F = $Z_{36}$.

Now we will proceed onto define the notion of semigroup linearly independent linearly dependent interval subset of a semigroup interval vector space.

**DEFINITION 2.2.3:** *Let S be a semigroup interval vector space over the semigroup F. A set of interval elements B = {$s_1$, $s_2$, ..., $s_n$} of S is a said to be a semigroup linearly independent interval subset if $s_i \neq cs_j$; for all c ∈ F and $s_i$, $s_j$ ∈ B; i ≠ j; 1 ≤ i, j ≤ n.*

*If for some $s_i$ = $cs_j$, c ∈ F; i≠j; $s_i$, $s_j$ ∈ B then we say the semigroup interval subset is linearly dependent or not linearly independent.*

*If B is a semigroup linearly independent interval subset of S and B generates S, the semigroup interval vector space over F; that is if every element s ∈ B can be got as s = $cs_i$, c ∈ F and $s_i$ ∈ S; 1 ≤ i ≤ n.*

We will illustrate this by some examples.

*Example 2.2.12*: Let S = {[0, n] | n = 0, 1, 2, ..., ∞} be a semigroup interval vector space over the semigroup F = $Z^+$ ∪ {0}. Take B = {[0, 1]} ⊆ S, B generates S as a semigroup interval vector space over F.

*Example 2.2.13*: Let S = {[0, n] | n ∈ $Z_{12}$} be a semigroup interval vector space over the semigroup F = {0, 6} ⊆ $Z_{12}$ semigroup under addition modulo 12. Take B = {[0, 1], [0, 2], [0, 3], [0, 4], [0, 5], [0, 7], [0, 8], [0, 9], [0, 10], [0, 11]} ⊆ S, B generates S over F = {0, 6}.

*Example 2.2.14*: Let S = {[0, n] | n ∈ $Z_{24}$} be a semigroup interval vector space over the semigroup F = {0, 2, 4, 6, 8, 10, ..., 22} ⊆ $Z_{24}$, semigroup under addition modulo 24.

Take B = {[0, 2], [0, 0], [0, 4], [0, 8]} ⊆ S; B is a linearly dependent interval subset of S over F.

*Example 2.2.15*: Let S = {[0, n] | n = 1, 2, …, 11} be a semigroup interval vector space over the semigroup F = {0, 2, 4, …, 10} ⊆ $Z_{12}$, semigroup under addition modulo 12. Take B = {[0, 1], [0, 3], [0, 5], [0, 7]} ⊆ S, B is a linearly independent interval subset of S over F but is not a generating interval subset of S over F.

We will now proceed onto define the notion of semigroup interval linear algebra.

**DEFINITION 2.2.4:** *Let S be a semigroup interval vector space over the semigroup F. If S is also an interval semigroup under addition then we define S to be semigroup interval linear algebra over the semigroup F if c ($s_1$ + $s_2$) = $cs_1$ + $cs_2$ for s ∈ S and $c_1$, $c_2$ ∈ F.*

We will illustrate this situation by some examples.

*Example 2.2.16*: Let S = {[0, n] | n ∈ $Z^+$ ∪ {0}} be a semigroup interval linear algebra over the semigroup $Z^+$ ∪ {0} = F. Take P = {[0, 5n] | n ∈ $Z^+$ ∪ {0}} ⊆ S; P is a semigroup interval linear subalgebra of S over the semigroup F = $Z^+$ ∪ {0}.

*Example 2.2.17*: Let S = {[na, (n + 5) a] | a ∈ $Q^+$, n ∈ $Z^+$ ∪ {0}} be a semigroup interval linear algebra over F = $Z^+$ ∪ {0}. Take P {[na, (n + 5)a | a, n ∈ $Z^+$ ∪ {0}} ⊆ S; P is a semigroup interval linear subalgebra of S over F.

*Example 2.2.18*: Let S = {[0, na] | a ∈ $R^+$, n = 0, 1, 2, …, ∞} be a semigroup interval linear algebra over the semigroup F = $Z^+$ ∪ {0}. Consider P = {[0, na] | a ∈ $Q^+$; n = 0,1,2,..., ∞} ⊆ S, P is a semigroup interval linear subalgebra of S over the set F = $Z^+$ ∪ {0}. If the semigroup interval linear algebra S over the set F has no proper semigroup interval linear subalgebras then we define S to be a simple semigroup interval linear algebra.

We will give some examples of simple semigroup linear algebras.

***Example 2.2.19***: Let $S = \{[0, n] \mid n \in Z_7\} \subseteq Z_7^I$ be a semigroup interval linear algebra over the semigroup $F = Z_7$. Clearly S is a simple semigroup interval linear algebra over F.

***Example 2.2.20***: Let $S = \{[0, n] / n \in Z_p$, p any prime$\} \subseteq Z_p^I$ be a semigroup interval linear algebra over the set $F = Z_p$.

It is easily verified S is a simple semigroup interval linear algebra over the set F.

Now we define new concepts of substructures in these new algebraic structures.

**DEFINITION 2.2.5**: *Let S be a semigroup interval linear algebra over the semigroup F. If $P \subseteq S$ (P = {0} or P ≠ S) be a proper subsemigroup of S. If T be a proper subsemigroup of F and P is a semigroup interval linear algebra over the semigroup T then we call P to be a subsemigroup interval linear subalgebra of S over the subsemigroup T of F.*

*If S has no subsemigroup interval linear subalgebras then we define S to be a pseudo simple semigroup interval linear algebra over F.*

We will first illustrate this situation by some simple examples.

***Example 2.2.21***: Let $S = \{[0, n] \mid n \in Z_{29}\}$ be a semigroup interval linear algebra over the semigroup $F = Z_{29}$. S is a pseudo simple semigroup interval linear algebra over F.

***Example 2.2.22***: Let $S = \{[0, 3n] \mid n \in Z_{19}\} \subseteq Z_{19}^I$; be a semigroup interval linear algebra over the semigroup $F = Z_{19}$. It is easily verified that S is a pseudo simple interval linear algebra over F.

Now we define a semigroup interval linear algebra which is both simple and pseudo simple as a doubly simple semigroup interval linear algebra over F.

We will illustrate this situation by some simple examples.

***Example 2.2.23***: Let $S = \{[0, n] \mid n \in Z_5\} \subseteq Z_5^I$ be a semigroup interval linear algebra over the semigroup $F = Z_5$. S is a doubly simple semigroup interval linear algebra of over F.

***Example 2.2.24***: Let $S = \{[0, n] \mid n \in Z_{11}\} \subseteq Z_{11}^I$ be a semigroup interval linear algebra over the semigroup $F = Z_{11}$. S is a doubly simple semigroup interval linear algebra over the semigroup F.

In view of this we give a class of semigroup interval linear algebras which are doubly simple semigroup interval linear algebras.

**THEOREM 2.2.3**: *Let $S = \{[0, n] \mid n \in Z_p, p \text{ a prime}\} \subseteq Z_p^I$ be a semigroup interval linear algebra over the semigroup $Z_p$. S is a doubly simple semigroup interval linear algebra over $Z_p$.*

The proof is left as an exercise to the reader.

**THEOREM 2.2.4:** *Let $S = \{[0, n] \mid n \in Z^+ \cup \{0\}\} \subseteq Z_I$ be a semigroup interval linear algebra over the semigroup $F = Z^+ \cup \{0\}$. S has both subsemigroup interval linear subalgebras and semigroup interval linear subalgebras.*

*Proof:* All $T_p = \{[0, np] \mid p \in Z^+ \cup \{0\}\} \subseteq S$ is an interval semigroup under addition $T_p \subseteq S$ are semigroup interval linear subalgebras of S over the semigroup $F = Z^+ \cup \{0\}$. Consider $T_p \subseteq S$, $T_p$ is also a subsemigroup interval linear subalgebra of S over the subsemigroup $T = pZ^+ \cup \{0\} \subseteq F = Z^+ \cup \{0\}$.

Hence the claim.

**DEFINITION 2.2.6**: *Let R and S be two semigroup interval linear algebras defined over the same semigroup F. Let T be a mapping from R to S such that $T(c\alpha + \beta) = cT(\alpha) + T(\beta)$ for all $c \in F$ and $\alpha, \beta \in R$, then we define T to be a semigroup interval linear transformation from R to S.*

*If R = S we define T to be a semigroup interval linear operator on R.*

We will illustrate this by some simple examples

*Example 2.2.25*: Let $R = \{[0, n] \mid n \in Z^+ \cup \{0\}\}$ and $S = \{[0, n] / n \in Q^+ \cup \{0\}\}$ be two semigroup interval linear algebras over the semigroup $F = Z^+ \cup \{0\}$. The map $T: R \to S$ is defined by $T([0, n]) = [0, n]$, $n \in Z^+ \cup \{0\}$ is a semigroup interval linear transformation.

*Example 2.2.26*: Let $R = \{[n, 5n] \mid n \in Z^+ \cup \{0\}\}$ and $S = \{[n, 5n] \mid n \in R^+ \cup \{0\}\}$ be two semigroup interval linear algebras defined over the semigroup $F = Z^+ \cup \{0\}$. Define $T: R \to S$ by $T\{[n, 5n]\} = [n, 5n]$, for all $[n, 5n] \in R$.

It is easily verified T is a semigroup interval linear transformation of R to S and infact T is an embedding.

We will give an example of a semigroup interval linear operator.

*Example 2.2.27*: Let $S = \{[n, 2n] \mid n \in Z^+ \cup \{0\}\}$ be a semigroup interval linear algebra on the semigroup $F = Z^+ \cup \{0\}$. Define a interval map $T: S \to S$ by $T\{[n, 2n]\} = [2n, 4n]$ for all $[n, 2n] \in S$. T is clearly a semigroup interval linear operator on S.

Now we proceed on to define the notion of semigroup interval linear projection of a semigroup interval linear algebra.

**DEFINITION 2.2.7**: *Let S be a semigroup interval linear algebra over the semigroup F. Let $P \subseteq S$ be a proper semigroup interval* *linear subalgebra of S over the semigroup F. Let T from V to V be a semigroup interval linear operator over F. T is said to be a semigroup interval linear projection on P if T(v) = ω, if ω ∈ P and T(α u + v) = αT(u) + T(v), T(u) and T(v) ∈ P for all α ∈ F and u, v ∈ S.*

We will illustrate this situation by some examples.

***Example 2.2.28***: Let $S = \{[n, 5n] | n \in Q^+ \cup \{0\}\}$ be a semigroup interval linear algebra over the semigroup $F = Z^+ \cup \{0\}$. Take $P = \{[n, 5n] | n \in Z^+ \cup \{0\}\} \subseteq S$; P is a semigroup interval linear algebra over F.
Define T: S → S by

$$T([n, 5n]) = \begin{cases} [n, 5n] & \text{if } n \in Z^+ \\ [0, 0] & \text{if } n \notin Z^+ \end{cases}$$

We see T is a semigroup interval linear projection.

***Example 2.2.29***: Let $S = \{[0, n] | n \in Z_{30}\}$ be a semigroup interval linear algebra over the semigroup $F = Z_{30}$. Take $P = \{[0, n] | n \in \{0, 5, 10, 15, 20, 25\} \subseteq Z_{30}\} \subseteq S$. P is a semigroup interval linear subalgebra of S over F.
   Define η: S → S by η{[0, n]} = [0, 5n]; η is clearly a semigroup interval projection of S on P.

Now we proceed on to define the notion of pseudo semigroup interval linear operator on V.

**DEFINITION 2.2.8**: *Let S be a semigroup interval linear algebra over the semigroup F. Let P ⊆ S be a subsemigroup interval linear subalgebra of S over a subsemigroup R of F. Let T: S → P be a map such that T(αv + u) = T(α) T(v) + T(u) for all u, v ∈ S and T(α) ∈ R and α ∈ F.*
   *We call T to be a pseudo semigroup interval linear algebra operator on S.*

Interested reader is expected to construct examples of pseudo semigroup interval linear operator on S.

**DEFINITION 2.2.9**: *Let S be a semigroup interval vector space over the semigroup F. Let $W_1, W_2, \ldots, W_n$ be semigroup interval vector subspaces of S over F.*

*If $S = \bigcup_i W_i$ and $W_i \cap W_j = \phi$ or $\{0\}$, if $i \neq j$ then we say S is the direct union of the semigroup interval vector subspaces of the semigroup interval vector space S over F.*

We will illustrate this by some examples.

*Example 2.2.30*: Let $S = \{[0, n] | n \in Z_4\}$ be a interval semigroup linear algebra over the semigroup $F = Z_4$. S cannot be written as a union of semigroup interval sublinear algebras over F.

*Example 2.2.31*: Let $S = \{[0, n] | n \in Z_6\}$ be a interval semigroup vector space over $F = \{0, 3\}$. Take $W_1 = \{[0, n] / n \in \{0, 1, 3, 5\} \subseteq Z_6\}$ and $W_2 = \{[0, n] | n \in \{0, 2, 4\} \subseteq Z_6\}$; $W_1$ and $W_2$ are interval semigroup vector subspace of V over $F = \{0, 3\}$. Clearly $V = W_1 \cup W_2$ and $W_1 \cap W_2 = \{0\}$. Thus W is a direct union of semigroup interval vector subspaces of S.

*Example 2.2.32*: Let $G = \{[0, n] | n \in Z_{10}\}$ be a interval semigroup vector space over the semigroup $S = \{0, 5\}$. Let $W_1 = \{[0, n] | n \in \{0, 2, 4, 6, 8\}\} \subseteq G$ and $W_2 = \{[0, n] | n \in \{0, 1, 3, 5, 7, 9\}\} \subseteq G$ be interval semigroup vector subspaces of V over the semigroup $S = \{0, 5\}$. Clearly $V = W_1 \cup W_2$ and $W_1 \cap W_2 = \{0\}$. Thus V is a direct sum of the interval semigroup vector subspaces $W_1$ and $W_2$.

**DEFINITION 2.2.10**: *Let $V = \{[0, n] | n \in Z_n$ or $Z^+ \cup \{0\}$ or $R^+ \cup \{0\}$ or $Q^+ \cup \{0\}\}$ be a interval semigroup linear algebra over the semigroup S. Suppose $W_1, W_2, \ldots, W_m$ be semigroup interval linear subalgebras of V such that $V = W_1 + \ldots + W_m$ and $W_i \cap W_j = \{0\}$ or $\phi$ if $i \neq j$ $\{1 \leq i, j \leq m\}$ then we say V is a direct sum of interval semigroup linear subalgebras of V.*

We will illustrate this situation by some examples.

*Example 2.2.33*: Let

$$V = \left\{ \begin{bmatrix} [0, a_1] & [0, a_2] \\ [0, a_3] & [0, a_4] \\ [0, a_5] & [0, a_6] \end{bmatrix} \middle| \begin{array}{l} a_i \in Z^+ \cup \{0\} \\ 1 \le i \le 6 \end{array} \right\}$$

be an interval semigroup linear algebra over $F = 3Z^+ \cup \{0\}$.
Let

$$W_1 = \left\{ \begin{bmatrix} [0\, a_1] & [0\, a_2] \\ 0 & 0 \\ 0 & 0 \end{bmatrix} \middle| a_1, a_2 \in Z^+ \cup \{0\} \right\}$$

$$W_2 = \left\{ \begin{bmatrix} 0 & 0 \\ [0\, a_1] & 0 \\ [0\, a_2] & 0 \end{bmatrix} \middle| a_1, a_2 \in Z^+ \cup \{0\} \right\}$$

and

$$W_3 = \left\{ \begin{bmatrix} 0 & 0 \\ 0 & [0\, a_1] \\ 0 & [0\, a_2] \end{bmatrix} \middle| a_1, a_2 \in Z^+ \cup \{0\} \right\}$$

be interval semigroup linear subalgebras of V over $F = 3Z^+ \cup \{0\}$. Clearly $V = W_1 + W_2 + W_3$ and $W_i \cap W_j = (0)$; $i \ne j$ $1 \le i, j \le 3$. Thus V is the direct sum of interval linear semigroup subalgebras.

*Example 2.2.34*: Let

$$V = \left\{ \begin{bmatrix} a_1 & a_2 \\ a_3 & a_4 \end{bmatrix} \middle| a_i \in Z_8, 1 \le i \le 4 \right\}$$

be a semigroup interval linear algebra over $F = Z_8$.
Let

$$W_1 = \left\{ \begin{bmatrix} a_1 & 0 \\ 0 & 0 \end{bmatrix} \middle| a_1 \in Z_8 \right\},$$

$$W_2 = \left\{ \begin{bmatrix} 0 & a_2 \\ 0 & 0 \end{bmatrix} \middle| a_2 \in Z_8 \right\},$$

$$W_3 = \left\{ \begin{bmatrix} 0 & 0 \\ a_3 & 0 \end{bmatrix} \middle| a_3 \in Z_8 \right\}$$

and

$$W_4 = \left\{ \begin{bmatrix} 0 & 0 \\ 0 & a_4 \end{bmatrix} \middle| a_4 \in Z_8 \right\}$$

be semigroup interval linear subalgebras of V over $Z_8 = F$. We see $V = W_1 + W_2 + W_3 + W_4$ and $W_i \cap W_j = (0)$; $1 \leq i, j \leq 4$. Thus V is a direct sum of semigroup interval linear subalgebras.

Now we proceed on to define Group interval linear algebras.

**DEFINITION 2.2.11:** *Let V be a set of intervals with zero which is non empty. Let G be a group under addition. We call V to be a group interval vector space over G if the following conditions are true;*
   *(a) For every $v \in V$ and $g \in V$ gv and vg are in V*
   *(b) $0v = 0$ for every $v \in V$ and 0 is the additive identity of G.*

We will illustrate this situation by some examples.

*Example 2.2.35*: Let $V = \{[0, a_i] \mid a_i \in Z_9\}$ be a group interval a vector space over the group $Z_9 = G$ under addition modulo 9.

*Example 2.2.36*: Let $V = \{[0, a_i] \mid a_i \in Z_{25}\}$ be a group interval vector space over the additive group modulo 25.

We see $Z^+ \cup \{0\}$ is not a group likewise $Q^+ \cup \{0\}$, $R^+ \cup \{0\}$ and $C^+ \cup \{0\}$ are not groups under addition.

*Example 2.2.37*: Let

$$V = \left\{ \begin{bmatrix} [0, a_1] & [0, a_2] \\ [0, a_3] & [0, a_4] \end{bmatrix}, ([0, a_1],[0, a_2]), \begin{bmatrix} [0, a_1] \\ [0, a_2] \\ [0, a_3] \\ [0, a_4] \\ [0, a_5] \end{bmatrix} \middle| \begin{array}{l} a_i \in Z_{90} \\ 1 \leq i \leq 5 \end{array} \right\}$$

be a group interval vector space over the group $Z_{90} = G$, under addition modulo 90.

*Example 2.2.38*: Let

$$V = \left\{ \begin{bmatrix} [0, a_1] \\ [0, a_2] \end{bmatrix}, ([0, a_1],[0, a_2][0, a_3]) \middle| \begin{array}{l} a_i \in Z_{14} \\ 1 \leq i \leq 3 \end{array} \right\}$$

be the group interval vector space over the group $Z_{14} = G$ under addition modulo 14.

Now we proceed on to define substructures of group interval vector spaces.

**DEFINITION 2.2.12**: *Let V be a group interval vector space over the group G. Let $P \subseteq V$ be a proper subset of V and is a group interval vector space over G. We define P to be a group interval vector subspace over G.*

We will illustrate this situation by some examples.

*Example 2.2.39*: Let

$$W = \left\{ \begin{bmatrix} 0 & 0 \\ [0, a_1] & 0 \end{bmatrix}, \begin{bmatrix} 0 & 0 \\ 0 & [0, a_2] \end{bmatrix}, \begin{bmatrix} [0, a_1] & [0, a_2] & 0 \\ 0 & 0 & [0, a_3] \end{bmatrix} \middle| \begin{array}{l} a_i \in Z_{15} \\ 1 \leq i \leq 15 \end{array} \right\}$$

be a group interval vector space over the group $G = Z_{15}$. Let

$$P = \left\{ \begin{bmatrix} 0 & 0 \\ [0, a_1] & 0 \end{bmatrix} \middle| a_i \in Z_{15} \right\} \subseteq V.$$

P is a group interval vector subspace of V over the group $G = Z_{15}$.

*Example 2.2.40*: Let $V = \{[0, a_i] | a_i \in Z_{40}\}$ be a group interval vector space over the group $G = Z_{40}$. Take $P = \{[0, a_i] | a_i \in \{0, 2, 4, 6, 8, 10, ..., 38\} \subseteq Z_{40}\} \subseteq V$; P is a group interval vector subspace of V over G.

*Example 2.2.41*: Let

$$V = \left\{ \sum_{i=0}^{10} [0, a_i] x^i \middle| \begin{array}{l} a_i \in Z_7 \\ 0 \leq i \leq 10 \end{array} \right\}$$

be a group interval vector space over the additive group $G = Z_7$. Let

$$W = \left\{ \sum_{i=0}^{5} [0, a_i] x^i \middle| a_i \in Z_7 \right\} \subseteq V$$

be a group interval vector subspace of V over $G = Z_7$.

*Example 2.2.42*: Let

$$V = \left\{ \begin{bmatrix} [0, a_1] & [0, a_2] \\ [0, a_3] & [0, a_4] \\ [0, a_5] & [0, a_6] \end{bmatrix} \middle| a_i \in Z_{16}; 1 \leq i \leq 6 \right\}$$

be a group interval vector space over the group $G = Z_{16}$,

$$W = \left\{ \begin{bmatrix} [0, a_1] & 0 \\ [0, a_3] & 0 \\ [0, a_5] & 0 \end{bmatrix} \middle| a_1, a_2, a_4 \in Z_{16} \right\} \subseteq V;$$

is a group interval vector subspace of V over the group $G = Z_{16}$.

**DEFINITION 2.2.13:** *Let V be a group interval vector space over a group G. We say a proper subset P of V to be a linearly dependent subset of V if for any $p_1, p_2 \in P$ ($p_1 \neq p_2$) $p_1 = a p_2$ or $p_2 = a' p_1$ for some $a, a' \in G$.*

*If for no distinct pair of elements $p_1, p_2 \in P$ we have $a, a_1 \in G$ such that $p_1 = ap_2$ or $p_2 = a_1 p_1$ then we say the set P is a linearly independent set.*

*Example 2.2.43*: Let

$$V = \left\{ \begin{bmatrix} [0,a_1] & 0 \\ [0,a_2] & [0,a_3] \end{bmatrix}, \begin{bmatrix} [0,a_1] & 0 \\ [0,a_2] & 0 \end{bmatrix} \middle| a_i \in Z_{12}; 1 \leq i \leq 3 \right\}$$

be a group interval vector space over the group $G = Z_{12}$.
Consider

$$x = \begin{bmatrix} [0,1] & 0 \\ [0,2] & [0,4] \end{bmatrix}, y = \begin{bmatrix} [0,3] & 0 \\ [0,6] & 0 \end{bmatrix}$$

in V. Clearly x and y are linearly dependent as $3x = y$ for $3 \in G = Z_{12}$.

*Example 2.2.44*: Let

$$V = \left\{ \begin{bmatrix} [0,a_1] & [0,a_2] & [0,a_3] \\ [0,a_4] & [0,a_5] & [0,a_6] \end{bmatrix} \middle| a_i \in Z_{15}; 1 \leq i \leq 6 \right\}$$

be a group interval vector space over the group $G = Z_{15}$. Let

$$x = \begin{bmatrix} [0,1] & [0,2] & [0,3] \\ [0,4] & [0,1] & [0,2] \end{bmatrix} \text{ and } y = \begin{bmatrix} [0,4] & [0,8] & [0,12] \\ [0,1] & [0,4] & [0,8] \end{bmatrix}$$

be elements of V.

We see {x, y} forms a linearly dependent subset of V. For we see x = 4y where $4 \in Z_{15} = G$.

*Example 2.2.45*: Let V be a group interval vector space over a group G. Let H be a proper subgroup of G. If $W \subseteq V$ is such that W is a group interval vector space over the subgroup H of G then we define W to be a subgroup interval vector subspace of V over the subgroup H of G.

If W happens to be both a group interval vector subspace as well a subgroup interval vector subspace then we define W to be duo subgroup interval vector subspace. If V has no subgroup interval vector subspace then we define V to be a simple group interval vector space.

We will first illustrate this situation by some simple examples.

*Example 2.2.46*: Let

$$V = \left\{ \begin{bmatrix} [0,a_1] \\ [0,a_2] \\ [0,a_3] \end{bmatrix} \middle| a_i \in Z_{24}; 1 \le i \le 3 \right\}$$

be a group interval vector space over the group $G = Z_{24}$. Consider

$$W = \left\{ \begin{bmatrix} [0,a_1] \\ [0,a_2] \\ 0 \end{bmatrix} \middle| a_i \in Z_{24} \right\} \subseteq V.$$

It is easy to verify W is a subgroup interval vector subspace of V over the subgroup H = {0, 4, 8, 12, 16, 20} $\subseteq G = Z_{24}$.

It is further verified W is also a group interval vector subspace of V. Thus W is a duo subgroup interval subspace of V.

***Example 2.2.47***: Let $V = \{([0, a_1], [0, a_2], [0, a_3], [0, a_4], [0, a_5], [0, a_6], [0, a_7]) | a_i \in Z_{19}, 1 \leq i \leq 7\}$ be a group interval vector space over the group $G = Z_{19}$. It is easy to verify that V has no subgroup interval subspaces as $G = Z_{19}$ has no subgroups. However V has several group interval vector subspaces. For take $W_1 = \{([0, a_1], [0, a_2], [0, a_3], 0, 0, 0, 0) | a_i \in Z_{19}; 1 \leq i \leq 3\}$ $\subseteq$ V is a group interval vector subspace of V over the group G. $W_2 = \{([0, a], [0, a], \ldots, [0, a])$ where $a \in Z_{19}\} \subseteq$ V is a group interval vector subspace of V over the group G. $W_3 = \{([0, a_1], \ldots, 0 [0, a_7]) | a_1, a_7 \in Z_{19}\} \subseteq$ V is a group interval vector subspace of V.

***Example 2.2.48:*** Let

$$V = \left\{ \begin{bmatrix} [0,a] \\ [0,a] \\ [0,a] \\ [0,a] \\ [0,a] \end{bmatrix} \middle| a \in Z_{13} \right\}$$

be a group interval vector space over the group $G = Z_{13}$. It is easily verified V has no proper group interval vector subspace as well as subgroup interval vector subspace.

We cannot define the notion of pseudo semigroup interval vector subspace. However we can define the notion of pseudo set interval vector subspace of a group interval vector space.

**DEFINITION 2.2.14**: *Let V be a group interval vector space over the group G. S a proper subset of G. Let $W \subseteq V$; if W is a set interval vector subspace of V over the set $S \subseteq G$ then we define W to be a pseudo set interval vector subspace of V over the set S; $S \subseteq G$*

We will illustrate this situation by some examples.

***Example 2.2.49***: Let $V = \{[0, n] | n \in Z_{49}\}$ be a group interval vector space over the group $G = \{Z_{49}\}$. Consider $W = \{[0, 0],$ [0, 7], [0, 14], [0, 21], [0, 28], [0, 35], [0, 42]} ⊆ V; W is a pseudo set interval vector subspace of V over the set S = {0, 1, 7} ⊆ $Z_{49}$.

***Example 2.2.50***: Let V = {[0, n] / n ∈ $Z_{40}$} be a group interval vector space over the group G = $Z_{40}$. W = {[0, 0], [0, 10], [0, 20], [0, 30]} ⊆ V is pseudo set interval vector subspace of V over the set S = {0, 1, 2, 3} ⊆ $Z_{40}$.

Now we proceed onto define group interval linear algebras.

**DEFINITION 2.2.15**: *Let V be a group interval vector space over the group G. If V is a group under addition then we call V to be a group interval linear algebra.*

We will illustrate this by some examples.

***Example 2.2.51***: Let V = {[0, n] | n ∈ $Z_{25}$} be a group interval linear algebra over the group G = $Z_{25}$.

***Example 2.2.52***: Let V = {([0, $a_1$], [0, $a_2$], [0, $a_3$]) / $a_1, a_2, a_3$ ∈ $Z_{18}$} be a group interval linear algebra over the group $Z_{18}$ = G.

***Example 2.2.53***: Let
$$V = \left\{ \begin{bmatrix} [0,a_1] \\ [0,a_2] \\ [0,a_3] \\ [0,a_4] \end{bmatrix} \middle| a_1, a_2, a_3, a_4 \in Z_{143} \right\}$$

be a group interval linear algebra over the group G = $Z_{143}$.

***Example 2.2.54***: Let
$$V = \left\{ \sum_{i=0}^{27} [0, a_i] x^i \middle| \begin{array}{l} a_i \in Z_9 \\ 0 \le i \le 27 \end{array} \right\}$$

be a group interval linear algebra over the group G = $Z_9$.

Now having seen examples of group interval linear algebras we now proceed onto define group interval linear subalgebras.

**DEFINITION 2.2.16**: *Let V be a group interval linear algebra over the group G. Let W ⊆ V (W a proper subset of V), if W itself is a group interval linear algebra over the group G then we define W to be a group interval linear subalgebra of V over the group G.*

We will illustrate this situation by some examples.

***Example 2.2.55***: Let $V = \{[0, a] \mid a \in Z_{144}\}$ be a group interval linear algebra over the group $G = Z_{144}$. Consider $W = \{[0, a] / a \in \{2Z_{144}\}\} \subseteq V$; W is a group interval linear subalgebra of V over the group $G = Z_{144}$.

***Example 2.2.56***: Let $V = \{([0, a_1], [0, a_2], [0, a_3], [0, a_4], [0, a_5]) \mid a_1, a_2, a_3, a_4, a_5 \in Z_{48}$ be a group interval linear algebra over the group $G = Z_{48}$. Consider $W = \{([0, a_1], 0, 0, 0, [0, a_5]) \mid a_1, a_5 \in Z_{48}\} \subseteq V$; W is a group interval linear subalgebra of V over the group $G = Z_{48}$.

Now we proceed onto define the notion of direct sum of group interval linear algebras.

**DEFINITION 2.2.17:** *Let V be a group interval linear algebra over the group G. Let $W_1, W_2, …, W_n$ be a group interval linear subalgebras of V over the group G. We say V is a direct sum of the group interval linear subalgebras $W_1, W_2, …, W_n$ if*
    *(a) $V = W_1 + … + W_n$*
    *(b) $W_i \cap W_j = \{0\}$ if $i \neq j$; $1 \leq i, j \leq n$.*

We will illustrate this situation by some simple examples.

***Example 2.2.57***: Let V = {all 3 × 3 interval matrices with entries from $Z_{48}$} be a group interval linear algebra over the group $G = Z_{48}$.

Let

$$W_1 = \left\{ \begin{pmatrix} [0,a_1] & [0,a_2] & 0 \\ 0 & 0 & [0,a_3] \\ 0 & 0 & 0 \end{pmatrix} \middle| a_1, a_2, a_3 \in Z_{48} \right\},$$

$$W_2 = \left\{ \begin{pmatrix} 0 & 0 & [0,a_1] \\ [0,a_2] & 0 & 0 \\ 0 & [0,a_3] & 0 \end{pmatrix} \middle| a_1, a_2, a_3 \in Z_{48} \right\},$$

$$W_3 = \left\{ \begin{pmatrix} 0 & 0 & 0 \\ 0 & [0,a_1] & 0 \\ 0 & 0 & 0 \end{pmatrix} \middle| a_1 \in Z_{48} \right\}$$

and

$$W_4 = \left\{ \begin{pmatrix} 0 & 0 & 0 \\ 0 & 0 & 0 \\ [0,a_1] & 0 & [0,a_2] \end{pmatrix} \middle| a_1, a_2 \in Z_{48} \right\}$$

be group interval linear subalgebras of V over the group $G = Z_{48}$. Clearly $V = W_1 + W_2 + W_3 + W_4$ and

$$W_i \cap W_j = \begin{pmatrix} 0 & 0 & 0 \\ 0 & 0 & 0 \\ 0 & 0 & 0 \end{pmatrix}$$

if $i \neq j$; $1 \leq i, j \leq n$.

Thus V is the direct sum of group interval linear subalgebras $W_1$, $W_2$, $W_3$ and $W_4$.

*Example 2.2.58*: Let V = {Collection of all 4 × 2 interval matrices with entries from $Z_7$} be a group interval linear algebra over the group $G = Z_7$.
Choose

$$W_1 = \left\{ \begin{bmatrix} [0,a_1] & 0 & 0 & 0 \\ 0 & [0,a_2] & 0 & 0 \end{bmatrix} \middle| a_1, a_2 \in Z_7 \right\},$$

$$W_2 = \left\{ \begin{bmatrix} 0 & [0,a_1] & 0 & [0,a_2] \\ [0,a_3] & 0 & 0 & 0 \end{bmatrix} \middle| a_1, a_2 \in Z_7 \right\}$$

$$W_3 = \left\{ \begin{bmatrix} 0 & 0 & [0,a_1] & 0 \\ 0 & 0 & 0 & [0,a_2] \end{bmatrix} \middle| a_1, a_2 \in Z_7 \right\}$$

and

$$W_4 = \left\{ \begin{bmatrix} 0 & 0 & 0 & 0 \\ 0 & 0 & [0,a_1] & 0 \end{bmatrix} \middle| a_1 \in Z_7 \right\}$$

be group interval linear subalgebras of V over the group G. We see $V = W_1 + W_2 + W_3 + W_4$ and

$$W_i \cap W_j = \begin{pmatrix} 0 & 0 & 0 & 0 \\ 0 & 0 & 0 & 0 \end{pmatrix}; \ 1 \le i, j \le 4.$$

Thus V is a direct sum of group interval linear subalgebras.

Let

$$P_1 = \left\{ \begin{bmatrix} [0,a_1] & 0 & [0,a_2] & 0 \\ 0 & [0,a_3] & 0 & [0,a_4] \end{bmatrix} \middle| \begin{array}{l} a_i \in Z_7 \\ 1 \le i \le 4 \end{array} \right\},$$

$$P_2 = \left\{ \begin{bmatrix} 0 & 0 & [0,a_1] & [0,a_2] \\ 0 & [0,a_3] & 0 & [0,a_4] \end{bmatrix} \middle| \begin{array}{l} a_i \in Z_7 \\ 1 \le i \le 4 \end{array} \right\}$$

$$P_3 = \left\{ \begin{bmatrix} [0,a_1] & [0,a_2] & 0 & 0 \\ 0 & 0 & [0,a_3] & [0,a_4] \end{bmatrix} \middle| \begin{array}{l} a_i \in Z_7 \\ 1 \le i \le 4 \end{array} \right\}$$

and

$$P_4 = \left\{ \begin{bmatrix} 0 & 0 & [0,a_1] & [0,a_2] \\ [0,a_3] & [0,a_4] & 0 & 0 \end{bmatrix} \middle| \begin{array}{l} a_i \in Z_7 \\ 1 \le i \le 4 \end{array} \right\}$$

be group interval linear subalgebras of V over the group $G = Z_7$. We see $P_1 + P_2 + P_3 + P_4 = V$ but

$$P_i \cap P_j \neq \begin{pmatrix} 0 & 0 & 0 & 0 \\ 0 & 0 & 0 & 0 \end{pmatrix}$$

if $i \neq j$ ; $1 \leq i; j \leq 4$. Thus any collection of group interval linear subalgebras may not in general give a direct sum of V.

In view of this we have the following interesting definition.

**DEFINITION 2.2.18**: *Let V be a group interval linear algebra over the group G. Let $W_1$, $W_2$, ..., $W_n$ be n distinct group interval linear subalgebras of V over the group G.*

*We say V is a pseudo direct sum if*

(a) $V = W_1 + ... + W_n$
(b) $W_i \cap W_j \neq \{0\}$ *even if* $i \neq j$
(c) *We need $W_i$'s to be distinct that is $W_i \cap W_j \neq W_i$ or $W_i \cap W_j = W_j$ even if $i \neq j$ i.e., $W_i \cap W_j = W_p$ then $p \notin \{1, 2, ..., n\}$ that is $W_p$ does not belong to the collection of group interval linear subalgebras of V.*

We will illustrate this situation by some examples.

*Example 2.2.59*: Let V = {Collection of all $5 \times 2$ interval matrices with entries from $Z_{11}$} be the group interval linear algebra over the group $G = Z_{11}$.

Consider

$$W_1 = \left\{ \begin{bmatrix} [0,a_1] & 0 \\ 0 & [0,a_2] \\ [0,a_3] & 0 \\ 0 & [0,a_4] \\ [0,a_5] & 0 \end{bmatrix} \middle| \begin{array}{l} a_i \in Z_{11} \\ 1 \leq i \leq 5 \end{array} \right\},$$

$$W_2 = \left\{ \begin{bmatrix} [0,a_1] & [0,a_6] \\ [0,a_2] & 0 \\ [0,a_3] & 0 \\ [0,a_4] & 0 \\ [0,a_5] & 0 \end{bmatrix} \middle| \begin{array}{l} a_i \in Z_{11} \\ 1 \le i \le 6 \end{array} \right\},$$

$$W_3 = \left\{ \begin{bmatrix} 0 & [0,a_2] \\ 0 & [0,a_3] \\ 0 & [0,a_4] \\ 0 & [0,a_5] \\ [0,a_1] & [0,a_6] \end{bmatrix} \middle| \begin{array}{l} a_i \in Z_{11} \\ 1 \le i \le 6 \end{array} \right\}$$

and

$$W_4 = \left\{ \begin{bmatrix} 0 & [0,a_1] \\ [0,a_2] & 0 \\ 0 & [0,a_3] \\ [0,a_4] & 0 \\ 0 & [0,a_5] \end{bmatrix} \middle| \begin{array}{l} a_i \in Z_{11} \\ 1 \le i \le 5 \end{array} \right\}$$

be group interval linear subalgebras of V over $G = Z_{11}$. We see $V = W_1 + W_2 + W_3 + W_4$ and $W_i \cap W_j \ne 0$. If $i \ne j$. Further $W_1$, $W_2$, $W_3$ and $W_4$ are all distinct. Thus V is a pseudo direct sum of $W_1$, $W_2$, $W_3$ and $W_4$.

*Example 2.2.60*: Let

  V = {all 4 × 4 interval matrices with entries from $Z_3$}

be the group interval linear algebra over group $G = Z_3$.
  Consider

$$W_1 = \left\{ \begin{bmatrix} [0,a_1] & [0,a_2] & 0 & 0 \\ 0 & 0 & 0 & 0 \\ [0,a_3] & [0,a_4] & [0,a_5] & 0 \\ 0 & [0,a_6] & 0 & [0,a_7] \end{bmatrix} \,\bigg|\, a_1,a_2,\ldots,a_7 \in Z_3 \right\},$$

$$W_2 = \left\{ \begin{bmatrix} 0 & [0,a_2] & 0 & [0,a_3] \\ 0 & [0,a_1] & [0,a_6] & 0 \\ 0 & 0 & [0,a_4] & [0,a_5] \\ 0 & 0 & 0 & [0,a_7] \end{bmatrix} \,\bigg|\, a_i \in Z_3; 1 \le i \le 7 \right\},$$

$$W_3 = \left\{ \begin{bmatrix} [0,a_1] & [0,a_2] & 0 & 0 \\ 0 & 0 & 0 & 0 \\ [0,a_3] & 0 & 0 & [0,a_5] \\ [0,a_6] & 0 & 0 & [0,a_4] \end{bmatrix} \,\bigg|\, a_i \in Z_3; 1 \le i \le 6 \right\}$$

and

$$W_4 = \left\{ \begin{bmatrix} [0,a_1] & [0,a_2] & 0 & 0 \\ [0,a_3] & 0 & [0,a_6] & [0,a_4] \\ 0 & 0 & 0 & 0 \\ [0,a_8] & 0 & [0,a_5] & [0,a_7] \end{bmatrix} \,\bigg|\, a_i \in Z_3; 1 \le i \le 8 \right\}$$

be group interval linear subalgebras of V over the group $Z_3$.

We see $W_i \cap W_j \ne (0)$ if $i \ne j$ $1 \le i, j \le 4$

(a) $V = W_1 + W_2 + W_3 + W_4$
(b) $W_i \cap W_j \ne (0)$ if $i \ne j$; $1 \le i, j \le 4$.
(c) $W_1, W_2, W_3$ and $W_4$ are all distinct group interval linear subalgebras of V over G.

Thus V is a pseudo direct sum of group interval linear subalgebras $W_1, W_2, W_3$ and $W_4$.

We now define linear independence in group interval linear algebras.

**DEFINITION 2.2.19**: *Let V be a group interval linear algebra over the group G. Let $X \subset V$ be a proper subset of V, we say X is a linearly independent subset of V if $X = \{x_1, ..., x_n\}$ (where $x_i = [0, a_i]$, $1 \leq i \leq n$) and for some $n_i \in G$; $1 \leq i \leq n$; $\alpha_1 x_1 + \alpha_2 x_2 + ... + \alpha_n x_n = 0$ if and only if each $\alpha_i = 0$.*

*A linearly independent subset X of V is said to generate V if every element of $v \in V$ can be represented as*

$$v = \sum_{i=1}^{n} \alpha_i x_i ; \alpha_i \in G; 1 \leq i \leq n.$$

We will illustrate this situation by some examples.

***Example 2.2.61***: Let $V = \{([0, a_1], [0, a_2], [0, a_3], [0, a_4]) \mid a_i \in Z_5, 1 \leq i \leq 4\}$ be a group interval linear algebra over the group $G = Z_5$. Consider $X = \{x_1 = ([0, 1], 0, 0, 0), x_2 = (0, [0, 1], 0, 0), x_3 = (0, 0, [0, 1], 0)$ and $x_4 = (0, 0, 0, [0, 1])\} \subseteq V$. X is a linearly independent set and generates V over G so X is a basis of V over G.

***Example 2.2.62***: Let V = {set all $4 \times 2$ interval matrices with entries from $Z_{12}$} be a group interval linear algebra over the group G.

Consider

$$X = \left\{ \begin{bmatrix} [0,1] & 0 \\ 0 & 0 \\ 0 & 0 \\ 0 & 0 \end{bmatrix}, \begin{bmatrix} 0 & [0,1] \\ 0 & 0 \\ 0 & 0 \\ 0 & 0 \end{bmatrix}, \begin{bmatrix} 0 & 0 \\ [0,1] & 0 \\ 0 & 0 \\ 0 & 0 \end{bmatrix}, \right.$$

$$\left.\begin{bmatrix} 0 & 0 \\ 0 & [0,1] \\ 0 & 0 \\ 0 & 0 \end{bmatrix}, \begin{bmatrix} 0 & 0 \\ 0 & 0 \\ [0,1] & 0 \\ 0 & 0 \end{bmatrix},\right.$$

$$\left.\begin{bmatrix} 0 & 0 \\ 0 & 0 \\ 0 & [0,1] \\ 0 & 0 \end{bmatrix}, \begin{bmatrix} 0 & 0 \\ 0 & 0 \\ 0 & 0 \\ [0,1] & 0 \end{bmatrix}, \begin{bmatrix} 0 & 0 \\ 0 & 0 \\ 0 & 0 \\ 0 & [0,1] \end{bmatrix}\right\} \subseteq V;$$

X is a linearly independent set and generates V; hence X is a basis of V.

Here also we cannot define the notion of pseudo semigroup interval linear subalgebras of a group interval linear algebra.

However we can define the notion of pseudo group interval vector subspace of a group interval linear algebra.

**DEFINITION 2.2.20:** *Let V be a group interval linear algebra over the group G. If P is just a subset of V and is not a closed structure but is a group interval vector space over the group G, then we call P to be a pseudo group interval vector subspace of V.*

We will illustrate this situation by an example.

*Example 2.2.63*: Let

$$V = \left\{ \begin{bmatrix} [0,a_1] & [0,a_2] \\ [0,a_3] & [0,a_4] \\ [0,a_5] & [0,a_6] \end{bmatrix} \middle| a_i \in Z_7; 1 \le i \le 6 \right\}$$

be a group interval linear algebra over the group $Z_7$.

Consider

$$W = \left\{ \begin{bmatrix} [0,a_1] & 0 \\ [0,a_2] & 0 \\ [0,a_3] & [0,a_4] \end{bmatrix}, \begin{bmatrix} 0 & [0,a_1] \\ 0 & [0,a_2] \\ 0 & 0 \end{bmatrix}, \begin{bmatrix} [0,a_1] & [0,a_2] \\ 0 & 0 \\ [0,a_3] & 0 \end{bmatrix} \right\} \subseteq V;$$

W is a pseudo group interval vector subspace of V over the group G.

Now we will define in the next chapter the notion of fuzzy interval linear algebras.

**Chapter Three**

# SET FUZZY INTERVAL ALGEBRAS AND THEIR PROPERTIES

In this chapter we introduce the notion of set fuzzy interval linear algebras, semigroup fuzzy interval linear algebras and group fuzzy interval linear algebras and study their properties. This chapter has two sections. First section introduces set fuzzy set vector spaces and discusses their properties. Section two introduces the notion of set fuzzy interval vector spaces of type II and analyses their properties.

**3.1 Set Fuzzy Interval Vector Spaces and Their Properties**

In this section we introduce the notion of set fuzzy interval vector spaces and give a few properties associated with them.

**DEFINITION 3.1.1**: *A fuzzy vector space (V, $\eta$) or $\eta$ V or V$\eta$ is an ordinary vector space V defined over the field F with a map $\eta : V \to [0, 1]$ satisfying the following conditions.*

   *(a)   $\eta(a + b) \geq \min \{\eta(a), \eta(b)\}$*
   *(b)   $\eta(-a) = \eta(a)$*
   *(c)   $\eta(0) = 1$*
   *(d)   $\eta(ra) \geq \eta(a)$*

   *for all a, b $\in$ V and r $\in$ F, where F is the field. V$\eta$ or V$\eta$ or $\eta$V will denote the fuzzy vector space.*

For more about these notions refer [53].

**DEFINITION 3.1.2**: *Let V be a set vector space over the set S. We say V with the map $\eta$ is a fuzzy set vector space or set fuzzy vector space if $\eta: V \to [0, 1]$ and $\eta(ra) \geq \eta(a)$ for all $a \in V$ and $r \in S$. We call V$\eta$ or $\eta$V or V$\eta$ to be the fuzzy set vector space over the set S.*

For more about these notions please refer [52].

Likewise we define a set fuzzy linear algebra (or fuzzy set linear algebra) (V, η) or Vη or ηV to be an ordinary set linear algebra V with a map η : V → [0, 1] such that η(a + b) > min (η(a), η(b)) for a, b ∈ V.

*Notation:* We say an interval [0, a] to be a fuzzy interval if $0 \leq a \leq 1$. $[0, 0] = (0)$ and $[0, 1]$ is the fuzzy set. We include both in the fuzzy interval. $[0, ½]$, $[0, 0.3]$ $[0, \frac{1}{\sqrt{2}}]$, $[0, 0.0031]$ etc are fuzzy intervals.

We will denote the collection of all fuzzy intervals by $I [0, 1] = \{[0, a] \mid 0 \leq a \leq 1\}$. Clearly the cardinality of $I [0, 1]$ is infinite.

Now we proceed onto define fuzzy set interval vector space or set fuzzy interval vector space over the set S.

**DEFINITION 3.1.3**: *Let V be a set interval vector space over the set S. We say V with the map η is a fuzzy set interval set vector space or set fuzzy interval vector space if $I_\eta : V \to I[0, 1]$ and $I_\eta (r[0, a]) > I_\eta ([0, a])$ for all $[0, a] \in V$ and $r \in S$. We call $VI_\eta$ or $I_\eta V$ to be the fuzzy set interval vector space over the set S.*

We will illustrate this situation by some examples.

*Example 3.1.1*: Let $V = \{[0, a_1], [0, a_2], [0, a_3], [0, a_4], [0, a_5])] \mid a_i \in Z_5, 1 \leq i \leq 5\}$ be a set interval vector space over the set $S = \{0, 1, 2, 3\}$. $I\eta: V \to I[0, 1]$ is defined as follows.

$$I_\eta ([0, a_i]) = \begin{cases} [0, \dfrac{1}{a_i}] & \text{if } a_i \neq 0 \\ [0,1] & \text{if } a_i = 0. \end{cases}$$

$V_{I\eta}$ is a set fuzzy interval vector space.

*Example 3.1.2*: Let $V = \{[0, a_i] \mid a_i \in Z^+ \cup \{0\}\}$ be a set interval vector space over the set $S = \{0, 1, 2, 3, 4, 5, 8\}$. Define $I_\eta : V \to I[0, 1]$ as follows:

$$I_\eta [0, a_i] = \begin{cases} [0, \dfrac{1}{a_i}] & \text{if } a_i \neq 0 \\ [0,1] & \text{if } a_i = 0. \end{cases}$$

$I_\eta V$ is a fuzzy set interval vector space.

*Example 3.1.3:* Let

$$V = \left\{ \begin{bmatrix} [0, a_1] \\ [0, a_2] \\ \vdots \\ [0, a_8] \end{bmatrix} \middle| a_i \in Z_{15}; 1 \leq i \leq 9 \right\}$$

be a set interval vector space over the set $S = \{0, 1, 3, 5, 7\} \subseteq Z_{15}$. Define $I\eta : V \to I[0, 1]$ by

$$I\eta\,[0, a_i] = \begin{cases} [0, \dfrac{1}{a_i}] & \text{if } a_i \neq 0 \\ [0,1] & \text{if } a_i = 0. \end{cases} ;$$

$I\eta V$ is a set fuzzy interval vector space.

**DEFINITION 3.1.4**: *Let V be a set interval linear algebra over the set S. A set fuzzy interval linear algebra (or fuzzy set interval linear algebra) (V, ηI) or VηI is a map ηI: V → I [0, 1] such that ηI(a + b) ≥ min( ηI(a), ηI(b)) for every a, b ∈ V.*

We will illustrate this situation by some examples.

*Example 3.1.4*: Let

$$V = \left\{ \sum_{i=0}^{\infty} [0, a_i] x^i \,\middle|\, a_i \in Z^+ \cup \{0\} \right\}$$

be a set interval linear algebra over the set S = {0, 1, 2, 5, 7, 13, 16}.

Define $\eta I : V \to I\,[0, 1]$ as

$$\eta I\,(p(x) = \sum_{i=0}^{n}[0, a_i]x^i\,)$$

$$= \begin{cases} [0, \dfrac{1}{\deg p(x)}] & \text{if } p(x) \text{ is not a constant} \\ [0,1] & \text{if } p(x) \text{ is a constant} \end{cases}$$

VηI is a set fuzzy interval linear algebra.

*Example 3.1.5*: Let V = {([0, $a_1$], [0, $a_2$], [0, $a_3$], [0, $a_4$]) | $a_i \in Z_{18}$; $1 \leq i \leq 4$} be a set interval linear algebra over the set S = {0, 4, 5, 9} $\subseteq Z_{18}$. Define ηI: V → I [0, 1] by

$$\eta I([0, a_i]) = \begin{cases} [0, \dfrac{1}{a_i}] & \text{if } a_i \neq 0 \\ [0,1] & \text{if } a_i = 0. \end{cases}$$

$V_{\eta I}$ is a set fuzzy interval linear algebra.

***Example 3.1.6***: Let

$$V = \left\{ \begin{pmatrix} [0,a_1] & [0,a_3] \\ [0,a_3] & [0,a_4] \end{pmatrix} \right.$$

where $a_i \in Z^+ \cup \{0\}\}$ be a set interval linear algebra over the set $S = \{3Z^+, 2Z^+, 0\} \subseteq Z^+ \cup \{0\}$.

Define $\eta I : V \to I[0, 1]$ by

$$\eta I \begin{pmatrix} [0,a_1] & [0,a_3] \\ [0,a_3] & [0,a_4] \end{pmatrix} = \begin{cases} [0,\dfrac{1}{a_1}] & \text{if } a_1 \neq 0 \\[6pt] [0,\dfrac{1}{a_2}] & \text{if } a_2 \neq 0 \text{ and } a_1 = 0 \\[6pt] [0,\dfrac{1}{a_3}] & \text{if } a_3 \neq 0 \text{ and } a_1 = 0 = a_2 \\[6pt] [0,\dfrac{1}{a_4}] & \text{if } a_4 \neq 0 \text{ and } a_1 = a_2 = a_3 = 0 \\[6pt] [0,1] & \text{if } a_1 = a_2 = a_3 = a_4 = 0 \end{cases}$$

$V_{\eta I}$ is a fuzzy set interval linear algebra.

Now we proceed onto define set fuzzy interval substructures.

**DEFINITION 3.1.5**: *Let V be a set interval vector space over the set S. Let W be a set interval vector subspace of V over the set S.*

*The map $I\eta : W \to I[0, 1]$ such that $W\eta I$ is a set fuzzy interval vector space, is called the set fuzzy interval vector subspace of V and is denoted by $I_{\eta w}$ or $\eta_{wI}$.*

We will illustrate this situation by examples.

***Example 3.1.7***: Let $V = \{([0, a_1], [0, a_2], \ldots, [0, a_{12}]) \mid a_i \in Z_{18}; 1 \leq i \leq 12\}$ be a set interval vector space over the set $S = \{0, 2, 4, 16\} \subseteq Z_{18}$. Let $W = \{([0, a_1], [0, a_2], \ldots, [0, a_{12}])\ a_i \in 2Z_{18};$

$\{0, 2, 4, 6, 8, 10, 12, 14, 16\} \subseteq Z_{18}$; $1 \le i \le 12\}$ be a set interval vector subspace of V over S.

Define $\eta I : W \to I[0, 1]$ by

$$\eta I ([0, a_1], [0, a_2], \ldots, [0, a_{12}]) = \begin{cases} [0, \frac{1}{12}] & \text{if no } a_i = 0 \\ [0, \frac{1}{10}] & \text{if some } a_i = 0 \\ [0, 1] & \text{if all } a_i = 0 \end{cases}$$

(W, $\eta I$) is the set fuzzy interval vector subsubspace of V.

*Note:* It is important and interesting to note that W$\eta$I need not be extendable to V$\eta$I in general.

***Example 3.1.8***: Let $V = \{([0, a_1], [0, a_2], \ldots, [0, a_8])\ |\ a_i \in Z^+ \cup \{0\}; 1 \le i \le 8\}$ be a set interval vector space over the set $S = \{0, 5, 12, 13, 90, 184, 249, 1000\} \subseteq Z^+ \cup \{0\}$. Choose $W = \{([0, a_1], [0, a_2], \ldots, [0, a_8]\ |\ a_i \in 5Z^+ \cup \{0\}\} \subseteq V$ be a set interval vector subspace of V over the set S.
Define $\eta I : W \to I [0, 1]$ by

$$\eta I (x) = \begin{cases} [0, \dfrac{1}{\sum_{i=1}^{} a_i}] & \text{if } \sum a_i \ne 0 \\ [0, 1] & \text{if } \sum a_i = 0 \end{cases}$$

(W, $\eta I$) is a set fuzzy interval vector subspace.

***Example 3.1.9***: Let

$$V = \{[0, a_i], \begin{bmatrix} [0, a_1] \\ [0, a_2] \\ [0, a_3] \end{bmatrix}, ([0, a_1], [0, a_2], [0, a_3], [0, a_4])$$

where $a_i \in Z_{24}\}$ be a set interval vector space over the set $S = \{0, 1, 2, 5, 6, 20, 21\} \subseteq Z_{24}$. $W = \{[0, a_i]\ |\ a_i \in Z_{24}\} \subseteq V$ be a set interval vector subspace of V over $Z_{24}$. Define $\eta I : W \to I [0, 1]$ as follows.

$$\eta I ([0, a_i]) = \begin{cases} [0, \dfrac{1}{a_i}] & \text{if } a_i \neq 0 \\ [0,1] & \text{if } a_i = 0 \end{cases}$$

(W, $\eta I$) is a set fuzzy interval subspace of V. Clearly $\eta I$ cannot be extended to whole of V.
Suppose

$$T = \left\{ \begin{bmatrix} [0, a_1] \\ [0, a_2] \\ [0, a_3] \end{bmatrix} \,\middle|\, \begin{array}{l} a_i \in Z_{24} \\ 1 \leq i \leq 3 \end{array} \right\} \subseteq V$$

be a set interval vector subspace of V.
Define $\eta I : T \to I [0, 1]$ by

$$\eta I \begin{bmatrix} [0, a_1] \\ [0, a_2] \\ [0, a_3] \end{bmatrix} = \begin{cases} [0, \dfrac{1}{3}] & \text{if } a_i \neq 0, i = 1, 2, 3 \\ [0, \dfrac{1}{2}] & \text{if atleast one of } a_i \neq 0, 1 \leq i \leq 3 \\ [0,1] & \text{if } a_i = 0, 1 \leq i \leq 3 \end{cases}$$

(T, $\eta I$) is a fuzzy set interval vector subspace of V. Clearly $\eta I$ cannot be extended to whole of V.

**DEFINITION 3.1.6**: *Let V be a set interval linear algebra over the set S. Suppose $W \subseteq V$ be a set interval linear algebra of V over the set S. Suppose $\eta I : W \to I [0, 1]$ is such that (W, $\eta I$) or $\eta I W$ is a fuzzy set interval linear algebra then we define (W, $\eta I$) to be a fuzzy set interval linear subalgebra of V.*

We will illustrate this situation by some examples.

***Example 3.1.10***: Let

$$V = \left\{ \begin{pmatrix} [0,a_1] & [0,a_3] \\ [0,a_3] & [0,a_4] \end{pmatrix} \right.$$

where $a_i \in Z^+ \cup \{0\}$; $1 \leq i \leq 4\}$ be a set interval linear algebra over the set $S = \{0, 2, 5, 8, 11, 16\} \subseteq Z^+ \cup \{0\}$.
Let

$$W = \left\{ \begin{pmatrix} [0,a_1] & [0,a_3] \\ 0 & [0,a_3] \end{pmatrix} \right.$$

where $a_i \in Z^+ \cup \{0\}$; $1 \leq i \leq 3\} \subseteq V$ be a set interval linear subalgebra of V over the set S.

Define $\eta I : W \to I [0, 1]$ as follows.

$$\eta I \begin{pmatrix} [0,a_1] & [0,a_3] \\ 0 & [0,a_3] \end{pmatrix} = \begin{cases} [0, \dfrac{1}{a_1 + a_2 + a_3}] & \text{if } a_1 + a_2 + a_3 \neq 0 \\ [0,1] & \text{if } 0 = a_1 = a_2 = a_3 \end{cases}$$

(W, ηI) or ηI W is a set fuzzy interval linear subalgebra of V.

***Example 3.1.11***: Let

$$V = \left\{ \begin{bmatrix} [0,a_1] \\ [0,a_2] \\ [0,a_3] \\ [0,a_4] \\ [0,a_5] \\ [0,a_6] \end{bmatrix} \middle| a_i \in Z_{12}; 1 \leq i \leq 6 \right\}$$

be a set interval linear algebra over the set $S = \{0, 2, 1, 3\} \subseteq V$.
Choose

$$W = \left\{ \begin{bmatrix} [0,a_1] \\ 0 \\ [0,a_2] \\ 0 \\ [0,a_3] \\ 0 \end{bmatrix} \middle| a_i \in Z_{12}; 1 \le i \le 3 \right\} \subseteq V$$

be a set interval linear subalgebra of V.
Define $\eta I : W \to I[0, 1]$ by

$$\eta I = \begin{bmatrix} [0,a_1] \\ 0 \\ [0,a_2] \\ 0 \\ [0,a_3] \\ 0 \end{bmatrix} = \begin{cases} [0,\frac{1}{a_1}] & \text{if } a_1 \ne 0; a_2 = a_3 = 0 \\ [0,\frac{1}{a_2}] & \text{if } a_2 \ne 0; a_1 = 0 = a_3 \\ [0,\frac{1}{a_3}] & \text{if } a_3 \ne 0; a_1 = 0 = a_2 \\ [0,\frac{1}{3}] & \text{if } a_i \ne 0; 1 \le i \le 3 \text{ or any two nonzero} \\ [0,1] & \text{if } a_1 = 0; i = 1,2,3 \end{cases}$$

(W, $\eta I$) is a fuzzy set interval linear subalgebra of V.

Now we proceed onto define the notion of semigroup fuzzy interval vector space.

**DEFINITION 3.1.7**: *Let V be a semigroup interval vector space defined over the semigroup S. (V, $\eta I$) or V$\eta I$, the semigroup fuzzy interval vector space is a map $\eta I : V \to I[0, 1]$ satisfying the following condition:*
$$\eta I(ra) \ge \eta I(a)$$
*for all $a \in V$ and $r \in S$.*

We will illustrate this situation by some simple examples.

**Example 3.1.12**: Let V = {[0, a]| a ∈ $Z_{124}$} be a semigroup interval vector space over the semigroup S = $Z_{124}$.
Define ηI: V → I [0, 1] as

$$\eta I([0, a]) = \begin{cases} [0, \frac{1}{a}] & \text{if } a \neq 0 \\ [0,1] & \text{if } a = 0 \end{cases}$$

(V, ηI) is a fuzzy semigroup interval vector space or semigroup fuzzy interval vector space.

**Example 3.1.13:** Let
$$V = \left\{ \begin{bmatrix} [0, a_1] \\ [0, a_2] \end{bmatrix}, \right.$$

([0, $a_1$], [0, $a_2$], [0, $a_3$], [0, $a_4$]) | $a_i$ ∈ $Z^+ \cup \{0\}$; 1 ≤ i ≤ 4} be a semigroup interval vector space over the semigroup S = $Z^+ \cup \{0\}$. Define ηI : V → I [0, 1] as

$$\eta I \left( \begin{bmatrix} [0, a_1] \\ [0, a_2] \end{bmatrix} \right) = \begin{cases} \left[ 0, \frac{1}{a_1 + a_2} \right] & \text{if } a_1 + a_2 \neq 0 \\ [0,1] & \text{if } a_1 = a_2 = 0 \end{cases}$$

and

ηI ([0, $a_1$], [0, $a_2$], [0, $a_3$], [0, $a_4$]) =

$$\begin{cases} \left[ 0, \frac{1}{a_1} \right] & \text{if } a_i \neq 0; \ i = 1, 2, 3, 4 \\ \left[ 0, \frac{1}{2} \right] & \text{if atleast one of } a_i \neq 0; \ 1 \leq i \leq 4 \\ [0,1] & \text{if } a_i = 0; \ i = 1, 2, 3, 4 \end{cases}$$

(V, ηI) or ηIV is a fuzzy semigroup interval vector space or semigroup fuzzy interval vector space.

*Example 3.1.14*: Let

$$V = \left\{ \begin{bmatrix} [0,a_1] & [0,a_2] \\ [0,a_3] & [0,a_4] \\ [0,a_5] & [0,a_6] \\ [0,a_7] & [0,a_8] \end{bmatrix}, ([0,a_1],[0,a_2]), \begin{bmatrix} [0,a_1] \\ [0,a_2] \\ [0,a_3] \end{bmatrix} \middle| \begin{array}{l} a_i \in Z_{23}; \\ 1 \le i \le 8 \end{array} \right\}$$

be a semigroup interval vector space over the semigroup $Z_{23}$.
Define $\eta I : V \to I[0, 1]$ as follows:

$$\eta I \left( \begin{bmatrix} [0,a_1] & [0,a_2] \\ [0,a_3] & [0,a_4] \\ [0,a_5] & [0,a_6] \\ [0,a_7] & [0,a_8] \end{bmatrix} \right) = \begin{cases} \left[0, \dfrac{1}{5}\right] & \text{if atleast one } a_i \ne 0 \\ [0,1] & \text{if all } a_i = 0;\ 1 \le i \le 8 \end{cases}$$

$$\eta I\,(([0, a_1], [[0, a_2]) = \begin{cases} \left[0, \dfrac{1}{a_1}\right] & \text{if } a_1 \ne 0 \text{ and } a_2 = 0 \\ \left[0, \dfrac{1}{a_2}\right] & \text{if } a_2 \ne 0 \text{ and } a_1 = 0 \\ \left[0, \dfrac{1}{10}\right] & \text{if } a_1 \ne 0 \text{ and } a_2 \ne 0 \\ [0,1] & \text{if } a_1 = a_2 = 0 \end{cases}$$

and

$$\eta \left( \begin{bmatrix} [0,a_1] \\ [0,a_2] \\ [0,a_3] \end{bmatrix} \right) = \begin{cases} \left[0, \dfrac{1}{9}\right] & \text{if atleast one of } a_i \ne 0\ \ 1 \le i \le 3 \\ [0,1] & \text{if } a_i = 0;\ i = 1, 2, 3 \end{cases}$$

(V, ηI) or ηIV is a fuzzy semigroup interval linear algebra or semigroup fuzzy interval linear algebra.

Now we proceed onto illustrate only by examples fuzzy semigroup interval linear algebras and leave the simple task of defining semigroup fuzzy interval linear algebras to the reader.

*Example 3.1.15:* Let

$$V = \left\{ \begin{bmatrix} [0,a_1] & [0,a_2] \\ [0,a_3] & [0,a_4] \\ [0,a_5] & [0,a_6] \end{bmatrix} \middle| a_i \in Z^+ \cup \{0\}; 1 \leq i \leq 6 \right\}$$

be a semigroup interval linear algebra over the semigroup $S = Z^+ \cup \{0\}$.

Define $\eta I : V \to I[0, 1]$ as follows:

$$\eta I = \left( \begin{bmatrix} [0,a_1] & [0,a_2] \\ [0,a_3] & [0,a_4] \\ [0,a_5] & [0,a_6] \end{bmatrix} \right) =$$

$$\begin{cases} \left[ 0, \dfrac{1}{a_1 + a_2 + \ldots + a_6} \right] & \text{if atleast one of } a_1 + a_2 + \ldots + a_6 \neq 0 \\ [0,1] & \text{if } a_1 = a_2 = \ldots = a_6 = 0 \end{cases}$$

$(V, \eta I)$ or $\eta I\, V$ is a fuzzy semigroup interval linear algebra or semigroup fuzzy interval linear algebra.

*Example 3.1.16*: Let

$$V = \left\{ \begin{bmatrix} [0,a_1] & [0,a_2] \\ [0,a_3] & [0,a_4] \end{bmatrix} \middle| a_i \in Z_7, 1 \leq i \leq 4 \right\}$$

be a semigroup interval linear algebra over the semigroup $Z_7$. Define $\eta I : V \to I[0, 1]$ as follows:

$$\eta I \left( \begin{bmatrix} [0,a_1] & [0,a_2] \\ [0,a_3] & [0,a_4] \end{bmatrix} \right) = \begin{cases} \left[0, \dfrac{1}{a_1}\right] & \text{if } a_1 \neq 0 \\ \left[0, \dfrac{1}{a_2}\right] & \text{if } a_1 = 0 \\ [0,1] & \text{if } a_1 = 0, i = 1,2 \end{cases}$$

V$\eta$I is a semigroup fuzzy interval linear algebra.

*Example 3.1.17*: Let

$$V = \left\{ \sum_{i=0}^{\infty} [0, a_i] x^i \,\middle|\, a_i \in Q^+ \cup \{0\} \right\}$$

be a semigroup interval linear algebra over the semigroup $S = Z^+ \cup \{0\}$.

Define $\eta I : V \to I [0, 1]$ as follows:

$$\eta I \left( \sum_{i=0}^{\infty} [0, a_i] x^i \right) =$$

$$\begin{cases} \left[0, \dfrac{1}{8}\right] & \text{if the degree of the interval polynomial is} \\ & \text{greater than or equal to three} \\ [0,1] & \text{if the degree of the polynomial is less than} \\ & \text{three this includes zero polynomial} \end{cases}$$

(V, $\eta$I) is a semigroup fuzzy interval linear algebra.

As in case of semigroup interval vector spaces we can define in case of semigroup interval linear algebras the concept of fuzzy semigroup linear subalgebras.

We just illustrate this situation by examples.

***Example 3.1.18***: Let

$$V = \left\{ \begin{bmatrix} [0,a_1] & [0,a_2] & [0,a_3] \\ [0,a_7] & [0,a_4] & [0,a_5] \\ [0,a_8] & [0,a_9] & [0,a_6] \end{bmatrix} \middle| a_i \in Z_8; 1 \le i \le 9 \right\}$$

be a semigroup interval linear algebra defined over the semigroup $S = \{0, 2, 4, 6\}$, under addition modulo 8.
    Consider

$$W = \left\{ \begin{bmatrix} [0,a_1] & [0,a_2] & [0,a_3] \\ 0 & [0,a_4] & [0,a_5] \\ 0 & 0 & [0,a_6] \end{bmatrix} \middle| a_i \in Z_8; 1 \le i \le 6 \right\} \subseteq V;$$

W is a semigroup interval linear subalgebra of V.
    Define $\eta I : W \to I\,[0, 1]$

$$\eta I \left( \begin{bmatrix} [0,a_1] & [0,a_2] & [0,a_3] \\ 0 & [0,a_4] & [0,a_5] \\ 0 & 0 & [0,a_6] \end{bmatrix} \right) = \begin{cases} \left[0, \dfrac{1}{a_1}\right] & \text{if } a_1 \ne 0 \\[4pt] \left[0, \dfrac{1}{a_2}\right] & \text{if } a_2 \ne 0 \text{ if } a_1 = 0 \\[4pt] \left[0, \dfrac{1}{a_3}\right] & \text{if } a_3 \ne 0 \text{ if } a_1 = a_2 = 0 \\[4pt] \left[0, \dfrac{1}{a_4}\right] & \text{if } a_3 \ne 0 \text{ if } a_1 = a_2 = 0 \\[4pt] [0,1] & \text{if } a_1 = 0, i = 1,2,...,6 \end{cases}$$

(W, $\eta I$) is a semigroup fuzzy interval linear subalgebra.

***Example 3.1.19***: Let

$$V = \left\{ \sum_{i=0}^{25} [0,a_i] x^i \,\middle|\, a_i \in Z_{40} \right\}$$

be a semigroup interval linear algebra defined over the semigroup $S = \{0, 10, 20, 30\} \subseteq Z_{40}$.

$$W = \left\{ \sum_{i=0}^{10} [0, a_i] x^i \,\middle|\, a_i \in Z_{40} \right\} \subseteq V$$

be a semigroup interval linear subalgebra of V over S.

Define $\eta I : W \to I[0, 1]$ as follows:

$$\eta I \left\{ p(x) = \sum_{i=0}^{10} [0, a_i] x^i \right\} =$$

$$\begin{cases} \left[ 0, \dfrac{1}{a_i} \right]; & \begin{array}{l} [0, a_i] \text{ corresponds to the coefficient interval} \\ \text{of the highest degree of x in p(x)} \end{array} \\ \\ [0,1] & \text{if p (x) is a constant polynomial} \end{cases}$$

(W, $\eta$I) is a fuzzy semigroup interval linear subalgebra.

*Example 3.1.20:* Let $V = \{[0, a_i] \mid a_i \in Z^+ \cup \{0\}\}$ be a semigroup interval linear algebra over the semigroup $S = \{3Z^+ \cup \{0\}\}$. Consider $W = \{[0, a_i] \mid a_i \in 5Z^+ \cup \{0\}\} \subseteq V$, W is a semigroup interval linear subalgebra over the semigroup $S = \{3Z^+ \cup \{0\}\}$.

Define $\eta I : W \to I[0, 1]$ as follows:

$$\eta I\,([0, a_i]) = \begin{cases} \left[ 0, \dfrac{1}{a_i} \right] & \text{if } a_i \neq 0 \\ [0,1] & \text{if } a_1 = 0 \end{cases}$$

(W, $\eta$I) is a fuzzy semigroup interval linear subalgebra.

Now we can define for group interval vector spaces the notion of group fuzzy interval vector spaces or fuzzy group interval vector spaces.

**DEFINITION 3.1.8**: *Let V be a group interval linear algebra defined over the group G.*
*Let $\eta I : V \to I [0, 1]$ such that*
$$\eta (a + b) \geq \min \{\eta(a), \eta(b)\}$$
$$\eta(-a) = \eta(a)$$
$$\eta(0) = 1$$
$$\eta(ra) \geq \eta(a)$$
*for all $a, b \in V$ and $r \in G$.*
*We call $V\eta I$ or $(V, \eta I)$ to be the group fuzzy interval linear algebra.*

We will illustrate this situation by some examples.

***Example 3.1.21:*** Let $V = \{([0, a_1], [0, a_2], [0, a_3], [0, a_4]) \mid a_i \in Z_{40}; 1 \leq i \leq 4\}$ be a group interval linear algebra over the group $G = \{0, 10, 20, 30\} \subseteq Z_{40}$.
Define $\eta I : V \to I [0, 1]$ as follows:

$$\eta I ([0, a_1], [0, a_2], [0, a_3], [0, a_4]) = \begin{cases} \left[0, \dfrac{1}{a_1}\right] & \text{if } a_1 \neq 0 \\ [0,1] & \text{if } a_1 = 0 \end{cases}$$

$(V, \eta I)$ is a group fuzzy interval vector space.

***Example 3.1.22***: Let
$$V = \left\{\sum_{i=0}^{\infty}[0, a_i]x^i \,\middle|\, a_i \in Z_{11}\right\}$$

be a group interval linear algebra over the $G = Z_{11}$. Define $\eta I : V \to I [0, 1]$ as follows:

$$\eta I\ (p(x)) = \begin{cases} \left[0, \dfrac{1}{\deg(p(x))}\right] \\ [0,1] \text{ if } \deg p(x) = 0 \end{cases}$$

$(V, \eta I)$ is a fuzzy group interval vector space; or group fuzzy interval vector space.

*Example 3.1.23*: Let

$$V = \left\{ \begin{bmatrix} [0, a_1] \\ [0, a_2] \\ [0, a_3] \\ [0, a_4] \end{bmatrix} \middle| a_i \in Z_{25}; 1 \le i \le 4 \right\}$$

be a group interval linear algebra over the group $G = Z_{25}$. Define $\eta I : V \to I\,[0, 1]$ as follows:

$$\eta I\left(\begin{bmatrix} [0, a_1] \\ [0, a_2] \\ [0, a_3] \\ [0, a_4] \end{bmatrix}\right) = \begin{cases} \left[0, \dfrac{1}{a_1}\right] \text{ if } a_1 \ne 0 \\ \left[0, \dfrac{1}{a_2}\right] \text{ if } a_2 \ne 0 \text{ if } a_1 = 0 \\ \left[0, \dfrac{1}{a_3}\right] \text{ if } a_3 \ne 0 \text{ if } a_1 = a_2 = 0 \\ \left[0, \dfrac{1}{a_4}\right] \text{ if } a_4 \ne 0 \text{ if } a_1 = a_2 = a_3 = 0 \\ [0,1] \text{ if } a_i = 0, 1 \le i \le 4 \end{cases}$$

$(V, \eta I)$ is a group fuzzy interval linear algebra.

The concept of group fuzzy interval linear subalgebra and group fuzzy interval vector subspaces is left for the reader to define in an analogous way. However we give examples of them.

*Example 3.1.24*: Let

$$V = \left\{ \begin{bmatrix} [0,a_1] & [0,a_2] \\ [0,a_3] & [0,a_4] \\ [0,a_5] & [0,a_6] \end{bmatrix} \middle| a_i \in Z_{21}, 1 \le i \le 6 \right\}$$

be group interval vector space.
Define $\eta I : V \to I[0, 1]$ as follows:

$$\eta I \left( \begin{bmatrix} [0,a_1] & [0,a_2] \\ [0,a_3] & [0,a_4] \\ [0,a_5] & [0,a_6] \end{bmatrix} \right) = \begin{cases} \left[0, \dfrac{1}{\max\{a_i\}}\right]; 1 \le i \le 6 \\ [0,1] \text{ if } a_1 = 0, i = 1, 2, ..., 6 \end{cases}$$

That is if

$$x = \begin{bmatrix} [0,8] & [0,17] \\ [0,4] & [0,1] \\ [,2] & [0,19] \end{bmatrix} \in V$$

then

$$\eta I (x) = \left\{ \left[0, \dfrac{1}{19}\right] \right\}.$$

Thus $(V, \eta I)$ is a group fuzzy interval vector space.
    Take

$$W = \left\{ \begin{bmatrix} [0,a_1] & 0 \\ [0,a_2] & 0 \\ [0,a_3] & 0 \end{bmatrix} \middle| a_i \in Z_{21}, 1 \le i \le 3 \right\} \subseteq V.$$

W is a group interval vector subspace of V over the group $G = Z_{21}$.
    $\eta I : W \to I[0, 1]$ is defined as follows:

$$\eta I \begin{bmatrix} [0,a_1] & 0 \\ [0,a_2] & 0 \\ [0,a_3] & 0 \end{bmatrix} = \begin{cases} \left[0, \dfrac{1}{a_2}\right] & \text{if } a_2 \neq 0 \\ [0,1] & \text{if } a_1 = 0 \end{cases}$$

Clearly (W, $\eta I$) is a fuzzy group interval vector subspace of V.

*Example 3.1.25*: Let V = {Collection of all $10 \times 12$ interval matrices; $([0, a_i])$ with entries from $Z_{36}$ that $a_i \in Z_{36}$; $1 \leq i \leq 20$} be a group interval vector space over the group $G = Z_{36}$.

W = {$([0, a_i])$ demotes all matrices with entries from $2Z_{36}$}.
Define

$$\eta I\,([0, a_i]) = \begin{cases} \left[0, \dfrac{1}{\max\{a_i\}}\right]; 1 \leq i \leq 36;\ a_i \in 2Z_{36} \\ [0,1] \text{ if } a_1 = 0, i = 1, 2, \ldots, 36 \end{cases}$$

(W, $\eta I$) is a group fuzzy interval vector subspace.

*Example 3.1.26*: Let V = {All upper triangular $4 \times 4$ interval matrices constructed using $Z_{13}$} be the group interval vector space over the group $G = Z_{13}$.

Let W = {all $4 \times 4$ diagonal interval matrices with entries from $Z_{13}$} $\subseteq$ V; W is a group interval vector subspace of V over the group $G = Z_{13}$.

Define $\eta I : W \to I[0, 1]$ as follows.

$$\eta I \begin{pmatrix} \begin{bmatrix} [0,a_1] & 0 & 0 & 0 \\ 0 & [0,a_2] & 0 & 0 \\ 0 & 0 & [0,a_3] & 0 \\ 0 & 0 & 0 & [0,a_4] \end{bmatrix} \end{pmatrix} = \begin{cases} \left[0, \dfrac{1}{\max\{a_i\}}\right] \\ [0,1] \text{ if } a_i = 0, 1 \leq i \leq 4 \end{cases}$$

(W, $\eta I$) is a group fuzzy interval vector subspace. Suppose

$$x = \begin{bmatrix} [0,3] & 0 & 0 & 0 \\ 0 & [0,7] & 0 & 0 \\ 0 & 0 & [0,11] & 0 \\ 0 & 0 & 0 & [0,1] \end{bmatrix} \in W$$

then

$$\eta I(x) = \left\{ \left[0, \frac{1}{11}\right] \right\}.$$

We see in case of group interval linear algebras or group interval vector spaces we cannot use groups other than $Z_n$, under addition modulo n. As Z or Q or R cannot be used since all the intervals we use are of the form $[0, a_i]$. $0 \le a_i$.

Now having seen fuzzy set interval vector spaces, fuzzy semigroup interval vector spaces and group fuzzy interval vector spaces we proceed onto define another type of fuzzy set interval vector spaces, fuzzy semigroup interval vector spaces and fuzzy group interval vector spaces by constructing directly and not using set interval vector spaces, semigroup interval vector spaces or group interval vector spaces. These we call as type II set fuzzy interval vector spaces and so on. Those fuzzy interval vector spaces constructed in section 3.1 will be known as type I spaces.

In the following section we define type II fuzzy interval spaces.

## 3.2 Set Fuzzy Interval Vector Spaces of Type II and Their Properties

In this section we proceed on to define set fuzzy interval vector spaces of type II, semigroup fuzzy interval vector spaces of type II and group fuzzy interval vector spaces of type II using fuzzy intervals recall $I[0,1] = \{$all intervals of the form $[0, a_i]$; $0 \le a_i \le 1\}$; known as fuzzy intervals.

**DEFINITION 3.2.1**: *Let $V = \{[0, a_i] / 0 \le a_i \le 1; [0, a_i] \in I[0,1]\}$. Let S be a set such that for each $s \in S$ and $v \in V$, sv and vs are* *in V. We then call V to be a set fuzzy interval vector space of type II.*

We will illustrate this by some examples.

**Example 3.2.1**: Let $V = \{[0, a_i] \mid 0 \leq a_i \leq 1\}$ be a set fuzzy interval vector space over the set $S = \{0, 1, \frac{1}{2}, 1/2^2, \ldots, 1/2^n\}$. Here for any $v = [0, a_i]$ and $s = \dfrac{1}{2^r}$ ($r \leq n$) we have

$$sv = \left[0, \frac{a_i}{2^r}\right] = vs$$

and $sv \in V$. Thus V is a fuzzy set interval vector space of type II over the set S.

**Example 3.2.2**: Let

$$V = \left\{ \begin{bmatrix} [0,a_1] \\ [0,a_2] \\ \vdots \\ [0,a_5] \end{bmatrix}, \left([0,a_1] \; [0,a_2] \; [0,a_3]\right) \;\middle|\; 0 \leq a_i \leq 1; 1 \leq i \leq 5 \right\}$$

be a fuzzy set interval vector space of type II over the set $S = \{0, 1, 1/5, 1/10, 1/121, 1/142\}$.

**Example 3.2.3**: Let

$$V = \left\{ \begin{bmatrix} [0,a_1] & [0,a_2] \\ [0,a_3] & [0,a_4] \\ [0,a_5] & [0,a_6] \\ [0,a_7] & [0,a_8] \\ [0,a_9] & [0,a_{10}] \end{bmatrix}, \left([0,a_1] \; [0,a_2] \; [0,a_3] \; [0,a_4] \; [0,a_5]\right) \right.$$

$$\left\{ \begin{bmatrix} [0,a_1] & [0,a_2] \\ [0,a_3] & [0,a_4] \end{bmatrix} \middle| 0 \le a_i \le 1; 1 \le i \le 10 \right\}$$

be a set fuzzy interval vector space over the set

$$S = \left\{ \frac{1}{3^n} \middle| n = 0,1,...,27 \right\}.$$

*Example 3.2.4:* Let

$$V = \left\{ \sum_{i=0}^{\infty} [0,a_i]x^i \middle| 0 \le a_i \le 1 \right\}$$

be a set fuzzy interval vector space over the set $S = \{0, 1, 1/3, 1/7, 1/5, 1/11, 1/13, 1/19, 1/17, 1/23\}$ of type II.

Now we define substructures of set fuzzy interval vector space.

**DEFINITION 3.2.2**: *Let V be a set fuzzy interval vector space over the set S of type II.*

*Let $W \subseteq V$; if W is a set fuzzy interval vector space over the set S of type II, then we define W to be a set fuzzy interval vector subspace of V over the set S of type II.*

We will illustrate this situation by examples.

*Example 3.2.5:* Let

$$V = \left\{ \begin{bmatrix} [0,a_1] & [0,a_2] \\ [0,a_3] & [0,a_4] \end{bmatrix}, \begin{bmatrix} [0,a_1] \\ [0,a_2] \\ [0,a_3] \\ [0,a_4] \end{bmatrix} \middle| 0 \le a_i \le 1; 1 \le i \le 4 \right\}$$

be a set fuzzy interval vector space of type II over the set

$$S = \left\{ \frac{1}{2^n}, 0 \middle| n = 0,1,2,...,41 \right\}.$$

Let

$$W = \left\{ \begin{bmatrix} [0,a_1] & [0,a_2] \\ [0,a_3] & [0,a_4] \end{bmatrix} \middle| 0 \le a_i \le 1; 1 \le i \le 4 \right\} \subseteq V;$$

W is a set fuzzy interval vector subspace of type II over the set

$$S = \left\{ 0, \frac{1}{2^n} \middle| n = 0,1,2,...,41 \right\}.$$

*Example 3.2.6*: Let

$$V = \left\{ \begin{bmatrix} [0,a_1] \\ [0,a_2] \\ [0,a_3] \\ [0,a_4] \\ [0,a_5] \end{bmatrix} \middle| 0 \le a_i \le 1; 1 \le i \le 5 \right\}$$

be a set fuzzy interval vector space of type II over the set $S = \{0, 1/3, 1/3^2, ..., 1/3^7\}$.

$$W = \left\{ \begin{bmatrix} 0 \\ [0,a_1] \\ 0 \\ [0,a_2] \\ 0 \end{bmatrix} \middle| a_1, a_2 \in [0,1] \right\} \subseteq V;$$

W is a set fuzzy interval vector subspace of type II of V over the set S.

We give the definition of yet another substructure.

**DEFINITION 3.2.3**: *Let V be a set fuzzy interval vector space over the set S of type II. Let $W \subseteq V$ be a proper subset of V and $P \subseteq S$ be a set fuzzy subset S. W be a set fuzzy interval vector* *space of type II over the set P, we call W to be a subset fuzzy interval vector subspace of V of type II over the subset P of S.*

We will illustrate this situation by some simple examples.

*Example 3.2.7:* Let

$$V = \left\{ \begin{bmatrix} [0,a_1] & [0,a_2] \\ [0,a_3] & [0,a_4] \end{bmatrix} \middle| a_i \in [0,1], 1 \leq i \leq 4 \right\}$$

be a set fuzzy interval vector space over the set $S = \{0, 1, 1/2^n, 1/3^n; n = 1, 2, \ldots, 12\}$ of type II. Let

$$W = \left\{ \begin{bmatrix} [0,a_1] & [0,a_2] \\ 0 & [0,a_4] \end{bmatrix} \middle| 0 \leq a_1, a_2, a_3 \leq 1 \right\} \subseteq V$$

and $P = \{0, 1, 1/2^n \mid n = 1, 2, \ldots, 12\} \subseteq S$. W is a subset fuzzy interval vector subspace of V of type two over the subset P of S.

*Example 3.2.8*: Let

$$V = \left\{ \begin{bmatrix} [0,a_1] \\ [0,a_2] \end{bmatrix}, ([0,a_1][0,a_2][0,a_3][0,a_4][0,a_5][0,a_6]) \middle| \begin{array}{l} 0 \leq a_i \leq 1; \\ i = 1,\ldots,6 \end{array} \right\}$$

be a set fuzzy interval vector space of type II over the set

$S = \{0, 1/2^n, 1, 1/3^m, 1/5^m, 1/7^n \mid 1 \leq m \leq 8, 1 \leq n \leq 12\}$.
Choose

$$W = \left\{ \begin{bmatrix} [0,a_1] \\ 0 \end{bmatrix}, ([0,a_1]\ 0\ [0,a_2]\ 0\ [0,a_3]\ 0) \middle| \begin{array}{l} 0 \leq a_i \leq 1; \\ 1 \leq i \leq 3 \end{array} \right\} \subseteq V$$

and $P = \{0, 1, 1/7^n \mid 1 \leq n \leq 12\} \subseteq S$. W is a subset fuzzy interval vector subspace of V of type II over the subset P of S.

*Example 3.2.9:* Let

$$V = \left\{ \begin{bmatrix} [0,a_1] & [0,a_2] \\ [0,a_3] & [0,a_4] \\ [0,a_5] & [0,a_6] \\ [0,a_7] & [0,a_8] \end{bmatrix}, \begin{bmatrix} [0,a_1] & [0,a_2] & [0,a_3] & [0,a_4] \\ [0,a_5] & [0,a_6] & [0,a_7] & [0,a_8] \end{bmatrix} \middle| \begin{array}{l} 0 \le a_i \le 1; \\ 1 \le i \le 8 \end{array} \right\}$$

be a set fuzzy interval vector space of type II over the set

$$S = \left\{ 0, \frac{1}{n} \middle| n \in Z^+ \right\}.$$

Choose

$$W = \left\{ \begin{bmatrix} [0,a_1] & [0,a_2] & [0,a_3] & [0,a_4] \\ [0,a_5] & [0,a_6] & [0,a_7] & [0,a_8] \end{bmatrix} \middle| 0 \le a_i \le 1; 1 \le i \le 8 \right\} \subseteq V$$

is a subset fuzzy interval vector space of type II over the subset

$$P = \left\{ 0, \frac{1}{4n} \middle| n \in Z^+ \right\} \subseteq S \text{ of } V.$$

*Example 3.2.10*: Let

$$V = \left\{ \begin{bmatrix} [0,a_1] & [0,a_2] \\ [0,a_3] & [0,a_4] \end{bmatrix} \middle| 0 \le a_i \le 1; 1 \le i \le 4 \right\}$$

be a set fuzzy interval linear algebra over the set

$$S = \left\{ \frac{1}{2^n}, 0 \middle| n = 0, 1, 2, \ldots, \infty \right\}$$

(the operation on V is max, i.e., if

$$x = \begin{bmatrix} [0,a_1] & [0,a_2] \\ [0,a_3] & [0,a_4] \end{bmatrix} \text{ and } y = \begin{bmatrix} [0,b_1] & [0,b_2] \\ [0,b_3] & [0,b_4] \end{bmatrix}$$

are in V then

$$x + y = \begin{bmatrix} [0, \max(a_1, b_1)] & [0, \max(a_2, b_2)] \\ [0, \max(a_3, b_3)] & [0, \max(a_4, b_4)] \end{bmatrix}.$$

Thus max (x, y) denoted by x + y is an associate closed commutative operation on V).

Choose

$$W = \left\{ \begin{bmatrix} [0, a_1] & [0, a_2] \\ 0 & [0, a_3] \end{bmatrix} \middle| 1 \le i \le 3; 0 \le a_i \le 1 \right\} \subseteq V;$$

W is a subset fuzzy interval vector subspace of V defined over the subset

$$P = \left\{ \frac{1}{2^{n+6}} \middle| n \in Z^+ \right\} \subseteq S.$$

Now we proceed onto define set fuzzy interval linear algebra formally.

**DEFINITION 3.2.4**: *Let V be a set fuzzy interval vector space over a set S. If on V is defined a closed associative binary operation denoted by '+' such that s (a + b) = sa + sb; for all s ∈ S and a, b ∈ V. Then we define V to be a set fuzzy interval linear algebra of type II.*

We will illustrate this by some simple examples.

***Example 3.2.11***: Let V = {[0, $a_i$] | 0 ≤ $a_i$ ≤ 1} be a set fuzzy interval linear algebra over the set

$$S = \left\{ \frac{1}{2^n}, 0, 1 \middle| n \in Z^+ \right\}.$$

V is closed under max operation that is for [0, a] and [0, b] in V we have [0, a] + [0, b] = [0, max {a, b}].

***Example 3.2.12***: Let V = {collection of all 3 × 5 interval fuzzy matrices with entries from I [0, 1]} be a set fuzzy interval linear algebra over the set

$$S = \left\{\frac{1}{n}, 0 \,\middle|\, n = 1, 2, 3, \ldots\right\}.$$

***Example 3.2.13***: Let

$$V = \left\{\sum_{i=0}^{\infty}[0, a_i]x^i \,\middle|\, 0 \leq a_i \leq 1\right\}$$

be a set fuzzy interval linear algebra over the set

$$S = \left\{0, \frac{1}{2^n} \,\middle|\, n = 1, 2, \ldots\right\}.$$

We will define set fuzzy interval linear subalgebra.

**DEFINITION 3.2.5**: *Let V be a set fuzzy interval linear algebra over the set*

$$S = \left\{\frac{1}{n}, 0 \,\middle|\, n = 1, 2, \ldots\right\}.$$

*Choose $W \subseteq V$; suppose W be a set fuzzy interval linear algebra over the set S; we define W to be set fuzzy interval linear subalgebra of V over S of type II.*

We will illustrate this situation by some examples.

***Example 3.2.14***: Let

$$V = \left\{\begin{bmatrix} [0, a_1] & [0, a_2] \\ [0, a_3] & [0, a_4] \end{bmatrix} \,\middle|\, 0 \leq a_i \leq 1; 1 \leq i \leq 4\right\}$$

be a set fuzzy interval linear algebra with max operation defined over the set

$$S = \left\{\frac{1}{2n+1}, 0 \,\middle|\, n = 1, 2, \ldots\right\}.$$

Choose

$$W = \left\{ \begin{bmatrix} [0,a_1] & [0,a_2] \\ 0 & [0,a_3] \end{bmatrix} \middle| 1 \leq i \leq 3; 0 \leq a_i \leq 1 \right\} \subseteq V;$$

W is a set fuzzy interval linear subalgebra of V over the set S.

*Example 3.2.15*: Let

$$V = \left\{ \begin{bmatrix} [0,a_1] \\ [0,a_2] \\ [0,a_3] \\ [0,a_4] \\ [0,a_5] \\ [0,a_6] \end{bmatrix} \middle| 0 \leq a_i \leq 1; 1 \leq i \leq 6 \right\}$$

be a set fuzzy interval linear algebra over the set

$$S = \left\{ \frac{1}{3n+1}, 0 \middle| n = 0,1,2,... \right\}.$$

Let

$$W = \left\{ \begin{bmatrix} [0,a_1] \\ 0 \\ [0,a_2] \\ 0 \\ [0,a_3] \\ 0 \end{bmatrix} \middle| 0 \leq a_i \leq 1; 1 \leq i \leq 3 \right\} \subseteq V$$

be the set fuzzy interval linear subalgebra over the set

$$S = \left\{ \frac{1}{3n+1}, 0 \middle| n = 0,1,2,... \right\} \text{ of } V.$$

**DEFINITION 3.2.6**: *Let V be a set fuzzy interval linear algebra over the set S. Suppose $W \subseteq V$; if W is a set fuzzy interval linear algebra over the subset P of S, then we define W to be subset fuzzy interval linear subalgebra of V of type II over the subset P of S.*

We will illustrate this situation by an example.

*Example 3.2.16*: Let

$$V = \left\{ \begin{bmatrix} [0,a_i] \\ [0,b_i] \end{bmatrix} \,\middle|\, \begin{array}{l} 0 \le a_i \le 1 \\ 0 \le b_i \le 1 \end{array} \right\}$$

be a set fuzzy interval linear algebra over the set

$$S = \left\{ 0, \frac{1}{n} \,\middle|\, n = 1, 2, \ldots \right\}$$

with min operation on V. That is min $\{[0, a_i], [0, b_i]\} = [0,$ min $\{a_i, b_i\}]$. Choose

$$W = \left\{ \begin{bmatrix} [0,a_i] \\ 0 \end{bmatrix} \,\middle|\, 0 \le a_i \le 1 \right\} \subseteq V$$

is a subset fuzzy interval linear subalgebra over the subset

$$P = \left\{ 0, \frac{1}{4n} \,\middle|\, n = 1, 2, \ldots, \infty \right\} \subseteq S \text{ of } S.$$

We as in case of usual set vector spaces and set interval vector spaces define set fuzzy interval vector space linear transformation.

We will illustrate this by some examples.

*Example 3.2.17*: Let

$$V = \left\{ \begin{bmatrix} [0,a_1] & [0,a_2] \\ [0,a_3] & [0,a_4] \end{bmatrix} \middle| 0 \le a_i \le 1; 1 \le i \le 4 \right\}$$

and

$$W = \left\{ \begin{bmatrix} [0,a_1] \\ [0,a_2] \\ [0,a_3] \\ [0,a_4] \end{bmatrix} \middle| 0 \le a_i \le 1; 1 \le i \le 4 \right\}$$

be set fuzzy interval linear algebras defined over the set

$$S = \left\{ \frac{1}{2n-1}, 0 \middle| n = 1, 2, \ldots \right\}.$$

Define $T_F : V \to W$ as

$$T_F(A) = T_F = \left( \begin{bmatrix} [0,a_1] & [0,a_2] \\ [0,a_3] & [0,a_4] \end{bmatrix} \right) = \begin{bmatrix} [0,a_1] \\ [0,a_2] \\ [0,a_3] \\ [0,a_4] \end{bmatrix}$$

for every A in V. $T_F$ is a set linear transformation of V into W.

Note as in case of vector spaces we see in case of set fuzzy interval vector spaces define linear transformation over the same set.

*Example 3.2.18*: Let

$$V = \left\{ \begin{bmatrix} [0,a_1] \\ [0,a_2] \\ [0,a_3] \end{bmatrix} \middle| 0 \le a_i \le 1; 1 \le i \le 3 \right\}$$

and

$$W = \left\{ \begin{bmatrix} [0,a_1] & [0,a_2] & [0,a_3] \\ [0,a_4] & [0,a_5] & [0,a_6] \end{bmatrix} \middle| 0 \le a_i \le 1; 1 \le i \le 6 \right\}$$

be set fuzzy interval linear algebra defined over the same set

$$S = \left\{ 0, \frac{1}{2^n} \middle| n = 1, 2, \ldots \right\}.$$

Define $T_F : V \to W$ by

$$T_F = \left( \begin{bmatrix} [0,a_1] \\ [0,a_2] \\ [0,a_3] \end{bmatrix} \right) = \begin{bmatrix} [0,a_1] & [0,a_2] & [0,a_3] \\ [0,a_1] & [0,a_2] & [0,a_3] \end{bmatrix}.$$

$T_F$ is a linear transformation of V into W.

We define a special fuzzy semigroup S as follows:

Let $S = I [0, 1] = \{[0, a_i] \mid 0 \le a_i \le 1\}$. On S define an associative closed binary operation * so that S is a semigroup. We call (S, *) a special fuzzy semigroup.

Throughout this chapter by a fuzzy interval semigroup we mean a semigroup constructed using fuzzy intervals.

**DEFINITION 3.2.7:** *Let V be a set whose elements are constructed using fuzzy intervals from I [0, 1]. S any additive semigroup with 1. We call V to be a fuzzy interval semigroup vector space of level two or type II if the following conditions hold good.*

*(a) $sv \in V$ for all $s \in S$ and $v \in V$.*
*(b) $0.v = 0 \in V$ for all $v \in V$ and $0 \in S$; 0 is the zero vector*
*(c) $(s_1 + s_2) v = s_1v + s_2v$ for all $s_1, s_2 \in S$ and $v \in V$.*

We illustrate this by the following examples. The terms type II and level two are used as synonym.

*Example 3.2.19*: Let $V = \{[0, a_i] \mid 0 \le a_i \le 1\}$ be a semigroup fuzzy interval vector space of level two over the semigroup

$$S = \left\{0, \frac{1}{2^n} \,\middle|\, n = 1, 2, \ldots \right\}$$

under multiplication.
Let

$$V = \left\{ \begin{bmatrix} [0, a_1] \\ [0, a_2] \\ [0, a_3] \end{bmatrix}, \left([0, a_1] \quad [0, a_2] \quad [0, a_3] \quad [0, a_4] \quad [0, a_5]\right) \right\}$$

be a semigroup fuzzy interval vector space of level two over the semigroup

$$(S, o) = \left\{ 0, 1, \frac{1}{2^n}, \frac{1}{3^m} \,\middle|\, m, n \in Z^+ \right\}$$

under 'o' the max operation that is

$$\left\{ \frac{1}{2^m} o \frac{1}{2^n} \right\} \max \left\{ \frac{1}{2^m}, \frac{1}{3^n} \right\} = \frac{1}{2^m} \text{ if } \frac{1}{2^m} > \frac{1}{3^n} \,;\, \frac{1}{3^n} \text{ if } \frac{1}{3^n} > \frac{1}{2^m}$$

*Example 3.2.20*: Let

$$W = \left\{ \begin{bmatrix} [0, a_1] \\ [0, a_2] \\ [0, a_3] \end{bmatrix}, \begin{pmatrix} [0, a_1] & [0, a_2] & [0, a_3] & [0, a_4] \\ [0, a_5] & [0, a_6] & [0, a_7] & [0, a_8] \end{pmatrix} \,\middle|\, \begin{array}{l} 1 \leq i \leq 8 \\ 0 \leq a_i \leq 1 \end{array} \right\}$$

be a semigroup fuzzy interval vector space of level two over the semigroup

$$S = \left\{ 0, 1, \frac{1}{5^n}, \frac{1}{8^m} \,\middle|\, m, n \in Z^+ \right\} \text{ under min operation.}$$

*Example 3.2.21*: Let V =

$$\left\{ \begin{bmatrix} [0,a_1] & [0,a_2] \\ [0,a_3] & [0,a_4] \\ [0,a_5] & [0,a_6] \\ [0,a_7] & [0,a_8] \end{bmatrix}, \begin{pmatrix} [0,a_1] & [0,a_2] & [0,a_3] & [0,a_4] & [0,a_5] \\ [0,a_6] & [0,a_7] & [0,a_8] & [0,a_9] & [0,a_{10}] \end{pmatrix} \middle| \begin{matrix} 0 \le a_i \le 1 \\ 1 \le i \le 10 \end{matrix} \right\}$$

be a semigroup fuzzy interval vector space of level two over the semigroup (S, o).

Now we proceed onto define semigroup fuzzy interval linear algebra V over the semigroup (S, o).

**DEFINITION 3.2.8:** *Let V be a fuzzy semigroup interval vector space over the semigroup (S, o) of type II. If V itself is a special fuzzy semigroup under some operation say '+' and so (a + b) = s o b + s o b for all s ∈ S and a, b ∈ V then we call V to be a fuzzy semigroup interval linear algebra over S of type II.*

We will illustrate this situation by some examples.

*Example 3.2.22*: Let

$$V = \left\{ \begin{bmatrix} [0,a_1] & [0,a_2] \\ [0,a_3] & [0,a_4] \\ [0,a_5] & [0,a_6] \end{bmatrix} \middle| \begin{matrix} 0 \le a_i \le 1 \\ 1 \le i \le 6 \end{matrix} \right\}$$

be a semigroup interval fuzzy linear algebra of type II over the semigroup

$$(S, o) = \left\{ \frac{1}{2^n}, 0, 1 \middle| n = 1, 2, ..., \infty \right\}$$

and 'o' is the min operation. On V we have min operation so that V is a semigroup.

***Example 3.2.23:*** Let

$$V = \left\{ \begin{bmatrix} [0,a_1] \\ [0,a_2] \end{bmatrix} \middle| \begin{array}{l} 0 \le a_i \le 1 \\ 1 \le i \le 2 \end{array} \right\}$$

be a fuzzy semigroup interval linear algebra of type II over the semigroup

$$(S, o) = \left\{ \frac{1}{2^n}, \frac{1}{12^m}, 0, 1 \middle| m, n \in Z^+ \right\}$$

and 'o' is the max operation in S}.

Now having seen examples of fuzzy semigroup interval linear algebras over a semigroup of type II we describe an interesting property related with them.

**THEOREM 3.2.1**: *Every fuzzy semigroup interval linear algebra is a fuzzy semigroup interval vector space over the semigroup.*

The proof is simple as one part follows immediately from the definition and other part is obvious by some examples given on fuzzy semigroup interval vector spaces.

**DEFINITION 3.2.9**: *Let V be a semigroup fuzzy interval linear algebra over the semigroup S. Let W ⊆ V; if W is a semigroup fuzzy interval linear algebra over S then we define W to be a semigroup fuzzy interval linear subalgebra of V over the semigroup S.*

We will illustrate this by some examples.

***Example 3.2.24***: Let V = {All 5 × 5 fuzzy interval matrices with entries from I [0, 1]} with min operation be a fuzzy semigroup interval linear algebra of type II over the semigroup

$$S = \left\{ 0, 1, \frac{1}{2^n} \middle| n = 1, 2, \ldots \right\}$$

with max operation.

Thus if $v = [0, a_i] \in V$ and $S = \dfrac{1}{2^r} \in S$ then

$$sv = vs = \dfrac{1}{2^r} [0, a_i] = \left[0, \dfrac{a_i}{2^r}\right].$$

Consider $M = \{$all upper triangular fuzzy interval $5 \times 5$ matrices with entries from $I[0, 1]\} \subseteq V$; $M$ is a fuzzy semigroup interval linear subalgebra of $V$ over $S$ of type II.

*Example 3.2.25*: Let

$$V = \left\{\sum_{i=0}^{\infty}[0, a_i]x^i \;\middle|\; 0 \le a_i \le 1\right\}$$

be a fuzzy semigroup interval linear algebra of type II with min operation (i.e., if

$$\sum_{i=0}^{\infty}[0, a_i]x^i = p(x)$$

and

$$q(x) = \sum_{i=0}^{\infty}[0, b_i]x^i$$

are in V then

$$p(x) + q(x) = \sum_{i=0}^{\infty}[0, \min\{a_i, b_i\}]x^i$$

over the semigroup

$$S = \left\{0, 1, \dfrac{1}{5^n} \;\middle|\; n = 1, 2, ...\right\}$$

and

$$sv = s\, p(x)\; (s = \dfrac{1}{5^{20}}\; p(x) = \sum_{i=0}^{\infty}[0, a_i]x^i\,)$$

$$= \sum_{i=0}^{\infty}\left[0, \dfrac{1}{5^{20}} a_i\right]x^i.$$

$$W = \left\{ \sum_{i=0}^{\infty} [0, a_i] x^{2i} \,\middle|\, 0 \leq a_i \leq 1 \right\} \subseteq V;$$

W is a semigroup fuzzy interval linear subalgebra of V over the semigroup S of type II.

Now we proceed onto define the notion of fuzzy subsemigroup interval sublinear algebra of V over the subsemigroup P of S.

**DEFINITION 3.2.10**: *Let V be a fuzzy semigroup interval linear algebra of type II over the semigroup S. Let $W \subseteq V$ and $P \subseteq S$ where W and P are proper subsets of V and S respectively. If W is a fuzzy semigroup interval linear algebra of type II over the semigroup P then we define W to be a fuzzy subsemigroup interval linear subalgebra of type II over the subsemigroup P of the semigroup S.*

We illustrate this situation by some examples.

*Example 3.2.26*: Let

$$V = \left\{ \begin{bmatrix} [0,a_1] & [0,a_2] \\ [0,a_3] & [0,a_4] \\ [0,a_5] & [0,a_6] \\ [0,a_7] & [0,a_8] \\ [0,a_9] & [0,a_{10}] \end{bmatrix} \,\middle|\, 0 \leq a_i \leq 1; 1 \leq i \leq 10 \right\}$$

be a fuzzy semigroup interval linear algebra of type II over the semigroup

$$S = \left\{ 1, 0, \frac{1}{2^n} \,\middle|\, n = 1, 2, \ldots \right\}.$$

Let

$$W = \left\{ \begin{bmatrix} [0,a_1] & 0 \\ 0 & [0,a_2] \\ [0,a_3] & 0 \\ 0 & [0,a_4] \\ [0,a_5] & 0 \end{bmatrix} \middle| 0 \leq a_i \leq 1; 1 \leq i \leq 5 \right\} \subseteq V;$$

and

$$P = \left\{ 1, 0, \frac{1}{2^{3n}} \middle| n = 1, 2, \ldots \right\} \subseteq S.$$

W is a fuzzy subsemigroup interval linear subalgebra of type II over the subsemigroup $P \subseteq S$.

*Example 3.2.27*: Let

$$V = \left\{ \begin{bmatrix} [0,a_1] & [0,a_2] \\ [0,a_3] & [0,a_4] \\ [0,a_5] & [0,a_6] \\ [0,a_7] & [0,a_8] \end{bmatrix} \middle| \begin{array}{l} 0 \leq a_i \leq 1 \\ 1 \leq i \leq 8 \end{array} \right\}$$

be a special fuzzy semigroup interval linear algebra over the semigroup

$$S = \left\{ 0, 1, \frac{1}{2^n}, \frac{1}{5^m} \middle| m, n \in Z^+ \right\}.$$

Choose

$$W = \left\{ \begin{bmatrix} [0,a_1] & 0 \\ [0,a_2] & 0 \\ [0,a_3] & 0 \\ [0,a_4] & 0 \end{bmatrix} \middle| \begin{array}{l} 0 \leq a_i \leq 4 \\ 1 \leq i \leq 4 \end{array} \right\} \subseteq V,$$

W is a fuzzy subsemigroup interval linear subalgebra of V over the subsemigroup

$$P = \left\{ 0, 1, \frac{1}{5^n} \,\middle|\, n \in Z^+ \right\} \subseteq S$$

of type II of V.

Now if V has no fuzzy semigroup interval linear subalgebra over F then we define V to be a simple fuzzy semigroup interval linear algebra of type II. We say V is said to be a pseudo simple fuzzy semigroup interval linear algebra over S of type II if V has no fuzzy subsemigroup interval linear algebra over S. We say V is doubly simple if V has no fuzzy semigroup interval linear subalgebras and fuzzy subsemigroup interval linear subalgebras.

We will illustrate this situation by examples.

*Example 3.2.28:* Let

$$V = \left\{ \begin{bmatrix} [0, a_1] & 0 \\ 0 & [0, a_1] \end{bmatrix} \,\middle|\, a_1 = 1 \right\}$$

be a semigroup fuzzy linear algebra over the semigroup $S = \{0, 1\}$ with min operation. It is easily verified V is a doubly simple semigroup fuzzy interval linear algebra of type II over S.

*Example 3.2.29*: Let

$$\left\{ \begin{bmatrix} [0,1/2] \\ [0,1/2] \\ [0,1/2] \end{bmatrix}, \begin{bmatrix} [0,1/4] \\ [0,1/4] \\ [0,1/4] \end{bmatrix}, \begin{bmatrix} [0,1] \\ [0,1] \\ [0,1] \end{bmatrix}, \begin{bmatrix} [0] \\ [0] \\ [0] \end{bmatrix} \right\} = V$$

be a fuzzy semigroup interval linear algebra over the semigroup $S = \{0, 1\}$. V is a pseudo simple fuzzy semigroup interval linear algebra as V has no proper fuzzy subsemigroup interval linear subalgebra of type II.
 However V has fuzzy semigroup interval linear subalgebra of type two over S.

**DEFINITION 3.2.11**: *Let V be a fuzzy semigroup interval linear algebra of type II over the semigroup S. Let W ⊆ V be a fuzzy semigroup interval linear subalgebra of V of type II.*

*Let T : V → W be such that T(v) = w for every v ∈ V and w ∈ W. T is a fuzzy semigroup interval projection of V into W.*

We will illustrate this by some examples.

*Example 3.2.30*: Let

$$V = \left\{ \begin{bmatrix} [0,a_1] & [0,a_2] \\ [0,a_3] & [0,a_4] \end{bmatrix} \middle| 0 \le a_i \le 1; 1 \le i \le 4 \right\}$$

be a fuzzy semigroup interval linear algebra over a semigroup

$$S = \left\{ 0, 1, \frac{1}{2^n} \middle| n \in Z^+ \right\}.$$

Choose

$$W = \left\{ \begin{bmatrix} [0,a_1] & [0,a_2] \\ 0 & [0,a_3] \end{bmatrix} \middle| 0 \le a_i \le 1; 1 \le i \le 3 \right\} \subseteq V;$$

W is a semigroup fuzzy interval linear subalgebra interval linear subalgebra of V over S.
  Define T : V → W by

$$T = \begin{bmatrix} [0,a_1] & [0,a_2] \\ [0,a_3] & [0,a_4] \end{bmatrix} = \begin{bmatrix} [0,a_1] & [0,a_2] \\ 0 & [0,a_3] \end{bmatrix}$$

T is a projection of V into W. Infact T is a linear operator on V.

*Example 3.2.31*: Let

$$V = \left\{ \sum_{i=0}^{\infty} [0, a_i] x^i \middle| 0 \le a_i \le 1 \right\}$$

be a semigroup fuzzy interval linear algebra over a semigroup S of type II.

Let
$$W = \left\{ \sum_{i=0}^{\infty} [0, a_i] x^{2i} \,\middle|\, 0 \leq a_i \leq 1 \right\} \subseteq V$$

be a semigroup fuzzy interval linear subalgebra of type II over a semigroup S of V.

Define $T : V \to W$ as
$$T\left( \sum_{i=0}^{\infty} [0, a_i] x^i \right) = \sum_{i=0}^{\infty} [0, a_i] x^{2i}$$

T is a linear operator on V which is a projection of V onto W.

We can as in case of linear transformation of vector spaces (linear algebras) define linear transformations of fuzzy semigroup interval linear algebras V (vector spaces) of type II to fuzzy semigroup interval linear algebras W (vector spaces) of type II provided V and W are defined over the same semigroup.

We will illustrate this situation by some examples.

*Example 3.2.32*: Let

$$V = \left\{ \begin{bmatrix} [0, a_1] \\ [0, a_2] \\ [0, a_3] \\ [0, a_4] \end{bmatrix} \,\middle|\, \begin{array}{l} 0 \leq a_i \leq 1 \\ 1 \leq i \leq 4 \end{array} \right\}$$

and

$$W = \left\{ \begin{bmatrix} [0, a_1] & [0, a_2] & [0, a_3] \\ [0, a_4] & [0, a_5] & [0, a_6] \end{bmatrix} \,\middle|\, \begin{array}{l} 0 \leq a_i \leq 1 \\ 1 \leq i \leq 6 \end{array} \right\}$$

be fuzzy semigroup interval linear algebras defined over the semigroup

$$S = \left\{0, 1, \frac{1}{2^n}, \frac{1}{10^m} \,\middle|\, m, n \in Z^+ \right\}.$$

Define $T : V \to W$ as follows:

$$T\left(\begin{bmatrix} [0,a_1] \\ [0,a_2] \\ [0,a_3] \\ [0,a_4] \end{bmatrix}\right) = \begin{bmatrix} [0,a_1] & 0 & [0,a_2] \\ 0 & [0,a_3] & [0,a_4] \end{bmatrix}.$$

It is easily verified that T is a linear transformation of V into W.

*Example 3.2.33*: Let

$$V = \left\{ \begin{bmatrix} [0,a_1] \\ [0,a_2] \end{bmatrix}, \begin{bmatrix} [0,a_1] & [0,a_2] & [0,a_3] \\ [0,a_4] & [0,a_5] & [0,a_6] \\ [0,a_7] & [0,a_8] & [0,a_9] \end{bmatrix} \,\middle|\, \begin{array}{l} 0 \le a_i \le 1 \\ 1 \le i \le 9 \end{array} \right\}$$

and

$$W = \left\{ \begin{bmatrix} [0,a_1] & [0,a_2] \\ [0,a_3] & [0,a_4] \end{bmatrix}, \begin{bmatrix} [0,a_1] & [0,a_2] \\ [0,a_3] & [0,a_4] \\ [0,a_5] & [0,a_6] \\ [0,a_7] & [0,a_8] \\ [0,a_9] & [0,a_{10}] \\ [0,a_{11}] & [0,a_{12}] \\ [0,a_{13}] & [0,a_{14}] \\ [0,a_{15}] & [0,a_{16}] \end{bmatrix} \,\middle|\, \begin{array}{l} 0 \le a_i \le 1; \\ 1 \le i \le 16 \end{array} \right\}$$

be fuzzy semigroup interval vector spaces defined over the semigroup

$$S = \left\{ \frac{1}{2^n}, 0, 1 \,\middle|\, n \in 32Z^+ \right\}$$

of type II.

Define $T: V \to W$ as follows

$$T\left(\begin{bmatrix}[0,a_1]\\ [0,a_2]\end{bmatrix}\right) = \begin{bmatrix}[0,a_1] & 0\\ 0 & [0,a_2]\end{bmatrix}$$

and

$$T\left(\begin{bmatrix}[0,a_1] & [0,a_2] & [0,a_3]\\ [0,a_4] & [0,a_5] & [0,a_6]\\ [0,a_7] & [0,a_8] & [0,a_9]\end{bmatrix}\right) = \begin{bmatrix}[0,a_1] & 0\\ 0 & [0,a_2]\\ [0,a_3] & 0\\ 0 & [0,a_4]\\ [0,a_5] & 0\\ 0 & [0,a_6]\\ [0,a_7] & 0\\ 0 & [0,a_8]\end{bmatrix}.$$

It is easily verified that T is a linear transformation of V to W.

Now we proceed onto define some more properties of semigroup fuzzy interval vector spaces and linear algebras of type II.

**DEFINITION 3.2.12**: *Let V be a fuzzy semigroup interval vector space of type II defined over the semigroup S. Let $W_1$, $W_2$, $W_3$, ..., $W_n$ be a semigroup interval subvector spaces of V over the semigroup S. If $V = \bigcup_{i=1}^{n} W_i$ but $W_i \cap W_j \neq \phi$ or {0} if $i \neq j$ then we call V to be the pseudo direct union of fuzzy semigroup vector spaces of V over semigroup S of type II.*

The reader is expected to give some examples of these vector spaces.

**DEFINITION 3.2.13**: *Let V be a fuzzy semigroup interval vector space of type II over the semigroup S. Let $W_1$, $W_2$, ..., $W_n$ be fuzzy semigroup interval vector subspaces of V of type II. We say $W_1$, $W_2$, ..., $W_n$ is a direct union of semigroup fuzzy interval* *vector subspaces of V. if* $V = \bigcup_{i=1}^{n} W_i$ *and* $W_i \cap W_j = \phi$ *or {0} if i* $\neq j;\ 1 \leq i, j \leq n.$

The reader is expected to give examples of direct union of semigroup fuzzy interval vector subspaces of type II.

Now we proceed onto define direct sum of fuzzy semigroup interval linear subalgebras of a fuzzy interval semigroup of type II.

**DEFINITION 3.2.14**: *Let V be a fuzzy semigroup interval linear algebra over a semigroup S of type II. We say V is a direct sum of semigroup fuzzy interval linear subalgebras $W_1, W_2, …, W_n$ of V if*

(a) $V = W_1 + … + W_n$
(b) $W_i \cap W_j = \{0\}$ *or* $\phi$ *if* $i \neq j$ ; $1 \leq j, j \leq n.$

We will illustrate this situation by an example.

*Example 3.2.34*: Let

$$V = \left\{ \begin{bmatrix} [0,a_1] & [0,a_2] \\ [0,a_3] & [0,a_4] \end{bmatrix} \middle| 0 \leq a_i \leq 1; 1 \leq i \leq 4 \right\}$$

be a fuzzy semigroup interval linear algebra of type II defined over the semigroup

$$S = \left\{ 0, 1, \frac{1}{5^n} \middle| n \in Z^+ \right\}.$$

Choose

$$W_1 = \left\{ \begin{bmatrix} [0,a_i] & 0 \\ 0 & 0 \end{bmatrix} \middle| 0 \leq a_i \leq 1 \right\}, \quad W_2 = \left\{ \begin{bmatrix} 0 & [0,a_i] \\ 0 & 0 \end{bmatrix} \middle| 0 \leq a_i \leq 1 \right\},$$

$$W_3 = \left\{ \begin{bmatrix} 0 & 0 \\ [0,a_i] & 0 \end{bmatrix} \middle| 0 \le a_i \le 1 \right\}$$

and

$$W_4 = \left\{ \begin{bmatrix} 0 & 0 \\ 0 & [0,a_i] \end{bmatrix} \middle| 0 \le a_i \le 1 \right\}$$

to be fuzzy semigroup interval linear subalgebras of V of type II over the semigroup S.

$$V = W_1 + W_2 + W_3 + W_4$$

and

$$W_i \cap W_j = \begin{pmatrix} 0 & 0 \\ 0 & 0 \end{pmatrix}$$

if $i \ne j$ and $1 \le i, j \le n$.

If in the definition we have $W_i$'s to be such that $W_i \cap W_j \ne (0)$ or $\phi$ and $W_i \subseteq W_j$; $1 \le i, j \le n$ then we define V to be a pseudo direct sum of fuzzy semigroup interval linear algebras.

We will illustrate this by an example.

*Example 3.2.35*: Let

$$V = \left\{ \begin{bmatrix} [0,a_1] & [0,a_2] \\ [0,a_3] & [0,a_4] \\ [0,a_5] & [0,a_6] \\ [0,a_7] & [0,a_8] \\ [0,a_9] & [0,a_{10}] \end{bmatrix} \middle| 0 \le a_i \le 1; 1 \le i \le 10 \right\}$$

be fuzzy semigroup interval linear algebra of type II over the semigroup

$$S = \left\{ 0, 1, \frac{1}{2^n} \middle| n \in Z^+ \right\}.$$

Choose

$$W_1 = \left\{ \begin{bmatrix} [0,a_1] & 0 \\ 0 & [0,a_2] \\ 0 & 0 \\ 0 & 0 \\ 0 & [0,a_3] \end{bmatrix} \middle| 0 \leq a_i \leq 1; 1 \leq i \leq 3 \right\},$$

$$W_2 = \left\{ \begin{bmatrix} [0,a_1] & 0 \\ 0 & 0 \\ 0 & 0 \\ [0,a_2] & [0,a_3] \\ [0,a_4] & 0 \end{bmatrix} \middle| 0 \leq a_i \leq 1; 1 \leq i \leq 4 \right\},$$

$$W_3 = \left\{ \begin{bmatrix} 0 & [0,a_1] \\ 0 & [0,a_3] \\ [0,a_2] & 0 \\ 0 & [0,a_4] \\ 0 & 0 \end{bmatrix} \middle| 0 \leq a_i \leq 1; 1 \leq i \leq 4 \right\},$$

$$W_4 = \left\{ \begin{bmatrix} [0,a_1] & 0 \\ [0,a_2] & [0,a_3] \\ [0,a_4] & [0,a_5] \\ 0 & 0 \\ [0,a_6] & [0,a_7] \end{bmatrix} \middle| 0 \leq a_i \leq 1; 1 \leq i \leq 7 \right\}$$

and

$$W_5 = \left\{ \begin{bmatrix} 0 & 0 \\ [0,a_1] & 0 \\ 0 & [0,a_2] \\ [0,a_3] & [0,a_4] \\ [0,a_5] & [0,a_6] \end{bmatrix} \middle| 0 \leq a_i \leq 1; 1 \leq i \leq 6 \right\}.$$

We see $V = W_1 + W_2 + W_3 + W_4 + W_5$

But

$$W_i \cap W_j \neq \begin{bmatrix} 0 & 0 \\ 0 & 0 \\ 0 & 0 \\ 0 & 0 \\ 0 & 0 \end{bmatrix}$$

if $i \neq j$. $1 \leq i, j \leq 5$. Thus V is a pseudo direct sum $W_1, \ldots, W_5$.

As it is not an easy task to define group fuzzy interval vector spaces, we proceed to work in different direction.

**Chapter Four**

# SET INTERVAL BIVECTOR SPACES AND THEIR GENERALIZATION

In this chapter we for the first time introduce the notion of set interval bivector spaces and generalize them to set interval n-vector spaces, n ≥ 3. We also define semigroup interval bivector spaces and group interval bivector spaces and generalize both these concepts to bisemigroup interval bivector spaces, bigroup interval bivector spaces, set group interval bivector spaces and so on. This chapter has four sections.

### 4.1 Set Interval Bivector Spaces and Their Properties

In this section we introduce the new notion of set interval bivector spaces and enumerate a few of their properties.

**DEFINITION 4.1.1**: *Let $V = V_1 \cup V_2$ where $V_1$ and $V_2$ are two distinct set interval vector spaces defined over the same set S. That is $V_1 \not\subset V_2$ and $V_2 \not\subset V_1$ we may have $V_1 \cap V_2 = \phi$ or non empty. We define V to be a set interval bivector space over the set S.*

We will illustrate this situation by some examples.

*Example 4.1.1*: Let $V = V_1 \cup V_2$

$$= \{[0, a_i] \mid a_i \in Z^+ \cup \{0\}\} \cup \left\{ \begin{bmatrix} [0,a_1] \\ [0,a_2] \\ [0,a_3] \end{bmatrix} \middle| \begin{array}{c} a_i \in Z^+ \cup \{0\} \\ 1 \le i \le 3 \end{array} \right\}$$

be set interval bivector space over the set $S = \{2, 4, 3, 5, 10, 12, 124, 149, 5021\}$.

*Example 4.1.2*: Let

$$V_1 \cup V_2 = \left\{ \begin{bmatrix} [0,a_1] & [0,a_2] & [0,a_3] \\ [0,a_4] & [0,a_5] & [0,a_6] \\ [0,a_7] & [0,a_8] & [0,a_9] \end{bmatrix}, [0, a_1], [0, a_2] \right.$$

where $a_i \in Z_{12}; 1 \le i \le 9\}$

$$\cup \left\{ \begin{bmatrix} [0,a_1] & [0,a_2] \\ [0,a_3] & [0,a_4] \end{bmatrix}, \begin{bmatrix} [0,a_1] \\ [0,a_2] \\ [0,a_3] \\ [0,a_4] \\ [0,a_5] \end{bmatrix} \middle| a_i \in Z_{12}; 1 \le i \le 5 \right\}$$

be a set interval bivector space over the set $S = \{0, 2, 6, 5, 8, 11\} \subseteq Z_{12}$.

*Example 4.1.3*: Let

$$V = V_1 \cup V_2 =$$

$$\left\{ \begin{bmatrix} [0,a_1] & [0,a_2] & [0,a_3] & [0,a_4] \\ [0,a_5] & [0,a_6] & [0,a_7] & [0,a_8] \end{bmatrix}, \begin{bmatrix} [0,a_1] \\ [0,a_2] \\ [0,a_3] \\ [0,a_4] \\ [0,a_5] \end{bmatrix}, [0,a_1] \;\middle|\; \begin{array}{l} a_i \in Q^+ \cup \{0\}; \\ 1 \le i \le 8 \end{array} \right\}$$

$$\cup \left\{ \begin{bmatrix} [0,a_1] & [0,a_2] & [0,a_3] & [0,a_4] \\ [0,a_5] & [0,a_6] & [0,a_7] & [0,a_8] \\ [0,a_9] & [0,a_{10}] & [0,a_{11}] & [0,a_{12}] \\ [0,a_{13}] & [0,a_{14}] & [0,a_{15}] & [0,a_{16}] \end{bmatrix} \;\middle|\; \begin{array}{l} a_i \in Q^+ \cup \{0\}; \\ 1 \le i \le 16 \end{array} \right\}$$

be a set interval bivector space defined over the set $S = \{0, \frac{1}{2}, 3/17, 25/4, 2, 4, 6, 21, 49\}$.

*Example 4.1.4*: Let $V = V_1 \cup V_2 = \{$All $10 \times 10$ interval matrices with intervals of the from $[0, a_i]$ with $a_i \in Z_7\} \cup$

$$\left\{ \sum_{i=0}^{\infty} [0,a_i] x^i \;\middle|\; a_i \in Z_7 \right\}$$

be a set interval bivector space over the set $S = \{0, 3, 5, 1\} \subseteq Z_7$.

We see examples 4.1.1, 4.1.3 and 4.1.4 give set interval bivector spaces of infinite order where as example 4.1.2 is of finite order.

Now we can define substructure in them.

**DEFINITION 4.1.2**: *Let $V = V_1 \cup V_2$ be a set interval bivector space over the set S. Suppose $W = W_1 \cup W_2 \subseteq V_1 \cup V_2$ be a proper biset of V and if $W = W_1 \cup W_2 \subseteq V$ is a set interval* *bivector space over the set S then we define W to be a set interval bivector subspace of V over the set S.*

We will illustrate this situation by some examples.

***Example 4.1.5***: Let $V = V_1 \cup V_2 =$

$$\left\{ \begin{bmatrix} [0,a_1] & [0,a_2] & [0,a_3] \\ [0,a_4] & [0,a_5] & [0,a_6] \end{bmatrix}, \begin{bmatrix} [0,a_1] \\ [0,a_2] \\ [0,a_3] \end{bmatrix} \middle| \begin{array}{l} a_i \in Z_{19}; \\ 1 \leq i \leq 6 \end{array} \right\} \cup$$

$$\left\{ \sum_{i=0}^{25} [0,a_i] x^i \middle| a_i \in Z_{19} \right\}$$

be a set interval bivector space over the set $S = \{0, 2, 5, 7, 9, 12, 17\} \subseteq Z_{19}$.

Choose

$$W = \left\{ \begin{bmatrix} [0,a_1] \\ [0,a_2] \\ [0,a_3] \end{bmatrix} \middle| \begin{array}{l} a_i \in Z_{19}; \\ 1 \leq i \leq 3 \end{array} \right\} \cup \left\{ \sum_{i=0}^{25} [0,a_i] x^{2i} \middle| a_i \in Z_{19} \right\}$$

$$= W_1 \cup W_2 \subseteq V_1 \cup V_2 = V$$

is a set interval bivector subspace of V over the set S.

***Example 4.1.6***: Let $V = V_1 \cup V_2 = \{[0, a_i] \mid a_i \in Z^+ \cup \{0\}\} \cup$

$$\left\{ \begin{bmatrix} [0,a_1] & [0,a_2] & [0,a_3] \\ [0,a_4] & [0,a_5] & [0,a_6] \\ [0,a_7] & [0,a_8] & [0,a_9] \end{bmatrix} \middle| a_i \in Z^+ \cup \{0\}; 1 \leq i \leq 9 \right\}$$

be a set interval bivector space over the set $S = 2Z^+ \cup 5Z^+ \cup \{0\}$. Take

$$W = W_1 \cup W_2 = \{[0, a_i] \mid a_i \in 7Z^+ \cup \{0\}\} \cup$$

$$\left\{ \begin{bmatrix} [0,a_1] & [0,a_2] & [0,a_3] \\ 0 & [0,a_4] & [0,a_5] \\ 0 & 0 & [0,a_6] \end{bmatrix} \middle| a_i \in 3Z^+ \cup \{0\} \right\}$$

$\subseteq V_1 \cup V_2$; $W = W_1 \cup W_2$ is a set interval bivector subspace of V over the set S.

*Example 4.1.7*: Let $V = V_1 \cup V_2 = \{$All $7 \times 7$ interval matrices with interval of the form $[0, a_i]$ $a_i \in Z_{18}\} \cup \{$All $1 \times 9$ interval row matrices with intervals of the form with $a_i \in Z_{18}$ be a set interval bivector space over the set $S = \{0, 1, 2, 4, 5, 7\} \subseteq Z_{18}$. We see $W = W_1 \cup W_2 = \{$All $7 \times 7$ diagonal interval matrices with intervals of the form $[0, a_i]$ with $a_i$ from $Z_{18}\} \cup \{([0, a_1], 0, [0, a_2], 0, [0, a_3], 0, [0, a_4], 0, [0, a_5]) / a_i \in Z_{18}; 1 \leq i \leq 5\} \subseteq V_1 \cup V_2 = V$ is a set interval bivector subspace of V over the set S.

Now having see examples of subspaces we now proceed on to define subset interval bivector subspaces.

**DEFINITION 4.1.3**: *Let $V = V_1 \cup V_2$ be a set interval bivector space over the set S. Let $W = W_1 \cup W_2 \subseteq V_1 \cup V_2 = V$ be a proper bisubset of V and $P \subseteq S$ be a proper subset of S. If W is a set interval bivector space over the set P then we define W to be a subset interval bivector subspace of V over the subset P of S.*

We will illustrate this by some simple examples.

*Example 4.1.8*: Let $V = V_1 \cup V_2 =$

$$\left\{ \begin{bmatrix} [0,a_1] \\ [0,a_2] \\ [0,a_3] \end{bmatrix}, ([0,a_1] \quad [0,a_2] \quad [0,a_3] \quad [0,a_4]) \middle| \begin{matrix} a_i \in Z_{24}; \\ 1 \leq i \leq 4 \end{matrix} \right\} \cup$$

{All 5 × 5 interval matrices with entries from $Z_{24}$} be a set interval bivector space over the set S = {0, 2, 3, 5, 6, 8, 9, 10, 14, 22} $\subseteq Z_{24}$. Choose

$$W = W_1 \cup W_2 = \left\{ \begin{bmatrix} [0,a_1] \\ [0,a_2] \\ [0,a_3] \end{bmatrix} \, \middle| \, \begin{array}{l} a_i \in Z_{24}; \\ 1 \le i \le 3 \end{array} \right\} \cup$$

{All 5 × 5 interval upper triangular matrices with entries from $Z_{24}$} $\subseteq V_1 \cup V_2 = V$. Choose P = {0, 2, 5, 10, 14, 22} $\subseteq S \subseteq Z_{24}$. $W = W_1 \cup W_2$ is a subset interval bivector subspace of V over the subset P of S.

Now we proceed onto define the notion of set interval linear bialgebra.

**DEFINITION 4.1.4:** *Let $V = V_1 \cup V_2$ be a set interval bivector space over the set S.*
*Suppose V is closed under addition and if $s(a + b) = sa + sb$ for all $s \in S$ and $a, b \in V$ then we call V to be a set interval bilinear algebra over S.*

We will illustrate this situation by some examples.

*Example 4.1.9*: Let $V = V_1 \cup V_2$ be a set interval bilinear algebra over the set S; where

$$V = V_1 \cup V_2$$
$$= \left\{ \begin{bmatrix} [0,a_1] & [0,a_2] \\ [0,a_3] & [0,a_4] \end{bmatrix} \, \middle| \, \begin{array}{l} a_i \in Z_{40}; \\ 1 \le i \le 4 \end{array} \right\} \cup \left\{ \sum_{i=0}^{\infty} [0,a_i] x^i \, \middle| \, a_i \in Z_{40} \right\}$$

and
$$S = \{0, 2, 5, 7, 10, 14, 32\} \subseteq Z_{40}.$$

*Example 4.1.10*: Let
$$V = V_1 \cup V_2$$

$$= \left\{ \begin{bmatrix} [0,a_1] \\ [0,a_2] \\ [0,a_3] \\ [0,a_4] \\ [0,a_5] \end{bmatrix} \middle| \begin{array}{l} a_i \in Z^+ \cup \{0\}; \\ 1 \le i \le 5 \end{array} \right\} \cup$$

$$\left\{ \begin{bmatrix} [0,a_1] & [0,a_2] & [0,a_3] \\ [0,a_4] & [0,a_5] & [0,a_6] \\ [0,a_7] & [0,a_8] & [0,a_9] \end{bmatrix} \middle| \begin{array}{l} a_i \in Z^+ \cup \{0\}; \\ 1 \le i \le 9 \end{array} \right\}$$

be a set interval bilinear algebra over the set $S = \{0, 3, 14, 27, 52, 75, 130\} \subseteq Z^+ \cup \{0\}$.

**Example 4.1.11**: Let $V = V_1 \cup V_2$

$$= \left\{ \begin{bmatrix} [0,a_1] & [0,a_2] & [0,a_3] & [0,a_4] & [0,a_5] \\ [0,a_6] & [0,a_7] & [0,a_8] & [0,a_9] & [0,a_{10}] \end{bmatrix} \middle| \begin{array}{l} a_i \in R^+ \cup \{0\} \\ 1 \le i \le 10 \end{array} \right\} \cup$$

{All $12 \times 11$ interval matrices with intervals from $Q^+ \cup \{0\}$ of the form $[0, a_i]$; $a_i \in Q^+ \cup \{0\}\}$ be a set interval bilinear algebra over the set $S = \{0, 12, \sqrt{3}, \sqrt{41}, \sqrt{5}/12, 412, \dfrac{\sqrt{53}}{150}\}$.

Now we see all the set interval bilinear algebras given in the examples 4.1.9, 4.1.10 and 4.1.11 are of infinite order.

We will give one example of a finite set interval bilinear algebra.

**Example 4.1.12**: Let

$$V = \{[0, a_i] \mid a_i \in Z_8\} \cup \left\{ \begin{bmatrix} [0,a_1] & [0,a_2] \\ [0,a_3] & [0,a_4] \end{bmatrix} \middle| \begin{array}{l} a_i \in Z_8; \\ 1 \le i \le 4 \end{array} \right\}$$

be a set interval linear bialgebra over the set $S = \{0, 1, 2, 3, 4\} \subseteq Z_8$. V is a finite order set interval linear bialgebra or finite order set interval bilinear algebra over the set S.

Now we proceed onto define the notion of set interval bilinear subalgebra of a set interval bilinear algebra over the set S.

**DEFINITION 4.1.5**: *Let $V = V_1 \cup V_2$ be a set interval linear bialgebra over the set S. Choose $W = W_1 \cup W_2 \subseteq V_1 \cup V_2 = V$; suppose W is a set interval linear bialgebra over the set S then we call W to be a set interval linear sub bialgebra of V over the set S.*

We will illustrate this situation by some examples.

*Example 4.1.13*: Let

$$V = V_1 \cup V_2 = \left\{ \begin{bmatrix} [0,a_1] & [0,a_2] & [0,a_3] \\ [0,a_4] & [0,a_5] & [0,a_6] \\ [0,a_7] & [0,a_8] & [0,a_9] \end{bmatrix} \middle| \begin{array}{c} a_i \in Z_{16}; \\ 1 \le i \le 9 \end{array} \right\} \cup$$

$$\left\{ \begin{bmatrix} [0,a_1] & [0,a_2] \\ [0,a_3] & [0,a_4] \\ [0,a_5] & [0,a_6] \\ [0,a_7] & [0,a_8] \\ [0,a_9] & [0,a_{10}] \end{bmatrix} \middle| \begin{array}{c} a_i \in Z_{16}; \\ 1 \le i \le 10 \end{array} \right\}$$

be a set interval linear bialgebra over the set $S = \{0, 3, 5, 8, 7, 10\} \subseteq Z_{16}$. Choose $W = W_1 \cup W_2$

$$= \left\{ \begin{bmatrix} [0,a_1] & [0,a_2] & [0,a_3] \\ 0 & 0 & 0 \\ [0,a_7] & [0,a_8] & [0,a_9] \end{bmatrix} \middle| \begin{array}{c} a_i \in Z_{16}; \\ i = 1,2,3,7,8,9 \end{array} \right\}$$

$$\cup \left\{ \begin{bmatrix} [0,a_1] & 0 \\ 0 & [0,a_2] \\ [0,a_3] & 0 \\ 0 & [0,a_4] \\ [0,a_5] & 0 \end{bmatrix} \,\middle|\, \begin{array}{l} a_i \in Z_{16}; \\ 1 \leq i \leq 5 \end{array} \right\}$$

$\subseteq V_1 \cup V_2 = V$; W is a set interval linear subbialgebra of V over the set S.

*Example 4.1.14*: Let

$$V = V_1 \cup V_2 = \left\{ \sum_{i=0}^{27} [0,a_i] x^i \,\middle|\, a_i \in Q^+ \cup \{0\} \right\} \cup$$

{All $10 \times 10$ interval matrices with entries from $Q^+ \cup \{0\}$} be a set interval linear bialgebra over the set $S = \{3Z^+, 0, 7Z^+\}$. Choose

$$W = W_1 \cup W_2 = \left\{ \sum_{i=0}^{20} [0,a_i] x^i \,\middle|\, a_i \in Z^+ \cup \{0\} \right\} \cup$$

{all $10 \times 10$ upper triangular interval matrices with entries from $Q^+ \cup \{0\}$} $\subseteq V_1 \cup V_2 = V$; W is a set interval linear bisubalgebra of V over the set S.

Now we proceed onto define other special type of substructures.

**DEFINITION 4.1.6**: *Let $V = V_1 \cup V_2$ be a set interval bilinear algebra over the set S.*

*Choose $W = W_1 \cup W_2 \subseteq V_1 \cup V_2$ and P a proper subset of S such that $W = W_1 \cup W_2$ is a set interval bilinear algebra over P. W is defined as a subset interval linear subbialgebra of V over the subset P of S. If V has no subset interval linear subalgebra then we define V to be a pseudo simple set linear bialgebra.*

First we will illustrate this by some simple examples.

*Example 4.1.15*: Let
$$V = V_1 \cup V_2$$
$$= \left\{ \begin{bmatrix} [0,a_1] & [0,a_2] \\ [0,a_3] & [0,a_4] \end{bmatrix} \middle| \begin{array}{l} a_i \in Z_{27}; \\ 1 \leq i \leq 4 \end{array} \right\} \cup \{[0, a_i] \mid a_i \in Z_{27}\}$$

be a set interval linear bialgebra over the set $S = Z_{27}$. Let

$$W = W_1 \cup W_2 = \left\{ \begin{bmatrix} [0,a_1] & 0 \\ [0,a_2] & 0 \end{bmatrix} \middle| \begin{array}{l} a_i \in Z_{27}; \\ 1 \leq i \leq 2 \end{array} \right\} \cup$$

$\{[0, a_i] \mid a_i \in \{0, 3, 6, 9, 12, 15, 18, 21, 24\} \subseteq Z_{27}\} \subseteq V_1 \cup V_2 = V$ be a subset interval linear subbialgebra of V over subset $P = \{0, 9, 18\} \subseteq S$.

*Example 4.1.16*: Let $V = V_1 \cup V_2 = \{$All $5 \times 5$ interval matrices with entries from $Q^+ \cup \{0\}\} \cup$

$$\left\{ \sum_{i=0}^{30}[0, a_i]x^i \middle| \begin{array}{l} a_i \in Q^+ \cup \{0\}; \\ 0 \leq i \leq 30 \end{array} \right\}$$

be a set interval linear bialgebra over the set $S = \{0, 3Z^+, 11Z^+, 17Z^+\}$. Choose $W = W_1 \cup W_2 = \{$all $5 \times 5$ interval upper matrices with entries from $Q^+ \cup \{0\}\} \cup$

$$\left\{ \sum_{i=0}^{20}[0, a_i]x^i \middle| \begin{array}{l} a_i \in Q^+ \cup \{0\}; \\ 0 \leq i \leq 20 \end{array} \right\}$$

$\subseteq V_1 \cup V_2 = V$, W is a subset interval linear subbialgebra over the subset $P = \{0, 3Z^+, 17Z^+\} \subseteq S$.

Now will give some examples of a pseudo simple set interval linear bialgebra.

***Example 4.1.17***: Let $V = V_1 \cup V_2 = \{[0, a], [0, a], [0, a], [0, a]) \mid a \in Z_5\} \cup$

$$\left\{ \begin{bmatrix} [0,a] \\ [0,a] \\ [0,a] \\ [0,a] \end{bmatrix} \middle| a \in Z_5 \right\}$$

be a set linear bialgebra of the set $S = \{0, 1\}$.
   Clearly V is a pseudo simple set linear bialgebra over S.

***Example 4.1.18***: Let

$$V = V_1 \cup V_2$$

$$= \left\{ \begin{bmatrix} [0,a] & [0,a] & [0,a] \\ [0,a] & [0,a] & [0,a] \\ [0,a] & [0,a] & [0,a] \end{bmatrix} \middle| a \in Z_3 \right\} \cup \left\{ \begin{bmatrix} [0,a] \\ [0,a] \\ [0,a] \end{bmatrix} \middle| a \in Z_3 \right\}$$

be a set interval linear bialgebra over the set $S = \{0, 1\}$. Clearly V is a pseudo simple set interval linear bialgebra over the set S.

We define pseudo set interval bivector space of a set interval linear bialgebra.

***Example 4.1.19***: Let $V = V_1 \cup V_2 =$

$$\left\{ \begin{bmatrix} [0,a_1] & [0,a_2] & [0,a_3] \\ [0,a_4] & [0,a_5] & [0,a_6] \end{bmatrix} \middle| \begin{array}{l} a_i \in Z_7; \\ 0 \le i \le 6 \end{array} \right\} \cup$$

$$\left\{ \begin{bmatrix} [0,a_1] & [0,a_2] \\ [0,a_3] & [0,a_4] \\ [0,a_5] & [0,a_6] \end{bmatrix} \middle| \begin{array}{l} a_i \in Z_7; \\ 0 \le i \le 6 \end{array} \right\}$$

be a set interval linear bialgebra of V over the set $S = \{0, 1, 2, 3\}$. Choose

$$W = W_1 \cup W_2$$

$$= \left\{ \begin{bmatrix} [0,a_1] & 0 & [0,a_2] \\ 0 & [0,a_3] & [0,a_4] \end{bmatrix}, \begin{bmatrix} [0,a_1] & 0 & 0 \\ [0,a_2] & 0 & 0 \end{bmatrix} \right\} \cup$$

$$\left\{ \begin{bmatrix} [0,a_1] & 0 \\ 0 & [0,a_2] \\ [0,a_3] & 0 \end{bmatrix}, \begin{bmatrix} 0 & [0,a_1] \\ [0,a_2] & 0 \\ 0 & [0,a_3] \end{bmatrix} \right\}$$

$\subseteq V_1 \cup V_2 = V$, W is a pseudo set interval bivector subspace of V over the set S.

*Example 4.1.20:* Let $V = V_1 \cup V_2$

$= \{([0, a_1], [0, a_2], [0, a_3], [0, a_4], [0, a_5]) \mid a_i \in Z_7; 1 \leq i \leq 5\} \cup$

$$\left\{ \begin{bmatrix} [0,a] & [0,a] \\ [0,a] & [0,a] \end{bmatrix} \Big| a \in Z_7 \right\}$$

be a set interval bilinear algebra over the set $S = \{0, 1\}$.

Choose $W = W_1 \cup W_2 =$

$\{([0, a], 0, [0, a], 0, [0, a]), ([0, a], [0, a], 0, [0, a], 0) \mid a \in Z_7\} \cup$

$$\left\{ \begin{bmatrix} [0,a] & 0 \\ [0,a] & 0 \end{bmatrix}, \begin{bmatrix} 0 & [0,a] \\ 0 & [0,a] \end{bmatrix} \, a \in Z_7 \right\}$$

$\subseteq V_1 \cup V_2 = V$, W is a pseudo set interval vector bisubspace of V over the set S.

As in case of usual set bivector spaces we define the bigenerations of set interval bivector spaces.

We will illustrate this by an simple example.

***Example 4.1.21***: Let $V = V_1 \cup V_2 = \{\{[0, a_i] \mid a_i \in Z_7\} \cup$

$$\left\{ \begin{bmatrix} [0,a_1] \\ [0,a_1] \end{bmatrix} \middle| a_i \in Z_7 \right\}$$

be a set interval bivector space over the set $S = Z_7$. The bigenerator of V is

$$X = \{[0, 1]\} \cup \left\{ \begin{bmatrix} [0,1] \\ [0,1] \end{bmatrix} \right\}.$$

Clearly the bidimension of V is finite and is (0, 1).

Interested reader is expected to derive more properties. However the concept of bilinearly independent set and bibasis can also be defined in an analogous way. We see the basis of the set interval bivector space given in example 4.1.21 is

$$\{[0, 1]\} \cup \begin{bmatrix} [0,1] \\ [0,1] \end{bmatrix}.$$

The bidimension is {1, 1}.

We will illustrate this by another example.

***Example 4.1.22***: Let
$$V = V_1 \cup V_2$$

$$= \left\{ \begin{bmatrix} [0,a_1] & [0,a_2] \\ [0,a_3] & [0,a_4] \\ [0,a_5] & [0,a_6] \end{bmatrix} \middle| a_i \in Z^+ \cup \{0\} \right\}$$

$$\cup \left\{ \sum_{i=0}^{6} [0,a_i] x^i \middle| a_i \in Z^+ \cup \{0\} \right\}$$

be a set interval bilinear algebra over the set $S = Z^+ \cup \{0\}$.

$$X = \left\{ \begin{bmatrix} [0,1] & 0 \\ 0 & 0 \\ 0 & 0 \end{bmatrix}, \begin{bmatrix} 0 & [0,1] \\ 0 & 0 \\ 0 & 0 \end{bmatrix}, \begin{bmatrix} 0 & 0 \\ [0,1] & 0 \\ 0 & 0 \end{bmatrix}, \right.$$

$$\left. \begin{bmatrix} 0 & 0 \\ 0 & [0,1] \\ 0 & 0 \end{bmatrix}, \begin{bmatrix} 0 & 0 \\ 0 & 0 \\ [0,1] & 0 \end{bmatrix}, \begin{bmatrix} 0 & 0 \\ 0 & 0 \\ 0 & [0,1] \end{bmatrix} \right\}$$

$$\cup \{1, x, x^2, x^3, x^4, x^5, x^6\} = X_1 \cup X_2$$

is a bilinearly independent bisubset of V and $X = X_1 \cup X_2$ bigenerates V thus X is a bibasis of V.

We define yet another set of interval bivector spaces called biset interval bivector spaces.

**DEFINITION 4.1.7**: *Let $V = V_1 \cup V_2$ where $V_1$ is a set interval vector space over the set $S_1$ and $V_2$ is another set interval vector space over the set $S_2$ where $V_1$ and $V_2$ are distinct that is $V_1 \not\subseteq V_2$ and $V_2 \not\subseteq V_1$ and $S_1$ and $S_2$ are distinct that is $S_1 \not\subseteq S_2$ and $S_2 \not\subseteq S_1$.*

*Then we define $V = V_1 \cup V_2$ to be a biset interval vector bispace over the biset $S = S_1 \cup S_2$ or V is a biset interval bivector space over the biset $S = S_1 \cup S_2$.*

We will illustrate this situation by some simple examples.

*Example 4.1.23*: Let $V = V_1 \cup V_2$

$$= \left\{ \begin{bmatrix} [0,a_1] & [0,a_2] & [0,a_3] \\ [0,a_4] & [0,a_5] & [0,a_6] \end{bmatrix}, \begin{bmatrix} [0,a_1] \\ [0,a_2] \end{bmatrix} \middle| \begin{array}{l} a_i \in Z_{12}; \\ 1 \leq i \leq 6 \end{array} \right\}$$

$$\cup \left\{ \begin{bmatrix} [0,a_1] \\ [0,a_2] \\ [0,a_3] \\ [0,a_4] \end{bmatrix} \middle| \begin{array}{l} a_i \in Z_{42}; \\ 1 \leq i \leq 4 \end{array} \right\}$$

be a biset interval bivector space over the biset $S = S_1 \cup S_2 = Z_{12} \cup Z_{42}$.

***Example 4.1.24:*** Let $V = V_1 \cup V_2$

$$= \left\{ \begin{bmatrix} [0,a_1] \\ [0,a_2] \\ [0,a_3] \\ [0,a_4] \end{bmatrix} \middle| \begin{array}{c} a_i \in Z^+ \cup \{0\} \\ 1 \leq i \leq 4 \end{array} \right\} \cup \{[0, a_i] \mid a_i \in Z_7\}$$

be a biset interval bivector space over the biset $S = S_1 \cup S_2 = (Z^+ \cup \{0\}) \cup Z_7$.

***Example 4.1.25:*** Let $V = V_1 \cup V_2$

$$= \left\{ \sum_{i=0}^{24} [0,a_i] x^i \middle| a_i \in Z_{45} \right\} \cup$$

{all $10 \times 10$ interval matrices with entries from $Z^+ \cup \{0\}$} be a biset interval bivector space over the biset $S = S_1 \cup S_2 = Z_{45} \cup 3Z^+ \cup \{0\}$.

Now we proceed onto define substructure in this bispace.

**DEFINITION 4.1.8:** *Let $V = V_1 \cup V_2$ be a biset interval bivector space over the biset $S = S_1 \cup S_2$. Let $W = W_1 \cup W_2 \subseteq V_1 \cup V_2 = V$ be a proper subset of V.*

*If $W = W_1 \cup W_2$ is a biset interval bivector space over the biset $S = S_1 \cup S_2$, then we define W to be a biset interval bivector sub bispace of V over the biset S.*

We will first illustrate this situation by some examples.

***Example 4.1.26***: Let $V = V_1 \cup V_2 =$

$$\left\{ \begin{bmatrix} [0,a_1] & [0,a_2] \\ [0,a_3] & [0,a_4] \end{bmatrix} \middle| \begin{array}{c} a_i \in Q^+ \cup \{0\}; \\ 1 \leq i \leq 4 \end{array} \right\} \cup$$

{All 17 × 17 upper triangular interval matrices with entries from $Z_{12}$} be a biset interval bivector space over biset the $S = (Q^+ \cup \{0\}) \cup Z_{12} = S_1 \cup S_2$.

Take $W = W_1 \cup W_2 =$

$$\left\{ \begin{bmatrix} [0,a_1] & [0,a_2] \\ 0 & [0,a_3] \end{bmatrix} \middle| \begin{array}{l} a_i \in Q^+ \cup \{0\}; \\ 1 \leq i \leq 3 \end{array} \right\} \cup$$

{All 17 × 17 diagonal interval matrices with entries from $Z_{12}$} $\subseteq V_1 \cup V_2 = V$, W is a biset interval bivector subspace of V over the biset $S = S_1 \cup S_2$.

***Example 4.1.27***: Let $V = V_1 \cup V_2 =$

$$\left\{ \begin{bmatrix} [0,a_1] & [0,a_2] & [0,a_3] \\ [0,a_4] & [0,a_5] & [0,a_6] \end{bmatrix} \middle| \begin{array}{l} a_i \in Z_{42} \\ 1 \leq i \leq 6 \end{array} \right\} \cup$$

$$\left\{ \sum_{i=0}^{25} [0,a_i]x^i \middle| \begin{array}{l} a_i \in Z_{25}; \\ 1 \leq i \leq 6 \end{array} \right\}$$

be a biset interval bivector space over the biset $S = S_1 \cup S_2 = Z_{42} \cup Z_{25}$. Choose

$$W = W_1 \cup W_2$$
$$= \left\{ \begin{bmatrix} [0,a_1] & [0,a_2] & [0,a_3] \\ [0,a_1] & [0,a_2] & [0,a_3] \end{bmatrix} \middle| \begin{array}{l} a_i \in Z_{42}; \\ 1 \leq i \leq 3 \end{array} \right\} \cup \quad \subseteq Z_{25}\}$$

$\subseteq V_1 \cup V_2 = V$; W is a biset interval bivector subspace of V over S.

Now we proceed onto define the notion of quasi set interval bivector spaces and quasi biset interval bivector spaces.

**DEFINITION 4.1.9**: *Let $V = V_1 \cup V_2$ be such that $V_1$ is a set vector space over the set S and $V_2$ a set interval vector space over the same set S then we define $V = V_1 \cup V_2$ to be a quasi set interval bivector space over the set S.*

We will first illustrate this situation by some examples.

***Example 4.1.28***: Let $V = V_1 \cup V_2$

$$= \left\{ \left[ \begin{pmatrix} a_1 & a_2 & a_3 \\ a_4 & a_5 & a_6 \\ a_7 & a_8 & a_9 \end{pmatrix} \right] \middle| \begin{array}{l} a_i \in Q^+ \cup \{0\}; \\ 1 \le i \le 9 \end{array} \right\} \cup$$

$$\left\{ \begin{bmatrix} [0,a_1] & [0,a_2] & [0,a_3] & [0,a_4] \\ [0,a_5] & [0,a_6] & [0,a_7] & [0,a_8] \end{bmatrix} \middle| \begin{array}{l} a_i \in Q^+ \cup \{0\}; \\ 1 \le i \le 8 \end{array} \right\}$$

be a quasi set interval bivector space over the set $S = Z^+ \cup \{0\}$.

***Example 4.1.29***: Let $V = V_1 \cup V_2$

$$= \left\{ \begin{bmatrix} [0,a_1] & [0,a_5] \\ [0,a_2] & [0,a_6] \\ [0,a_3] & [0,a_7] \\ [0,a_4] & [0,a_8] \end{bmatrix} \middle| \begin{array}{l} a_i \in Z_8; \\ 1 \le i \le 8 \end{array} \right\} \cup \left\{ \sum_{i=0}^{26} [0,a_i] x^i \middle| a_i \in Z_8 \right\}$$

be a quasi set interval bivector space over the set $S = Z_8$.

***Example 4.1.30:*** Let $V = V_1 \cup V_2$

$$= \left\{ \sum_{i=0}^{\infty} [0,a_i] x^i \middle| a_i \in Z^+ \cup \{0\} \right\} \cup$$

{All $8 \times 8$ interval matrices with entries from $Z^+ \cup \{0\}$} be a quasi set interval bivector space over the set $S = 3Z^+ \cup \{0\}$.

Now we define quasi set interval bivector subspace in an analogous way.

**DEFINITION 4.1.10:** *Let $V = V_1 \cup V_2$ be a quasi set interval bivector space over the set S. Let $W = W_1 \cup W_2 \subseteq V_1 \cup V_2$ ; W* *is a quasi set interval bivector subspace of V over the set S, if W is a quasi set interval bivector spaces over the sets.*

For instance if $V = V_1 \cup V_2$

$$= \left\{\sum_{i=0}^{40}[0,a_i]x^i \,\Big|\, a_i \in Z_{28}\right\} \cup \left\{\begin{bmatrix}[0,a_1] & [0,a_2]\\ [0,a_3] & [0,a_4]\\ [0,a_5] & [0,a_6]\end{bmatrix} \,\Bigg|\, \begin{array}{l}a_i \in Z_{28};\\ 1 \le i \le 6\end{array}\right\}$$

be a quasi set interval bivector space over the set $S = Z_{28}$.
 Let $W = W_1 \cup W_2$

$$= \left\{\sum_{i=0}^{20} a_i x^i \,\Big|\, a_i \in Z_{28}\right\} \cup \left\{\begin{bmatrix}[0,a_1] & 0\\ 0 & [0,a_2]\\ [0,a_3] & 0\end{bmatrix} \,\Bigg|\, \begin{array}{l}a_i \in Z_{28};\\ 1 \le i \le 3\end{array}\right\}$$

$\subseteq V_1 \cup V_2$; W is a quasi set interval bivector subspace of V over the set S.

***Example 4.1.31***: Let $V = V_1 \cup V_2$

$$= \left\{\begin{bmatrix}[0,a_1] & [0,a_2] & [0,a_3] & [0,a_4] & [0,a_5]\\ [0,a_6] & [0,a_7] & [0,a_8] & [0,a_9] & [0,a_{10}]\end{bmatrix} \,\Bigg|\, \begin{array}{l}a_i \in Z^+ \cup \{0\};\\ 1 \le i \le 10\end{array}\right\} \cup$$

$\{10 \times 10$ upper triangular matrices with entries from $Z^+ \cup \{0\}\}$ be a quasi set interval bivector space over the set $S = 3Z^+ \cup \{0\}$.
 Choose $W = W_1 \cup W_2$

$$= \left\{\begin{bmatrix}[0,a_1] & 0 & [0,a_2] & 0 & [0,a_3]\\ 0 & [0,a_4] & 0 & [0,a_5] & 0\end{bmatrix} \,\Bigg|\, \begin{array}{l}a_i \in Z^+ \cup \{0\};\\ 1 \le i \le 5\end{array}\right\} \cup$$

$\{$All $10 \times 10$ diagonal matrices with entries from $13Z^+ \cup \{0\}\}$
$\subseteq V_1 \cup V_2 = V$; W is a quasi set interval bivector subspace of V over the set S.

Now we proceed onto define the new notion of quasi biset interval bivector space.

**DEFINITION 4.1.11**: *Let $V = V_1 \cup V_2$ where $V_1$ is a set vector space over the set $S_1$ and $V_2$ is a set interval vector space over the set $S_2$. We call $V = V_1 \cup V_2$ to be a quasi biset interval bivector space over the biset $S = S_1 \cup S_2$.*

We will illustrate this situation by some examples.

*Example 4.1.32*: Let $V = V_1 \cup V_2$

$$= \left\{ \left[ \begin{pmatrix} a_1 & a_2 & a_3 \\ a_4 & a_5 & a_6 \\ a_7 & a_8 & a_9 \end{pmatrix} \right] \middle| a_i \in Z_{18}; 1 \leq i \leq 9 \right\} \cup \left\{ \begin{bmatrix} [0,a_1] & [0,a_2] \\ [0,a_3] & [0,a_4] \\ [0,a_5] & [0,a_6] \\ [0,a_7] & [0,a_8] \end{bmatrix} \middle| a_i \in Z_{11}; 1 \leq i \leq 8 \right\}$$

be a quasi biset interval vector bispace over the biset $S = S_1 \cup S_2 = Z_{18} \cup Z_{11}$.

*Example 4.1.33*: Let $V = V_1 \cup V_2 = \{$all $12 \times 12$ matrices with entries from $Z^+ \cup \{0\}\} \cup$

$$\left\{ \sum_{i=0}^{\infty} [0,a_i] x^i \middle| a_i \in Z_{29} \right\}$$

be a quasi biset interval bivector space over the biset $S = S_1 \cup S_2 = (13Z^+ \cup \{0\}) \cup Z_{29}$.

*Example 4.1.34*: Let $V = V_1 \cup V_2 =$

$$\left\{ \sum_{i=0}^{15} a_i x^i \middle| a_i \in 3Z^+ \cup \{0\} \right\} \cup$$

$$\left\{ \begin{bmatrix} [0,a_1] & [0,a_2] & [0,a_3] \\ [0,a_4] & [0,a_5] & [0,a_6] \\ [0,a_7] & [0,a_8] & [0,a_9] \end{bmatrix} \middle| \begin{array}{c} a_i \in 5Z^+ \cup \{0\}; \\ 1 \le i \le 9 \end{array} \right\}$$

be a quasi biset interval bivector space over the biset $S = S_1 \cup S_2 = (13Z^+ \cup \{0\}) \cup \{15Z^+ \cup \{0\})$.

Now we give examples of quasi biset interval bivector spaces their substructures.

***Example 4.1.35***: Let $V = V_1 \cup V_2 = \{$All $6 \times 6$ interval matrices with entries from $Z_7\} \cup$

$$\left\{ \sum_{i=0}^{\infty} a_i x^i \,\middle|\, a_i \in Q^+ \cup \{0\} \right\}$$

be a quasi biset interval bivector space over the biset $S = Z_7 \cup Q^+ \cup \{0\}$. Choose $W = W_1 \cup W_2 = \{$all $6 \times 6$ upper triangular interval matrices with entries from $Z_7\} \cup$

$$\left\{ \sum_{i=0}^{\infty} a_i x^{2i} \,\middle|\, a_i \in Q^+ \cup \{0\} \right\}$$

$\subseteq V_1 \cup V_2 = V$ be a quasi biset interval bivector subspace of V over the biset $Z_7 \cup Q^+ \cup \{0\}$.

***Example 4.1.36***: Let $V = V_1 \cup V_2 =$

$$\left\{ \sum_{i=0}^{\infty} [0, a_i] x^i \,\middle|\, a_i \in Z_{49} \right\}$$

$\cup \{$All $16 \times 16$ matrices with entries from $Z_{81}\}$ be a quasi biset interval bivector space over the biset $S = S_1 \cup S_2 = Z_{49} \cup Z_{81}$. Choose $W = W_1 \cup W_2$

$$= \left\{ \sum_{i=0}^{\infty} [0, a_i] x^{2i} \,\middle|\, a_i \in \{0, 7, 14, 21, 28, 35, 42\} \subseteq Z_{49} \right\}$$

∪ {all 16 × 16 diagonal matrices with entries from $Z_{81}$} ⊆ $V_1$ ∪ $V_2$ = V. W is a quasi biset interval bivector subspace of V over the biset S = $S_1$ ∪ $S_2$.

***Example 4.1.37***: Let V = $V_1$ ∪ $V_2$ = {[0, $a_i$] | $a_i$ ∈ {0, 1, 2, 3, 4} = $Z_5$} ∪

$$\left\{ \begin{bmatrix} a & a \\ a & a \end{bmatrix} \,\middle|\, a \in Z_7 \right\}$$

be a quasi biset interval bivector space over the biset S = $S_1$ ∪ $S_2$ = $Z_5$ ∪ $Z_7$. It is easy to verify V has no quasi biset interval bivector subspace over S.

We define those quasi biset interval bivector spaces which has no subspace to be a simple quasi biset interval bivector space.

Vector space given an example 4.1.37 is a simple quasi biset interval bivector space.

Now we proceed onto define the notion of quasi subbiset interval bivector space.

**DEFINITION 4.1.12**: *Let V = $V_1$ ∪ $V_2$ be a quasi biset interval bivector space over the biset S = $S_1$ ∪ $S_2$. Let W = $W_1$ ∪ $W_2$ ⊆ $V_1$ ∪ $V_2$; where W = $W_1$ ∪ $W_2$ is a quasi biset interval bivector space over the biset P = $P_1$ ∪ $P_2$ ⊆ $S_1$ ∪ $S_2$ = S (where P is a proper subbiset of S) then we call W to be a quasi subbiset interval bivector subspace of V over the subbiset P of S.*

We will illustrate this situation by some examples. If V has no proper quasi subbiset interval bivector subspace then we call V to be a pseudo simple quasi biset interval bivector space. If V is both simple and pseudo simple then we call V to be a doubly simple interval space.

We will illustrate this by some examples.

***Example 4.1.38***: Let $V = V_1 \cup V_2 = \{[0, a] \mid a \in Z_7\} \cup$

$$\left\{ \begin{bmatrix} a \\ a \\ a \end{bmatrix} \middle| a \in Z_5 \right\}$$

be a quasi biset interval bivector space over the biset $S = \{0, 1\} \cup Z_5 = S_1 \cup S_2$. $V$ is a doubly simple quasi biset interval bivector space over the biset $S$.

***Example 4.1.39***: Let $V = V_1 \cup V_2 = \{(a, a, a, a, a, a, a, a) \mid a \in Z_2\} \cup$

$$\left\{ \begin{bmatrix} [0,a] & [0,a] \\ [0,a] & [0,a] \\ [0,a] & [0,a] \end{bmatrix} \middle| a \in Z_3 \right\}$$

be a quasi biset interval bivector space over the biset $S = Z_2 \cup Z_3$. Clearly $V$ is a quasi doubly simple interval bivector space over the biset $S = Z_2 \cup Z_3$.

Now we proceed onto define the notion of quasi set interval linear algebra semiquasi set interval linear algebra.

**DEFINITION 4.1.13**: *Let $V = V_1 \cup V_2$ be a quasi set interval bivector space over the set $S$. Suppose each $V_i$ is closed under the operation, addition for i=1, 2, then we define $V = V_1 \cup V_2$ to be a quasi set interval linear bialgebra over the set $S$.*

We will illustrate this situation by some examples.

***Example 4.1.40***: Let $V = V_1 \cup V_2$
$= \left\{ \sum_i a_i x^i \mid a_i \in Z^+ \cup \{0\} \right\} \cup \left\{ \sum_{i=0}^{\infty} [0, a_i] x^i \middle| a_i \in Z^+ \cup \{0\} \right\}$

be a quasi set interval linear bialgebra over the set $S = 3Z^+ \cup \{0\}$.

*Example 4.1.41*: Let $V = V_1 \cup V_2$

$$= \left\{ \sum_{i=0}^{9} [0,a_i]x^i \;\middle|\; \begin{array}{l} a_i \in Z_7; \\ 0 \le i \le 9 \end{array} \right\} \cup \left\{ \begin{bmatrix} [0,a_1] & [0,a_2] \\ [0,a_3] & [0,a_4] \\ [0,a_5] & [0,a_6] \end{bmatrix} \;\middle|\; \begin{array}{l} a_i \in Z_7; \\ 1 \le i \le 6 \end{array} \right\}$$

be a quasi set interval linear bialgebra over the set $S = Z_7$. Clearly V is of finite order where as V given in example 4.1.40 of infinite order.

*Example 4.1.42*: Let $V = V_1 \cup V_2$ = {All 5 × 5 interval matrices with entries from $Q^+ \cup \{0\}$} $\cup$ {all 3 × 7 matrices with entries from $Q^+ \cup \{0\}$} be a quasi set interval bilinear algebra over the set $S = 13Z^+ \cup \{0\}$.

Now we proceed onto define semi quasi interval bilinear algebra (linear bialgebra) over the set S.

**DEFINITION 4.1.14**: *Let $V = V_1 \cup V_2$ where $V_1$ is a set interval linear algebra over the set S and $V_2$ is a set vector space over the same set S (or $V_1$ is a set interval vector space over the set S and $V_2$ is a set linear algebra over the set S). We define V to be a semi quasi set interval bilinear algebra over the set S.*

We will illustrate this situation by some examples.

*Example 4.1.43*: Let $V = V_1 \cup V_2$

$$= \left\{ [0,a], \begin{bmatrix} [0,a_1] & [0,a_2] \\ [0,a_3] & [0,a_4] \end{bmatrix} \;\middle|\; \begin{array}{l} a_i, a \in Z_{45}; \\ 1 \le i \le 4 \end{array} \right\} \cup$$

{All $9 \times 3$ matrices with entries from $Z_{45}$} be a semi quasi set interval bilinear algebra over the set $S = \{0, 1, 5, 7, 14, 27, 35, 42\} \subseteq Z_{45}$.

*Example 4.1.44:* Let $V = V_1 \cup V_2$

$$= \left\{ \begin{bmatrix} [0,a_1] & [0,a_2] & [0,a_3] & [0,a_4] & [0,a_5] \\ [0,a_6] & [0,a_7] & [0,a_8] & [0,a_9] & [0,a_{10}] \\ [0,a_{11}] & [0,a_{12}] & [0,a_{13}] & [0,a_{14}] & [0,a_{15}] \end{bmatrix} \middle| \begin{array}{l} a_i \in R^+ \cup \{0\} \\ 1 \leq i \leq 15 \end{array} \right\}$$

$$\cup \left\{ \begin{bmatrix} a \\ b \\ c \end{bmatrix}, [a,b,c,d,e,f] \middle| a,b,c,d,e,f \in R^+ \cup \{0\} \right\}$$

be a semi quasi set interval bilinear algebra over the set $S = \{0, 1, \sqrt{2}, \sqrt{7}/5, \frac{\sqrt{13}}{\sqrt{19}}, \sqrt{43}, 52, 75, 1031\} \subseteq R^+ \cup \{0\}$.

*Example 4.1.45*: Let $V = V_1 \cup V_2 = \{[0, a] \mid a \in Z_3\} \cup$

$$\left\{ \begin{bmatrix} a \\ b \end{bmatrix}, [a,b,c,d] \middle| a,b,c,d \in Z_3 \right\}$$

be a semi quasi set interval bilinear algebra over the set $S = \{0, 1, 2\} = Z_3$.

Now we proceed onto define quasi biset interval bilinear algebra defined over the biset $S = S_1 \cup S_2$.

**DEFINITION 4.1.15**: *Let $V = V_1 \cup V_2$ where $V_1$ is a set linear algebra over the set $S_1$ and $V_2$ is a set interval vector space over the set $S_2$, we define $V$ to be a quasi biset interval bilinear algebra over the biset $S = S_1 \cup S_2$.*

We will illustrate this situation by some examples.

***Example 4.1.46:*** Let $V = V_1 \cup V_2$

$$= \left\{ \begin{bmatrix} a & b \\ c & d \end{bmatrix}, \begin{bmatrix} a & b & c \\ d & e & f \\ g & h & i \end{bmatrix} \middle| a,b,c,d,e,f,g,h,i \in Z_{49} \right\} \cup$$

$$\left\{ \begin{bmatrix} [0,a_1] \\ [0,a_2] \\ [0,a_3] \end{bmatrix} \middle| \begin{array}{l} a_i \in Z_7; \\ 1 \le i \le 3 \end{array} \right\}$$

be a quasi biset interval bilinear algebra over the biset $S = Z_{19} \cup Z_7$.

***Example 4.1.47:*** Let $V = V_1 \cup V_2$

$$= \left\{ \begin{bmatrix} [0,a_1] \\ [0,a_2] \\ [0,a_3] \end{bmatrix}, ([0,a_1],[0,a_2]) \middle| \begin{array}{l} a_i \in Z^+ \cup \{0\}; \\ 1 \le i \le 3 \end{array} \right\} \cup$$

{all $3 \times 5$ matrices with entries from $Z_{49}$} be a quasi biset interval bilinear algebra over the biset $S = 5Z^+ \cup \{0\} \cup Z_{49}$.

We see the quasi biset interval bilinear algebra given in example 4.1.46 is of finite order where as the quasi biset interval bilinear algebra given in example 4.1.47 is of infinite order.

***Example 4.1.48:*** Let $V = V_1 \cup V_2 =$

$$\left\{ \sum_{i=0}^{\infty} [0,a_i] x^i \middle| a_i \in Z_5 \right\}$$

$$\cup \left\{ \begin{bmatrix} a & b \\ c & d \end{bmatrix}, [(a,b,c,d,e,f,g,h)] \middle| a,b,c,d,e,f,g,h \in Z^+ \cup \{0\} \right\}$$

be a quasi biset interval bilinear algebra over the biset $S = S_1 \cup S_2 = \{0, 1\} \cup \{1, 2, 5, 0, 7\}$.

Clearly the quasi biset interval bilinear algebra given in example 4.1.48 is of infinite order.

Now we will proceed onto give examples of substructures and the reader is given the simple task of defining these substructures.

**Example 4.1.49**: Let $V = V_1 \cup V_2 =$

$$\left\{ \begin{bmatrix} [0,a_1] & [0,a_3] & [0,a_5] \\ [0,a_2] & [0,a_4] & [0,a_6] \end{bmatrix}, \begin{bmatrix} [0,a_1] \\ [0,a_2] \\ [0,a_3] \\ [0,a_4] \end{bmatrix} \middle| a_i \in Z_{29}; 1 \le i \le 6 \right\} \cup$$

$$\left\{ \sum_{i=0}^{25} a_i x^i \middle| a_i \in Z_{29} \right\}$$

be a quasi set interval bilinear algebra over the set $S = Z_{29}$.

Choose $W = W_1 \cup W_2$

$$= \left\{ \begin{bmatrix} [0,a_1] \\ [0,a_2] \\ [0,a_3] \\ [0,a_4] \end{bmatrix} \middle| a_i \in Z_{29}; 1 \le i \le 4 \right\} \cup \left\{ \sum_{i=0}^{15} a_i x^i \middle| a_i \in Z_{29} \right\}$$

$\subseteq V_1 \cup V_2 = V$; W is a quasi set interval bilinear subalgebra of V over the set $S = Z_{29}$.

**Example 4.1.50**: Let $V = V_1 \cup V_2 =$

$$\left\{ (a,b,c), \begin{bmatrix} a & b \\ c & d \\ e & f \end{bmatrix} \middle| a,b,c,d,e,f \in Z^+ \cup \{0\} \right\} \cup$$

$\{[0, a] \mid a \in Z^+ \cup \{0\}\}$ be a quasi set interval bilinear algebra over the set $5Z^+ \cup \{0\} = S$. Take $W = W_1 \cup W_2 = \{(a, b, c) \mid a, b, c \in Z^+ \cup \{0\}\} \cup \{[0, a] \mid a \in 15Z^+ \cup \{0\}\} \subseteq V_1 \cup V_2 = V$; $W$ is a quasi set interval bilinear subalgebra of $V$ over the set $S = 5Z^+ \cup \{0\}$.

*Example 4.1.51*: Let $V = V_1 \cup V_2 =$

$$\left\{ \begin{bmatrix} a & b & e \\ c & d & f \end{bmatrix}, \begin{bmatrix} a & b \\ c & d \\ e & f \\ g & h \\ i & j \end{bmatrix} \middle| a,b,c,d,e,f,g,h,i,j \in Z_{17} \right\} \cup$$

$$\left\{ \begin{bmatrix} [0,a_1] & [0,a_2] & [0,a_3] \\ [0,a_4] & [0,a_5] & [0,a_6] \\ [0,a_7] & [0,a_8] & [0,a_9] \end{bmatrix} \middle| \begin{array}{l} a_i \in Z_{47}; \\ 1 \leq i \leq 9 \end{array} \right\}$$

be a quasi biset interval bilinear algebra over the biset $S = Z_{17} \cup Z_{47}$. Choose $W = W_1 \cup W_2 =$

$$\left\{ \begin{bmatrix} a & b & e \\ c & d & f \end{bmatrix} \middle| a,b,c,d,e,f \in Z_{17} \right\} \cup$$

$$\left\{ \begin{bmatrix} [0,a_1] & [0,a_2] & [0,a_3] \\ 0 & [0,a_4] & [0,a_5] \\ 0 & 0 & [0,a_6] \end{bmatrix} \middle| \begin{array}{l} a_i \in Z_{47}; \\ 1 \leq i \leq 6 \end{array} \right\}$$

$\subseteq V_1 \cup V_2 = V$, $W$ is a quasi biset interval bilinear subalgebra of $V$ over the biset $S = Z_{17} \cup Z_{47}$.

*Example 4.1.52*: Let $V = V_1 \cup V_2 = \{$All $7 \times 7$ matrices with entries from $Z_{42}\} \cup$

$$\left\{ ([0,a_1],[0,a_2],[0,a_3],[0,a_4]), \begin{bmatrix} [0,a_1] & [0,a_3] & [0,a_5] \\ [0,a_2] & [0,a_4] & [0,a_6] \end{bmatrix} \middle| \begin{array}{l} a_i \in Z^+ \cup \{0\}; \\ 1 \le i \le 6 \end{array} \right\}$$

be a quasi biset interval linear bialgebra over the biset $S = S_1 \cup S_2 = Z_{42} \cup Z^+ \cup \{0\}$. Choose $W = W_1 \cup W_2 = \{$all $7 \times 7$ upper triangular matrices with entries from $Z_{42}\} \cup \{([0, a_1], [0, a_2], [0, a_3], [0, a_4])| a_i \in Z^+ \cup \{0\}; 1 \le i \le 4\} \subseteq V_1 \cup V_2$ be a quasi biset interval bilinear subalgebra of V over the biset $S = S_1 \cup S_2 = Z_{42} \cup Z^+ \cup \{0\}$.

Now as in case of set interval bivector spaces we can define the notion of quasi subset interval bilinear subalgebras. As the definition is a matter of routine the reader is given that task. However we illustrate this situation by some examples.

*Example 4.1.53*: Let $V = V_1 \cup V_2 =$

$$\left\{ (a_1, a_2, a_3, a_4, a_5, a_6, a_7), \begin{bmatrix} a & b & c \\ d & e & f \\ g & h & i \end{bmatrix} \middle| a_i, a, b, c, d, e, \right.$$

$$\left. f, g, h, i \in Z_9; 1 \le i \le 6 \right\} \cup \left\{ \begin{bmatrix} [0,a_1] & [0,a_2] \\ [0,a_3] & [0,a_4] \\ [0,a_5] & [0,a_6] \end{bmatrix} \middle| \begin{array}{l} a_i \in Z_9 \\ 1 \le i \le 6 \end{array} \right\}$$

be a quasi set interval bilinear algebra over the set $S = Z_9$.
Choose $W = W_1 \cup W_2 =$
$\{(a_1, a_2, a_3, a_4, a_5, a_6, a_7) \mid a_i \in Z_9; 1 \le i \le 7\} \cup$

$$\left\{ \begin{bmatrix} [0,a_1] & 0 \\ 0 & [0,a_2] \\ [0,a_3] & 0 \end{bmatrix} \middle| \begin{array}{l} a_i \in Z_9; \\ 1 \le i \le 3 \end{array} \right\}$$

$\subseteq V = V_1 \cup V_2$ and $P = \{0, 3, 6\} \subseteq Z_9 = S$.

W is a quasi subset interval linear subalgebra of V over the subset $P = \{0, 3, 6\}$ of S.

***Example 4.1.54***: Let $V = V_1 \cup V_2 = \{$all $5 \times 5$ matrices with entries from $Z^+ \cup \{0\}\} \cup$

$$\left\{ \begin{bmatrix} [0,a_1] \\ [0,a_2] \\ [0,a_3] \\ [0,a_4] \end{bmatrix}, \begin{bmatrix} [0,a_1] & [0,a_2] & [0,a_3] \\ 0 & [0,a_4] & [0,a_5] \\ 0 & 0 & [0,a_6] \end{bmatrix} \middle| a_i \in 3Z^+ \cup \{0\}; 1 \le i \le 6 \right\}$$

be a quasi set interval bilinear algebra over the set $S = 3Z^+ \cup \{0\}$. Choose $W = W_1 \cup W_2 = \{$all $5 \times 5$ upper triangular matrices with entries from $Z^+ \cup \{0\}\} \cup$

$$\left\{ \begin{bmatrix} [0,a] & [0,a] & [0,a] \\ 0 & [0,a] & [0,a] \\ 0 & 0 & [0,a] \end{bmatrix} \middle| a_i \in 3Z^+ \cup \{0\} \right\} \subseteq V_1 \cup V_2 = V$$

and $P = 33Z^+ \cup \{0\} \subseteq S = 3Z^+ \cup \{0\}$. Clearly $W = W_1 \cup W_2$ is a quasi subset interval bilinear subalgebra of V over the subset P $\subseteq$ S. However it is possible that V has no quasi subset interval bilinear subalgebra in such cases we call V to be a pseudo simple quasi set interval bilinear algebra.

We will illustrate this situation by some examples.

***Example 4.1.55***: Let $V = V_1 \cup V_2 = \{[0, a] | a \in Z_3\} \cup$

$$\left\{ \begin{bmatrix} a & a \\ a & a \\ a & a \end{bmatrix} \middle| a \in Z_3 \right\}$$

be a quasi set interval bilinear algebra over the set $S = \{0, 1\} \subseteq Z_3$. Clearly V has no quasi subset interval bilinear algebra as well as V has no quasi set interval bilinear subalgebra. Thus V is a pseudo simple quasi set interval bilinear algebra as well as simple quasi set interval bilinear algebra which we choose to call as doubly simple quasi set bilinear algebra.

*Example 4.1.56*: Let

$$V = V_1 \cup V_2$$

$$= \left\{ \begin{bmatrix} [0,a] & [0,b] \\ [0,c] & [0,d] \end{bmatrix} \middle| a,b,c,d \in Z_7 \right\} \cup$$

$$\left\{ \sum_{i=0}^{\infty} a_i x^i, (a,b,c,d) \middle| a_i, a, b, c, d \in Z_7 \right\}$$

be a quasi set interval bilinear algebra over the set $S = \{0, 1\}$.

Since S cannot have proper subsets of order greater than or equal two we see V is a pseudo simple quasi set interval bilinear algebra. However V is not a simple quasi set interval bilinear algebras as

$$W = \left\{ \begin{bmatrix} [0,a] & [0,a] \\ [0,a] & [0,a] \end{bmatrix} a \in Z_7 \right\} \cup \left\{ \sum_{i=0}^{\infty} a_i x^i \middle| a_i \in Z_7 \right\}$$

$\subseteq V_1 \cup V_2 = V$ is a quasi set interval bilinear subalgebra of V over the set $S = \{0, 1\}$.

Thus V is not a doubly simple quasi set interval linear algebra over the set S.

Now we will give yet another example to show the different possibilities.

*Example 4.1.57*: Let $V = V_1 \cup V_2 =$

$$\left\{ \begin{bmatrix} a & a \\ a & a \end{bmatrix}, (a,a) \middle| a \in Z_{19} \right\} \cup \left\{ \begin{bmatrix} [0,a] \\ [0,a] \\ [0,a] \end{bmatrix} \middle| a \in Z_{19} \right\}$$

be a quasi set interval bilinear algebra over the set $S = Z_{19}$. Clearly V has no quasi set interval bilinear subalgebra over $S = Z_{19}$.

However we see S can have several subsets but V cannot have any proper pseudo quasi subset interval bilinear subalgebras.

Hence V is a doubly simple quasi set interval bilinear algebra over $Z_{19} = S$.

Now we proceed onto define the notion of quasi subbiset interval bilinear subalgebras.

**DEFINITION 4.1.16**: *Let $V = V_1 \cup V_2$ be a quasi biset interval bilinear algebra over the biset $S = S_1 \cup S_2$. Let $W = W_1 \cup W_2 \subseteq V_1 \cup V_2$ and $P = P_1 \cup P_2 \subseteq S_1 \cup S_2$ both W and P are proper bisubsets of V and S respectively. Suppose W is a quasi biset interval bilinear algebra over the biset $P = P_1 \cup P_2$ then we call W to be a quasi bisubset interval bilinear subalgebra of V over the bisubset P of S.*

*We will say V is pseudo simple quasi biset interval bilinear algebra over the bisubset interval bilinear subalgebra over a bisubset $P = P_1 \cup P_2$ of $S = S_1 \cup S_2$.*

We will illustrate these situations by some simple examples.

***Example 4.1.58***: Let $V = V_1 \cup V_2 =$

$$\left\{ \begin{bmatrix} a & b \\ c & d \end{bmatrix}, \begin{bmatrix} a & b & a & b & a & b \\ b & a & b & a & b & a \end{bmatrix} \middle| a,b,c,d \in Z_6 \right\} \cup$$

$$\left\{ \begin{bmatrix} [0,a_1] & [0,a_2] \\ [0,a_3] & [0,a_4] \\ [0,a_5] & [0,a_6] \\ [0,a_7] & [0,a_8] \end{bmatrix} \middle| a_i \in Z_8; 1 \leq i \leq 8 \right\}$$

be a quasi biset interval bilinear algebra over the biset $S = Z_6 \cup Z_8$.

Choose $W = W_1 \cup W_2$

$$= \left\{ \begin{bmatrix} a & b \\ c & d \end{bmatrix} \middle| a,b,c,d \in Z_6 \right\} \cup \left\{ \begin{bmatrix} [0,a_1] & 0 \\ 0 & [0,a_2] \\ [0,a_3] & 0 \\ 0 & [0,a_4] \end{bmatrix} \middle| a_i \in Z_8; 1 \le i \le 8 \right\}$$

$\subseteq V_1 \cup V_2 = V$ such that W is a quasi bisubset interval bilinear algebra over the subbiset $P = \{0, 3\} \cup \{0, 2, 4, 6\} \subseteq S_1 \cup S_2$.

*Example 4.1.59*: Let $V = V_1 \cup V_2 =$

$$\left\{ \begin{bmatrix} [0,a_1] & [0,a_2] & [0,a_3] & [0,a_4] \\ [0,a_5] & [0,a_6] & [0,a_7] & [0,a_8] \end{bmatrix}, \begin{bmatrix} [0,a_1] \\ [0,a_2] \\ [0,a_3] \\ [0,a_4] \end{bmatrix} \middle| a_i \in Z_3; 1 \le i \le 8 \right\} \cup$$

$$\left\{ \begin{bmatrix} a & b \\ c & d \end{bmatrix} \middle| a,b,c,d \in Z_2 \right\}$$

be a quasi biset interval bilinear algebra over the biset $S = Z_3 \cup \{0, 1\}$.

We see V has no quasi subbiset interval bilinear subalgebra over S as S does not contain any proper subbiset.

Thus V is a pseudo simple quasi biset interval bilinear algebra over the set $S = Z_3 \cup \{0, 1\}$.

*Example 4.1.60*: Let $V = V_1 \cup V_2$

$$= \left\{ \begin{bmatrix} a & a & a \\ a & a & a \\ a & a & a \end{bmatrix} \middle| a \in Z_5 \right\} \cup$$

$$\left\{ \left([0,a_1] \quad [0,a_2] \quad [0,a_3] \quad [0,a_4]\right), \begin{bmatrix} [0,a_1] \\ [0,a_2] \\ [0,a_3] \\ [0,a_4] \end{bmatrix} \middle| \begin{array}{l} a_i \in Z_2; \\ 1 \le i \le 4 \end{array} \right\}$$

be a quasi biset interval bilinear algebra over the biset $S = Z_5 \cup Z_2$. Clearly V has no quasi biset interval bilinear subalgebras. Also V does not contain any quasi subbiset interval bilinear subalgebras. Thus V is both a pseudo simple quasi biset bilinear algebra as well as simple quasi biset interval bilinear algebra. We call a quasi biset interval bilinear algebra which is both a simple quasi biset interval bilinear algebra as well as pseudo simple quasi biset interval bilinear algebra as doubly simple quasi biset interval bilinear algebra.

We have given examples of all types of quasi biset interval bilinear algebras. Now we proceed onto give some properties and the reader is expected to prove them.

**THEOREM 4.1.1**: *Every quasi set interval bilinear algebra over a set S is a quasi set interval bivector space and the converse in general is not true.*

**THEOREM 4.1.2:** *Every set interval bilinear algebra is a set interval bivector space and not conversely.*

**THEOREM 4.1.3:** *Every set interval bilinear algebra is a quasi set interval bilinear algebra and not conversely.*

**THEOREM 4.1.4:** *Every doubly simple quasi interval bilinear algebra is a simple quasi interval bilinear algebra.*

**THEOREM 4.1.5:** *Every doubly simple quasi set interval bilinear algebra is a pseudo simple quasi set interval bilinear algebra.*

In the next section we proceed onto define semigroup interval bivector spaces and bilinear algebras.

## 4.2 Semigroup Interval Bilinear Algebras and Their Properties

In this section we define semigroup interval bilinear algebras and several related structures and substructures associated with them. Main properties about them are discussed in this section.

**DEFINITION 4.2.1**: *Let $V = V_1 \cup V_2$ be such that $V_1$ is semigroup interval vector space over the semigroup S and $V_2$ is also a semigroup interval vector space over the same semigroup S; where $V_1$ and $V_2$ are distinct with $V_1 \not\subseteq V_2$ or $V_2 \not\subseteq V_1$.*

*We define $V = V_1 \cup V_2$ to be a semigroup interval bivector space over the semigroup S.*

We will illustrate this situation by some examples.

***Example 4.2.1***: Let $V = V_1 \cup V_2 =$

$$\left\{ \begin{bmatrix} [0,a_1] & [0,a_2] & [0,a_3] \\ [0,a_4] & [0,a_5] & [0,a_6] \end{bmatrix} \middle| \begin{array}{l} a_i \in Z_{19}; \\ 1 \leq i \leq 6 \end{array} \right\} \cup \left\{ \sum_{i=0}^{\infty} [0,a_i]x^i \middle| a_i \in Z_{19} \right\}$$

be a semigroup interval bivector space over the semigroup $S = Z_{19}$.

***Example 4.2.2***: Let $V = V_1 \cup V_2 =$

$$\left\{ \begin{bmatrix} [0,a_1] & [0,a_2] \\ [0,a_3] & [0,a_4] \\ [0,a_5] & [0,a_6] \\ [0,a_7] & [0,a_8] \\ [0,a_9] & [0,a_{10}] \end{bmatrix} \middle| a_i \in Z^+ \cup \{0\} \right\} \cup$$

$\{([0, a_1], [0, a_2], …, [0, a_{24}]) \mid a_i \in 3Z^+ \cup \{0\}\}$ be a semigroup interval bivector space over the semigroup $S = Z^+ \cup \{0\}$.

*Example 4.2.3*: Let $V = V_1 \cup V_2 =$

$$\left\{ \begin{bmatrix} [0,a_1] & 0 & [0,a_3] & 0 \\ 0 & [0,a_2] & 0 & [0,a_4] \\ [0,a_5] & 0 & [0,a_6] & 0 \end{bmatrix} \middle| \begin{matrix} a_i \in Z_{12}; \\ 1 \le i \le 6 \end{matrix} \right\} \cup$$

$$\left\{ \begin{bmatrix} [0,a_1] & [0,a_2] & [0,a_3] \\ [0,a_4] & [0,a_5] & [0,a_6] \\ [0,a_7] & [0,a_8] & [0,a_9] \end{bmatrix} \middle| \begin{matrix} a_i \in Z_{12}; \\ 1 \le i \le 9 \end{matrix} \right\}$$

be a semigroup interval bivector space over the semigroup $S = \{0, 3, 6, 9\}$.

Now we define two substructures in them.

**DEFINITION 4.2.2**: *Let $V = V_1 \cup V_2$ be a semigroup interval bivector space over the semigroup S. Choose $W = W_1 \cup W_2 \subseteq V_1 \cup V_2 = V$; W a proper subset of V; if W itself a semigroup interval bivector space over the semigroup S then we define W to be a semigroup interval bivector subspace of V over the semigroup S. If V has no proper semigroup interval bivector subspace then we call V to be a simple semigroup interval bivector space.*

We will illustrate this situation by some examples.

*Example 4.2.4*: Let $V = V_1 \cup V_2 =$

$$\left\{ \begin{bmatrix} [0,a_1] & [0,a_2] \\ [0,a_3] & [0,a_4] \end{bmatrix}, ([0,a_1] \ [0,a_2] \ [0,a_3]) \middle| \begin{matrix} a_i \in Z_5; \\ 1 \le i \le 4 \end{matrix} \right\} \cup$$

$$\left\{ \begin{bmatrix} [0,a_1] \\ [0,a_2] \\ [0,a_3] \\ [0,a_4] \end{bmatrix} \middle| \begin{matrix} a_i \in Z_5; \\ 1 \le i \le 4 \end{matrix} \right\}$$

be a semigroup interval bivector space over the semigroup $S = Z_5$. Choose $W = W_1 \cup W_2 =$

$$\left\{ \begin{bmatrix} [0,a_1] & [0,a_2] \\ [0,a_3] & [0,a_4] \end{bmatrix} \middle| \begin{array}{l} a_i \in Z_5; \\ 1 \le i \le 4 \end{array} \right\} \cup \left\{ \begin{bmatrix} [0,a_1] \\ 0 \\ [0,a_2] \\ 0 \end{bmatrix} \middle| \begin{array}{l} a_i \in Z_5; \\ 1 \le i \le 2 \end{array} \right\}$$

$\subseteq V_1 \cup V_2 = V$; $W$ is a semigroup interval bivector subspace of $V$ over the semigroup $Z_5$.

*Example 4.2.5*: Let $V = V_1 \cup V_2 =$

$$\left\{ [0,a], \begin{bmatrix} [0,a] \\ [0,a] \end{bmatrix} \middle| a \in Z^+ \cup \{0\} \right\} \cup$$

$$\left\{ \begin{bmatrix} [0,a] \\ [0,a] \\ [0,a] \\ [0,a] \end{bmatrix}, \sum_{i=0}^{7} [0,a_i]x^i \middle| a, a_i \in 3Z^+ \cup \{0\} \right\}$$

be a semigroup interval bivector space over the semigroup $S = 4Z^+ \cup \{0\}$. Take $W = W_1 \cup W_2 =$

$$\left\{ [0,a] \middle| a \in 3Z^+ \cup \{0\} \right\} \cup \left\{ \sum_{i=0}^{7} [0,a_i]x^i \middle| a_i \in 3Z^+ \cup \{0\} \right\}$$

$\subseteq V_1 \cup V_2 = V$; $W$ is a semigroup interval bivector subspace of $V$ over the semigroup $4Z^+ \cup \{0\} = S$.

We see the semigroup interval bivector space can be of finite order or infinite order. Clearly the semigroup interval bivector space given in example 4.2.4 is of finite order where as the semigroup interval bivector space given in example 4.2.5 is of infinite order.

***Example 4.2.6***: Let $V = V_1 \cup V_2 =$

$$\{[0, a] \mid a \in Z_{29}\} \cup \left\{ \begin{bmatrix} [0,a] \\ [0,a] \\ [0,a] \end{bmatrix} \middle| a \in Z_{29} \right\}$$

be a semigroup interval bivector space over the semigroup $S = Z_{29}$. V is a simple semigroup interval bivector space as V has no semigroup interval bivector subspaces.

***Example 4.2.7***: Let $V = V_1 \cup V_2$

$$= \left\{ \begin{bmatrix} [0,a] & [0,a] & [0,a] \\ [0,a] & [0,a] & [0,a] \\ [0,a] & [0,a] & [0,a] \end{bmatrix} \middle| a \in Z_{13} \right\} \cup$$

$\{([0, a], [0, a], [0, a], [0, a], [0, a]) \mid a \in Z_{13}\}$ be a semigroup interval bivector space over the semigroup $S = Z_{13}$. V is a simple semigroup interval bivector space over $S = Z_{13}$.

**DEFINITION 4.2.3**: *Let $V = V_1 \cup V_2$ be a semigroup interval bivector space over the semigroup S. Let $W = W_1 \cup W_2 \subseteq V_1 \cup V_2 = V$ and $P \subseteq S$ (W and P are proper bisubset and subsemigroup of V and S respectively).*
   *If $W = W_1 \cup W_2$ is a semigroup interval bivector space over the semigroup P then we define W to be a subsemigroup interval bivector subspace of V. If V has no subsemigroup interval bivector subspaces then we call V to be a pseudo simple semigroup interval bivector space over the semigroup S.*

We will illustrate this situation by some examples.

**Example 4.2.8**: Let $V = V_1 \cup V_2 =$

$$\left\{ \begin{bmatrix} [0,a_1] & [0,a_2] \\ [0,a_3] & [0,a_4] \end{bmatrix}, ([0,a_1] \ [0,a_2]) \,\middle|\, \begin{matrix} a_i \in Z_{12}; \\ 1 \leq i \leq 4 \end{matrix} \right\} \cup$$

$$\left\{ \begin{bmatrix} [0,a_1] & [0,a_2] \\ [0,a_3] & [0,a_4] \\ [0,a_5] & [0,a_6] \\ [0,a_7] & [0,a_8] \end{bmatrix} \,\middle|\, a_i \in Z_{12} \right\}$$

be a semigroup interval bivector space defined over the semigroup $S = Z_{12}$. Choose

$$W = \left\{ \begin{bmatrix} [0,a_1] & [0,a_2] \\ [0,a_3] & [0,a_4] \end{bmatrix} \,\middle|\, \begin{matrix} a_i \in Z_{12}; \\ 1 \leq i \leq 4 \end{matrix} \right\} \cup$$

$$\left\{ \begin{bmatrix} [0,a_1] & [0,a_2] \\ 0 & 0 \\ [0,a_3] & [0,a_4] \\ 0 & 0 \end{bmatrix} \,\middle|\, \begin{matrix} a_i \in Z_{12}; \\ 1 \leq i \leq 4 \end{matrix} \right\}$$

$\subseteq V_1 \cup V_2 = V$ and $P = \{0, 4, 8\} \subseteq Z_{12}$. $W = W_1 \cup W_2$ is a subsemigroup interval bivector subspace of V over the subsemigroup P of $S = Z_{12}$.

**Example 4.2.9**: Let $V = V_1 \cup V_2 =$

$$\left\{ \begin{bmatrix} [0,a_1] \\ [0,a_2] \\ [0,a_3] \end{bmatrix} \,\middle|\, \begin{matrix} a_i \in Z^+ \cup \{0\} \\ 1 \leq i \leq 3 \end{matrix} \right\} \cup$$

$$\left\{ ([0,a_1][0,a_2][0,a_3][0,a_4]), \begin{bmatrix} [0,a_1] & [0,a_2] & 0 \\ 0 & [0,a_3] & [0,a_4] \\ [0,a_6] & 0 & [0,a_5] \end{bmatrix} \middle| \begin{array}{l} a_i \in 3Z^+ \cup \{0\}; \\ 1 \le i \le 6 \end{array} \right\}$$

be a semigroup interval bivector space over the semigroup $S = 5Z^+ \cup \{0\}$.

Choose $W = W_1 \cup W_2 =$

$$\left\{ \begin{bmatrix} [0,a_1] \\ 0 \\ [0,a_3] \end{bmatrix} \middle| a_1, a_3 \in Z^+ \cup \{0\} \right\} \cup$$

$$\left\{ \begin{bmatrix} 0 & [0,a_1] & 0 \\ 0 & 0 & [0,a_4] \\ [0,a_6] & 0 & 0 \end{bmatrix} \middle| a_1, a_6, a_4 \text{ are in } Z^+ \cup \{0\} \right\}$$

$\subseteq V_1 \cup V_2 = V$; and $P = \{125Z^+ \cup \{0\}\} \subseteq S$. $W = W_1 \cup W_2$ is a subsemigroup interval bivector subspace of V over the subsemigroup P of S.

*Example 4.2.10*: Let $V = V_1 \cup V_2$

$$= \left\{ \begin{bmatrix} [0,a_1] \\ [0,a_2] \\ [0,a_3] \end{bmatrix} \middle| \begin{array}{l} a_i \in Z_5; \\ 1 \le i \le 3 \end{array} \right\} \cup$$

$\{[[0, a_1] \ [0, a_1] \ [0, a_1] \ [0, a_1]] \mid a_i \in Z_5; 1 \le i \le 4\}$

be a semigroup interval bivector space over the semigroup $S = Z_5$. Clearly V has no subsemigroup interval bivector subspace as S has no subsemigroups. Thus V is a pseudo simple semigroup interval bivector space over the semigroup $S = Z_5$.

***Example 4.2.11***: Let $V = V_1 \cup V_2 =$

$$\left\{ \begin{bmatrix} [0,a] & [0,a] \\ [0,a] & [0,a] \end{bmatrix} \middle| a \in Z_7 \right\} \cup$$

$\{([0, a], [0, a], [0, a], [0, a], [0, a], [0, a]) \mid a \in Z_7\}$ be a semigroup interval bivector space over the semigroup $S = Z_7$.

Clearly V is a doubly simple semigroup interval bivector space over the semigroup $S = Z_7$.

***Example 4.2.12***: Let $V = V_1 \cup V_2 =$

$$\left\{ \begin{bmatrix} [0,a] & [0,a] & [0,a] & [0,a] & [0,a] \\ [0,a] & [0,a] & [0,a] & [0,a] & [0,a] \\ [0,a] & [0,a] & [0,a] & [0,a] & [0,a] \end{bmatrix} \middle| a \in Z_5 \right\}$$

$$\cup \left\{ \sum_{n=0}^{\infty} [0, a_i] x^i \middle| a_i \in Z_5 \right\}$$

be a semigroup interval bivector space over the semigroup $S = Z_5$. Clearly V is a doubly simple semigroup interval bivector space over the semigroup $S = Z_5$.

We see there is difference between the semigroup interval bivector space described in example 4.2.11 and 4.2.12 for we see in example 4.2.11 both $V_1$ and $V_2$ are doubly simple where as we see in example 4.2.12 only $V_1$ is doubly simple and $V_2$ infact has a semigroup interval bivector subspace viz.,

$$W_2 = \left\{ \sum_{i=0}^{\infty} [0, a_i] x^{2i} \middle| a_i \in Z_5 \right\} \subseteq V_2;$$

so in view of this we are forced to define yet another new notion.

**DEFINITION 4.2.4**: *Let $V = V_1 \cup V_2$ be a semigroup interval bivector space over the semigroup S. If only one of $V_1$ or $V_2$ is doubly simple and one of $V_i$ is not simple or pseudo simple (or* *not used in the mutually exclusive sense) then we call V to be a semi simple semigroup interval bivector space.*

***Example 4.2.13:*** Let $V = V_1 \cup V_2 =$

$$\left\{ \sum_{i=0}^{\infty} a_i x^i \,\middle|\, a_i \in Z_{19} \right\} \cup \left\{ \begin{bmatrix} [0,a] & [0,a] \\ [0,a] & [0,a] \end{bmatrix} \,\middle|\, a \in Z_{19} \right\}$$

be a semigroup interval bivector space over the semigroup $S = Z_{19}$. Clearly V is a semi simple semigroup interval bivector space.

**THEOREM 4.2.1**: *Let $V = V_1 \cup V_2$ be a semigroup interval bivector space defined over the semigroup $S = Z_p$; p a prime. V can be either a doubly simple semigroup bviector space or a semi simple semigroup bivector space.*

Proof is left as an exercise to the reader.

Now we proceed onto define the notion of semigroup interval bilinear algebra.

**DEFINITION 4.2.5**: *Let $V = V_1 \cup V_2$ be a semigroup interval bivector space over the semigroup S. If both $V_1$ and $V_2$ are closed under addition that is they are semigroups under addition then we call V to be a semigroup interval bilinear algebra over the semigroup S.*

We will illustrate this situation by some examples.

***Example 4.2.14***: Let $V = V_1 \cup V_2$

$$= \left\{ \begin{bmatrix} [0,a_1] & [0,a_2] & [0,a_3] \\ [0,a_4] & [0,a_5] & [0,a_6] \end{bmatrix} \,\middle|\, \begin{array}{l} a_i \in Z_{12}; \\ 1 \leq i \leq 6 \end{array} \right\} \cup$$

{All $10 \times 10$ interval matrices with entries from $Z_{12}$} be a semigroup interval bilinear algebra over the semigroup $S = Z_{12}$.

***Example 4.2.15***: Let $V = V_1 \cup V_2$

$$= \left\{ \sum_{i=0}^{\infty} [0, a_i] x^i \,\middle|\, a_i \in Z^+ \cup \{0\} \right\} \cup$$

$\{\{([0, a_i] [0, a_i] [0, a_i])\}| a_i \in SZ^+ \cup \{0\}\}$ be a semigroup interval bilinear algebra over the semigroup $3Z^+ \cup \{0\} = S$.

We have an interesting related result.

**THEOREM 4.2.2**: *Let $V = V_1 \cup V_2$ be a semigroup interval bilinear algebra over the semigroup S then V is a semigroup interval bivector space over the semigroup S but the converse however is not true.*

The proof is left as an exercise to the reader.

Now we proceed onto define substructures of these structures.

**DEFINITION 4.2.6**: *Let $V = V_1 \cup V_2$ be a semigroup interval bilinear algebra over the semigroup S. Let $W = W_1 \cup W_2 \subseteq V_1 \cup V_2 = V$; suppose W is a semigroup interval bilinear algebra over the semigroup S then we call W to be a semigroup interval bilinear subalgebra of V over the semigroup S. If V has no semigroup interval bilinear subalgebra then we define V to be a simple semigroup interval bilinear algebra over the semigroup S.*

We will illustrate this situation by some examples.

***Example 4.2.16***: Let V =

$$\left\{ \begin{bmatrix} [0,a_1] & [0,a_2] & [0,a_3] & [0,a_4] & [0,a_5] \\ [0,a_6] & [0,a_7] & [0,a_8] & [0,a_9] & [0,a_{10}] \end{bmatrix} \,\middle|\, a_i \in Z_{12}; 1 \le i \le 10 \right\}$$

be a semigroup interval bilinear algebra over the semigroup $S = Z_{12}$. Choose $W = W_1 \cup W_2$

$$= \left\{ \begin{bmatrix} 0 & [0,a_2] & 0 & [0,a_1] & 0 \\ [0,a_3] & 0 & [0,a_4] & 0 & [0,a_5] \end{bmatrix} \middle| a_i \in Z_{12}; 1 \leq i \leq 5 \right\}$$

$$\cup \left\{ \sum_{i=0}^{\infty} [0,a_i] x^{2i} \middle| a_i \in Z_{12} \right\}$$

$\subseteq V_1 \cup V_2 = V$; W is a semigroup interval bilinear subalgebra over the semigroup $S = Z_{12}$.

*Example 4.2.17*: Let $V = V_1 \cup V_2 =$

$$\left\{ \begin{bmatrix} [0,a_1] & [0,a_2] & [0,a_3] \\ [0,a_4] & [0,a_5] & [0,a_6] \\ [0,a_7] & [0,a_8] & [0,a_9] \end{bmatrix} \middle| a_i \in Z^+ \cup \{0\}; 1 \leq i \leq 9 \right\} \cup$$

$$\left\{ \begin{bmatrix} [0,a_1] \\ [0,a_2] \\ [0,a_3] \\ [0,a_4] \\ [0,a_5] \end{bmatrix} \middle| a_i \in Z^+ \cup \{0\}; 1 \leq i \leq 5 \right\}$$

be a semigroup interval bilinear algebra over the semigroup $S = 3Z^+ \cup \{0\}$. Take $W = W_1 \cup W_2 =$

$$\left\{ \begin{bmatrix} [0,a_1] & [0,a_2] & [0,a_3] \\ 0 & [0,a_4] & [0,a_5] \\ 0 & 0 & [0,a_6] \end{bmatrix} \middle| a_i \in 3Z^+ \cup \{0\}; 1 \leq i \leq 6 \right\} \cup$$

$$\left\{ \begin{bmatrix} [0,a] \\ [0,a] \\ [0,a] \\ [0,a] \\ [0,a] \end{bmatrix} \middle| a_i \in 3Z^+ \cup \{0\} \right\}$$

⊆ V₁ ∪ V₂ = V; W is a semigroup interval bilinear subalgebra of V over the semigroup S = 3Z⁺ ∪ {0}.

***Example 4.2.18.*** Let V = V₁ ∪ V₂ =

$$\left\{ \begin{bmatrix} [0,a] & [0,a] \\ [0,a] & [0,a] \end{bmatrix} \bigg| a \in Z_7 \right\} \cup \left\{ \begin{bmatrix} [0,a] \\ [0,a] \\ [0,a] \\ [0,a] \\ [0,a] \end{bmatrix} \bigg| a \in Z_7 \right\}$$

be a semigroup interval bilinear algebra over the semigroup S = $Z_7$. We see V has no semigroup interval bilinear subalgebra; hence V is a simple semigroup interval bilinear algebra over the semigroup S = $Z_7$.

***Example 4.2.19.*** Let V = V₁ ∪ V₂ =

$$\left\{ \begin{bmatrix} [0,a] & [0,a] \\ [0,a] & [0,a] \\ [0,a] & [0,a] \\ [0,a] & [0,a] \\ [0,a] & [0,a] \end{bmatrix} \bigg| a \in Z_{11} \right\} \cup$$

{([0, a], [0, a], [0, a], [0, a], [0, a])|a ∈ $Z_{11}$} be a semigroup interval bilinear algebra over the semigroup S = $Z_{11}$. V is a simple semigroup interval bilinear algebra over the semigroup S = $Z_{11}$.

**DEFINITION 4.2.7:** *Let V = V₁ ∪ V₂ be a semigroup interval bilinear algebra over the semigroup S. Let W = W₁ ∪ W₂ ⊆ V₁ ∪ V₂ = V; be such that W is a semigroup interval bilinear algebra over a subsemigroup P of S. We define W = W₁ ∪ W₂ to be a subsemigroup interval bilinear subalgebra of V over the subsemigroup P of S. If V has no subsemigroup interval bilinear subalgebra then we define V to be a pseudo simple semigroup* *interval bilinear algebra over the semigroup S. If $V = V_1 \cup V_2$ is both a simple semigroup interval bilinear algebra as well as pseudo simple semigroup interval bilinear algebra over the semigroup S then we call V to be a doubly simple semigroup interval bilinear algebra over the semigroup S.*

***Example 4.2.20:*** Let $V = V_1 \cup V_2 =$

$$\left\{\begin{bmatrix} [0,a_1] & [0,a_2] \\ [0,a_3] & [0,a_4] \end{bmatrix} \middle| a_i \in Z_{12}; 1 \leq i \leq 4 \right\} \cup \left\{ \sum_{i=0}^{\infty} [0,a_i]x^i \middle| a_i \in Z_{12} \right\}$$

be a semigroup interval bilinear algebra over the semigroup $S = Z_{12}$ under addition modulo 12.
Choose $W = W_1 \cup W_2 =$

$$\left\{\begin{bmatrix} [0,a_1] & [0,a_2] \\ 0 & [0,a_3] \end{bmatrix} \text{ where } a_i \in Z_{12} \atop 1 \leq i \leq 3 \right\} \cup$$

$$\left[ \sum_{i=0}^{\infty} [0,a_i]x^{2i} \middle| a_i \in \{0,2,4,6,8,10\} \right]$$

$\subseteq V_1 \cup V_2 = V$ be a subsemigroup interval bilinear subalgebra of V over the subsemigroup $P = \{0, 6\} \subseteq Z_{12} = S$.

***Example 4.2.21***: Let $V = V_1 \cup V_2 =$

$$\left\{ \begin{bmatrix} [0,a_1] & [0,a_2] & [0,a_3] & [0,a_4] \\ [0,a_5] & [0,a_6] & [0,a_7] & [0,a_8] \end{bmatrix} \middle| a_i \in Z_5; 1 \leq i \leq 8 \right\} \cup$$

$$\left\{ \begin{bmatrix} [0,a_1] & [0,a_2] & [0,a_3] \\ [0,a_4] & [0,a_5] & [0,a_6] \\ [0,a_7] & [0,a_8] & [0,a_9] \\ [0,a_{10}] & [0,a_{11}] & [0,a_{12}] \\ [0,a_{13}] & [0,a_{14}] & [0,a_{15}] \end{bmatrix} \middle| a_i \in Z_5; 1 \leq i \leq 15 \right\}$$

be a semigroup interval bilinear algebra over the semigroup $S = Z_5$. Clearly S has no proper subsemigroups.

Take $W = W_1 \cup W_2 =$

$$\left\{ \begin{bmatrix} [0,a_1] & 0 & [0,a_2] & 0 \\ 0 & [0,a_3] & 0 & [0,a_4] \end{bmatrix} \middle| a_i \in Z_5; 1 \leq i \leq 4 \right\} \cup$$

$$\left\{ \begin{bmatrix} [0,a_1] & 0 & [0,a_2] \\ 0 & [0,a_3] & 0 \\ [0,a_4] & 0 & [0,a_5] \\ 0 & [0,a_6] & 0 \\ [0,a_7] & 0 & [0,a_8] \end{bmatrix} \middle| a_i \in Z_5; 1 \leq i \leq 8 \right\}$$

$\subseteq V_1 \cup V_2$ be a semigroup interval bilinear subalgebra of V over the semigroup $S = Z_5$.

However V has no proper subsemigroup interval bilinear subalgebra as S has no proper subsemigroups in $S = Z_5$ under addition modulo 5.

*Example 4.2.22:* Let $V = V_1 \cup V_2 =$

$$\left\{ \begin{bmatrix} [0,a] & [0,a] & [0,a] \\ [0,a] & [0,a] & [0,a] \end{bmatrix} \middle| a \in Z_{17} \right\} \cup \left\{ \begin{bmatrix} [0,a] & [0,a] \\ [0,a] & [0,a] \\ [0,a] & [0,a] \\ [0,a] & [0,a] \\ [0,a] & [0,a] \end{bmatrix} \middle| a \in Z_{17} \right\}$$

be a semigroup interval bilinear algebra over the semigroup $S = Z_{17}$. Clearly V has no semigroup interval bilinear subalgebra as well as V has no subsemigroup interval bilinear subalgebra. Thus V is a pseudo simple semigroup interval bilinear algebra over S.

Now having seen some of the substructures of the semigroup interval bilinear algebra we now proceed on to define more properties about them.

**DEFINITION 4.2.8:** *Let $V = V_1 \cup V_2$ be such that $V_1$ is a semigroup interval linear algebra over the semigroup S and $V_2$ is only a semigroup interval vector space over the same semigroup S and $V_2$ is not a linear algebra then we define $V = V_1 \cup V_2$ to be a quasi semigroup interval bilinear algebra over S.*

We will first illustrate this situation by some simple examples.

*Example 4.2.23*: Let $V = V_1 \cup V_2 =$

$$\left\{ \begin{bmatrix} [0,a_1] & [0,a_2] & [0,a_3] \\ [0,a_4] & [0,a_5] & [0,a_6] \\ [0,a_7] & [0,a_8] & [0,a_9] \end{bmatrix} \middle| a_i \in Z^+ \cup \{0\}; 1 \leq i \leq 9 \right\} \cup$$

$$\left\{ \begin{bmatrix} [0,a_1] \\ [0,a_2] \\ [0,a_3] \end{bmatrix}, \big[[0,a_1], [0,a_2], [0,a_3], [0,a_4], [0,a_5]\big] \middle| \begin{matrix} a_i \in 2Z^+ \cup \{0\}; \\ 1 \leq i \leq 5 \end{matrix} \right\}$$

be a quasi semigroup interval linear bialgebra over the semigroup $S = 6Z^+ \cup \{0\}$.

*Example 4.2.24*: Let $V = V_1 \cup V_2 =$

$$\left\{ \sum_{i=0}^{\infty} [0,a_i] x^i \middle| a_i \in Z_{47} \right\} \cup$$

$\{([0, a_1], 0, [0, a_2], 0, [0, a_3]), (0, [0, a_1], 0, [0, a_2], 0) \mid a_i \in Z_{47}, 1 \leq i \leq 3\}$ be a quasi semigroup interval bilinear algebra over the semigroup $S = Z_{47}$.

We can as in case of semigroup interval bilinear algebras define substructures. The definition is a matter of routine and is left as an exercise for the reader.

How ever we will illustrate this situation by some examples.

***Example 4.2.25***: Let $V = V_1 \cup V_2 =$

$$\left\{ \begin{bmatrix} [0,a] & [0,a] \\ [0,a] & [0,a] \end{bmatrix}, [[0,a] \quad [0,a] \quad [0,a] \quad [0,a] \quad [0,a]] \, \bigg| \, a_i \in Z_{421} \right\}$$

$$\cup \left\{ \sum_{i=0}^{\infty} [0, a_i] x^{2i} \, \bigg| \, a_i \in Z_{421} \right\}$$

$\subseteq V = V_1 \cup V_2$; W is a quasi semigroup interval bilinear subalgebra of V over the semigroup $S = Z_{421}$.

***Example 4.2.26***: Let $V = V_1 \cup V_2 =$

$$\left\{ \begin{bmatrix} [0,a] \\ [0,a] \end{bmatrix}, [0,a] \, \bigg| \, a \in Z_3 \right\} \cup \left\{ \begin{bmatrix} [0,a] & [0,a] & [0,a] \\ [0,b] & [0,b] & [0,b] \end{bmatrix} \, \bigg| \, a, b \in Z_3 \right\}$$

be the quasi semigroup interval bilinear algebra over the semigroup $S = Z_3$.

Consider $W = W_1 \cup W_2 =$

$$\{[0, a] \mid a \in Z_3\} \cup \left\{ \begin{bmatrix} [0,a] & [0,a] & [0,a] \\ 0 & 0 & 0 \end{bmatrix} \, \bigg| \, a \in Z_3 \right\}$$

$\subseteq V_1 \cup V_2$; W is a quasi semigroup interval bilinear subalgebra of V over $Z_3$.

If the quasi semigroup interval bilinear algebra V has no quasi semigroup interval bilinear subalgebras then we call V to be a simple quasi semigroup interval bilinear algebra.

We will first illustrate this situation by some examples.

**Example 4.2.27**: Let $V = V_1 \cup V_2 =$

$$\left\{ \begin{bmatrix} [0,a] & [0,a] \\ [0,b] & [0,b] \end{bmatrix}, \begin{bmatrix} [0,a] \\ [0,b] \end{bmatrix} \middle| a, b \in Z_5 \right\} \cup$$

$$\left\{ \begin{bmatrix} [0,a] & [0,a] & [0,a] & [0,a] & [0,a] \\ [0,a] & [0,a] & [0,a] & [0,a] & [0,a] \\ [0,a] & [0,a] & [0,a] & [0,a] & [0,a] \end{bmatrix} \middle| a \in Z_5 \right\}$$

be a quasi semigroup interval bilinear algebra over the semigroup $S = Z_5$. Clearly V has no quasi semigroup interval bilinear subalgebras. Hence V is a simple quasi semigroup interval bilinear algebra over $S = Z_5$.

We will now proceed on to give examples of quasi subsemigroup interval bilinear algebras. If the quasi semigroup interval bilinear algebra V has no quasi subsemigroup interval bilinear subalgebra then we call V to be a pseudo simple quasi semigroup interval bilinear algebra. If V is both a simple and a pseudo simple quasi semigroup interval bilinear algebra then we define V to be a doubly simple quasi semigroup interval bilinear algebra.

**Example 4.2.28:** Let $V = V_1 \cup V_2 =$

$$\left\{ \begin{bmatrix} [0,a_1] & [0,a_2] \\ [0,a_4] & [0,a_5] \end{bmatrix}, \begin{bmatrix} [0,a_3] \\ [0,a_6] \end{bmatrix} \middle| \begin{array}{l} a_i \in Z_{18}; \\ 1 \le i \le 6 \end{array} \right\} \cup$$

$$\left\{ \begin{bmatrix} [0,a_1] \\ [0,a_2] \\ [0,a_3] \\ [0,a_4] \end{bmatrix}, \begin{bmatrix} [0,a_1] & [0,a_2] & [0,a_3] \end{bmatrix} \middle| \begin{array}{l} a_i \in Z_{18}; \\ 1 \le i \le 4 \end{array} \right\}$$

be a quasi semigroup interval bilinear algebra over the semigroup $S = Z_{18}$.
Take $W = W_1 \cup W_2 =$

$$\left\{ \begin{bmatrix} [0,a] & [0,a] \\ [0,b] & [0,b] \end{bmatrix}, \begin{bmatrix} [0,a] \\ [0,b] \end{bmatrix} \middle| a,b \in \{0,2,4,6,8,10,12,14,16\} \subseteq Z_{18} \right\}$$
$$\cup \left\{ [[0,a_1] \quad [0,a_2] \quad [0,a_3]] \middle| a_i \in Z_{18}; 1 \le i \le 3 \right\}$$

$\subseteq V_1 \cup V_2 = V$ and $P = P\{0, 9\} \subseteq Z_{18}$ (P is a subsemigroup under addition modulo 18 of the semigroup $Z_{18}$).

W is a subsemigroup interval bilinear subalgebra of V over the subsemigroup $P \subseteq S = Z_{18}$.

*Example 4.2.29*: Let $V = V_1 \cup V_2 = \{$All $5 \times 5$ interval matrices with intervals of the form $[0, a_i]$ where $a_i \in Z^+ \cup \{0\}\} \cup$

$$\left\{ \begin{bmatrix} [0,a_1] & [0,a_2] \\ [0,a_3] & [0,a_4] \end{bmatrix}, \begin{bmatrix} [0,a_1] \\ [0,a_2] \\ [0,a_3] \\ [0,a_4] \\ [0,a_5] \\ [0,a_6] \end{bmatrix} \middle| a_i \in Z^+ \cup \{0\}; 1 \le i \le 6 \right\}$$

be a quasi semigroup interval bilinear algebra over the semigroup $S = Z^+ \cup \{0\}$. Let $W = W_1 \cup W_2 = \{$all $5 \times 5$ interval upper triangular matrices with intervals of the form $[0, a_i]$ where $a_i \in Z^+ \cup \{0\}\}$

$$\cup \left\{ \begin{bmatrix} [0,a_1] & [0,a_3] \\ [0,a_3] & [0,a_1] \end{bmatrix}, \begin{bmatrix} 0 \\ [0,a] \\ 0 \\ [0,a] \\ 0 \\ 0 \end{bmatrix} \middle| a, a_1, a_3 \in Z^+ \cup \{0\} \right\}$$

⊆ $V_1 \cup V_2$. W is a quasi subsemigroup interval bilinear subalgebra of V over the subsemigroup $P = 3Z^+ \cup \{0\} \subseteq Z^+ \cup \{0\} = S$.

*Example 4.2.30*: Let $V = V_1 \cup V_2 =$

$$\left\{ \begin{bmatrix} [0,a] & [0,a] \\ [0,a] & [0,a] \end{bmatrix}, \begin{bmatrix} [0,a] \\ [0,a] \end{bmatrix} \,\middle|\, a \in Z_7 \right\} \cup$$

$\{([0, a], [0, a], [0, a], [0, a], [0, a]) \mid a \in Z_7\}$ be a quasi semigroup interval bilinear algebra over the semigroup $S = Z_7$. Since S has no proper subsemigroups we see V is a pseudo simple quasi semigroup interval bilinear algebra over $S = Z_7$. Further as V has no proper semigroup interval bilinear algebras we see V is a simple quasi semigroup interval bilinear algebra. Thus V is a doubly simple quasi semigroup interval bilinear algebra over the semigroup $S = Z_7$.

Now we can define bilinear transformation of quasi semigroup interval bilinear algebras V to W also the notion of bilinear operator of a quasi semigroup interval bilinear algebra V.
    This task is left as an exercise for the reader.

### 4.3 Group Interval Bilinear Algebras and their Properties

In this section we proceed on to define the notion of group interval bivector spaces and describe a few of their properties associated with them.

**DEFINITION 4.3.1**: *Let $V = V_1 \cup V_2$ be such that each $V_i$ is a group interval vector space over the additive group G for i = 1, 2; such that*
    *(1) $V_1 \not\subseteq V_2$ and $V_2 \not\subseteq V_1$*
        *$V_1 \cap V_2 = \phi$ or non empty*
    *(2) For every $v = v_1 \cup v_2 \in V_1 \cup V_2 = V$ and $g \in G$ $gv = gv_1 \cup gv_2$ and $vg = v_1g \cup v_2g$ belong to $V = V_1 \cup V_2$.*

(3) $0.v = 0.v_1 \cup 0.v_2$
$= 0 \cup 0 \in V_1 \cup V_2 = V$

*0 is the additive identity of G.*

*We call V to be a group interval bivector space over the group G.*

We will illustrate this situation by some examples.

**Example 4.3.1**: Let $V = V_1 \cup V_2 =$

$$\left\{ \begin{bmatrix} [0,a] & [0,a] \\ [0,b] & [0,b] \end{bmatrix}, [0,b], \begin{bmatrix} [0,a] \\ [0,b] \\ [0,c] \\ [0,d] \\ [0,e] \end{bmatrix} \middle| a,b,c,d,e \in Z_{19} \right\}$$

$\cup \{([0, a_1], [0, a_2], [0, a_3], [0, a_4], [0, a_5]) \mid a_i \in Z_{19}; i = 1, 2, 3, 4, 5\}$ be a group interval bivector space over the group $G = Z_{19}$ (G is a group under addition modulo 19).

**Example 4.3.2:** Let $V = V_1 \cup V_2 =$

$$\left\{ \sum_{i=0}^{5} [0, a_i] x^i \middle| a_i \in Z_{12} \right\} \cup \left\{ \begin{bmatrix} [0, a_1] \\ [0, a_2] \\ [0, a_3] \\ [0, a_4] \\ [0, a_5] \\ [0, a_6] \end{bmatrix} \middle| a_i \in Z_{12}; 1 \leq i \leq 6 \right\}$$

be a group interval bivector space over the group $G = Z_{12}$ under addition modulo 12. We see both the group interval bivector spaces are of finite order. Further it is important at this juncture to state that we cannot built in this manner group interval bivector spaces using $Z^+ \cup \{0\}$ or $R^+ \cup \{0\}$ or $Q^+ \cup \{0\}$ or $C^+ \cup \{0\}$; as they are not groups under addition.

Thus we have our own limitations in dealing with them.

However we have infinite group interval bivector spaces using $Z_n$.

*Example 4.3.3:* Let $V = V_1 \cup V_2 =$

$$\left\{\sum_{i=0}^{\infty}[0,a]x^{2i} \,\Big|\, a_i \in Z_{42}\right\} \cup \left\{\sum_{i=0}^{\infty}[0,a_i]x^i, \sum_{i=0}^{\infty}[0,a_i]x^{3i} \,\Big|\, a, a_i \in Z_{42}\right\}$$

be a group interval bivector space over the group $G = Z_{42}$. Clearly V is of infinite order.

We now proceed onto define substructures related with these structures.

**DEFINITION 4.3.2**: *Let $V = V_1 \cup V_2$ be a group interval bivector space over the group G. Let $W = W_1 \cup W_2 \subseteq V_1 \cup V_2$; such that W is a group interval bivector space over the group G; then we define W to be a group interval bivector subspace of V over the group G.*

*We say V is a simple group interval bivector space if V has no proper group interval bivector subspace.*

We will illustrate this situation by some examples.

*Example 4.3.4*: Let $V = V_1 \cup V_2 =$

$$\left\{\begin{bmatrix}[0,a_1] & [0,a_2] \\ [0,a_3] & [0,a_4]\end{bmatrix}, \bigl([0,a_1] \quad [0,a_2] \quad [0,a_3]\bigr) \,\Bigg|\, a_i \in Z_{15}; 1 \le i \le 4\right\} \cup$$

$$\left\{\begin{bmatrix}[0,a_1] & [0,a_6] \\ [0,a_2] & [0,a_7] \\ [0,a_3] & [0,a_8] \\ [0,a_4] & [0,a_9] \\ [0,a_5] & [0,a_{10}]\end{bmatrix} \,\Bigg|\, a_i \in Z_{15}; 1 \le i \le 10\right\}$$

be a group interval bivector space over the group $G = Z_{15}$.
Take $W = W_1 \cup W_2 =$

$$\left\{ ([0,a_1] \quad [0,a_2] \quad [0,a_3]) \mid a_i \in Z_{15}; 1 \leq i \leq 3 \right\} \cup$$

$$\left\{ \begin{bmatrix} [0,a] & [0,a] \\ 0 & 0 \\ [0,a] & [0,a] \\ 0 & 0 \\ [0,a] & [0,a] \end{bmatrix} \mid a \in Z_{15} \right\}$$

$\subseteq V_1 \cup V_2 = V$; W is a group interval bivector subspace of V over the group G.

*Example 4.3.5*: Let $V = V_1 \cup V_2 =$

$$\left\{ \sum_{i=0}^{\infty} [0,a_i] x^i \mid a_i \in Z_{248} \right\}$$

$\cup$ {all $10 \times 10$ square interval matrices with entries from $I(Z_{248})$} be a group interval bivector space over the group $G = Z_{248}$.

Choose $W = W_1 \cup W_2 =$

$$\left\{ \sum_{i=0}^{\infty} [0,a_i] x^{2i} \mid a_i \in Z_{248} \right\}$$

$\cup$ {All $10 \times 10$ upper triangular matrices with entries from $I(Z_{248}) = \{[0, a_i] \mid a_i \in Z_{248}\} \subseteq V_1 \cup V_2$; W is a group interval bivector subspace of V over the group $G = Z_{248}$.

*Example 4.3.6:* Let $V = V_1 \cup V_2 =$

$$\left\{ \begin{bmatrix} [0,a] & [0,a] \\ [0,a] & [0,a] \end{bmatrix} \mid a \in Z_{15} \right\} \cup$$

$$\left\{ \begin{bmatrix} [0,a] & [0,a] & [0,a] & [0,a] \\ [0,a] & [0,a] & [0,a] & [0,a] \\ [0,a] & [0,a] & [0,a] & [0,a] \\ [0,a] & [0,a] & [0,a] & [0,a] \\ [0,a] & [0,a] & [0,a] & [0,a] \end{bmatrix} \,\bigg|\, a \in Z_{19} \right\}$$

be a group interval bivector space over the group $G = Z_{19}$. Clearly V is a simple group interval bivector space over G.

*Example 4.3.7*: Let $V = V_1 \cup V_2 =$

$$\left\{ \begin{bmatrix} [0,a] & [0,a] & [0,a] & [0,a] \\ [0,a] & [0,a] & [0,a] & [0,a] \\ [0,a] & [0,a] & [0,a] & [0,a] \end{bmatrix} \,\bigg|\, a \in Z_5 \right\} \cup$$

$$\left\{ \begin{bmatrix} [0,a] & [0,a] \\ [0,a] & [0,a] \\ [0,a] & [0,a] \\ [0,a] & [0,a] \end{bmatrix} \,\bigg|\, a \in Z_5 \right\}$$

be a group interval bivector space over the group $G = Z_5$. It is easily verified $V = V_1 \cup V_2$ simple group interval bivector space over the group $G = Z_5$.

Now we proceed onto define the notion of subgroup interval bivector subspace of a group interval bivector space.

**DEFINITION 4.3.3**: *Let $V = V_1 \cup V_2$ be a group interval bivector space over the group G. Let $W = W_1 \cup W_2 \subseteq V_1 \cup V_2$ and (e) ≠ $H \subseteq G$ be a subgroup of G. If $W = W_1 \cup W_2$ is a group interval bivector space over the group H then we define W to be a subgroup interval bivector subspace of V over the subgroup H of G. If V has no subgroup interval bivector subspace then we* *define V to be a pseudo simple group interval bivector space. If V is both a simple and pseudo simple group interval bivector space then we define V to be a doubly simple group interval bivector space over the group G.*

We will illustrate this situation by some examples.

***Example 4.3.8***: Let $V = V_1 \cup V_2 = \{([0, a_1], [0, a_2], [0, a_3], [0, a_4]) \mid a_i \in Z_{48}; 1 \leq i \leq 4\} \cup$

$$\left\{\sum_{i=0}^{8}[0,a_i]x^i \,\bigg|\, a_i \in Z_{48}\right\}$$

be a group interval bivector space over the group $G = Z_{48}$. Take $W = W_1 \cup W_2 = \{([0, a_1], 0, [0, a_2], 0) \mid a_i \in Z_{48}; 1 \leq i \leq 2\} \cup$

$$\left\{\sum_{i=0}^{8}[0,a_i]x^i \,\bigg|\, a_i \in \{0,2,4,6,8,...,44,46\}\right\}$$

$\subseteq V_1 \cup V_2$ and $H = \{0, 4, 8, 12, 16, 20, 24, 28, 32, 36, 40, 44\} \subseteq G$ a subgroup of $Z_{48}$ under addition modulo 48.

W is a subgroup interval bivector subspace of V over H the subgroup G.

***Example 4.3.9***: Let $V = V_1 \cup V_2 =$

$$\left\{\sum_{i=0}^{\infty}[0,a_i]x^i \,\bigg|\, a_i \in Z_{18}\right\}$$

$\cup$ {All $6 \times 6$ interval matrices with entries from $I(Z_{18})$} be a group interval bivector space over the group $G = Z_{18}$.

Take $W = W_1 \cup W_2 =$

$$\left\{\sum_{i=0}^{\infty}[0,a_i]x^{zi} \,\bigg|\, a_i \in Z_{18}\right\}$$

∪ {6 × 6 upper triangular interval matrices with entries from I $(Z_{18})$} ⊆ $V_1 \cup V_2$; W is a subgroup interval bivector subspace of V over the subgroup H = {0, 9} ⊆ $Z_{18}$.

*Example 4.3.10*: Let V = $V_1 \cup V_2$ =
$$\left\{ \sum_{i=0}^{\infty} [0, a_i] x^i \,\middle|\, a_i \in Z_{11} \right\}$$

∪ {set of all 11 × 15 interval matrices with entries from I $(Z_{11})$ = {[0, $a_i$] | $a_i \in Z_{11}$}} be a group interval bivector space over the group G = $Z_{11}$. We see G = $Z_{11}$ is a simple group under addition modulo 11. Hence V is a pseudo simple group interval bivector space over G.

However V has group interval bivector subspace so V is not a simple group interval bivector space over G.

*Example 4.3.11:* Let V = $V_1 \cup V_2$ = {([0, $a_1$], [0, $a_2$], [0, $a_3$], [0, $a_4$]) | $a_i \in Z_{43}$} ∪

$$\left\{ \begin{bmatrix} [0,a] & [0,b] \\ [0,a] & [0,b] \\ [0,a] & [0,b] \\ [0,a] & [0,b] \\ [0,a] & [0,b] \\ [0,a] & [0,b] \end{bmatrix} \,\middle|\, a, b \in Z_{43} \right\}$$

be a group interval bivector space over the group G = $Z_{43}$. G has no proper subgroups hence V is a pseudo simple group interval bivector space over the group G = $Z_{43}$.

In view of this we have the following theorems.

**THEOREM 4.3.1**: *Let V = $V_1 \cup V_2$ be a group interval bivector space over the group G = $Z_p$; p a prime then V is a pseudo simple group interval bivector space over the group G = $Z_p$.*

The proof is left as an exercise to the reader.

**THEOREM 4.3.2**: *Let $V = V_1 \cup V_2$ be a group interval bivector space over the group $G = Z_n$, n not a prime,*
1. *V in general is not a pseudo simple group interval bivector space over the group G*
2. *V is not a simple group interval bivector space over $G = Z_n$.*

This proof is also straight forward and hence left as an exercise for the reader to prove.

Now one can as in case of set interval bivector spaces define the notion of bilinear transformation of group interval bivector spaces. This task is also left as an exercise for the reader. Now we proceed onto define the notion of group interval bilinear algebras.

**DEFINITION 4.3.4**: *Let $V = V_1 \cup V_2$ be a group interval bivector space over the group G. We say V is a group interval bilinear algebra over the group G that is if both $V_1$ and $V_2$ are groups under addition.*

We will illustrate this situation by some examples.

*Example 4.3.12*: Let $V = V_1 \cup V_2 =$

$$\left\{ \begin{bmatrix} [0,a_1] & [0,a_2] & [0,a_3] & [0,a_4] \\ [0,a_5] & [0,a_6] & [0,a_7] & [0,a_8] \end{bmatrix} \middle| a_i \in Z_{12}; 1 \leq i \leq 8 \right\} \cup$$

$$\left\{ \sum_{i=0}^{\infty} [0,a_i] x^i \middle| a_i \in Z_{12} \right\}$$

be a group interval bilinear algebra over the group $G = Z_{12}$. V is of infinite order.

***Example 4.3.13***: Let $V = V_1 \cup V_2 =$

$$\left\{ \begin{bmatrix} [0,a_1] & [0,a_2] \\ [0,a_3] & [0,a_4] \\ [0,a_5] & [0,a_6] \\ [0,a_7] & [0,a_8] \end{bmatrix} \middle| a_i \in Z_7; 1 \le i \le 7 \right\} \cup$$

$$\left\{ \begin{bmatrix} [0,a_1] & [0,a_2] & [0,a_3] & [0,a_4] & [0,a_5] \\ [0,a_6] & [0,a_7] & [0,a_8] & [0,a_9] & [0,a_{10}] \end{bmatrix} \middle| a_i \in Z_7; 1 \le i \le 10 \right\}$$

be a group interval bilinear algebra over the group $G = Z_7$.
This V is of finite order.

Now we proceed onto give some properties enjoyed by them and define some substructures associated with them.

**THEOREM 4.3.3**: *Let $V = V_1 \cup V_2$ be a group interval bivector space over the group G; then in general V need not be a group interval bilinear algebra over the group G.*

The proof can be given by an appropriate example.

**THEOREM 4.3.4**: *Let $V = V_1 \cup V_2$ be a group interval bilinear algebra over a group G then V is a group interval bivector space over the group G.*

The proof directly follows from the definition of group interval bilinear algebras.

**DEFINITION 4.3.5:** *Let $V = V_1 \cup V_2$ be a group interval bilinear algebra over a group G. Let $W = W_1 \cup W_2 \subseteq V_1 \cup V_2 = V$; if W itself is a group interval bilinear algebra over the same group G then we define W to be a group interval bilinear subalgebra of V over G.*

*If V has no proper group interval bilinear subalgebra then we call V to be a simple group interval bilinear algebra.*

We will illustrate this situation by some examples.

*Example 4.3.14:* Let $V = V_1 \cup V_2 =$

$$\left\{ \begin{bmatrix} [0,a_1] & [0,a_2] & [0,a_3] & [0,a_4] \\ [0,a_5] & [0,a_6] & [0,a_7] & [0,a_8] \\ [0,a_9] & [0,a_{10}] & [0,a_{11}] & [0,a_{12}] \end{bmatrix} \,\middle|\, a_i \in Z_{11}; 1 \leq i \leq 12 \right\}$$

$$\cup \left\{ \sum_{i=0}^{\infty} [0, a_i] x^i \,\middle|\, a_i \in Z_{11} \right\}$$

be a group interval bilinear algebra over the group $G = Z_{11}$.
Take $W = W_1 \cup W_2 =$

$$\left\{ \begin{bmatrix} [0,a_1] & 0 & [0,a_3] & 0 \\ 0 & [0,a_4] & 0 & [0,a_2] \\ [0,a_5] & 0 & [0,a_6] & 0 \end{bmatrix} \,\middle|\, a_i \in Z_{11}; 1 \leq i \leq 6 \right\} \cup$$

$$\left\{ \sum_{i=0}^{\infty} [0, a_i] x^{2i} \,\middle|\, a_i \in Z_{11} \right\}$$

$\subseteq V_1 \cup V_2 = V$, W is a group interval bilinear subalgebra of V over the group $G = Z_{11}$.

*Example 4.3.15*: Let $V = V_1 \cup V_2 =$ {Collection of all $10 \times 10$ interval matrices with intervals of the form $[0, a_i]$ with $a_i \in Z_{20}$} $\cup$ {set of all $5 \times 5$ interval matrices with intervals of the form $[0, a_i]$ where $a_i \in Z_{20}$} be a group interval bilinear algebra over the group $G = Z_{20}$. Choose $W = W_1 \cup W_2 =$ {all $10 \times 10$ diagonal interval matrices with entries from $Z_{20}$ where intervals are of the form $[0, a_i]$} $\cup$ {all $5 \times 5$ upper triangular interval matrices with intervals of the form $[0, a_i]$; $a_i \in Z_{20}$} $\subseteq V_1 \cup V_2$; W is a group interval bilinear subalgebra of V over the group G.

*Example 4.3.16*: Let $V = V_1 \cup V_2 =$

$$\left\{ \begin{bmatrix} [0,a] & [0,a] & [0,a] \\ [0,a] & [0,a] & [0,a] \\ [0,a] & [0,a] & [0,a] \end{bmatrix} \middle| a \in Z_{13} \right\} \cup \left\{ \begin{bmatrix} [0,a] & [0,a] \\ [0,a] & [0,a] \\ [0,a] & [0,a] \\ [0,a] & [0,a] \\ [0,a] & [0,a] \end{bmatrix} \middle| a \in Z_{13} \right\}$$

be a group interval bilinear algebra over the group $G = Z_{13}$. We see V has no proper group interval bilinear subalgebras hence V is a simple group interval bilinear algebra over the group $G = Z_{13}$.

*Example 4.3.17*: Let $V = V_1 \cup V_2 =$

$$\left\{ \begin{bmatrix} [0,a] & [0,a] \\ [0,a] & [0,a] \end{bmatrix} \middle| a \in Z_3 \right\} \cup \left\{ \begin{bmatrix} [0,a] & [0,a] \\ [0,a] & [0,a] \\ [0,a] & [0,a] \\ [0,a] & [0,a] \end{bmatrix} \middle| a \in Z_3 \right\}$$

be a group interval bilinear algebra over the group $G = Z_3$. V is a simple group interval bilinear algebra over the group $G = Z_3$.

**DEFINITION 4.3.6**: *Let $V = V_1 \cup V_2$ be a group interval bilinear algebra over a group G. Let $W = W_1 \cup W_2 \subseteq V_1 \cup V_2$ be a proper bisubset of V and H a proper subgroup of G. If W is a group interval bilinear algebra over the group H then we define W to be a subgroup interval bilinear subalgebra of V over the subgroup H of G.*

We will illustrate this situation by some examples.

*Example 4.3.18*: Let $V = V_1 \cup V_2 =$

$$\left\{ \begin{bmatrix} [0,a_1] & [0,a_2] & [0,a_3] & [0,a_4] & [0,a_5] \\ [0,a_6] & [0,a_7] & [0,a_8] & [0,a_9] & [0,a_{10}] \end{bmatrix} \middle| a_i \in Z_{24}; 1 \le i \le 10 \right\}$$

$$\cup \left\{ \begin{bmatrix} [0,a_1] & [0,a_2] \\ [0,a_3] & [0,a_4] \\ [0,a_5] & [0,a_6] \\ [0,a_7] & [0,a_8] \\ [0,a_9] & [0,a_{10}] \\ [0,a_{10}] & [0,a_{12}] \end{bmatrix} \middle| a_i \in Z_{24}; 1 \le i \le 12 \right\}$$

be a group interval bilinear algebra over the group $G = Z_{24}$.

Choose $W = W_1 \cup W_2 =$

$$\left\{ \begin{bmatrix} [0,a] & 0 & [0,a] & 0 & [0,a] \\ [0,a] & 0 & [0,a] & 0 & [0,a] \end{bmatrix} \middle| a_i \in Z_{24} \right\} \cup$$

$$\left\{ \begin{bmatrix} [0,a] & 0 \\ 0 & [0,a] \\ [0,a] & 0 \\ 0 & [0,a] \\ [0,a] & 0 \\ 0 & [0,a] \end{bmatrix} \middle| a_i \in Z_{24} \right\}$$

$\subseteq V_1 \cup V_2$ and $H = \{0, 3, 6, 9, 12, 15, 18, 21\} \subseteq Z_{24} = G$ a proper subgroup of G under addition modulo 24. W is a subgroup interval bilinear subalgebra of V over the subgroup H of G.

***Example 4.3.19:*** Let $V = V_1 \cup V_2 =$

$$\left\{ \begin{bmatrix} [0,a_1] & [0,a_2] & [0,a_3] \\ [0,a_4] & [0,a_5] & [0,a_6] \\ [0,a_7] & [0,a_8] & [0,a_9] \end{bmatrix} \middle| a_i \in Z_{35}; 1 \le i \le 9 \right\} \cup$$

$$\left\{ \begin{bmatrix} [0,a_1] & [0,a_2] & [0,a_3] & [0,a_4] \\ [0,a_5] & [0,a_6] & [0,a_7] & [0,a_8] \\ [0,a_9] & [0,a_{10}] & [0,a_{11}] & [0,a_{12}] \\ [0,a_{13}] & [0,a_{14}] & [0,a_{15}] & [0,a_{16}] \end{bmatrix} \middle| a_i \in Z_{35}; 1 \le i \le 16 \right\}$$

be a group interval bilinear algebra over the group $G = Z_{35}$. Take $W = W_1 \cup W_2 =$

$$\left\{ \begin{bmatrix} [0,a] & [0,a] & [0,a] \\ 0 & [0,a] & [0,a] \\ 0 & 0 & [0,a] \end{bmatrix} \middle| a_i \in Z_{35} \right\} \cup$$

$$\left\{ \begin{bmatrix} [0,a_1] & [0,a_2] & [0,a_3] & [0,a_4] \\ [0,a_5] & [0,a_6] & [0,a_7] & [0,a_8] \\ [0,a_9] & [0,a_{10}] & [0,a_{11}] & [0,a_{12}] \\ [0,a_{13}] & [0,a_{14}] & [0,a_{15}] & [0,a_{16}] \end{bmatrix} \middle| a_i \in \{0,5,10,15,20,25,30\} \subseteq Z_{35} \right\}$$

$\subseteq V_1 \cup V_2 = V$ and $H = \{0, 7, 14, 21, 28\} \subseteq Z_{35}$ a proper subgroup of $Z_{35}$ under addition modulo 35. $W = W_1 \cup W_2$ is a subgroup interval bilinear subalgebra of V over the subgroup H of G.

*Example 4.3.20*: Let $V = V_1 \cup V_2 =$

$$\left\{ \begin{bmatrix} [0,a_1] & [0,a_2] \\ [0,a_3] & [0,a_4] \end{bmatrix} \middle| a_i \in Z_5; 1 \le i \le 4 \right\}$$

$\cup \{([0, a_1], [0, a_2], [0, a_3], [0, a_4], [0, a_5]) \mid a_i \in Z_5; 1 \le i \le 5\}$ be a group interval bilinear algebra over the group $G = Z_5$. As G has no proper subgroups we see V is a pseudo simple group interval bilinear algebra over the group $G = Z_5$. We see V has proper group interval bilinear subalgebras over G.

For take $W = W_1 \cup W_2 =$

$$\left\{ \begin{bmatrix} [0,a] & [0,a] \\ [0,a] & [0,a] \end{bmatrix} \middle| a \in Z_5 \right\} \cup$$

$\{([0, a], [0, a], [0, a], [0, a], [0, a]) \mid a \in Z_5\} \subseteq V_1 \cup V_2 = V$; W is a group interval bilinear subalgebra of V over the group $G = Z_5$. So V is not a simple group interval bilinear algebra. Thus V is not a doubly simple group interval bilinear algebra over the group G.

*Example 4.3.21*: Let $V = V_1 \cup V_2 =$

$$\left\{ \begin{bmatrix} [0,a] & [0,a] & [0,a] & [0,a] \\ [0,a] & [0,a] & [0,a] & [0,a] \end{bmatrix} \middle| a \in Z_{23} \right\} \cup$$

$$\left\{ \begin{bmatrix} [0,a] & [0,a] & [0,a] & [0,a] \\ [0,a] & [0,a] & [0,a] & [0,a] \\ [0,a] & [0,a] & [0,a] & [0,a] \\ [0,a] & [0,a] & [0,a] & [0,a] \\ [0,a] & [0,a] & [0,a] & [0,a] \end{bmatrix} \middle| a \in Z_{23} \right\}$$

be a group interval bilinear algebra over the group $G = Z_{23}$. V is a simple group interval bilinear algebra over the group $G = Z_{23}$, as V has no group interval bilinear subalgebras. Further V is a pseudo simple group interval bilinear algebra as V has no subgroup interval bilinear subalgebra. Thus V is a doubly simple group bilinear algebra over $G = Z_{23}$.

Now we can as in case of set interval bilinear algebras and semigroup interval bilinear algebras develop several properties about group interval bilinear algebras.

We just show the existence of a class of pseudo simple group interval bilinear algebras.

**THEOREM 4.3.5**: *Let $V = V_1 \cup V_2$ be a group interval bilinear algebra over the group $G = Z_p$; p a prime. V is a pseudo simple group interval bilinear algebra over the group G.*

*Proof:* Follows from the fact that G is a group which has no proper subgroups.

**THEOREM 4.3.6**: *Let $V = V_1 \cup V_2$ be a group interval linear algebra over the group $G = Z_n$, n not a prime and $V_1$ and $V_2$ take entries from $Z_n$. Then G is not a simple group interval bilinear algebra as well as G is not a pseudo simple group interval bilinear algebra.*

The proof is obvious from the fact that $Z_n$ has subgroups when n is not a prime and $V_1$ and $V_2$ constructed over $Z_n$ will certainly yield sub bispaces or sub bilinear algebras.

*Example 4.3.22:* Let $V = V_1 \cup V_2 =$

$$\left\{ \begin{bmatrix} [0,a] & [0,a] & [0,a] \\ [0,a] & [0,a] & [0,a] \end{bmatrix} \middle| a \in Z_{20} \right\} \cup$$

$$\left\{ \begin{bmatrix} [0,a] & [0,a] & [0,a] \\ [0,a] & [0,a] & [0,a] \\ [0,a] & [0,a] & [0,a] \\ [0,a] & [0,a] & [0,a] \end{bmatrix} \middle| a \in Z_{20} \right\}$$

be a group interval bilinear algebra over the group $G = Z_{20}$. Take $W = W_1 \cup W_2 =$

$$\left\{ \begin{bmatrix} [0,a] & [0,a] & [0,a] \\ [0,a] & [0,a] & [0,a] \end{bmatrix} \middle| a \in \{0,5,10,15\} \subseteq Z_{20} \right\} \cup$$

$$\left\{ \begin{bmatrix} [0,a] & [0,a] & [0,a] \\ [0,a] & [0,a] & [0,a] \\ [0,a] & [0,a] & [0,a] \\ [0,a] & [0,a] & [0,a] \end{bmatrix} \middle| a \in \{0,5,10,15\} \subseteq Z_{20} \right\}$$

$\subseteq V_1 \cup V_2$; W is a subgroup interval bilinear subalgebra over the subgroup H = {0, 5, 10, 15} $\subseteq Z_{20}$.

Now having seen some of the basic properties of group interval bilinear algebras, we can as in case of other bilinear algebras define bilinear transformation and bilinear operator.

We can define some more properties like quasi group bilinear algebra.

**DEFINITION 4.3.7:** *Let $V = V_1 \cup V_2$ be a group interval bivector space over the group G, if one of $V_1$ or $V_2$ (or in the mutually exclusive sense) is a group interval linear algebra then we define V to be a quasi group interval bilinear algebra over the group G.*

We will illustrate this situation by some examples.

*Example 4.3.23:* Let $V = V_1 \cup V_2 =$

$$\left\{ \begin{bmatrix} [0,a_1] & [0,a_2] & [0,a_3] \\ [0,a_4] & [0,a_5] & [0,a_6] \end{bmatrix}, \begin{bmatrix} [0,a] & [0,b] \\ [0,b] & [0,a] \\ [0,a] & [0,b] \\ [0,b] & [0,a] \\ [0,a] & [0,b] \\ [0,b] & [0,a] \\ [0,a] & [0,b] \\ [0,b] & [0,a] \\ [0,a] & [0,b] \end{bmatrix} \middle| \begin{array}{l} a_i, a, b \in Z_{45}; \\ 1 \leq i \leq 6 \end{array} \right\} \cup$$

$$\left\{ \sum_{i=0}^{\infty} [0, a_i] x^i \,\bigg|\, a_i \in Z_{45} \right\}$$

be a quasi group interval bilinear algebra over the group $G = Z_{45}$.

*Example 4.3.24:* Let $V = V_1 \cup V_2 =$

$$\left\{ \begin{bmatrix} [0,a] & [0,b] \\ [0,c] & [0,d] \\ [0,e] & [0,f] \end{bmatrix} \,\bigg|\, a,b,c,d,e,f \in Z_{240} \right\} \cup$$

$$\left\{ \sum_{i=0}^{\infty} [0, a_i] x^i, \left( [0,a_1] \quad [0,a_2] \quad \cdots \quad [0,a_9] \right) \,\bigg|\, a_1, a_2, \ldots, a_9 \in Z_{240} \right\}$$

be a quasi group interval bilinear algebra over the group $G = Z_{240}$.

Now we can as in case of group interval bilinear algebras define two types of substructures. We will however illustrate this situation by some examples.

*Example 4.3.25*: Let $V = V_1 \cup V_2 =$

$$\left\{ \sum_{i=0}^{\infty} [0, a_i] x^i \,\bigg|\, a_i \in Z_{14} \right\} \cup$$

$$\left\{ \begin{bmatrix} [0,a_1] & [0,a_2] \\ [0,a_3] & [0,a_4] \\ [0,a_5] & [0,a_6] \\ [0,a_7] & [0,a_8] \\ [0,a_9] & [0,a_{10}] \end{bmatrix}, \begin{bmatrix} [0,a_1] & [0,a_2] \\ [0,a_3] & [0,a_4] \end{bmatrix} \,\bigg|\, a_i \in Z_{14}; 1 \le i \le 10 \right\}$$

be a quasi group interval bilinear algebra over the group $G = Z_{14}$. Take $W = W_1 \cup W_2 =$

$$\left\{ \sum_{i=0}^{\infty} [0, a_i] x^{2i} \,\middle|\, a_i \in Z_{14} \right\} \cup$$

$$\left\{ \begin{bmatrix} [0,a] & [0,a] \\ [0,a] & [0,a] \\ [0,a] & [0,a] \\ [0,a] & [0,a] \\ [0,a] & [0,a] \end{bmatrix}, \begin{bmatrix} [0,a] & [0,a] \\ [0,a] & [0,a] \end{bmatrix} \,\middle|\, a_i \in Z_{14} \right\}$$

$\subseteq V_1 \cup V_2$ be a quasi group interval bilinear subalgebra of V over the group $G = Z_{14}$.

***Example 4.3.26***: Let $V = V_1 \cup V_2 =$

$$\left\{ \begin{bmatrix} [0,a_1] & [0,a_2] & [0,a_3] & [0,a_4] \\ [0,a_5] & [0,a_6] & [0,a_7] & [0,a_8] \\ [0,a_9] & [0,a_{10}] & [0,a_{11}] & [0,a_{12}] \\ [0,a_{13}] & [0,a_{14}] & [0,a_{15}] & [0,a_{16}] \end{bmatrix} \,\middle|\, a_i \in Z_{18}; 1 \leq i \leq 16 \right\} \cup$$

$$\left\{ \begin{bmatrix} [0,a_1] \\ [0,a_2] \\ [0,a_3] \\ [0,a_4] \end{bmatrix} ([0,a_1],[0,a_2],[0,a_3]) \,\middle|\, a_i \in Z_{18}; 1 \leq i \leq 4 \right\}$$

be a quasi group interval bilinear algebra over the group $G = Z_{18}$. Take $W = W_1 \cup W_2 =$

$$\left\{ \begin{bmatrix} [0,a_1] & [0,a_2] & [0,a_3] & [0,a_4] \\ [0,a_1] & [0,a_2] & [0,a_3] & [0,a_4] \\ [0,a_1] & [0,a_2] & [0,a_3] & [0,a_4] \\ [0,a_1] & [0,a_2] & [0,a_3] & [0,a_4] \end{bmatrix} \,\middle|\, a_i \in Z_{18}; 1 \leq i \leq 4 \right\} \cup$$

$$\left\{ \begin{bmatrix} [0,a] \\ [0,a] \\ [0,a] \\ [0,a] \end{bmatrix}, \big([0,a],[0,a],[0,a]\big) \,\Big|\, a_i \in Z_{18} \right\}$$

$\subseteq V_1 \cup V_2$ be a quasi group interval bilinear subalgebra of V over the group $G = Z_{18}$.

*Example 4.3.27*: Let $V = V_1 \cup V_2 =$

$$\left\{ \begin{bmatrix} [0,a] & [0,a] & [0,a] \\ [0,a] & [0,a] & [0,a] \\ [0,a] & [0,a] & [0,a] \end{bmatrix} \,\Big|\, a \in Z_7 \right\} \cup$$

$$\left\{ [0,a], \begin{bmatrix} [0,a] & [0,a] \\ [0,a] & [0,a] \\ [0,a] & [0,a] \\ [0,a] & [0,a] \end{bmatrix} \,\Big|\, a \in Z_7 \right\}$$

is a quasi group interval bilinear algebra over the group G. We see V is a simple quasi group interval bilinear algebra as V has no quasi group interval bilinear subalgebras.

*Example 4.3.28:* Let $V = V_1 \cup V_2 =$

$$\left\{ \begin{bmatrix} [0,a] & [0,a] & [0,a] & [0,a] \\ [0,a] & [0,a] & [0,a] & [0,a] \\ [0,a] & [0,a] & [0,a] & [0,a] \end{bmatrix} \,\Big|\, a \in Z_{47} \right\} \cup$$

$$\left\{ \begin{bmatrix} [0,a] \\ [0,a] \\ [0,a] \\ [0,a] \\ [0,a] \\ [0,a] \\ [0,a] \end{bmatrix}, ([0,a] \quad [0,a]) \;\middle|\; a \in Z_{47} \right\}$$

be a quasi group interval bilinear algebra over the group $G = Z_{47}$. Clearly V is a simple quasi group interval bilinear algebra over the group $G = Z_{47}$.

Next as in case of group interval bilinear algebras define the same notion in case of quasi group interval bilinear algebras. We will only illustrate this situation by some examples and the task of giving the definition is left as an exercise to the reader.

*Example 4.3.29*: Let $V = V_1 \cup V_2 =$

$$\left\{ \begin{bmatrix} [0,a_1] & [0,a_2] & [0,a_3] & [0,a_4] \\ [0,a_5] & [0,a_6] & [0,a_7] & [0,a_8] \end{bmatrix} \;\middle|\; a \in Z_{48}; 1 \le i \le 8 \right\} \cup$$

$$\left\{ \begin{bmatrix} [0,a_1] \\ [0,a_2] \\ [0,a_3] \\ [0,a_4] \\ [0,a_5] \\ [0,a_6] \end{bmatrix}, \sum_{i=0}^{\infty} [0,a_i] x^i \;\middle|\; a_i \in Z_{48}; 1 \le i \le 6 \right\}$$

be a quasi group interval bilinear algebra over the group $G = Z_{48}$. Take $H = \{0, 2, 4, 6, 8, 10, 12, \ldots, 46\} \subseteq Z_{48}$ to be a proper subgroup of $G = Z_{48}$. Choose $W = W_1 \cup W_2 =$

$$\left\{ \begin{bmatrix} [0,a_1] & [0,a_2] & [0,a_3] & [0,a_4] \\ [0,a_5] & [0,a_6] & [0,a_7] & [0,a_8] \end{bmatrix} \,\middle|\, a_i \in H; 1 \le i \le 8 \right\} \cup$$

$$\left\{ \sum_{i=0}^{\infty} [0,a_i] x^i \,\middle|\, a_i \in Z_{48} \right\}$$

$\subseteq V_1 \cup V_2$, W is a quasi subgroup interval bilinear subalgebra of V over the subgroup $H \subseteq Z_{48}$.

We call a quasi group interval bilinear algebra which has no quasi subgroup bilinear subalgebra to be pseudo quasi simple group interval bilinear algebra.

*Example 4.3.30*: Let $V = V_1 \cup V_2 =$

$$\left\{ \begin{bmatrix} [0,a_1] & [0,a_4] \\ [0,a_2] & [0,a_5] \\ [0,a_3] & [0,a_6] \end{bmatrix}, \big([0,a_1] \;\; [0,a_2] \;\; [0,a_3]\big) \,\middle|\, a_i \in Z_{13}; 1 \le i \le 6 \right\} \cup$$

$$\left\{ \sum_{i=0}^{\infty} [0,a_i] x^i \,\middle|\, a_i \in Z_{13} \right\}$$

be a quasi group interval bilinear algebra over the group $G = Z_{13}$. G has no proper subgroups, hence V is a pseudo simple quasi group interval bilinear algebra over the group $G = Z_{13}$. However V is not a simple quasi group interval bilinear algebra over the group $G = Z_{13}$. For take $W =$

$$\left\{ \begin{bmatrix} [0,a] & [0,a] \\ [0,a] & [0,a] \\ [0,a] & [0,a] \end{bmatrix}, \big([0,a] \;\; [0,a] \;\; [0,a]\big) \,\middle|\, a \in Z_{13} \right\} \cup$$

$$\left\{ \sum_{i=0}^{\infty} [0,a_i] x^{2i} \,\middle|\, a_i \in Z_{13} \right\}$$

$\subseteq V_1 \cup V_2$ is a quasi subgroup interval bilinear subalgebra of V over the group $G = Z_{13}$.

Thus V is not a doubly simple quasi group interval bilinear algebra over the group $G = Z_{13}$.

*Example 4.3.31*: Let $V = V_1 \cup V_2 = \{[0, a] \mid a \in Z_{43}\} \cup$

$$\{([0, a], [0, a]), \begin{bmatrix} [0,a] \\ [0,a] \\ [0,a] \end{bmatrix} \mid a \in Z_{43}\}$$

be a quasi group interval bilinear algebra over the group $G = Z_{43}$. V is a doubly simple quasi group interval bilinear algebra over the group $G = Z_{43}$.

*Example 4.3.32*: Let $V = V_1 \cup V_2 = \{([0, a], [0, a], [0, a], [0, a], [0, a]) \mid a \in Z_{47}\} \cup$

$$\{\begin{bmatrix} [0,a] & [0,a] \\ [0,a] & [0,a] \end{bmatrix}, ([0, a] [0, a]) \mid a \in Z_{47}\}$$

be a quasi group interval bilinear algebra over the group $G = Z_{47}$. V is a doubly simple quasi group interval bilinear algebra over the group $G = Z_{47}$.

**THEOREM 4.3.7**: *Let $V = V_1 \cup V_2$ be a quasi group interval bilinear algebra over the group $G = Z_p$; p a prime. Then V is a pseudo simple quasi group interval bilinear algebra over the group $G = Z_p$.*

Proof is straight forward and hence left as an exercise for the reader.

We see in general all quasi group interval bilinear algebra over $Z_p$, p a prime need not be a simple quasi group interval bilinear algebra over $Z_p$.

***Example 4.3.33***: Let $V = V_1 \cup V_2 = \{([0, a_1], [0, a_2], [0, a_3], [0, a_4], [0, a_5]) \mid a_i \in Z_7, 1 \le i \le 5\} \cup$

$$\left\{ \begin{bmatrix} [0,a_1] & [0,a_2] \\ [0,a_3] & [0,a_4] \end{bmatrix}, \begin{bmatrix} [0,a_1] \\ [0,a_2] \\ [0,a_3] \\ [0,a_4] \\ [0,a_5] \end{bmatrix} \middle| a_i \in Z_7; 1 \le i \le 5 \right\}$$

be a quasi group interval bilinear algebra over the group $G = Z_7$. Take $W = \{([0, a_1], [0, a_2], 0, 0, [0, a_3]) \mid a_i \in Z_7; 1 \le i \le 3\} \cup$

$$\left\{ \begin{bmatrix} [0,a] & [0,a] \\ 0 & [0,a] \end{bmatrix}, \begin{bmatrix} [0,a] \\ 0 \\ [0,a] \\ 0 \\ [0,a] \end{bmatrix} \middle| a \in Z_7 \right\}$$

$= W_1 \cup W_2 \subseteq V_1 \cup V_2$, W is a quasi group interval bilinear subalgebra of V over the group $G = Z_7$. Thus V is not a simple quasi group interval bilinear algebra over the group $G = Z_7$. Infact V has several such quasi group interval bilinear subalgebras.

Now we have a class of quasi group interval bilinear algebras which are not simple or pseudo simple. We illustrate this by the following theorem.

**THEOREM 4.3.8:** *Let $V = V_1 \cup V_2$ be a quasi group interval bilinear algebra over the group $G = Z_n$; n not a prime. V in general is not a doubly simple quasi group interval bilinear algebra over $Z_n = G$.*

Proof is straight forward and is left as an exercise for the reader.
    Now we proceed onto define yet another new structure.

## 4.4 Bisemigroup Interval Bilinear Algebras and Their Generalization

In this section we define the notion of bisemigroup interval bilinear algebras, bigroup interval bilinear algebras, set-semigroup interval bilinear algebras set group interval bilinear algebras and semigroup group interval bilinear algebras and describe some of their properties.

**DEFINITION 4.4.1:** *Let $V = V_1 \cup V_2$ where $V_i$ is a semigroup interval vector space over the semigroup $S_i$, $i = 1, 2$, such that $V_1 \not\subseteq V_2$, $V_2 \not\subseteq V_1$ and $S_1 \neq S_2$ $S_1 \neq S_2$ and $S_2 \not\subseteq S_1$. We define $V = V_1 \cup V_2$ to be a bisemigroup interval bivector space over the bisemigroup $S = S_1 \cup S_2$.*

We will illustrate this by some examples.

*Example 4.4.1*: Let $V = V_1 \cup V_2 =$

$$\left\{ \begin{bmatrix} [0,a_1] & [0,a_2] & [0,a_3] \\ [0,a_4] & [0,a_5] & [0,a_6] \end{bmatrix}, \begin{bmatrix} [0,a_1] \\ [0,a_2] \\ [0,a_3] \\ [0,a_4] \end{bmatrix} \middle| a_i \in Z^+ \cup \{0\} \right\} \cup$$

$$\left\{ ([0,a_1] \ [0,a_2] \ [0,a_3] \ [0,a_4]), \begin{bmatrix} [0,a_1] & [0,a_2] & [0,a_3] \\ [0,a_4] & [0,a_5] & [0,a_6] \\ [0,a_7] & [0,a_8] & [0,a_9] \end{bmatrix} \middle| \begin{matrix} a_i \in Z_{10}; \\ 1 \leq i \leq 9 \end{matrix} \right\}$$

be a bisemigroup interval bivector space over the bisemigroup $S = S_1 \cup S_2 = Z^+ \cup \{0\} \cup Z_{10}$.

*Example 4.4.2*: Let $V = V_1 \cup V_2 =$

$$\left\{ \begin{bmatrix} [0,a_1] & [0,a_2] \\ [0,a_3] & [0,a_4] \\ [0,a_5] & [0,a_6] \end{bmatrix}, ([0,a_1] \quad [0,a_2]) \,\middle|\, a_i \in Z_7; 1 \le i \le 6 \right\} \cup$$

$$\left\{ \begin{bmatrix} [0,a_1] \\ [0,a_2] \\ \vdots \\ [0,a_7] \end{bmatrix}, ([0,a_1] \quad [0,a_2] \quad [0,a_3] \quad [0,a_4]) \,\middle|\, \begin{array}{l} a_i \in Z_{92}; \\ 1 \le i \le 7 \end{array} \right\}$$

be a bisemigroup interval bivector space over the bisemigroup $S = S_1 \cup S_2 = Z_7 \cup Z_{92}$.

Now if in the definition 4.4.1 each $V_i$ is a semigroup interval linear algebra over $S_i$, $i = 1, 2$ then we call V to be a bisemigroup interval bilinear algebra over the bisemigroup $S = S_1 \cup S_2$.

We will illustrate this situation by some examples.

*Example 4.4.3*: Let $V = V_1 \cup V_2 =$

$$\left\{ \begin{bmatrix} [0,a_1] & [0,a_2] & [0,a_3] & [0,a_4] \\ [0,a_5] & [0,a_6] & [0,a_7] & [0,a_8] \\ [0,a_9] & [0,a_{10}] & [0,a_{11}] & [0,a_{12}] \end{bmatrix} \,\middle|\, a_i \in Q^+ \cup \{0\} \right\} \cup$$

$$\left\{ \sum_{i=0}^{\infty} [0,a_i] x^i \,\middle|\, a_i \in Z_{243} \right\}$$

be a semigroup interval bilinear algebra over the bisemigroup $S = S_1 \cup S_2 = (Q^+ \cup \{0\}) \cup Z_{243}$.

***Example 4.4.4:*** Let $V = V_1 \cup V_2 = \{$all $10 \times 10$ interval matrices with intervals of the form $[0, a_i]$ with $a_i \in R^+ \cup \{0\}\} \cup$

$$\left\{ \begin{bmatrix} [0,a] & [0,a] \\ [0,b] & [0,b] \end{bmatrix} \middle| a,b \in Z_8 \right\}$$

be a bisemigroup interval bilinear algebra with the bisemigroup $S = R^+ \cup \{0\} \cup Z_8$.

Both the bisemigroup interval bilinear algebras in example 4.4.3 and 4.4.4 are of infinite cardinality.

***Example 4.4.5:*** Let $V = V_1 \cup V_2 =$

$$\left\{ \begin{bmatrix} [0,a_1] & [0,a_2] \\ [0,a_3] & [0,a_4] \\ [0,a_5] & [0,a_6] \end{bmatrix} \middle| a_i \in Z_{17}; 1 \le i \le 6 \right\} \cup$$

$$\left\{ \begin{bmatrix} [0,a_1] & [0,a_2] & [0,a_3] & [0,a_4] & [0,a_5] \\ [0,a_6] & [0,a_7] & [0,a_8] & [0,a_9] & [0,a_{10}] \end{bmatrix} \middle| \begin{array}{l} a_i \in Z_{102}; \\ 1 \le i \le 10 \end{array} \right\}$$

be a bisemigroup interval bilinear algebra over the bisemigroup $S = Z_{17} \cup Z_{102}$. We see the bisemigroup interval bilinear algebra given in example 4.4.5 is of finite order.

We have as in case of other bilinear algebras the following theorem to be true.

**THEOREM 4.4.1:** *Let $V = V_1 \cup V_2$ be a bisemigroup interval bilinear algebra over the bisemigroup $S = S_1 \cup S_2$. V is clearly a bisemigroup interval bivector space over the bisemigroup $S = S_1 \cup S_2$; however if V is a bisemigroup interval bivector space over the bisemigroup $S = S_1 \cup S_2$; V need not in general be a bisemigroup interval bilinear algebra over the bisemigroup $S = S_1 \cup S_2$.*

The proof is simple and straight forward so left as an exercise to the reader.

Now as in case of other bisemigroup linear algebras and bisemigroup vector spaces we can define substructures.

Here we only illustrate these situations by some examples.

*Example 4.4.6*: Let $V = V_1 \cup V_2 =$

$$\left\{ \begin{bmatrix} [0,a_1] & [0,a_2] \\ [0,a_3] & [0,a_4] \end{bmatrix}, \begin{bmatrix} [0,a_1] \\ [0,a_2] \\ [0,a_3] \\ [0,a_4] \end{bmatrix} \,\middle|\, a_i \in Z_{13}; 1 \leq i \leq 4 \right\} \cup$$

$$\left\{ ([0,a_1]\ [0,a_2]\ [0,a_3]\ [0,a_4]), \begin{bmatrix} [0,a_1] & [0,a_2] & [0,a_3] \\ [0,a_4] & [0,a_5] & [0,a_6] \\ [0,a_7] & [0,a_8] & [0,a_9] \end{bmatrix} \,\middle|\, \begin{array}{l} a_i \in Z_{18}; \\ 1 \leq i \leq 9 \end{array} \right\}$$

be a bisemigroup interval bivector space over the bisemigroup $S = S_1 \cup S_2 = Z_{13} \cup Z_{18}$.

Take $W = W_1 \cup W_2 =$

$$\left\{ \begin{bmatrix} [0,a] & [0,a] \\ [0,a] & [0,a] \end{bmatrix}, \begin{bmatrix} [0,a_1] \\ 0 \\ [0,a_2] \\ 0 \end{bmatrix} \,\middle|\, a_1, a_2 \in Z_{13} \right\} \cup$$

$$\left\{ ([0,a_1]\ [0,a_2]\ 0\ 0), \begin{bmatrix} [0,a] & [0,a] & 0 \\ 0 & 0 & [0,a] \\ 0 & 0 & 0 \end{bmatrix} \,\middle|\, a, a_1, a_2 \in Z_{18} \right\}$$

$\subseteq V = V_1 \cup V_2$ is a bisemigroup interval bivector subspace of V over the bisemigroup $S = S_1 \cup S_2 = Z_{13} \cup Z_{18}$.

***Example 4.4.7:*** Let $V = V_1 \cup V_2 =$

$$\left\{ \begin{bmatrix} [0,a_1] & [0,a_2] \\ [0,a_3] & [0,a_4] \\ [0,a_5] & [0,a_6] \\ [0,a_7] & [0,a_8] \end{bmatrix}, [0,a_i] \middle| \begin{array}{l} a_i, a_j \in Z^+ \cup \{0\}; \\ 1 \le j \le 8; \quad i, j \in Z^+ \cup \{0\} \end{array} \right\} \cup$$

$$\left\{ \sum_{i=0}^{\infty} [0,a_i]x^i; \begin{bmatrix} [0,a_1] \\ [0,a_2] \\ [0,a_3] \\ [0,a_4] \\ [0,a_5] \end{bmatrix} \middle| \begin{array}{l} a_j, a_o \in Z_{11}; \\ 1 \le j \le 5 \end{array} \right\}$$

be a bisemigroup interval bivector space over the bisemigroup $S = S_1 \cup S_2 = (Z^+ \cup \{0\}) \cup \{Z_{11}\}$.

Let $W = W_1 \cup W_2 =$

$$\left\{ \begin{bmatrix} [0,a] & [0,a] \\ 0 & 0 \\ [0,a] & [0,a] \\ 0 & 0 \end{bmatrix}, [0,a] \middle| a \in Z^+ \cup \{0\} \right\} \cup$$

$$\left\{ \sum_{i=0}^{\infty} [0,a_i]x^{2i}; \begin{bmatrix} [0,a_1] \\ 0 \\ [0,a_2] \\ 0 \\ [0,a_3] \end{bmatrix} \middle| a_i, a_1, a_2, a_3 \in Z_{11} \right\}$$

$\subseteq V_1 \cup V_2$ is a bisemigroup interval bivector subspace of V over the bisemigroup $S = (Z^+ \cup \{0\}) \cup Z_{11}$.

***Example 4.4.8***: Let $V = V_1 \cup V_2 =$

$$\left\{ \begin{bmatrix} [0,a_1] & [0,a_2] & [0,a_3] \\ [0,a_4] & [0,a_5] & [0,a_6] \\ [0,a_7] & [0,a_8] & [0,a_9] \end{bmatrix} \middle| a_i \in Z_{15}; 1 \leq i \leq 9 \right\} \cup$$

$$\left\{ \begin{bmatrix} [0,a_1] \\ [0,a_2] \\ [0,a_3] \\ [0,a_4] \\ [0,a_5] \\ [0,a_6] \end{bmatrix} \middle| a_i \in Z_{18}; 1 \leq i \leq 6 \right\}$$

be a bisemigroup interval bilinear algebra over the bisemigroup $S = S_1 \cup S_2 = Z_{15} \cup Z_{18}$.

Take $W = W_1 \cup W_2 =$

$$\left\{ \begin{bmatrix} [0,a_1] & 0 & 0 \\ 0 & [0,a_2] & 0 \\ 0 & 0 & [0,a_3] \end{bmatrix} \middle| a_i \in Z_{15}; 1 \leq i \leq 3 \right\} \cup$$

$$\left\{ \begin{bmatrix} [0,a_1] \\ 0 \\ [0,a_2] \\ 0 \\ [0,a_3] \\ 0 \end{bmatrix} \middle| a_i \in Z_{18}; 1 \leq i \leq 3 \right\}$$

$\subseteq V_1 \cup V_2 = V$, W is a bisemigroup interval bilinear subalgebra of V over the bisemigroup $S = S_1 \cup S_2 = Z_{15} \cup Z_{18}$.

***Example 4.4.9***: Let $V = V_1 \cup V_2 = \{$All $12 \times 12$ interval matrices with intervals of the form $[0, a_i]$; $a_i \in Z_{48}\} \cup$

$$\left\{\sum_{i=0}^{\infty}[0,a_i]x^i \,\middle|\, a_i \in Z_{40}\right\}$$

be a bisemigroup interval bilinear algebra over the bisemigroup $S = S_1 \cup S_2 = Z_{48} \cup Z_{40}$. Choose $W = W_1 \cup W_2 = \{$All $12 \times 12$ upper triangular interval matrices with intervals of the form $[0, a_i]$; $a_i \in Z_{48}\} \cup$

$$\left\{\sum_{i=0}^{\infty}[0,a_i]x^{2i} \,\middle|\, a_i \in Z_{40}\right\} \subseteq V_1 \cup V_2;$$

W is a bisemigroup interval bilinear subalgebra of V over the bisemigroup S.

Now one can define bisubsemigroup interval bivector subspaces and bisubsemigroup interval bilinear subalgebras.

The task of defining these notions are left as an exercise to the reader.

We will however illustrate these situations by some examples.

***Example 4.4.10***: Let $V = V_1 \cup V_2 =$

$$\left\{ \begin{bmatrix} [0,a_1] & [0,a_2] & [0,a_3] & [0,a_4] \\ [0,a_5] & [0,a_6] & [0,a_7] & [0,a_8] \end{bmatrix}, \begin{bmatrix} [0,a_1] \\ [0,a_2] \\ [0,a_3] \\ [0,a_4] \\ [0,a_5] \end{bmatrix} \,\middle|\, \begin{array}{l} a_i \in Z_{24}; \\ 1 \le i \le 8 \end{array} \right\}$$

$$\cup \left\{ \sum_{i=0}^{\infty}[0,a_i]x^i, \begin{bmatrix} [0,a_1] & [0,a_2] \\ [0,a_3] & [0,a_4] \\ [0,a_5] & [0,a_6] \\ [0,a_7] & [0,a_8] \end{bmatrix} \,\middle|\, \begin{array}{l} a_j, a_i \in Z_{150}; \\ 1 \le i \le \infty, \\ 1 \le j \le 8 \end{array} \right\}$$

be a bisemigroup interval bivector space over the bisemigroup $S = S_1 \cup S_2 = Z_{24} \cup Z_{150}$.

Take $W = W_1 \cup W_2 =$

$$\left\{ \begin{bmatrix} [0,a] & 0 & [0,a] & 0 \\ [0,a] & 0 & [0,a] & 0 \end{bmatrix}, \begin{bmatrix} [0,a_1] \\ 0 \\ [0,a_2] \\ 0 \\ [0,a_3] \end{bmatrix} \middle| a, a_1, a_2, a_3 \in Z_{24} \right\} \cup$$

$$\left\{ \sum_{i=0}^{\infty} [0,a_i]x^{2i}, \begin{bmatrix} [0,a] & 0 \\ 0 & [0,a] \\ [0,a] & 0 \\ 0 & [0,a] \end{bmatrix} \middle| a_i, a \in Z_{150}; 1 \le i \le \infty \right\}$$

$\subseteq V_1 \cup V_2 = V$ and $T = T_1 \cup T_2 = \{0, 3, 6, 9, 12, 15, 18, 21\} \cup \{0, 10, 20, \ldots, 140\} \subseteq S_1 \cup S_2$.

It is easily verified W is a bisubsemigroup interval bivector subspace of V over the bisubsemigroup $T = T_1 \cup T_2$ of $S = S_1 \cup S_2$.

***Example 4.4.11***: Let $V = V_1 \cup V_2 = $ {All $9 \times 9$ interval matrices with intervals of the form $[0, a_i]$; $a_i \in Z^+ \cup \{0\}$} $\cup$ {All $4 \times 5$ interval matrices with intervals of the form $[0, a_i]$; $a_i \in Z_{14}$} be a bisemigroup interval bilinear algebra over the bisemigroup $S = S_1 \cup S_2 = Z^+ \cup \{0\} \cup Z_{14}$.

Take $T = T_1 \cup T_2 = \{3Z^+ \cup \{0\}\} \cup \{0, 2, 4, 6, 8, 10, 12\} \subseteq S_1 \cup S_2 = Z^+ \cup \{0\} \cup Z_{14}$ and $W = W_1 \cup W_2 = $ {All $9 \times 9$ interval matrices with intervals of the form $[0, a_i]$ with $a_i \in 7Z^+ \cup \{0\}\} \cup$

$$\left\{ \begin{bmatrix} [0,a] & 0 & [0,a] & 0 & [0,a] \\ 0 & [0,a] & 0 & [0,a] & 0 \\ [0,a] & 0 & [0,a] & 0 & [0,a] \\ 0 & [0,a] & 0 & [0,a] & 0 \end{bmatrix} \,\middle|\, a \in Z_{14} \right\}$$

$\subseteq V_1 \cup V_2$.

Clearly W is a bisubsemigroup interval bilinear subalgebra of V over the bisubsemigroup $T = T_1 \cup T_2 = 3Z^+ \cup \{0\} \cup 2Z_{14}$ $\subseteq V_1 \cup V_2$.

If V has no proper bisubsemigroup interval bilinear subalgebras we call V to be a simple bisemigroup interval bilinear algebra. If V has no proper bisubsemigroup interval bilinear subalgebras we call V to be a pseudo simple bisemigroup interval bilinear algebra. If V is both simple and pseudo simple then we call V to be a doubly simple bisemigroup interval bilinear algebra.

We will illustrate all these three situations by examples.

***Example 4.4.12:*** Let $V = V_1 \cup V_2 =$ {All $8 \times 8$ upper triangular interval matrices with entries from $[0, a_i]$ with $a_i \in Z_7$} $\cup$

$$\left\{ \sum_{i=0}^{\infty} [0, a_i] x^i \,\middle|\, a_i \in Z_{19} \right\}$$

be a bisemigroup interval bilinear algebra over the bisemigroup $S = S_1 \cup S_2 = Z_7 \cup Z_{19}$. V is a pseudo simple bisemigroup interval bilinear algebra but is not a simple bisemigroup interval bilinear algebra as $W = W_1 \cup W_2 =$ {all diagonal interval matrices with intervals of the form $[0, a_i] \mid a_i \in Z_7$} $\cup$

$$\left\{ \sum_{i=0}^{\infty} [0, a_i] x^{2i} \,\middle|\, a_i \in Z_{19} \right\}$$

⊆ $V_1 \cup V_2$ is bisemigroup interval bilinear algebra over the bisemigroup $S = Z_7 \cup Z_{19}$ so V is not simple. As $S = Z_7 \cup Z_{19}$ has no proper subsemigroups V is pseudo simple. Thus V is not doubly simple.

*Example 4.4.13*: Let $V = V_1 \cup V_2 =$

$$\left\{ \begin{bmatrix} [0,a] & [0,a] \\ [0,a] & [0,a] \end{bmatrix} \middle| a \in Z_7 \right\} \cup$$

$$\left\{ \begin{bmatrix} [0,a] & [0,a] & [0,a] & [0,a] \\ [0,a] & [0,a] & [0,a] & [0,a] \\ [0,a] & [0,a] & [0,a] & [0,a] \end{bmatrix} \middle| a \in Z_{11} \right\}$$

be a bisemigroup interval bilinear algebra over the bisemigroup $S = S_1 \cup S_2 = \{Z_7\} \cup \{Z_{11}\}$. We see V is a doubly simple bisemigroup interval bilinear algebra over the bisemigroup $S = S_1 \cup S_2 = Z_7 \cup Z_{11}$.

We have a class of pseudo simple bisemigroup interval bilinear algebras over a bisemigroup $S = S_1 \cup S_2$.

**THEOREM 4.4.2**: *Let $V = V_1 \cup V_2$ be a bisemigroup interval bilinear algebra over the bisemigroup $S = S_1 \cup S_2 = Z_p \cup Z_q$ where p and q two distinct primes. Then V is a pseudo simple bisemigroup interval bilinear algebra over S. V need not in general be simple.*

The proof is left as an exercise to the reader.

**THEOREM 4.4.3**: *Let $V = V_1 \cup V_2$ be a bisemigroup interval bilinear algebra over the bisemigroup $S = S_1 \cup S_2 = Z_m \cup Z_n$ where $m \neq n$ and m and n are not primes. Then V not in general a doubly simple semigroup interval bilinear algebra over the bisemigroup $S = S_1 \cup S_2 = Z_m \cup Z_n$.*

The proof is left as an exercise for the reader.

**THEOREM 4.4.4:** *Let $V = V_1 \cup V_2$ be a bisemigroup interval bilinear algebra over the bisemigroup $S = S_1 \cup S_2$ where one of $S_1$ is $Z^+ \cup \{0\}$ or $Q^+ \cup \{0\}$ or $R^+ \cup \{0\}$ and other $S_2$ is $Z_n$ or some subsemigroup of $Z^+ \cup \{0\}$ or $Q^+ \cup \{0\}$ or $R^+ \cup \{0\}$ such that $S_1 \not\subseteq S_2$ and $S_2 \not\subseteq S_1$ then $V$ is not a doubly simple bisemigroup interval bilinear algebra over the bisemigroup.*

This proof is also left for the reader.

We will give some illustrative examples.

*Example 4.4.14*: Let $V = V_1 \cup V_2 =$

$$\left\{ \begin{bmatrix} [0,a] & [0,a] & [0,a] \\ [0,a] & [0,a] & [0,a] \\ [0,a] & [0,a] & [0,a] \end{bmatrix} \middle| a \in Z^+ \cup \{0\} \right\} \cup$$

$$\left\{ \begin{bmatrix} [0,a] \\ [0,a] \\ [0,a] \\ [0,a] \\ [0,a] \\ [0,a] \\ [0,a] \end{bmatrix} \middle| a \in Z_{45} \right\}$$

be a bisemigroup interval bilinear algebra over the bisemigroup $S = S_1 \cup S_2 = Z^+ \cup \{0\} \cup Z_{45}$.
Take $W = W_1 \cup W_2$

$$= \left\{ \begin{bmatrix} [0,a] & [0,a] & [0,a] \\ [0,a] & [0,a] & [0,a] \\ [0,a] & [0,a] & [0,a] \end{bmatrix} \middle| a \in 3Z^+ \cup \{0\} \right\} \cup$$

$$\left\{ \begin{bmatrix} [0,a] \\ [0,a] \\ [0,a] \\ [0,a] \\ [0,a] \\ [0,a] \\ [0,a] \end{bmatrix} \,\middle|\, a \in \{0,5,10,15,20,25,30,35,40\} \subseteq Z_{45} \right\}$$

$\subseteq V_1 \cup V_2$ is bisubsemigroup interval bilinear subalgebra of V over the bisubsemigroup $T = T_1 \cup T_2 = 3Z^+ \cup \{0\} \cup \{0, 5, 10, 15, 20, 25, 30, 35, 40\} \subseteq S_1 \cup S_2 = Z^+ \cup \{0\} \cup Z_{45}$. Thus V is not a pseudo simple bisemigroup interval bilinear subalgebra of V over the bisubsemigroup $T = T_1 \cup T_2 \subseteq S_1 \cup S_2$.

Also W is a bisemigroup interval bilinear subalgebra of V over the bisemigroup $S = S_1 \cup S_2$ so, V is not a simple bisemigroup interval bilinear algebra over the bisemigroup $S = S_1 \cup S_2 = Z^+ \cup \{0\} \cup Z_{45}$.

***Example 4.4.15***: Let $V = V_1 \cup V_2 =$

$$\left\{ \begin{bmatrix} [0,a_1] & [0,a_2] \\ [0,a_3] & [0,a_4] \\ [0,a_5] & [0,a_6] \end{bmatrix} \,\middle|\, a_i \in 3Z^+ \cup \{0\} \right\} \cup$$

$\{([0, a_1], [0, a_2], \ldots, [0, a_{10}]) \mid a_i \in 7Z^+ \cup \{0\}\}$ be a bisemigroup interval bilinear algebra over the bisemigroup $S = S_1 \cup S_2 = \{3Z^+ \cup \{0\}\} \cup \{7Z^+ \cup \{0\}\}$. It is easy to verify V is not a doubly simple bisemigroup interval bilinear algebra over the bisemigroup S.

Now we define set- semigroup interval bivector space over the biset $S = S_1 \cup S_2$ where one of $S_i$ is a semigroup and other is a set, $i = 1, 2$.

**DEFINITION 4.4.2**: *Let $V = V_1 \cup V_2$ be such that $V_1$ is a set interval vector space over the set $S_1$ and $V_2$ is a semigroup interval vector space over the semigroup $S_2$. Then we define $V = V_1 \cup V_2$ to be a set-semigroup interval bivector space over the set-semigroup $S = S_1 \cup S_2$.*

We will illustrate this situation by some examples.

*Example 4.4.16:* Let $V = V_1 \cup V_2 =$

$$\left\{ \begin{bmatrix} [0,a_1] & [0,a_2] \\ [0,a_3] & [0,a_4] \\ [0,a_5] & [0,a_6] \end{bmatrix}, ([0,a_1] \; [0,a_2] \; [0,a_3]) \;\middle|\; a_i \in \{0,1,2,4,8,12,15\} \right\}$$

$$\cup \left\{ \begin{bmatrix} [0,a_1] \\ [0,a_2] \\ [0,a_3] \\ [0,a_4] \end{bmatrix}, \begin{bmatrix} [0,a_1] & [0,a_2] & [0,a_3] \\ [0,a_4] & [0,a_5] & [0,a_6] \\ [0,a_7] & [0,a_8] & [0,a_9] \\ [0,a_{10}] & [0,a_{11}] & [0,a_{12}] \end{bmatrix} \;\middle|\; \begin{array}{l} a_i \in Z^+ \cup \{0\}, \\ 1 \leq i \leq 12 \end{array} \right\}$$

be a set-semigroup interval bivector space over the set-semigroup $S = S_1 \cup S_2 = \{0, 1\} \cup \{3Z^+ \cup \{0\}\}$.

*Example 4.4.17:* Let $V = V_1 \cup V_2 =$

$$\left\{ \begin{bmatrix} [0,a_1] \\ [0,a_2] \\ [0,a_3] \end{bmatrix}, ([0,a_1] \; [0,a_2] \; ... \; [0,a_9]) \;\middle|\; a_i \in \{0,1,2,...,\infty\} \right\} \cup$$

$$\left\{ \begin{bmatrix} [0,a_1] & [0,a_2] \\ [0,a_3] & [0,a_4] \end{bmatrix}, ([0,a]) \;\middle|\; a, a_1, a_2, a_3, a_4 \in \{1, 1+i, i, 0, 2, 4\} \right\}$$

be a semigroup-set interval bivector space over the semigroup set $S = S_1 \cup S_2 = 3Z^+ \cup \{0\} \cup \{1, 0\}$.

We can as in case of other interval algebraic structures define substructures. We will leave the task of giving formal definitions to the reader but, however we will give some illustrative examples.

*Example 4.4.18*: Let $V = V_1 \cup V_2 =$

$$\left\{ \begin{bmatrix} [0,a_1] & [0,a_2] \\ [0,a_3] & [0,a_4] \end{bmatrix}, \sum_{i=0}^{\infty} [0,a_i] x^{2i} \,\middle|\, a_i \in Q^+ \cup \{0\} \right\} \cup$$

$$\left\{ [0,a_i], \begin{bmatrix} [0,a_1] \\ [0,a_2] \\ [0,a_3] \\ [0,a_4] \end{bmatrix} \,\middle|\, a_i, a_j \in \{\sqrt{3}, \sqrt{5}, \sqrt{15}, Z^+ \cup \{0\}\} \right\}$$

be a semigroup-set interval bivector space over the semigroup set $S = S_1 \cup S_2 = Q^+ \cup \{0\} \cup \{0, \sqrt{3}, 1, \sqrt{5}, \sqrt{15}\}$.

Take $W = W_1 \cup W_2 =$

$$\left\{ \begin{bmatrix} [0,a_1] & [0,a_2] \\ 0 & [0,a_3] \end{bmatrix}, \left\{ \sum_{i=0}^{\infty} [0,a_i] x^{16i} \,\middle|\, a_i \in Z^+ \cup \{0\} \right\} \right\} \cup$$

$$\left\{ [0,a_i], \begin{bmatrix} [0,a] \\ 0 \\ 0 \\ [0,a] \end{bmatrix} \,\middle|\, a, a_i \in \{0, \sqrt{3}, \sqrt{5}, \sqrt{15}, 4Z^+ \cup \{0\}\} \right\}$$

$\subseteq V = V_1 \cup V_2$; W is semigroup set interval bivector subspace of V over the semigroup-set $S = S_1 \cup S_2$.

***Example 4.4.19:*** Let $V = V_1 \cup V_2 =$

$$\left\{ \begin{bmatrix} [0,a_1] & [0,a_2] & [0,a_3] & [0,a_4] \\ [0,a_5] & [0,a_6] & [0,a_7] & [0,a_8] \end{bmatrix}, \begin{bmatrix} [0,a_1] \\ [0,a_2] \\ [0,a_3] \\ [0,a_4] \end{bmatrix} \middle| a_i \in \{0,1,2,\sqrt{3},\sqrt{41},7,9\} \right\}$$

$$\cup \left\{ \begin{bmatrix} [0,a_1] & [0,a_2] \\ [0,a_3] & [0,a_4] \\ [0,a_5] & [0,a_6] \\ [0,a_7] & [0,a_8] \\ [0,a_9] & [0,a_{10}] \\ [0,a_{11}] & [0,a_{12}] \end{bmatrix}, \bigl([0,a_1]\ [0,a_2]\ [0,a_3]\bigr) \middle| \begin{array}{l} a_i \in Z_{48}; \\ 1 \le i \le 12 \end{array} \right\}$$

be a set-semigroup interval bivector space over the set-semigroup $S = S_1 \cup S_2 = \{0, 1\} \cup \{Z_{48}\}$.

Take $W = W_1 \cup W_2 =$

$$\left\{ \begin{bmatrix} [0,a] & [0,b] & [0,a] & [0,b] \\ [0,b] & [0,a] & [0,b] & [0,a] \end{bmatrix}, \begin{bmatrix} [0,a] \\ 0 \\ [0,a] \\ 0 \end{bmatrix} \middle| a,b \in \{0,1,2,\sqrt{3},\sqrt{41},7,9\} \right\}$$

$$\cup \left\{ \begin{bmatrix} [0,a] & 0 \\ 0 & [0,b] \\ [0,a] & 0 \\ 0 & [0,b] \\ [0,a] & 0 \\ 0 & [0,b] \end{bmatrix}, \bigl([0,a]\ [0,b]\ 0\bigr) \middle| a,b \in Z_{48} \right\}$$

$\subseteq V_1 \cup V_2$. W is a set-semigroup interval bivector subspace of V over the set-semigroup $S = S_1 \cup S_2$.

***Example 4.4.20***: Let $V = V_1 \cup V_2$ = {All $10 \times 10$ interval matrices with entries of the form $[0, a_i]$ with $a_i \in Q^+ \cup \{0\}\} \cup$

$$\left\{ \sum_{i=0}^{\infty} a_i x^i; \begin{bmatrix} [0,a_1] & [0,a_2] & [0,a_3] \\ [0,a_4] & [0,a_5] & [0,a_6] \end{bmatrix} \middle| \begin{array}{l} a_i, a_j \in \{5Z^+ \cup \{0\}, \sqrt{19}, \sqrt{2}, \sqrt{5}\}; \\ 1 \le j \le 6 \end{array} \right\}$$

be a semigroup –set interval bivector space over the semigroup-set $S = S_1 \cup S_2 = Q^+ \cup \{0\} \cup \{5Z^+ \cup \{0\}, \sqrt{19}, \sqrt{2}, \sqrt{5}\}$.

Take $W = W_1 \cup W_2$ = {All $10 \times 10$ upper triangular interval matrices with entries of the form $[0, a_i]$ with $a_i \in Q^+ \cup \{0\}\} \cup$

$$\left\{ \sum_{i=0}^{\infty} a_i x^i; \begin{bmatrix} [0,a] & [0,a] & [0,a] \\ [0,a] & [0,a] & [0,a] \end{bmatrix} \middle| a \in \{5Z^+ \cup \{0\}, \sqrt{19}, \sqrt{2}, \sqrt{5}\} \right\}$$

$\subseteq V_1 \cup V_2$ and $T = T_1 \cup T_2 = (7Z^+ \cup \{0\}) \cup \{10Z^+ \cup \{0\}, \sqrt{2}\} \subseteq S_1 \cup S_2 = Q^+ \cup \{0\} \cup \{5Z^+ \cup \{0\}, \sqrt{19}, \sqrt{2}, \sqrt{5}\}$. $W = W_1 \cup W_2$ is a subsemigroup-subset interval bivector subspace of V over the subsemigroup-subset $T = T_1 \cup T_2$ of $S = S_1 \cup S_2$.

***Example 4.4.21***: Let $V = V_1 \cup V_2$ = {collection of all $5 \times 5$ interval matrices with entries of the form $[0, a_i]$; and $1 \times 8$ interval matrices of the form $[0, b_j]$ with $b_j, a_i \in \{11Z^+ \cup \{0\}, \sqrt{3}, \sqrt{7}, \sqrt{19}, \sqrt{41}, \sqrt{23}, \sqrt{43}, \sqrt{101}\} = S_1\} \cup$ {Collection of all $16 \times 16$ interval matrices with intervals of the form $[0, a_i]$ and all $7 \times 1$ interval matrices with intervals of the form $[0, b_j]$; $a_i, b_j \in 7Z^+ \cup \{0\} = S_2$} be a set-semigroup interval bivector space defined over the set-semigroup $S = S_1 \cup S_2$.

Take $W = W_1 \cup W_2$ = {Collection of all $5 \times 5$ diagonal interval matrices with intervals of the form $[0, a_i]$, $([0, a], 0, [0, a], 0, [0, a], 0, [0, a], 0) / a_i, a \in S_1\} \cup$ {Collection of all $16 \times 16$ upper triangular interval matrices with intervals of the form $[0, a_i]$ and

$$\left.\begin{bmatrix} [0,a_1] \\ 0 \\ [0,a_2] \\ 0 \\ [0,a_3] \\ 0 \\ [0,a_4] \end{bmatrix} \right| a_i, a_j \in 7Z^+ \cup \{0\}; 1 \le j \le 4\}$$

$\subseteq V_1 \cup V_2$ and $T = T_1 \cup T_2 = \{33Z^+ \cup \{0\}, \sqrt{3}, \sqrt{19}, \sqrt{101}\} \cup \{21Z^+ \cup \{0\}\} \subseteq S_1 \cup S_2$. W is a subset-subsemigroup interval bivector subspace of V over the subset-subsemigroup $T = T_1 \cup T_2 \subseteq S_1 \cup S_2 = S$.

If in the definition of the set-semigroup (semigroup-set) interval bivector space $V = V_1 \cup V_2$ over $S = S_1 \cup S_2$ if $V_1$ is a set interval linear algebra and $V_2$ is a semigroup interval linear algebra then we define V to be a set-semigroup interval bilinear algebra over $S = S_1 \cup S_2$.

We will illustrate this situation by some examples.

***Example 4.4.22***: Let $V = V_1 \cup V_2 =$

$$\left\{ \sum_{i=0}^{\infty} [0, a_i] x^i \,\middle|\, \begin{array}{l} a_i, b_i \in \{13Z^+ \cup \{0\}, a_i + b_i\sqrt{13}, \\ a_i + b_i\sqrt{2}, a_i + b_i\sqrt{3}\} = S_1 \end{array} \right\}$$

$\cup$ {All $5 \times 5$ interval matrices with intervals of the form $[0, a_i]$ and $a_i \in Z_{28}$} be a set-semigroup interval bilinear algebra over the set-semigroup, $S = S_1 \cup Z_{48} = S_1 \cup S_2$.

***Example 4.4.23***: Let $V = V_1 \cup V_2 =$

$$\left\{ \begin{bmatrix} [0,a_1] & [0,a_2] & [0,a_3] & [0,a_4] & [0,a_5] \\ [0,a_6] & [0,a_7] & [0,a_8] & [0,a_9] & [0,a_{10}] \end{bmatrix} \,\middle|\, a_i \in Q^+ \cup \{0\} \right\}$$

∪ {All 3 × 3 interval matrices with intervals of the form [0, $a_i$] / $a_i \in Z_{27}$} be a semigroup-set interval bilinear algebra over the semigroup set $S = S_1 \cup S_2 = Z^+ \cup \{0\} \cup \{\{0, 1, 2, 3, 9, 14, 20\} \subseteq Z_{27}\}$.

Now we give examples of substructures.

*Example 4.4.24*: Let $V = V_1 \cup V_2 = $ {all 9 × 9 interval matrices with intervals of the form [0, $a_i$] with $a_i \in Z_{240}$} ∪

$$\left\{ \sum_{i=0}^{\infty} [0, a_i] x^i \,\middle|\, a_i \in Z^+ \cup \{0\} \right\}$$

be a semigroup set interval bilinear algebra over the semigroup-set $S = S_1 \cup S_2 = Z_{240} \cup \{8Z^+ \cup \{0\}, 5Z^+ \cup \{0\}\}$.

Take $W = W_1 \cup W_2 = $ {All 9 × 9 interval upper triangular matrices with entries from $Z_{240}$} ∪

$$\left\{ \sum_{i=0}^{\infty} [0, a_i] x^{2i} \,\middle|\, a_i \in Z^+ \cup \{0\} \right\} \subseteq V_1 \cup V_2;$$

W is a semigroup set interval bilinear subalgebra of V over the semigroup-set $S = S_1 \cup S_2$.

*Example 4.4.25*: Let $V = V_1 \cup V_2 = $

$$\left\{ \begin{bmatrix} [0, a_1] & [0, a_2] & [0, a_3] & [0, a_4] \\ [0, a_5] & [0, a_6] & [0, a_7] & [0, a_8] \\ [0, a_9] & [0, a_{10}] & [0, a_{11}] & [0, a_{12}] \end{bmatrix} \,\middle|\, \begin{array}{l} a_i \in 5Z^+ \cup \{0\}; \\ 1 \leq i \leq 12 \end{array} \right\} \cup$$

$$\left\{ \begin{bmatrix} [0, a_1] & [0, a_2] \\ [0, a_3] & [0, a_4] \\ [0, a_5] & [0, a_6] \\ [0, a_7] & [0, a_8] \\ [0, a_9] & [0, a_{10}] \end{bmatrix} \,\middle|\, \begin{array}{l} a_i \in Z_{12}; \\ 1 \leq i \leq 10 \end{array} \right\}$$

be a set-semigroup interval bilinear algebra over the set semigroup $S = S_1 \cup S_2 = \{15Z^+ \cup \{0\}, 40Z^+ \cup \{0\}\} \cup Z_{12}$.

Take $W = W_1 \cup W_2 =$

$$\left\{ \begin{bmatrix} [0,a_1] & 0 & [0,a_2] & 0 \\ 0 & [0,a_3] & 0 & [0,a_4] \\ [0,a_5] & 0 & [0,a_6] & 0 \end{bmatrix} \middle| \begin{array}{l} a_i \in 5Z^+ \cup \{0\}; \\ 1 \leq i \leq 6 \end{array} \right\} \cup$$

$$\left\{ \begin{bmatrix} [0,a_1] & 0 \\ [0,a_2] & 0 \\ 0 & [0,a_1] \\ 0 & [0,a_2] \\ [0,a_1] & [0,a_2] \end{bmatrix} \middle| a_1, a_2 \in Z_{12} \right\}$$

$\subseteq V_1 \cup V_2$; W is a set- semigroup interval bilinear subalgebra of V over the set-semigroup $S = S_1 \cup S_2$.

*Example 4.4.26*: Let $V = V_1 \cup V_2 = $ {Collection of all $6 \times 6$ interval matrices with intervals of the form $[0, a_i]$, $a_i \in Z^+ \cup \{0\}\} \cup$

$$\left\{ \sum_{i=0}^{\infty} [0, a_i] x^i \middle| a_i \in Z_{36} \right\}$$

be a semigroup-set interval bilinear algebra over the semigroup-set $S = S_1 \cup S_2 = 3Z^+ \cup \{0\} \cup \{\{0, 2, 5, 1, 6, 7, 9, 14, 32, 30, 35\} \subseteq Z_{36}\}$. Choose $W = W_1 \cup W_2 = $ {Collection of all $6 \times 6$ interval upper triangular matrices with intervals of the form $[0, a_i]$ with $a_i \in Z^+ \cup \{0\}\} \cup$

$$\left\{ \sum_{i=0}^{\infty} [0, a_i] x^{2i} \middle| a_i \in Z_{36} \right\}$$

$\subseteq V_1 \cup V_2$ and $T = T_1 \cup T_2 = \{15Z^+ \cup \{0\}\} \cup \{\{0, 1, 2, 5, 30, 35\} \subseteq S_2 \subseteq Z_{36}\} \subseteq S_1 \cup S_2$.

Clearly W is a subsemigroup-subset interval bilinear subalgebra of V over the subsemigroup-subset $T = T_1 \cup T_2$ of $S_1 \cup S_2$.

*Example 4.4.27:* Let $V = V_1 \cup V_2 =$

$$\left\{ \sum_{i=0}^{\infty} [0, a_i] x^i \,\middle|\, a_i \in Z^+ \cup \{0\} \right\} \cup$$

$$\left\{ \begin{bmatrix} [0,a_1] & [0,a_2] & [0,a_3] & [0,a_4] & [0,a_5] \\ [0,a_6] & [0,a_7] & [0,a_8] & [0,a_9] & [0,a_{10}] \end{bmatrix} \,\middle|\, \begin{array}{l} a_i \in Z_{412}; \\ 1 \le i \le 10 \end{array} \right\}$$

be a set-semigroup interval bilinear algebra over the set-semigroup $S = S_1 \cup S_2 = \{2Z^+ \cup \{0\}, 5Z^+ \cup \{0\}, 3Z^+ \cup \{0\}\} \cup Z_{412}$. Choose $W = W_1 \cup W_2 =$

$$\left\{ \sum_{i=0}^{\infty} [0, a_i] x^{2i} \,\middle|\, a_i \in 30Z^+ \cup \{0\} \right\} \cup$$

$$\left\{ \begin{bmatrix} [0,a] & [0,a] & [0,a] & [0,a] & [0,a] \\ [0,b] & [0,b] & [0,b] & [0,b] & [0,b] \end{bmatrix} \,\middle|\, a, b \in Z_{412} \right\} ;$$

$T = T_1 \cup T_2 = \{4Z^+ \cup \{0\}, 15Z^+ \cup \{0\}\} \cup \{2Z_{412} = \{0, 2, 4, \ldots, 410\} \subseteq Z_{412}\} \subseteq S_1 \cup S_2$.

W is a subset-subsemigroup interval bilinear subalgebra of V over the subset-subsemigroup $T = T_1 \cup T_2 \subseteq S_1 \cup S_2$.

Now having seen examples of substructure we can define set-semigroup (semigroup-set) interval bilinear transformation provided the spaces are defined on the same set-semigroup (semigroup-set).

Let $V = V_1 \cup V_2$ and $P = P_1 \cup P_2$ be any two set-semigroup interval bivector space over the set-semigroup $S = S_1 \cup S_2$; that is $V_1$ and $P_1$ are set interval vector spaces over the same set $S_1$ and $V_2$ and $P_2$ are semigroup interval vector spaces over the same semigroup $S_2$. The bimap $T = T_1 \cup T_2 : V_1 \cup V_2 \to P_1 \cup$ $P_2$ where $T_1 : V_1 \to P_1$ and $T_2 : V_2 \to P_2$ are such that $T_1$ is a set linear interval vector space transformation and $T_2$ is a semigroup linear interval vector space transformation, then the bimap $T = T_1 \cup T_2$ is defined as the set-semigroup interval linear bitransformation of V in to P.

Interested reader can define properties analogous to usual linear transformations.

If $V = P$ that is $V_1 = P_1$ and $V_2 = P_2$ then we define T to be a set-semigroup interval linear bioperator. The transformations for set-semigroup interval bilinear algebra can be defined using, some simple and appropriate modifications.

Now we can derive almost all properties of these algebraic structures in an analogous way.

Now we can also define quasi set-semigroup interval linear algebras and their substructures in an analogous way.

Now we proceed on to define bigroup interval bivector spaces, set group (group-set) interval bivector spaces and semigroup-group (group-semigroup) interval bivector spaces and derive a few properties associated with them.

**DEFINITION 4.4.3**: *Let $V = V_1 \cup V_2$ be such that $V_i$ is a group interval vector space over the group $G_i$; $i = 1, 2$ and $V_i \not\subseteq V_j$, $V_j \not\subseteq V_i$ if if $i \neq j$ and $G_i \not\subseteq G_j$, $G_j \neq G_i$ if $i \neq j$; , $1 \leq i , j \leq 2$.*

*Then we define $V = V_1 \cup V_2$ to be a bigroup interval bivector space over the bigroup $G = G_1 \cup G_2$.*

***Example 4.4.28***: Let $V = V_1 \cup V_2 =$

$$\left\{ \begin{bmatrix} [0,a_1] & [0,a_2] & [0,a_3] \\ [0,a_4] & [0,a_5] & [0,a_6] \end{bmatrix}, \begin{bmatrix} [0,a_1] \\ [0,a_2] \\ [0,a_3] \\ [0,a_4] \\ [0,a_5] \\ [0,a_6] \\ [0,a_7] \end{bmatrix} \begin{matrix} a_i \in Z_{42}; \\ 1 \leq i \leq 7 \end{matrix} \right\} \cup$$

$$\left\{ \begin{bmatrix} [0,a_1] & [0,a_2] \\ [0,a_3] & [0,a_4] \\ [0,a_5] & [0,a_6] \\ [0,a_7] & [0,a_8] \\ [0,a_9] & [0,a_{10}] \end{bmatrix}, \bigl([0,a_1] \quad [0,a_2] \quad \ldots \quad [0,a_6]\bigr) \,\middle|\, \begin{array}{l} a_i \in Z_{30}; \\ 1 \le i \le 10 \end{array} \right\}$$

be a bigroup interval bivector space over the bigroup $G = G_1 \cup G_2 = Z_{42} \cup Z_{30}$.

*Example 4.4.29*: Let $V = V_1 \cup V_2 =$

$$\left\{ \bigl([0,a_1] \quad [0,a_2] \quad \ldots \quad [0,a_{15}]\bigr), \sum_{i=0}^{\infty}[0,a_i]x^i \,\middle|\, a_i \in Z_{12} \right\} \cup$$

$$\left\{ \begin{bmatrix} [0,a_1] \\ [0,a_2] \\ \vdots \\ [0,a_{15}] \end{bmatrix}, \begin{pmatrix} [0,a_1] & [0,a_2] & \ldots & [0,a_9] \\ [0,a_{10}] & [0,a_{11}] & \ldots & [0,a_{18}] \end{pmatrix} \,\middle|\, \begin{array}{l} a_i \in Z_{29}; \\ 1 \le i \le 18 \end{array} \right\}$$

be a bigroup interval bivector space over the bigroup $G = G_1 \cup G_2 = Z_{12} \cup Z_{29}$.

Now we will give examples of their substructure and the task of giving definition is left as an exercise for the reader.

*Example 4.4.30:* Let $V = V_1 \cup V_2 = \{$all $5 \times 5$ interval matrices with intervals of the form $[0, a_i]; a_i \in Z_{310}\} \cup$

$$\left\{ \begin{bmatrix} [0,a_1] \\ [0,a_2] \\ [0,a_3] \\ [0,a_4] \end{bmatrix} \,\middle|\, a_1, a_2, a_3, a_4 \in Z_{310} \right\}, \cup$$

$$\left. \sum_{i=0}^{\infty}[0,a_i]x^i, ([0,a_1] \quad [0,a_2] \quad \ldots \quad [0,a_{17}]) \right| a_i, a_j \in Z_{46}; 1 \leq j \leq 17 \right\}$$

be a bigroup interval bivector space over the bigroup $G = G_1 \cup G_2 = Z_{310} \cup Z_{46}$.

Take $W = W_1 \cup W_2 = \{$all $5 \times 5$ upper triangular interval matrices with intervals of the form $[0, a_i]$; $a_i \in Z_{310}\} \cup$

$$\left\{ \begin{bmatrix} [0,a] \\ [0,a] \\ [0,a] \\ [0,a] \end{bmatrix} \middle| a \in Z_{310} \right\} \cup$$

$$\left\{ \sum_{i=0}^{\infty}[0,a_i]x^{2i}, ([0,a] \quad [0,a] \quad \ldots \quad [0,a]) \middle| a_i, a \in Z_{46} \right\}$$

$\subseteq V_1 \cup V_2$ is a bigroup interval bivector subspace of V over the bigroup $G = G_1 \cup G_2 = Z_{310} \cup Z_{46}$.

***Example 4.4.31***: Let $V = V_1 \cup V_2 =$

$$\left\{ [0,a_i], \begin{bmatrix} [0,a_1] & [0,a_6] & [0,a_{11}] \\ [0,a_2] & [0,a_7] & [0,a_{12}] \\ [0,a_3] & [0,a_8] & [0,a_{13}] \\ [0,a_4] & [0,a_9] & [0,a_{14}] \\ [0,a_5] & [0,a_{10}] & [0,a_{15}] \end{bmatrix} \middle| a_i, a_j \in Z_{19}; 1 \leq j \leq 5 \right\} \cup$$

$$\left\{ \sum_{i=0}^{\infty}[0,a_i]x^{4i}, \begin{bmatrix} [0,a_1] & [0,a_5] & [0,a_9] \\ [0,a_2] & [0,a_6] & [0,a_{10}] \\ [0,a_3] & [0,a_7] & [0,a_{11}] \\ [0,a_4] & [0,a_8] & [0,a_{12}] \end{bmatrix} \middle| a_i, a_j \in Z_{24}; 1 \leq j \leq 12 \right\}$$

be a bigroup interval bivector space over the bigroup $G = Z_{19} \cup Z_{24} = G_1 \cup G_2$.

Take $W = W_1 \cup W_2 = $

$$\left\{ [0,a], \begin{bmatrix} [0,a] & 0 & [0,a] \\ 0 & [0,a] & 0 \\ [0,a] & 0 & [0,a] \\ 0 & [0,a] & 0 \\ [0,a] & 0 & [0,a] \end{bmatrix} \middle| a \in Z_{19} \right\} \cup$$

$$\left\{ \sum_{i=0}^{\infty} [0,a_i]x^{4i}, \begin{bmatrix} [0,a] & [0,a] & [0,a] \\ 0 & 0 & 0 \\ [0,a] & [0,a] & [0,a] \\ 0 & 0 & 0 \end{bmatrix} \middle| a_i, a \in \{2Z_{24} = \{0,2,\ldots,22\} \subseteq Z_{24} \right\}$$

$\subseteq V_1 \cup V_2$; W is a bigroup interval bivector subspace of V over the bigroup $G = Z_{19} \cup Z_{24}$.

*Example 4.4.32*: Let $V = V_1 \cup V_2 =$

$$\left\{ \begin{bmatrix} [0,a_1] & [0,a_2] \\ [0,a_3] & [0,a_4] \\ [0,a_5] & [0,a_6] \end{bmatrix}, \begin{bmatrix} [0,a_1] \\ [0,a_2] \\ \vdots \\ [0,a_{16}] \end{bmatrix} \middle| \begin{array}{l} a_i \in Z_{45}; \\ 1 \leq i \leq 16 \end{array} \right\} \cup$$

$$\left\{ ([0,a_1] \quad [0,a_2] \quad \ldots \quad [0,a_{10}]), \sum_{i=0}^{\infty}[0,a_i]x^i \middle| \begin{array}{l} a_j, a_i \in Z_{248}; \\ 1 \leq j \leq 10 \end{array} \right\}$$

be a bigroup interval bivector space over the bigroup $G = G_1 \cup G_2 = Z_{45} \cup Z_{248}$.

Take $W = W_1 \cup W_2 =$

$$\left\{ \begin{bmatrix} [0,a] & [0,a] \\ 0 & 0 \\ [0,b] & [0,b] \end{bmatrix}, \begin{bmatrix} 0 \\ [0,a_1] \\ 0 \\ [0,a_2] \\ 0 \\ 0 \\ 0 \\ 0 \\ 0 \\ 0 \\ 0 \\ 0 \\ [0,a_3] \\ 0 \\ 0 \\ [0,a_4] \end{bmatrix} \;\Bigg|\; a,b,a_1,a_2,a_3,a_4 \in Z_{45} \right\}$$

$$\cup \left\{ \left([0,a_1] \;\; 0 \;\; [0,a_2] \;\; 0 \;\; [0,a_3] \;\; 0 \;\; [0,a_4] \;\; 0 \;\; [0,a_5] \;\; 0\right), \sum_{i=0}^{\infty}[0,a_i]x^{2i} \;\Bigg|\; \begin{array}{l} a_i, a_j \in Z_{248}; \\ 1 \le j \le 5 \end{array} \right\}$$

$\subseteq V_1 \cup V_2$; be a bigroup interval bivector subspace of V over the bigroup G.

***Example 4.4.33***: Let $V = V_1 \cup V_2 =$

$$\left\{ \begin{bmatrix} [0,a_1] \\ [0,a_2] \\ \vdots \\ [0,a_{12}] \end{bmatrix}, \sum_{i=0}^{\infty}[0,a_i]x^i \;\Bigg|\; a_i, a_j \in Z_8; 1 \le j \le 12 \right\} \cup$$

{All 8 × 8 interval matrices with interval entries of the form [0, $a_i$] ; $a_i \in Z_{15}$, ( [0, $a_1$], [0, $a_2$], [0, $a_3$], [0, $a_4$] ) | $a_i$, $a_j \in Z_{15}$; $1 \le j \le 4$} be a bigroup interval bivector space over the bigroup $G = G_1 \cup G_2 = Z_8 \cup Z_{15}$.

Take $W = W_1 \cup W_2 =$

$$\left\{ \begin{bmatrix} [0,a] \\ [0,a] \\ \vdots \\ [0,a] \end{bmatrix}, \sum_{i=0}^{\infty}[0,a_i]x^{2i} \,\middle|\, a, a_i \in Z_8 \right\}$$

$\cup$ {All 8 × 8 upper triangular interval matrices with intervals of the form [0, $a_i$]; ([0, $a_1$], 0, [0, $a_2$], 0) | $a_i$, $a_1$, $a_2 \in Z_{15}$} $\subseteq V_1 \cup V_2$ be a subbigroup interval bivector subspace of V over the subbigroup $T = T_1 \cup T_2 = \{0, 2, 4, 6\} \cup \{0, 5, 10\} \subseteq Z_8 \cup Z_{15} = G_1 \cup G_2$.

If a bigroup interval bivector space V over the bigroup $G = G_1 \cup G_2$ has no subbigroup interval bivector subspace over the bigroup $G = G_1 \cup G_2$ then we say V to be a pseudo simple bigroup interval bivector subspace over the bigroup G. If V has no bigroup interval bivector subspace then we call V to be a simple bigroup interval bivector space. If V is both simple and pseudo simple then we call V to be a doubly simple bigroup interval bivector space.

We will give some illustrative examples of them.

*Example 4.4.34*: Let $V = V_1 \cup V_2 =$

$$\left\{ \sum_{i=0}^{\infty}[0,a_i]x^i ; \left([0,a_1] \;\; [0,a_2] \;\; [0,a_3]\right) \,\middle|\, a_i, a_1, a_2, a_3 \in Z_7 \right\} \cup$$

$$\left\{ \begin{bmatrix} [0,a_1] \\ [0,a_2] \\ [0,a_3] \\ [0,a_4] \end{bmatrix}, \right.$$

all 10 × 10 interval matrices with intervals of the form $[0, a_i]$; $a_i \in Z_5$} be a bigroup interval bivector space over the bigroup $G = G_1 \cup G_2 = Z_7 \cup Z_5$. We see the bigroup $G = Z_7 \cup Z_5$ is bisimple as it has no subgroups. Thus V is a pseudo simple bigroup interval bivector space over G.

However V is not doubly simple for take $W = W_1 \cup W_2 =$

$$\left\{ \sum_{i=0}^{\infty} [0, a_i] x^{2i}; ([0,a] \ [0,a] \ [0,a]) \middle| a_i, a \in Z_7 \right\} \cup \left\{ \begin{bmatrix} [0,a] \\ [0,a] \\ 0 \\ 0 \end{bmatrix} \right.,$$

all 10 × 10 upper triangular interval matrices with entries from $Z_5$ with intervals of the form $[0, a_i]$} $\subseteq V_1 \cup V_2$ ; W is a bigroup interval bivector subspace of V over the bigroup G. Thus V is not doubly simple.

*Example 4.4.35*: Let $V = V_1 \cup V_2 =$

$$\left\{ \begin{bmatrix} [0,a] & [0,a] \\ [0,a] & [0,a] \end{bmatrix}, [0,a] \middle| a \in Z_3 \right\} \cup \left\{ \begin{bmatrix} [0,a] \\ [0,a] \\ [0,a] \\ [0,a] \end{bmatrix} \middle| a \in Z_5 \right\}$$

be a bigroup interval bivector space over bigroup $G = G_1 \cup G_2 = Z_3 \cup Z_5$. We see V is a doubly simple bigroup interval bivector space over the bigroup G.

In view of this we have the following theorem.

**THEOREM 4.4.5:** *Let $V = V_1 \cup V_2$ be a bigroup interval bivector space over the bigroup $G = Z_p \cup Z_q$, p and q are two distinct primes. Then*
  1. *V is a pseudo simple bigroup interval bivector space over the bigroup G.*

2. *V is not a doubly simple bigroup interval bivector space over the bigroup in general.*

The proof is left as an exercise for the reader to prove.

**THEOREM 4.4.6**: *Let $V = V_1 \cup V_2$ be a bigroup interval bivector space over the bigroup $G = Z_n \cup Z_m$ where m and n are non primes; V is not simple or pseudo simple.*

This proof is also left as an exercise to the reader.

We see bigroup interval bivector spaces can be built only on finite bigroups of the form $G = Z_m \cup Z_n$, we cannot use $Z^+$ or $R^+$ or $Q^+$ or $C$.

Now we can define bigroup interval bilinear algebras. We give only examples of them.

*Example 4.4.36*: Let $V = V_1 \cup V_2 =$

$$\left\{ \begin{bmatrix} [0,a_1] & [0,a_2] & [0,a_3] & [0,a_4] & [0,a_5] \\ [0,a_6] & [0,a_7] & [0,a_8] & [0,a_9] & [0,a_{10}] \end{bmatrix} \middle| \begin{array}{l} a_i \in Z_{24}; \\ 1 \le i \le 10 \end{array} \right\}$$

$$\cup \left\{ \begin{bmatrix} [0,a_1] \\ [0,a_2] \\ \vdots \\ [0,a_{14}] \\ [0,a_{15}] \end{bmatrix} \middle| a_i \in Z_{32}; 1 \le i \le 15 \right\}$$

be a bigroup interval bilinear algebra over the bigroup $G = G_1 \cup G_2 = Z_{24} \cup Z_{32}$.

*Example 4.4.37*: Let $V = V_1 \cup V_2 = $ {All $9 \times 9$ interval matrices with intervals of the form $[0, a_i]$, $a_i \in Z_{38}$} $\cup$

$$\left\{ \sum_{i=0}^{\infty} [0, a_i] x^i \;\middle|\; a_i \in Z_{17} \right\}$$

be a bigroup interval bilinear algebra over the bigroup $G = G_1 \cup G_2 = Z_{38} \cup Z_{17}$.

We will illustrate the substructures by some examples.

***Example 4.4.38:*** Let $V = V_1 \cup V_2 = \{([0, a_1], [0, a_2], \ldots, [0, a_9]) \mid a_i \in Z_{28}; 1 \le i \le 9\} \cup$

$$\left\{ \begin{bmatrix} [0, a_1] \\ [0, a_2] \\ \vdots \\ [0, a_{11}] \\ [0, a_{12}] \end{bmatrix} \;\middle|\; a_i \in Z_{15}; 1 \le i \le 12 \right\}$$

be a bigroup interval bilinear algebra over the bigroup $G = G_1 \cup G_2 = Z_{28} \cup Z_{15}$.

Take $W = W_1 \cup W_2 = \{([0, a], [0, a], \ldots, [0, a]) \mid a \in Z_{28}\} \cup$

$$\left\{ \begin{bmatrix} [0, a_1] \\ [0, a_2] \\ \vdots \\ [0, a_{11}] \\ [0, a_{12}] \end{bmatrix} \;\middle|\; a_i \in \{0, 3, 6, 9, 12\} \subseteq Z_{15} \right\}$$

$\subseteq V_1 \cup V_2$ ; W is a bigroup interval bilinear subalgebra of V over the bigroup $G = Z_{28} \cup Z_{15}$.

***Example 4.4.39***: Let $V = V_1 \cup V_2 =$

$$\left\{ \begin{bmatrix} [0, a_1] & [0, a_2] \\ [0, a_3] & [0, a_4] \end{bmatrix} \;\middle|\; \begin{array}{l} a_i \in Z_3; \\ 1 \le i \le 4 \end{array} \right\}$$

$\cup \{([0, a_1], [0, a_2], \ldots, [0, a_7]) | a_i \in Z_{11}; 1 \le i \le 7\}$ be a bigroup interval bilinear algebra over the bigroup $G = Z_3 \cup Z_7$.

Take $W = W_1 \cup W_2 =$

$$\left\{ \begin{bmatrix} [0,a_1] & [0,a_2] \\ 0 & [0,a_3] \end{bmatrix} \middle| \begin{array}{l} a_i \in Z_3; \\ 1 \le i \le 3 \end{array} \right\} \cup$$

$\{([0, a], [0, a], \ldots, [0, a]) \mid a \in Z_{11}\} \subseteq V_1 \cup V_2$ is a bigroup interval bilinear subalgebra of V over the bigroup $G = Z_3 \cup Z_7$.

*Example 4.4.40*: Let $V = V_1 \cup V_2 =$

$$\left\{ \begin{bmatrix} [0,a_1] & [0,a_2] & [0,a_3] \\ [0,a_4] & [0,a_5] & [0,a_6] \end{bmatrix} \middle| \begin{array}{l} a_i \in Z_{18}; \\ 1 \le i \le 6 \end{array} \right\} \cup$$

$$\left\{ \begin{bmatrix} [0,a_1] & [0,a_2] \\ [0,a_3] & [0,a_4] \\ [0,a_5] & [0,a_6] \\ [0,a_7] & [0,a_8] \end{bmatrix} \middle| \begin{array}{l} a_i \in Z_{40}; \\ 1 \le i \le 8 \end{array} \right\}$$

be a bigroup interval bilinear algebra over the bigroup $G = G_1 \cup G_2 = Z_{18} \cup Z_{40}$. Take $H = H_1 \cup H_2 = \{0, 6, 12\} \cup \{0, 10, 20, 30\} \subseteq G_1 \cup G_2 = Z_{18} \cup Z_{40}$ and $W = W_1 \cup W_2 =$

$$\left\{ \begin{bmatrix} [0,a_1] & 0 & [0,a_2] \\ 0 & [0,a_3] & 0 \end{bmatrix} \middle| a_1, a_2, a_3 \in Z_{18} \right\} \cup$$

$$\left\{ \begin{bmatrix} [0,a_1] & 0 \\ 0 & [0,a_2] \\ [0,a_3] & 0 \\ 0 & [0,a_4] \end{bmatrix} \middle| \begin{array}{l} a_i \in Z_{40}; \\ 1 \le i \le 4 \end{array} \right\}$$

$\subseteq V_1 \cup V_2$ is a bisubgroup interval bilinear subalgebra of V over the bisubgroup $H = H_1 \cup H_2$ of $G = G_1 \cup G_2$.

**DEFINITION 4.4.4**: *Let $V = V_1 \cup V_2$ be where $V_1$ is a group interval vector space over the group $G_1$ and $V_2$ is a semigroup interval vector space over the semigroup $S_2$. $V$ is a group semigroup interval bivector space over the group-semigroup $G_1 \cup S_2$.*

We will illustrate this situation by some examples.

*Example 4.4.41:* Let $V = V_1 \cup V_2 =$

$$\left\{ \begin{bmatrix} [0,a_1] & [0,a_2] & [0,a_3] & [0,a_4] \\ [0,a_5] & [0,a_6] & [0,a_7] & [0,a_8] \end{bmatrix}, [0,a_i] \;\middle|\; \begin{array}{l} a_i, a_j \in Z^+ \cup \{0\}; \\ 1 \le j \le 8 \end{array} \right\}$$

$$\cup \left\{ \begin{bmatrix} [0,a_1] & [0,a_2] \\ [0,a_3] & [0,a_4] \\ [0,a_5] & [0,a_6] \\ [0,a_7] & [0,a_8] \end{bmatrix}, \sum_{i=0}^{\infty} [0,a_i]x^i \;\middle|\; \begin{array}{l} a_i, a_j \in Z_{17}; \\ 1 \le j \le 8 \end{array} \right\}$$

be a semigroup-group interval bivector space over the semigroup-group $G = (Z^+ \cup \{0\}) \cup Z_{17}$.

*Example 4.4.42:* Let $V = V_1 \cup V_2 =$

$$\left\{ \begin{bmatrix} [0,a_1] & [0,a_2] & [0,a_3] & [0,a_4] & [0,a_5] \\ [0,a_6] & [0,a_7] & [0,a_8] & [0,a_9] & [0,a_{10}] \end{bmatrix}, [0,a_i] \;\middle|\; \begin{array}{l} a_i \in Z_{48}; \\ 1 \le i \le 10 \end{array} \right\}$$

$$\cup \left\{ \begin{bmatrix} [0,a_1] & [0,a_2] \\ [0,a_3] & [0,a_4] \\ [0,a_5] & [0,a_6] \\ [0,a_7] & [0,a_8] \\ [0,a_9] & [0,a_{10}] \\ [0,a_{11}] & [0,a_{12}] \\ [0,a_{13}] & [0,a_{14}] \end{bmatrix}, \sum_{i=0}^{\infty} [0,a_i]x^i \;\middle|\; \begin{array}{l} a_i, a_j \in Q^+ \cup \{0\}; \\ 1 \le j \le 14 \end{array} \right\}$$

be a group-semigroup interval bivector space over the group-semigroup $G = Z_{48} \cup Q^+ \cup \{0\}$.

We can define substructures in a analogous way. We give only examples of them.

*Example 4.4.43*: Let $V = V_1 \cup V_2 =$

$$\left\{ \sum_{i=0}^{\infty} [0, a_i]x^i, \begin{bmatrix} [0, a_1] \\ [0, a_2] \\ [0, a_3] \\ [0, a_4] \end{bmatrix} \middle| a_i, a_j \in Q^+ \cup \{0\}; \; 1 \leq j \leq 4 \right\} \cup$$

$$\left\{ ([0, a_1] \; [0, a_2] \; ... \; [0, a_{11}]), \begin{bmatrix} [0, a_1] & [0, a_2] \\ [0, a_3] & [0, a_4] \\ [0, a_5] & [0, a_6] \end{bmatrix} \middle| \begin{array}{l} 1 \leq i \leq 11 \\ a_i \in Z_{19} \end{array} \right\}$$

be a semigroup-group interval bivector space over the semigroup-group $G = Q^+ \cup \{0\} \cup Z_{19}$.

Take $W = W_1 \cup W_2 =$

$$\left\{ \sum_{i=0}^{\infty} [0, a_i]x^{2i}, \begin{bmatrix} [0, a] \\ 0 \\ [0, a] \\ 0 \end{bmatrix} \middle| a_i, a \in Q^+ \cup \{0\} \right\} \cup$$

$$\left\{ \left([0, a_1], 0, [0, a_2], 0, [0, a_3], 0, [0, a_4], 0, [0, a_5], 0, [0, a_6]\right), \begin{bmatrix} [0, a_1] & 0 \\ 0 & [0, a_2] \\ [0, a_3] & 0 \end{bmatrix} \middle| \begin{array}{l} 1 \leq i \leq 6 \\ a_i \in Z_{19} \end{array} \right\}$$

$\subseteq V_1 \cup V_2$; W is a semigroup-group interval bivector subspace of V over the semigroup-group $Q^+ \cup \{0\} \cup Z_{19}$.

***Example 4.4.44***: Let $V = V_1 \cup V_2 =$

$$\left\{ \sum_{i=0}^{\infty} [0, a_i] x^{2i}, \sum_{i=0}^{\infty} [0, a_i] x^{3i} \;\middle|\; a_i \in Z_{320} \right\} \cup$$

{All $5 \times 5$ interval matrices with intervals of the form $[0, a_i]$; $a_i \in 3Z^+ \cup \{0\}$, $([0, a_1], [0, a_2], [0, a_3], [0, a_4], [0, a_5]) \mid a_j \in 4Z^+ \cup \{0\}\}$ be a group-semigroup interval bivector space over the group-semigroup $Z_{320} \cup 12Z^+ \cup \{0\}$.

Take $W = W_1 \cup W_2 =$

$$\left\{ \sum_{i=0}^{\infty} [0, a_i] x^{4i}, \sum_{i=0}^{\infty} [0, a_i] x^{9i} \;\middle|\; a_i \in Z_{320} \right\}$$

$\cup$ {All $5 \times 5$ interval upper triangular matrices with intervals of the form $[0, a_i]$, $a_i \in 3Z^+ \cup \{0\}$, $([0, a_1]\, [0, a_2]\, [0, a_3]) \mid a_1, a_2, a_3 \in 4Z^+ \cup \{0\}\} \subseteq V_1 \cup V_2$; $W$ is a group-semigroup interval bivector subspace of $V$ over the group-semigroup, $Z_{320} \cup 12Z^+ \cup \{0\}$.

***Example 4.4.45:*** Let $V = V_1 \cup V_2 =$

$$\left\{ \sum_{i=0}^{\infty} [0, a_i] x^{2i}, \sum_{i=0}^{\infty} [0, a_i] x^{5i} \;\middle|\; a_i \in Z_{196} \right\} \cup$$

$$\left\{ \begin{bmatrix} [0, a_1] & [0, a_2] \\ [0, a_3] & [0, a_4] \\ [0, a_5] & [0, a_6] \end{bmatrix}, ([0, a_1]\ \ [0, a_2]\ \ \ldots\ \ [0, a_{18}]) \;\middle|\; \begin{array}{l} a_i \in Z^+ \cup \{0\} \\ 1 \le i \le 18 \end{array} \right\}$$

be a group-semigroup interval bivector space over the group semigroup $Z_{196} \cup Z^+ \cup \{0\}$. Choose $W = W_1 \cup W_2 =$

$$\left\{ \sum_{i=0}^{\infty} [0, a_i] x^{4i}, \sum_{i=0}^{\infty} [0, a_i] x^{10i} \;\middle|\; a_i \in Z_{196} \right\} \cup$$

$$\left\{ \begin{bmatrix} [0,a] & [0,a] \\ [0,b] & [0,b] \\ [0,c] & [0,c] \end{bmatrix}, ([0,a] \quad [0,a] \quad ... \quad [0,a]) \,\middle|\, a, b, c \in Z^+ \cup \{0\} \right\}$$

$\subseteq V_1 \cup V_2$; W is a subgroup-subsemigroup interval bivector subspace over the subgroup-subsemigroup $\{0, 2, 4, 6, 8, …, 194\} \cup 3Z^+ \cup \{0\} \subseteq Z_{196} \cup Z^+ \cup \{0\}$.

***Example 4.4.46***: Let $V = V_1 \cup V_2 =$ {all $8 \times 8$ interval matrices with intervals of the form $[0, a_i]$ ; $a_i \in Q^+ \cup \{0\}$, all $3 \times 2$ interval matrices with intervals of the form $[0, a_i]$ with $a_i \in Z^+ \cup \{0\}\} \cup$

$$\left\{ \sum_{i=0}^{\infty} [0, a_i]x^i, ([0, a_1] \quad [0, a_2] \quad ... \quad [0, a_8]) \,\middle|\, \begin{array}{l} a_i, a_j \in Z_{48}; \\ 1 \leq j \leq 8 \end{array} \right\}$$

be a semigroup-group interval bivector space over the semigroup-group $Z^+ \cup \{0\} \cup Z_{48}$.

Let $W = W_1 \cup W_2 = $ {all $8 \times 8$ interval upper triangular matrices with intervals of the form $[0, a_i]$, $a_i \in Q^+ \cup \{0\}$,

$$\begin{bmatrix} [0,a] & [0,a] \\ [0,a] & [0,a] \\ [0,a] & [0,a] \end{bmatrix} \,\middle|\, a \in Z^+ \cup \{0\} \right\} \cup$$

$$\left\{ \sum_{i=0}^{\infty} [0, a_i]x^{3i}, ([0, a_1]\, 0\, [0, a_2]\, 0\, [0, a_3]\, 0\, [0, a_4]\, 0) \,\middle|\, \begin{array}{l} a_i, a_j \in Z_{48}; \\ 1 \leq j \leq 4 \end{array} \right\}$$

$\subseteq V_1 \cup V_2$ and $T = T_1 \cup T_2 = \{4Z^+ \cup \{0\}\} \cup \{0, 2, 4, 6, 8, …, 46\} \subseteq Z^+ \cup \{0\} \cup Z_{48}$. W is a subsemigroup-subgroup interval bivector subspace of V over the subsemigroup-subgroup $T = T_1 \cup T_2$.

Now we can in an analogous way define group-semigroup (semigroup group) interval bilinear algebra and their substructures.

We will illustrate these situations only by examples.

*Example 4.4.47:* Let $V = V_1 \cup V_2 =$

$$\left\{ \begin{bmatrix} [0,a_1] & [0,a_2] & [0,a_3] \\ [0,a_4] & [0,a_5] & [0,a_6] \\ [0,a_7] & [0,a_8] & [0,a_9] \end{bmatrix} \middle| \begin{array}{l} a_i \in Z_{49}; \\ 1 \leq i \leq 9 \end{array} \right\} \cup$$

$$\left\{ \begin{bmatrix} [0,a_1] & [0,a_2] \\ [0,a_3] & [0,a_4] \\ [0,a_5] & [0,a_6] \\ [0,a_7] & [0,a_8] \\ [0,a_9] & [0,a_{10}] \\ [0,a_{11}] & [0,a_{12}] \end{bmatrix} \middle| \begin{array}{l} a_i \in Z^+ \cup \{0\}; \\ 1 \leq i \leq 12 \end{array} \right\}$$

be a group-semigroup interval bilinear algebra over the group-semigroup $Z_{49} \cup Z^+ \cup \{0\}$.

*Example 4.4.48*: Let $V = V_1 \cup V_2 =$

$$\left\{ \sum_{i=0}^{\infty} [0,a_i]x^i \middle| a_i \in Q^+ \cup \{0\} \right\}$$

$\cup$ {all $7 \times 7$ interval matrices with intervals of the form $[0, a_i]$, $a_i \in Z_{412}$} be a semigroup-group interval bilinear algebra over the semigroup-group $Q^+ \cup \{0\} \cup Z_{412}$.

*Example 4.4.49:* Let $V = V_1 \cup V_2 =$ {all $8 \times 8$ interval matrices with intervals of the form $[0, a_i]$; $a_i \in Z_{480}$} $\cup$

$$\left\{ \sum_{i=0}^{\infty} [0, a_i] x^i \;\middle|\; a_i \in 3Z^+ \cup \{0\} \right\}$$

be a group-semigroup interval bilinear algebra over the group-semigroup $S = Z_{480} \cup 3Z^+ \cup \{0\}$. Choose $W = W_1 \cup W_2 = \{$all $8 \times 8$ interval upper triangular matrices with intervals $[0, a_i]$; $a_i \in Z_{480}\} \cup$

$$\left\{ \sum_{i=0}^{\infty} [0, a_i] x^{2i} \;\middle|\; a_i \in 3Z^+ \cup \{0\} \right\}$$

$\subseteq V_1 \cup V_2$ be a group-semigroup interval bilinear subalgebra of V over the group-semigroup S.

*Example 4.4.50:* Let $V = V_1 \cup V_2 =$

$$\left\{ \begin{bmatrix} [0,a_1] & [0,a_3] & [0,a_5] \\ [0,a_2] & [0,a_4] & [0,a_6] \end{bmatrix} \;\middle|\; \begin{array}{l} a_i \in Z^+ \cup \{0\}; \\ 1 \leq i \leq 6 \end{array} \right\} \cup$$

$$\left\{ \sum_{i=0}^{\infty} [0, a_i] x^i \;\middle|\; a_i \in Z_{15} \right\}$$

be a semigroup group interval bilinear algebra over the semigroup-group $Z^+ \cup \{0\} \cup Z_{15}$.

Take $W = W_1 \cup W_2 =$

$$\left\{ \begin{bmatrix} [0,a_1] & 0 & [0,a_2] \\ 0 & [0,a_3] & 0 \end{bmatrix} \;\middle|\; \begin{array}{l} a_i \in Z^+ \cup \{0\}; \\ 1 \leq i \leq 3 \end{array} \right\} \cup$$

$$\left\{ \sum_{i=0}^{\infty} [0, a_i] x^{2i} \;\middle|\; a_i \in Z_{15} \right\}$$

$\subseteq V_1 \cup V_2$, W is a semigroup-group interval bilinear subalgebra over the semigroup-group.

***Example 4.4.51***: Let $V = V_1 \cup V_2 =$

$$\left\{ \sum_{i=0}^{\infty} [0, a_i] x^i \,\middle|\, a_i \in Q^+ \cup \{0\} \right\} \cup$$

$$\left\{ \begin{bmatrix} [0, a_1] & [0, a_2] & [0, a_3] & [0, a_4] & [0, a_5] \\ [0, a_6] & [0, a_7] & [0, a_8] & [0, a_9] & [0, a_{10}] \\ [0, a_{11}] & [0, a_{12}] & [0, a_{13}] & [0, a_{14}] & [0, a_{15}] \end{bmatrix} \,\middle|\, \begin{array}{l} a_i \in Z_{30}; \\ 1 \le i \le 15 \end{array} \right\}$$

be a semigroup-group interval bilinear algebra over the semigroup group, $G = Q^+ \cup \{0\} \cup Z_{30}$.

Take $W = W_1 \cup W_2 =$

$$\left\{ \sum_{i=0}^{\infty} [0, a_i] x^{2i} \,\middle|\, a_i \in Q^+ \cup \{0\} \right\} \cup$$

$$\left\{ \begin{bmatrix} [0, a_1] & 0 & [0, a_2] & 0 & [0, a_3] \\ 0 & [0, a_4] & 0 & [0, a_5] & 0 \\ [0, a_6] & 0 & [0, a_4] & 0 & [0, a_8] \end{bmatrix} \,\middle|\, \begin{array}{l} a_i \in Z_{30}; \\ 1 \le i \le 8 \end{array} \right\}$$

$\subseteq V_1 \cup V_2$ and $H = 3Z^+ \cup \{0\} \cup \{0, 5, 10, 15, 20, 25\} \subseteq Q^+ \cup \{0\} \cup Z_{30}$. W is a subsemigroup-subgroup interval bilinear subalgebra of V over the subsemigroup-subgroup H of G.

***Example 4.4.52:*** Let $V = V_1 \cup V_2 = \{$all $16 \times 16$ interval matrices with intervals of the form $[0, a_i] \mid a_i \in Z_{12}\} \cup$

$$\left\{ \sum_{i=0}^{\infty} [0, a_i] x^i \,\middle|\, a_i \in Z^+ \cup \{0\} \right\}$$

be a group-semigroup interval bilinear algebra over the group-semigroup $S = Z_{12} \cup Z^+ \cup \{0\}$. Choose $W = W_1 \cup W_2 = \{$all $16 \times 16$ upper triangular interval of the form $[0, a_i] \mid a_i \in Z_{12}\} \cup$

$$\left\{ \sum_{i=0}^{\infty} [0, a_i] x^{2i} \;\middle|\; a_i \in Z^+ \cup \{0\} \right\}$$

and $T = T_1 \cup T_2 = \{0, 2, 4, 6, 8, 10\} \cup 5Z^+ \cup \{0\} \subseteq Z_{12} \cup Z^+ \cup \{0\}$. W is a subgroup-subsemigroup interval bilinear subalgebra of V over the subgroup-subsemigroup T of S.

We say a group-semigroup (semigroup-group) interval bilinear algebra V (bivector space) is pseudo simple if V has no subgroup subsemigroup (subsemigroup-subgroup) interval bilinear subalgebras (or bivector subspaces). A group-semigroup (semigroup-group) interval bilinear algebra (bivector space) is simple if V has no group-semigroup (semigroup-group) interval bilinear subalgebra (bivector subspace).

We will illustrate this situation by some examples.

*Example 4.4.53:* Let $V = V_1 \cup V_2 =$

$$\left\{ \begin{bmatrix} [0,a] & [0,a] \\ [0,a] & [0,a] \end{bmatrix} \;\middle|\; a \in Z_3 \right\} \cup \left\{ \begin{bmatrix} [0,a] \\ [0,a] \\ [0,a] \end{bmatrix} \;\middle|\; a \in Z^+ \cup \{0\} \right\}$$

be a group-semigroup bilinear algebra over the group-semigroup $Z_3 \cup Z^+ \cup \{0\}$. Clearly V is pseudo simple as well simple. Hence V is a doubly simple group-semigroup bilinear algebra over the group semigroup.

**THEOREM 4.4.7**: *Let $V = V_1 \cup V_2$ be a group-semigroup (semigroup-group) interval bilinear algebra over the group-semigroup $Z_n \cup Z^+ \cup \{0\}$ (semigroup-group $Z^+ \cup \{0\} \cup Z_n$). Then V is not simple or pseudo simple provided n is a compositive number.*

The proof is left as an exercise to the reader.

**THEOREM 4.4.8:** *Let $V = V_1 \cup V_2$ be a group semigroup (semigroup-group) bilinear algebra (bivector space) over the group-semigroup $(Z_p - Z^+ \cup \{0\})$ p a prime. Then V is pseudo simple and need not in general be doubly simple.*

This proof is also left as an exercise for the reader.
Now we proceed onto define set-group (group-set) interval bilinear algebra (bivector space) over the set-group (group-set).

**DEFINITION 4.4.5**: *Let $V = V_1 \cup V_2$ be such that $V_1$ is a set interval vector space over the set $S_1$ and $V_2$ is a group interval vector space over the group $G_2$. We define $V = V_1 \cup V_2$ to be a set-group interval bivector space over the set-group.*

We will illustrate this situation by some examples.

*Example 4.4.54:* Let $V = V_1 \cup V_2 =$

$$\left\{ \sum_{i=0}^{\infty} [0, a_i] x^{3i}, \sum_{i=0}^{\infty} [0, a_i] x^{2i} \;\middle|\; a_i \in Z^+ \cup \{0\} \right\} \cup$$

$$\left\{ \begin{bmatrix} [0, a_1] \\ [0, a_2] \\ [0, a_3] \\ [0, a_4] \end{bmatrix}, \bigl([0, a_1], [0, a_2], [0, a_3]\bigr) \;\middle|\; \begin{array}{l} a_i \in Z_9; \\ 1 \le i \le 4 \end{array} \right\}$$

be a set group interval bivector space over the set- group, $\{0, 2, 17, 41, 142, 250\} \cup Z_9$.

*Example 4.4.55:* Let $V = V_1 \cup V_2 =$

$$\left\{ \begin{bmatrix} [0, a_1] \\ [0, a_2] \\ [0, a_3] \\ [0, a_4] \end{bmatrix}, \bigl([0, a_1], ..., [0, a_{12}]\bigr) \;\middle|\; a_i \in Z_{48}; 1 \le i \le 12 \right\} \cup$$

$$\left\{ \sum_{i=0}^{\infty} [0, a_i] x^i, [0, a_i] \,\middle|\, a_i \in 5Z^+ \cup \{0\} \right\}$$

be a group-set interval bivector space over the group set $Z_{48} \cup \{0, 1, 2, 7, 15Z^+\}$.

*Example 4.4.56:* Let $V = V_1 \cup V_2 =$

$$\left\{ \sum_{i=0}^{\infty} [0, a_i] x^{2i}, \sum_{i=0}^{\infty} [0, a_i] x^{3i} \,\middle|\, a_i \in 4Z^+ \cup \{0\} \right\} \cup$$

$$\left\{ \begin{bmatrix} [0, a_1] \\ [0, a_2] \\ \vdots \\ [0, a_8] \end{bmatrix}, ([0, a_1], [0, a_2]) \,\middle|\, a_i \in Z_{49}; 1 \le i \le 8 \right\}$$

be a set-group interval bivector space over the set-group $S \cup G = \{16 Z^+ \cup \{0\}, 4, 8\} \cup Z_{49}$.

Take $W = W_1 \cup W_2 =$

$$\left\{ \sum_{i=0}^{\infty} [0, a_i] x^{4i}, \sum_{i=0}^{\infty} [0, a_i] x^{3i} \,\middle|\, a_i \in 16Z^+ \cup \{0\} \right\} \cup$$

$$\left\{ \begin{bmatrix} [0, a_1] \\ 0 \\ [0, a_1] \\ 0 \\ [0, a_1] \\ 0 \\ [0, a_1] \\ 0 \end{bmatrix}, ([0, a], [0, a]) \,\middle|\, a_1, a \in Z_{49} \right\}$$

⊆ $V_1 \cup V_2$; W is a set - group interval bivector subspace of V over the set-group $S \cup G$.

***Example 4.4.57:*** Let $V = V_1 \cup V_2 =$

$$\left\{ \begin{bmatrix} [0,a_1] \\ [0,a_2] \\ [0,a_3] \end{bmatrix}, ([0,a_1],[0,a_2]...,[0,a_7]) \,\middle|\, a_i \in Z^+ \cup \{0\} \right\} \cup$$

$$\left\{ \begin{bmatrix} [0,a_1] & [0,a_2] \\ [0,a_3] & [0,a_4] \\ [0,a_5] & [0,a_6] \\ [0,a_7] & [0,a_8] \end{bmatrix}, \begin{bmatrix} [0,a_1] \\ [0,a_2] \\ \vdots \\ [0,a_{11}] \end{bmatrix} \,\middle|\, \begin{array}{l} a_i \in Z_{24}; \\ 1 \le i \le 11 \end{array} \right\}$$

be a set-group interval bivector space over the set-group $\{3Z^+, 4Z^+, 17Z^+, 13Z^+, 0\} \cup Z_{24} = S_1 \cup G_2$.

Choose $W = W_1 \cup W_2 =$

$$\left\{ ([0,a_1],...,[0,a_7]) \,\middle|\, \begin{array}{l} a, a_i \in 3Z^+ \cup \{0\}; \\ 1 \le i \le 7 \end{array} \right\} \cup$$

$$\left\{ \begin{bmatrix} [0,a_1] & 0 \\ 0 & [0,a_2] \\ [0,a_3] & 0 \\ 0 & [0,a_4] \end{bmatrix}, \begin{bmatrix} [0,a_1] \\ [0,a_2] \\ \vdots \\ [0,a_{11}] \end{bmatrix} \,\middle|\, a_i \in 2Z_{24} = \{0,...,22\} \subseteq Z_{24} \right\}$$

⊆ $V_1 \cup V_2$ and $T = T_1 \cup T_2 = \{3Z^+, 13Z^+, 4Z^+, 0\} \cup \{0, 4, 8, 12, 16, 20\} \subseteq S_1 \cup G_2$. W is a subset-subgroup interval bivector subspace of V over the subset-subgroup T of $S_1 \cup G_2$.

***Example 4.4.58:*** Let $V = V_1 \cup V_2 = \left\{ \sum_{i=0}^{\infty} [0,a_i]x^i \,\middle|\, a_i \in Z_{48} \right\}$

∪ {All 6 × 6 interval matrices with intervals of the form $[0, a_i]$; $a_i \in Z^+ \cup \{0\}$} be a group-set interval bilinear algebra over the group-set $G = G_1 \cup S_2 = Z_{48} \cup \{30Z^+ \cup \{0\}, 2Z^+\}$.

*Example 4.4.59:* Let $V = V_1 \cup V_2 = $ {set of all 9 × 5 interval matrices with intervals of the form $[0, a_i]$; $a_i \in Q^+ \cup \{0\}$} ∪ {all 3 × 8 interval matrices with intervals of the form $[0, a_i]$; $a_i \in Z_{41}$} be a set - group interval bilinear algebra over the set - group $S_1 \cup G_2 = \{2Z^+, 5Z^+, 7Z^+, 0\} \cup Z_{41}$.

We now state the theorem the proof of which is direct.

**THEOREM 4.4.9**: *Every set-group (group-set) interval bilinear algebra over a set-group (group-set) is a set-group (group-set) interval bivector space but not conversely.*

*Example 4.4.60:* Let $V = V_1 \cup V_2 =$

$$\left\{ \sum_{i=0}^{\infty} [0, a_i] x^i \,\middle|\, a_i \in Z^+ \cup \{0\} \right\}$$

∪ {all 4 × 4 interval matrices with interval entries of the form $[0, a_i]$; $a_i \in Z_{43}$} be a set - group interval bilinear algebra over the set-group $S_1 \cup G_2 = \{3Z^+, 2Z^+, 7Z^+, 0\} \cup \{Z_{43}\}$.
Choose $W = W_1 \cup W_2 =$

$$\left\{ \sum_{i=0}^{\infty} [0, a_i] x^{2i} \,\middle|\, a_i \in Z^+ \cup \{0\} \right\}$$

∪ {All 4 × 4 upper triangular interval matrices of intervals of the form $[0, a_i]$; $a_i \in Z_{43}$} $\subseteq V_1 \cup V_2$; W is a set-group interval bilinear subalgebra of V over $S_1 \cup G_2$.

*Example 4.4.61:* Let $V = V_1 \cup V_2 =$
$$\left\{ \sum_{i=0}^{\infty} [0, a_i] x^i \,\middle|\, a_i \in Q^+ \cup \{0\} \right\}$$

∪ {all 5 × 2 interval matrices with intervals of the form $[0, a_i]$; $a_i \in Z_{48}$} be a set-group interval bilinear algebra over the set - group $S_1 \cup G_2 = \{7Z^+ \cup \{0\}, 3Z^+ \cup \{0\}, 4Z^+\} \cup Z_{48}$.

Choose $W = W_1 \cup W_2 =$

$$\left\{\sum_{i=0}^{\infty}[0,a_i]x^{2i} \,\middle|\, a_i \in Z^+ \cup \{0\}\right\} \cup$$

$$\left\{\begin{bmatrix}[0,a_1] & 0 & [0,a_2] & 0 & [0,a_3] \\ 0 & [0,a_1] & 0 & [0,a_2] & 0\end{bmatrix} \,\middle|\, a_1, a_2, a_3 \in Z\right\}_{48}$$

$\subseteq V_1 \cup V_2$ and $P_1 \cup P_2 = \{3Z^+, 4Z^+, 0\} \cup \{0, 12, 24, 36\} \subseteq S_1 \cup G_2$. W is a subset - subgroup interval bilinear subalgebra of V over the subset - subgroup $P_1 \cup H_2$ of $S_1 \cup G_2$.

*Example 4.4.62*: Let $V = V_1 \cup V_2 =$

$$\left\{\sum_{i=0}^{\infty}[0,a_i]x^i \,\middle|\, a_i \in Z^+ \cup \{0\}\right\} \cup$$

$$\left\{\begin{bmatrix}[0,a] & [0,a] & [0,a] & [0,a] \\ [0,b] & [0,b] & [0,b] & [0,b] \\ [0,c] & [0,c] & [0,c] & [0,c]\end{bmatrix} \,\middle|\, a, b, c \in Z_7\right\}$$

be a set - group interval bilinear algebra over the set-group $S = S_1 \cup G_2 = \{2Z^+, 5Z^+, 0\} \cup Z_7$. Clearly V is a pseudo simple set - group interval bilinear algebra over S.

We can derive several interesting properties related with them as in case of usual bigroup-linear algebras.

**Chapter Five**

# APPLICATIONS OF THE SPECIAL CLASSES OF INTERVAL LINEAR ALGEBRAS

These new classes of interval linear algebras find their applications in fields, which demand the solution to be in intervals and in finite element methods. The present day trend is scientists, technologists and medical experts seek interval solutions to single valued solutions. For interval solutions give them more freedom to work and also one can choose the best-suited solution from that interval.

These structures can be best utilized in the study of finite element analysis. These structures are well suited for the analysis of stiffness matrices when interval solutions are in demand. These new structures have several limitations for they cannot be built using intervals of the form [–a, b] where a and b are in Z.

These structures can be used in all mathematical models and fuzzy models, which demand interval solutions. For more about interval algebraic structures refer [52].

**Chapter Six**

# SUGGESTED PROBLEMS

In this chapter we propose over 100 problems, which will be a challenge to the reader.

1. Find some interesting properties about set complex interval vector spaces.

2. Give an example of a order 21 set modulo integer vector space built using $Z_{40}$.

3. Does their exists a set modulo integer vector space of cardinality 12 built using $Z_7$? Justify your claim.

4. Obtain some interesting properties enjoyed by the set real interval spaces.

5. Does there exists a set modulo integer interval vector space of order 149? Justify your claim.

6. Let $V = \left\{ \begin{bmatrix} [0, a_1] & [0, a_2] \\ [0, a_3] & [0, a_4] \\ [0, a_5] & [0, a_6] \end{bmatrix} \middle| a_i \in Z_{17} \right\}$, be a set modulo integer linear algebra over the set $S = \{0, 1, 2, 5\} \subseteq Z_{17}$. Find set modulo integer interval linear subalgebras of V.

7. Obtain some interesting properties about set rational interval vector spaces.

8. Let $S = \{[a, b] \mid a, b \in Q^+ \cup \{0\}; a \leq b\}$ be a set rational interval vector space over the set $S = \{0, 1\}$. Find set rational interval vector subspaces of V. Can S be generated finitely? Justify your claim.

9. Obtain some interesting properties about set complex interval linear algebras.

10. Give an example of a doubly simple set interval integer linear algebra.

11. Give an example of a semigroup interval vector space which is not a semigroup interval linear algebra.

12. Give some interesting properties of semigroup interval vector spaces.

13. Give an example of a finitely generated semigroup interval linear algebra.

14. Give an example of a simple semigroup interval vector space.

15. Give an example of a pseudo simple semigroup interval vector space which is not simple.

16. Give an example of a doubly simple semigroup interval vector space.

17. Give an example of a pseudo semigroup interval linear algebra.

18. Does there exists a semigroup interval linear algebra which cannot be written as a direct sum? Justify your claim!

19. Obtain some interesting properties about group interval linear algebras.

20. Does there exists an infinite group interval linear algebra?

21. Does their exists a group interval linear algebra of order 43?

22. Give an example of a group linear algebra of order 25.

23. Let $X = \left\{ \sum_{i=0}^{\infty} [0, a_i] x^i \mid a_i \in Z_7 \right\}$. Is X a group interval linear algebra over the group $Z_7$? Is X finite or infinite?

24. Give an example of a set fuzzy interval vector space.

25. Obtain some interesting properties about set fuzzy interval linear algebras.

26. Give an example of a semigroup fuzzy interval linear algebra.

27. Let V = {All 5 × 5 interval matrices with intervals of the form 
$\{[0, a_i], \begin{bmatrix} [0, a_1] \\ [0, a_2] \\ \vdots \\ [0, a_9] \end{bmatrix} \mid a_i \in Z^+ \cup \{0\}; 1 \leq i \leq 25\}$ be a semigroup interval vector space over the semigroup $S = Z^+ \cup \{0\}$. Obtain atleast 5 fuzzy semigroup interval vector spaces or semigroup fuzzy interval vector spaces.

28. Obtain some interesting properties about group interval vector spaces.

29. Let V = {all 6 × 6 interval matrices with intervals of the form [0, $a_i$]; $a_i \in Z_{18}$} be a group interval linear algebra over the group G = $Z_{18}$.
   i. Obtain group fuzzy interval linear algebras.
   ii. Does V have subgroup interval linear subalgebras?
   iii. Find at least 3 group interval linear subalgebras.
   iv. Define a linear operator on V with non trivial kernel.

30. Bring out the difference between type I and type II semigroup fuzzy interval linear algebras.

31. Obtain some interesting properties enjoyed by type II group fuzzy interval linear algebras?

32. Let V = $V_1 \cup V_2$ = $\left\{ \sum_{i=0}^{\infty}[0,a_i]x^{2i}, \sum_{i=0}^{\infty}[0,a_i]x^{3i}; a_i \in Z^+ \cup \{0\} \right\}$

   $\cup \left\{ \begin{bmatrix} [0,a_1] \\ [0,a_2] \\ \vdots \\ [0,a_7] \end{bmatrix}, ([0,a_1] \cdots [0,a_8]) \Big| a_i \in 3Z^+ \cup \{0\} \right\}$ be a set

   interval bivector space over the set S = {$2Z^+$, $3Z^+$, 0}.
   i. Find set interval bivector subspaces of V.
   ii. Find subset interval bivector subspace of V.
   iii. Define a bilinear operator on V.
   iv. Find a generating set of V.

33. Give an example of a set interval bivector space which is not a set interval linear algebra of finite dimension.

34. Give an example of a pseudo simple set interval linear algebra.

35. Let V = $V_1 \cup V_2$ = {all 7 × 7 interval matrices with intervals of the form [0, $a_i$]| $a_i \in Z_7$} $\cup$ {[0, $a_i$]| $a_i \in Z_7$} be a set interval bilinear algebra over the set S = $Z_7$.

i. Find a linear bioperator on V which has a non trivial bikernel.
  ii. Is V simple?
  iii. Can V have subset interval bilinear algebras?

36. Give an example of a pseudo simple set interval linear bialgebra which is not simple.

37. Obtain some important properties about set interval linear bialgebras.

38. Give an example of a set interval linear bialgebra of bidimension (5, 9).

39. What is the difference between a set interval bilinear algebra and a biset interval bilinear algebra?

40. Let $V = V_1 \cup V_2 =$ {all $10 \times 10$ interval matrices with intervals of the form $[0, a_i]$, $a_i \in 5 Z^+ \cup \{0\}\} \cup \{3 \times 7$ interval matrices with intervals of the form $[0, a_i] \in 3Z^+ \cup \{0\}\}$ be a biset interval bivector space of V over the biset $S = 10Z^+ \cup \{0\} \cup 6Z^+ \cup \{0\} = S_1 \cup S_2$.
   i. Find a bigenerating bisubset of V.
   ii. Is V finite bidimensional?
   iii. Find biset interval bivector subspaces of V.
   iv. Is V pseudo simple? Justify your answer.
   v. Define a nontrivial one to one bilinear operator on V.

41. Give an example of a quasi biset interval bivector space.

42. Let
$$V = V_1 \cup V_2 = \left\{ ([0,a_1] \quad \cdots \quad [0,a_8]) \big| a_i \in Z_7 \right\} \cup \left\{ \begin{bmatrix} a_1 & a_2 & a_3 \\ a_4 & a_5 & a_6 \\ a_7 & a_8 & a_9 \end{bmatrix} \bigg| a_i \in Z_{11} \right\}$$
be a quasi biset interval bivector space over the biset $S = Z_7 \cup Z_{11}$. Is V bisimple? What is the bigenerating subbiset of V?

43. Give an example of a doubly simple quasi biset interval bivector space over the biset S.

44. Let $V = V_1 \cup V_2 = \left\{ \sum_{i=0}^{\infty} [0, a_i] x^i \mid a_i \in Z_{13} \right\} \cup \{$All $3 \times 8$

    interval matrices with intervals of the form $[0, a_i]$ / $a_i \in Z_{13}\}$ be a quasi set interval linear bialgebra over the set $S = Z_{13}$. Is V simple? Justify!

45. Give an example of a semi quasi set interval bilinear algebra.

46. Determine some special properties enjoyed by semigroup interval bilinear algebras.

47. Let
$$V = V_1 \cup V_2$$
$$= \left\{ \begin{bmatrix} [0,a_1] & [0,a_2] & [0,a_3] & [0,a_4] \\ [0,a_5] & [0,a_6] & [0,a_7] & [0,a_8] \end{bmatrix} \middle| a_i \in Z_{17}; 1 \le i \le 8 \right\} \cup$$
$$\left\{ \sum_{i=0}^{\infty} [0, a_i] x^i \middle| a_i \in Z_{17} \right\}$$

    be a semigroup interval bilinear algebra over the semigroup $Z_{17}$.
    i. Is V pseudo simple?
    ii. Is V simple? Justify
    iii. Find a generating bisubset of V.

48. Prove there exists an infinite class of semigroup interval bilinear algebras which are pseudo simple.

49. Give an example of a simple semigroup interval bivector space.

50. Give an example of a double simple semigroup interval bivector space.

51. Give some interesting applications of semigroup interval bilinear algebras.

52. Let
$$V = V_1 \cup V_2$$
$$= \left\{ \begin{bmatrix} [0,a_1] \\ [0,a_2] \\ \vdots \\ [0,a_9] \end{bmatrix}, \left([0,a_1] \cdots [0,a_{14}]\right) \middle| \begin{array}{l} a_i \in Z^+ \cup \{0\}; \\ 1 \leq i \leq 14 \end{array} \right\} \cup$$

{All 7 × 7 interval matrices with intervals of the form [0, $a_i$] / $a_i \in 3Z^+ \cup \{0\}$} be a quasi semigroup interval bilinear algebra over the semigroup $S = 2Z^+ \cup \{0\}$.
   i. Find substructures of V.
   ii. Find a bilinear operator on V
   iii. Can V be a made into a quasi fuzzy semigroup interval bilinear algebra?
   iv. Is V pseudo simple?
   v. Is V doubly simple?

53. Give an example of a doubly simple quasi semigroup interval bilinear algebra.

54. Describe some important properties enjoyed by group interval bivector space.

55. Give an example of a simple group interval bivector space.

56. Let $V = V_1 \cup V_2 = \left\{ \begin{bmatrix} [0,a_1] \\ [0,a_2] \\ \vdots \\ [0,a_9] \end{bmatrix}, \left([0,a][0,b]\right) \middle| a,b,a_i \in Z_{23} \right\} \cup$

$$\left\{ \begin{bmatrix} [0,a_1] & [0,a_2] & [0,a_3] & [0,a_4] \\ [0,a_5] & [0,a_6] & [0,a_7] & [0,a_8] \\ [0,a_9] & [0,a_{10}] & [0,a_{11}] & [0,a_{12}] \end{bmatrix}, \sum_{i=0}^{\infty}[0,a_i]x^i \left| \begin{array}{l} a_j, a_i \in Z_{23}; \\ 1 \le j \le 12 \end{array} \right. \right\}$$

be a group interval bivector space over the group $G = Z_{23}$.
   i. What is the bidimension of V?
   ii. Find group interval bivector subspaces of V.
   iii. Is V pseudo simple? Justify.
   iv. Find all generating bisubset of V.

57. Let $V = V_1 \cup V_2 = \left\{ \sum_{i=0}^{7}[0,a_i]x^i \,\middle|\, a_i \in Z_{40}; 0 \le i \le 7 \right\} \cup$

$$\left\{ \begin{bmatrix} [0,a_1] \\ [0,a_2] \\ \vdots \\ [0,a_{10}] \end{bmatrix} \middle| \begin{array}{l} a_i \in Z_{40}; \\ 1 \le i \le 10 \end{array} \right\}$$ be a group interval bivector space over

the group $G = Z_{40}$.
   i. Find at least two group interval bivector subspaces of V.
   ii. Is V simple?
   iii. Prove V is not pseudo simple.
   iv. Find a bibasis of V.
   v. Find the bidimension of V.

58. Give an example of a simple group interval bivector space.

59. Give an example of a doubly simple group interval bivector space.

60. Give an example of a pseudo simple group interval bivector space which is not simple.

61. Let $V = V_1 \cup V_2 = $ {all $10 \times 10$ intervals matrices with intervals of the form $[0, a_i]$; $a_i \in Z_{43}$} $\cup$ {All $5 \times 8$ interval matrices with intervals of the form $[0, a_i]$; $a_i \in Z_{43}$} be a group interval bilinear algebra over the group $G = Z_{43}$.

i. Prove V is pseudo simple.
ii. Prove V is not simple
iii. Find atleast 3 group interval bilinear subalgebras of V.
iv. What is the bidimension of V?
v. Find a bigenerating bisubset of V.
vi. Find a bilinear operator on V.

62. Give an example of a quasi group interval bilinear algebra over the group $Z_n$.

63. Is $V = V_1 \cup V_2 = \left\{ \begin{bmatrix} [0,a] \\ [0,a] \\ [0,a] \end{bmatrix} \middle| a \in Z_3 \right\} \cup \{([0, a]\ [0, a]\ [0, a])\ /$

    $a \in Z_3\}$ a doubly simple group interval bilinear algebra over the group $G = Z_3$? Justify your claim.

64. Let $V = V_1 \cup V_2 = \left\{ \sum_{i=0}^{7} [0, a_i]x^i \middle| a_i \in Z_5 \right\} \cup$

    $\left\{ \begin{bmatrix} [0,a] & [0,a] \\ [0,a] & [0,a] \\ [0,a] & [0,a] \end{bmatrix} \middle| a \in Z_5 \right\}$ be a quasi group interval bilinear

    algebra over the group $G = Z_5$.
    i. Is V simple?
    ii. Is V doubly simple?
    iii. Is V pseudo simple?

65. Give some interesting properties about bigroup interval bilinear algebra.

66. Give an example of a simple bigroup interval bilinear algebra.

67. Prove all bigroup interval algebras built using the bigroups $Z_p \cup Z_p$ ( p and q two distinct primes) are always pseudo simple.

68. Can one have bigroup interval bivector spaces using positive reals or positive rationals? Substantiate your answer.

69. Let

$$V = V_1 \cup V_2$$

$$= \left\{ \begin{bmatrix} [0,a_1] \\ [0,a_2] \\ \vdots \\ [0,a_{12}] \end{bmatrix} \middle| \begin{array}{l} a_i \in Z_{50}; \\ 1 \leq i \leq 12 \end{array} \right\} \cup$$

$$\left\{ \begin{bmatrix} [0,a_1] & [0,a_2] & \ldots & [0,a_{10}] \\ [0,a_{11}] & [0,a_{12}] & \ldots & [0,a_{20}] \\ [0,a_{21}] & [0,a_{22}] & \ldots & [0,a_{30}] \end{bmatrix} \middle| \begin{array}{l} a_i \in Z_{28}; \\ 1 \leq i \leq 30 \end{array} \right\}$$

be a bigroup interval bilinear algebra over the bigroup $G = G_1 \cup G_2 = Z_{50} \cup Z_{28}$.
  i. Find atleast two subbigroup interval bilinear subalgebras.
  ii. Find atleast three bigroup interval bilinear subalgebras.
  iii. Find a generating biset of V.
  iv. Find a bilinear operator on V.

70. Let

$$V = V_1 \cup V_2$$

$$= \left\{ \sum [0,a_i]x^i \middle| a_i \in Z_7 \right\} \cup \left\{ \begin{bmatrix} [0,a_1] \\ [0,a_2] \\ \vdots \\ [0,a_{19}] \end{bmatrix} \middle| \begin{array}{l} a_i \in Z_{11}; \\ 1 \leq i \leq 19 \end{array} \right\}$$

be a bigroup interval bilinear algebra over the bigroup $G = G_1 \cup G_2 = Z_7 \cup Z_{11}$.
  i. Prove V is pseudo simple.
  ii. Find bigroup interval bilinear subalgebras of V.
  iii. Prove V is not doubly simple.
  iv. Define a bilinear operator $T = T_1 \cup T_2 : V_1 \cup V_2 \rightarrow V_1 \cup V_2$ such that biker $T \neq \{0\} \cup \{0\}$.
  v. Find a generating biset of V.
  vi. Is bidimension of V finite?

71. Let $V = V_1 \cup V_2 = \left\{ \sum [0, a_i]x^i \mid a_i \in Z^+ \cup \{0\} \right\} \cup \{$All $6 \times 6$

    interval matrices with intervals of the form $[0, a_i]$ where $a_i \in Z_{420}\}$ be a semigroup - group interval bilinear algebra over the semigroup - group $Z^+ \cup \{0\} \cup Z_{420}$.
    i. Find substructures of V.
    ii. Prove V is not a doubly simple space.
    iii. Find a $T : T_1 \cup T_2 : V_1 \cup V_2 \to V_1 \cup V_2$ such that bikerT $= \{0\} \cup \{0\}$.
    iv. Find $T : T_1 \cup T_2 : V = V_1 \cup V_2 \to V = V_1 \cup V_2$ such that biker $T = \{0\} \cup \{S\}$. $S \neq 0$.

72. Let $V = V_1 \cup V_2 = \left\{ \begin{bmatrix} [0,a_1] & [0,a_2] \\ [0,a_3] & [0,a_4] \\ [0,a_5] & [0,a_6] \end{bmatrix} \middle| \begin{array}{c} a_i \in Z^+ \cup \{0\}; \\ 1 \leq i \leq 6 \end{array} \right\} \cup$

    $\left\{ \sum_{i=0}^{\infty} [0, a_i]x^i \mid a_i \in Z_{47} \right\}$ be a set semigroup interval bilinear

    algebra over the set - semigroup $3Z^+ \cup \{0\} \cup Z_{47}$.
    i. Find a set - semigroup interval bilinear subalgebras.
    ii. Find a subset - subsemigroup interval bilinear subalgebras.
    iii. Find a bilinear bioperator on V which is one to one.

73. Give some interesting results about biset interval bilinear algebras.

74. Give examples of infinite biset interval bilinear algebra which is
    i. simple biset interval bilinear algebra.
    ii. Pseudo simple biset interval bilinear algebra.
    iii. Doubly simple biset interval bilinear algebra.

75. Let $V = V_1 \cup V_2$ be a set - semigroup interval bilinear algebra. Define a bilinear operator on V, which is one to one; where $V = \{$all $3 \times 3$ interval matrices with intervals of the form $[0, a_i]$; $a_i \in Z_{17}\} \cup \{ \Sigma [0, a_i] x^i \mid a_i \in Z_{17})\} = V_1 \cup V_2$ is a set-semigroup interval bilinear algebra over the set-semigroup $S = \{3Z^+, \{0\}, 5Z^+\} \cup Z_{17} = S_1 \cup S_2$.
  i. Is V pseudo simple?
  ii. Does V have set semigroup interval bilinear subalgebra?
  iii. Is V doubly simple? Justify.

76. Give an example of a set-semigroup interval bivector space which is not a set-semigroup interval bilinear algebra of finite over which is doubly simple!

77. Let $V = V_1 \cup V_2$ = {Collection of $2 \times 10$ interval matrices with entries form $Z_{15}$} $\cup$ {$3 \times 5$ interval matrices with entries from $Z_{12}$} be a bisemigroup interval bilinear algebra over the bisemigroup $S = S_1 \cup S_2 = Z_{15} \cup Z_{12}$.
  i. Is V simple?
  ii. Is V pseudo simple?
  iii. Is V doubly simple?
  iv. Is V finite?
  v. Determine a bilinear operator which is not one to one.
  vi. What is the bidimension of V?

78. Determine some interesting properties about group-semigroup interval bivector spaces of finite dimension?

79. Is it possible to have a bigroup interval bilinear algebra of infinite dimension?

80. Is it possible to have a bigroup interval bivector space of finite dimension?

81. Give an example of a set - group interval bilinear algebra of bidimension (9, 8).

82. Let $V = V_1 \cup V_2$ = {all $5 \times 6$ interval matrices of the form $[0, a_i]$; $a_i \in Z_7$} $\cup$ {all $6 \times 5$ interval matrices over the set $Z_7$ with intervals of the form $[0, a_i]$} be a group interval bilinear algebra over the group $Z_7$.
  i. Is V simple?
  ii. Prove V is pseudo simple.

    iii. Prove V is not doubly simple.
    iv. What is the biorder of V?
    v. Define a one to one bilinear operator on V.

83. Prove if $V = V = V_1 \cup V_2$ is a group interval bilinear algebra over a group G; then the set of all bioperators on V is again a group interval bilinear algebra over the group G.

84. Obtain some interesting properties about bilinear transformations on group interval bivector spaces $V = V_1 \cup V_2$ and $W = W_1 \cup W_2$ defined over the group G.

85. Let $V = V_1 \cup V_2 = \{$all $3 \times 1$ be a set of all interval matrices with intervals of the form $[0, a_i]$; $a_i \in Z_{12}\} \cup \{1 \times 7$ be a set of all interval matrices with intervals of the form $[0, a_i]$; $a_i \in Z_{19}\}$ be a bisemigroup interval bilinear algebra over the bisemigroup $S = S_1 \cup S_2 = Z_{12} \cup Z_{19}$.
    i. Find all bisemigroup interval bilinear subalgebras of V over S.
    ii. Is V pseudo simple?
    iii. Find a bigenerating subset of V.
    iv. Find all bilinear operators on V.

86. Let $V = V_1 \cup V_2$ be as in problem 77.
    i. Find a bigenerating subset of V.
    ii. What is the bidimension of V over S?

87. Obtain some interesting properties on set - group interval bivector spaces of finite order.

88. Let $V = V_1 \cup V_2 = \left\{ \sum_{i=0}^{\infty}[0,a_i]x^i \,\middle|\, a_i \in Z_3 \right\} \cup \{$All $8 \times$ interval matrices with intervals of the form $[0, a_i]$ where $a_i \in Z_{42}\} \cup \left\{ \sum[0,a_i]x^i \,\middle|\, a_i \in Z^+ \cup \{0\} \right\}$ be a group-semigroup bilinear algebra defined over the group-semigroup $S = S_1 \cup S_2 = Z_{42} \cup 3Z^+ \cup \{0\}$.

     i. Find atleast 3 group-semigroup interval bilinear subalgebras of V over S.
     ii. Find atleast 3 pseudo subgroup-subsemigroup interval

89. Give an example of a doubly simple set - group interval bivector space.

90. Give an example of a pseudo simple bigroup interval bivector space which is not simple.

91. Let $V = V_1 \cup V_2 = \left\{ \begin{bmatrix} [0,a] \\ [0,a] \\ [0,a] \\ [0,a] \end{bmatrix} \,\middle|\, a_i \in Z_5 \right\} \cup \{([0, a], [0, a], [0,$ a], [0, a], [0, a], [0, a]) / $a \in Z_{11}\}$ be a bigroup interval algebra over the bigroup $G = G_1 \cup G_2 = Z_5 \cup Z_{11}$.
    i. Is V simple?
    ii. Is V doubly simple?
   iii. Is V pseudo simple? Justify your claim.

92. Prove there exists an infinite class of doubly simple bigroup interval bilinear algebras!

93. Prove their exists an infinite class of bigroup interval bilinear algebras which are not pseudo simple!

94. Does there exists an infinite classes of set-group interval bilinear algebras? Justify your claim.

95. Does there exist a bigroup interval bilinear algebras built over the bigroup $G = G_1 \cup G_2$, where both $G_1$ and $G_2$ are of infinite order?

96. Give some innovative results on biset interval bivector spaces of infinite order.

97. Let V = {1 × 9 interval matrices using $Z_7$} ∪ {9 × 1 interval matrices using $Z_7$} be a group interval bilinear algebra over the group $Z_7$.
   i. Prove V is not doubly simple!
   ii. Find all bilinear operators on V and show it is a group interval bilinear algebra over $Z_7$.

98. Let V = $V_1$ ∪ $V_2$ = {all 2 × 2 interval matrices using $3Z^+$ ∪ {0} and $5Z^+$ ∪ {0}} ∪ {3 × 3 interval matrices using $7Z^+$ ∪ {0}, $3Z^+$ ∪ {0}} be a bisemigroup interval bilinear algebra over the bisemigroup S = $S_1$ ∪ $S_2$ = $5Z^+$ ∪ {0} ∪ $7Z^+$ ∪ {0}.
   i. Is V simple?
   ii. Find subbisemigroup interval bilinear subalgebras of V.

99. Show V in problem (98) is not doubly simple.

100. For V in problem (98) prove set of all bilinear operators on V is again a bisemigroup bilinear algebra over S.

101. Give an example of set-semigroup interval bilinear algebra which is not a set-group interval bilinear algebra.

102. Is every set-group interval bilinear algebra a set - semigroup interval bilinear algebra?

103. Show a biset interval bilinear vector space in general not a bigroup or bisemigroup interval bivector space.

104. Obtain conditions on a bigroup interval bilinear algebra V = $V_1$ ∪ $V_2$ so that V is never a nontrivial bisemigroup interval bilinear algebra.

105. Derive some important and interesting properties related with bisemigroup interval bivector spaces.

106. Give an example of a bisemigroup interval bivector space which is not a bisemigroup interval bilinear algebra.

107. Let $V = V_1 \cup V_2 = \left\{ \sum_{i=0}^{\infty} [0,a_i]x^i \mid a_i \in 3Z^+ \cup \{0\} \right\} \cup \{\Sigma [0, a_i]$ $x^i \mid a_i \in 5Z^+ \cup \{0\}\}$ be a bisemigroup interval bilinear algebra defined over the bisemigroup $S = S_1 \cup S_2 = 3Z^+ \cup \{0\} \cup 5Z^+ \cup \{0\}$.
    i. Find a bigenerating subset of V.
    ii. Find atleast two bisemigroup interval bilinear subalgebras.
    iii. Find atleast two subbisemigroup interval bilinear subalgebras.

108. Give an example of a pseudo simple bisemigroup interval bilinear algebra.

109. Give an example of a simple bisemigroup interval bilinear algebra.

110. Give an example of a doubly simple bisemigroup interval bilinear algebra of finite order.

111. Give an example of a group - set interval bilinear algebra of infinite order.

112. Give an example of a group semigroup interval bilinear-algebra of finite order.

113. Prove a set-group interval bivector space in general is not a semigroup-group interval bivector space.

# FURTHER READING

41. VASANTHA KANDASAMY, W.B., On a new class of semivector spaces, *Varahmihir J. of Math. Sci.*, **1**, 23-30, 2003.

42. VASANTHA KANDASAMY and THIRUVEGADAM, N., Application of pseudo best approximation to coding theory, *Ultra Sci.,* **17**, 139-144, 2005.

43. VASANTHA KANDASAMY and RAJKUMAR, R. A New class of bicodes and its properties, (To appear).

44. VASANTHA KANDASAMY, W.B., On fuzzy semifields and fuzzy semivector spaces, *U. Sci. Phy. Sci.,* **7,** 115-116, 1995.

45. VASANTHA KANDASAMY, W.B., On semipotent linear operators and matrices, *U. Sci. Phy. Sci.*, **8**, 254-256, 1996.

46. VASANTHA KANDASAMY, W.B., Semivector spaces over semifields, *Zeszyty Nauwoke Politechniki*, **17,** 43-51, 1993.

47. VASANTHA KANDASAMY, W.B., *Smarandache Fuzzy Algebra,* American Research Press, Rehoboth, 2003.

48. VASANTHA KANDASAMY, W.B., *Smarandache rings,* American Research Press, Rehoboth, 2002.

49. VASANTHA KANDASAMY, W.B., Smarandache semirings and semifields, *Smarandache Notions Journal*, **7** 88-91, 2001**.**

50. VASANTHA KANDASAMY, W.B., *Smarandache Semirings, Semifields and Semivector spaces*, American Research Press, Rehoboth, 2002.

51. VASANTHA KANDASAMY, W.B., SMARANDACHE, Florentin and K. ILANTHENRAL, *Introduction to Linear Bialgebra,* Hexis, Phoenix, 2005.

52. VASANTHA KANDASAMY, W.B., and SMARANDACHE, Florentin, *Fuzzy Interval Matrices, Neutrosophic Interval Matrices and their Application,* Hexis, Phoenix, 2005.

53. VASANTHA KANDASAMY, W.B., SMARANDACHE, Florentin and K. ILANTHENRAL, *Set Linear Algebra and Set Fuzzy Linear Algebra,* InfoLearnQuest, Phoenix, 2008.
240

# INDEX

# ABOUT THE AUTHORS

**Dr.W.B.Vasantha Kandasamy** is an Associate Professor in the Department of Mathematics, Indian Institute of Technology Madras, Chennai. In the past decade she has guided 13 Ph.D. scholars in the different fields of non-associative algebras, algebraic coding theory, transportation theory, fuzzy groups, and applications of fuzzy theory of the problems faced in chemical industries and cement industries. She has to her credit 646 research papers. She has guided over 68 M.Sc. and M.Tech. projects. She has worked in collaboration projects with the Indian Space Research Organization and with the Tamil Nadu State AIDS Control Society. She is presently working on a research project funded by the Board of Research in Nuclear Sciences, Government of India. This is her $51^{st}$ book.

On India's 60th Independence Day, Dr.Vasantha was conferred the Kalpana Chawla Award for Courage and Daring Enterprise by the State Government of Tamil Nadu in recognition of her sustained fight for social justice in the Indian Institute of Technology (IIT) Madras and for her contribution to mathematics. The award, instituted in the memory of Indian-American astronaut Kalpana Chawla who died aboard Space Shuttle Columbia, carried a cash prize of five lakh rupees (the highest prize-money for any Indian award) and a gold medal.
She can be contacted at vasanthakandasamy@gmail.com
Web Site: http://mat.iitm.ac.in/home/wbv/public_html/

**Dr. Florentin Smarandache** is a Professor of Mathematics at the University of New Mexico in USA. He published over 75 books and 150 articles and notes in mathematics, physics, philosophy, psychology, rebus, literature.

In mathematics his research is in number theory, non-Euclidean geometry, synthetic geometry, algebraic structures, statistics, neutrosophic logic and set (generalizations of fuzzy logic and set respectively), neutrosophic probability (generalization of classical and imprecise probability). Also, small contributions to nuclear and particle physics, information fusion, neutrosophy (a generalization of dialectics), law of sensations and stimuli, etc. He can be contacted at smarand@unm.edu